\documentstyle{amsppt}
\magnification=1200
\catcode`\@=11
\redefine\logo@{}
\catcode`\@=13

\define \bn{\Bbb N}
\define \bz{\Bbb Z}
\define \bq{\Bbb Q}
\define \br{\Bbb R}
\define \bc{\Bbb C}

\define \bk{\Bbb K}
\define \bo{\Bbb O}
\define \M{{\Cal M}}
\define\Ha{{\Cal H}}
\define\La{{\Cal L}}

\define\T{{\Cal T}}
\define\N{{\Cal N}}

\define\rk{\text{rk}~}


\define\pr{\text{pr}}

\define\arc{\text{arc}}

\define\sh{\text{sh}}
\define\ch{\text{ch}}
\define\th{\text{th}}

\define\NE{\text{NE}}
\define\NS{\text{NS}}
\define\Exc{\text{Exc}}

\TagsOnRight

\document

\topmatter

\title
On the classification of hyperbolic root systems of
the rank three. Part III
\endtitle

\author
Viacheslav V. Nikulin \footnote{Supported by
Grant of Russian Fund of Fundamental Research
\hfill\hfill}
\endauthor

\address
Steklov Mathematical Institute,
ul. Gubkina 8, Moscow 117966, GSP-1, Russia
\endaddress
\email
slava\@nikulin.mian.su 
\endemail

\abstract
See Parts I and II in alg-geom/9711032 and 9712033.  

Here we classify maximal hyperbolic root systems of the 
rank three having restricted arithmetic type and a generalized 
lattice Weyl vector $\rho$ with $\rho^2<0$ (i. e. of the hyperbolic type). 
We give classification of all reflective of hyperbolic type elementary 
hyperbolic lattices of the rank three. For elliptic (when $\rho^2>0$) 
and parabolic (when $\rho^2=0$) types it was done in Parts I and II.

We apply the same arguments as for elliptic and parabolic types: 
the method of narrow parts of polyhedra in hyperbolic spaces,  
and class numbers of central 
symmetries. But we should say that for the hyperbolic type all 
considerations are much more complicated and required much more 
calculations and time.  

These results are important, for example, for Theory of 
Lorentzian Kac--Moody algebras and some aspects of Mirror Symmetry. 
We also apply these results to prove boundedness of families of algebraic 
surfaces with almost finite polyhedral Mori cone 
(see math.AG/9806047 about this subject). 
\endabstract

\rightheadtext
{Hyperbolic root systems, III}
\leftheadtext{V.V. Nikulin}
\endtopmatter

\head
Introduction
\endhead

See Introductions in Parts I and II. 

Here we give classification of maximal hyperbolic root systems of the 
rank three having restricted arithmetic type and a generalized 
lattice Weyl vector 
$\rho$ with $\rho^2<0$ (i. e. of the hyperbolic type). More exactly, 
we classify all 
reflective of hyperbolic type (or hyperbolically reflective) 
elementary hyperbolic lattices of the rank three. 
For elliptic (when $\rho^2>0$) or parabolic (when $\rho^2=0$) types it was 
done in Parts I and II. 

The main problem for the hyperbolic type is to develop an appropriate 
method of narrow parts of polyhedra of restricted hyperbolic type in 
hyperbolic spaces. This is much more complicated than for polyhedra of 
elliptic and restricted parabolic types when it 
is known long time ago (see \cite{N4}, \cite{N5}, \cite{N11}). 
In Sect. 8 we develop the method of narrow parts for restricted hyperbolic 
polygons in a hyperbolic plane. These results are crucial for further 
developing of this method in higher dimension. We mention that the first 
results about the method of narrow parts for the hyperbolic type were  
obtained in \cite{N13}. 

Except of hyperbolic root systems of hyperbolic type which is 
the main subject of the paper (see results in Sect. 7), 
we here apply the method of narrow parts of restricted hyperbolic 
polyhedra in hyperbolic spaces also to finiteness of the set of 
reflective hyperbolic lattices of hyperbolic 
type over purely real algebraic number fields (Theorem 8.4.1.1) 
and to boundedness of families of algebraic surfaces with almost finite 
polyhedral Mori cone (Theorem 8.4.2.1). 

We continue numeration of sections and statements 
from the Parts I and II.  

\smallpagebreak 

This paper was written during my stay in Steklov Mathematical Institute, 
Mos\-cow, and Max-Planck-Institut f\"ur Mathematik, Bonn. I am grateful 
to the Institutes for hospitality.  

\head
7. Classification of reflective hyperbolic lattices of the rank $3$
and of hyperbolic type: formulations
\endhead

Like for elliptic or parabolic type (see Sect. 2), 
we reduce (or restrict) the classification to hyperbolically 
reflective hyperbolic lattices $S$ of the rank three with 
square-free determinant $d$. All main 
(i. e. which are even when $d$ is even) these reflective 
lattices are given in Table 4 below. We give the invariants 
$d$, $\mu$ of $S$ which define the $S$ up to isomorphism. 
We also give the lattice $S$ itself for some standard bases 
(see Sect. 2.3.1). 
It may have one of three forms:  
$U\oplus \langle -n\rangle$ or  
$\langle n_1\rangle \oplus \langle n_2\rangle 
\oplus \langle n_3\rangle(\epsilon_1/2,\,\epsilon_2/2,\,\epsilon_3/2)$ 
where $\epsilon_i=0$ or $1$, or 
$\langle n \rangle \oplus A_2(1/3,\,1/3,\,-1/3)$ where 
$$
A_2=\pmatrix-2&1\\1&-2\endpmatrix .
\tag{7.1}
$$    
We also give the invariant $h$. For hyperbolic type it is equal to  
$0$ or $2$. In Table 4  we also describe the fundamental 
polygon $\M$ for $W(S)$ (equivalently, we describe the  
set $P(\M)_{\pr}$ in the standard bases of $S$). 
To describe $P(\M)_{\pr}$, we give the following 
data. We give the primitive generalized lattice Weyl vector $\rho$ 
(for the fixed $\M$, the $\rho$ is unique up to replacing 
by $-\rho$) and  its square $\rho^2$. The set $P(\M)_{\pr}$ is 
divided in two subsets: 
$$
{P(\M)_{\pr}}_+=\{\alpha \in P(\M)_{\pr}\ |\ (\alpha,\,\rho)>0\}
\tag{7.2} 
$$ 
and 
$$
{P(\M)_{\pr}}_-=\{\alpha \in P(\M)_{\pr}\ |\ (\alpha,\,\rho)<0\}. 
\tag{7.3} 
$$
Both these subsets give two infinite chains of consecutive sides of 
$\M$ and are naturally ordered (up to changing the order to 
the opposite one). To describe these chains, we give the part  
$\Gamma_+$ (and $\Gamma_-$ if necessary) of the chain 
${P(\M)_{\pr}}_+$ (respectively ${P(\M)_{\pr}}_-$) 
and the Gram matrix $G(\Gamma_+)$ (respectively 
$G(\Gamma_-)$). We also give generators of the group 
$A(P(\M)_{\pr})$ of symmetries of $\M$. These data are 
sufficient to find the set $P(\M)_{\pr}$ because we have    
$$
P(\M)_{\pr}=A(P(\M)_{\pr})(\Gamma_+)
\tag{7.4}
$$
(or we have 
$$
P(\M)_{\pr}=A(P(\M)_{\pr})(\Gamma_+\cup \Gamma_-)
\tag{7.5}
$$
if we also give the set $\Gamma_-$). If the invariant $h=2$, then 
the group $A(P(\M)_{\pr})$ is generated by two central symmetries 
$C_1$ and $C_2$, and it is sufficient to give only $\Gamma_+$ 
because $C_1(\rho)=C_2(\rho)=-\rho$.   
If $h=0$ (there are five these cases) 
and the group $A(P(\M)_{\pr})$ is generated by the  
skew-symmetry $S_1$ along the line orthogonal to $\rho$ (there are 
three these cases), then it is also sufficient to give only the 
set $\Gamma_+$ because $S_1(\rho)=-\rho$. 
Only for two cases (when $(d,\,\mu,\,h)=(65,\,0,\,0)$, $(74,\,1,\,0)$),     
the group $A(P(\M)_{\pr})$ is generated by the translation $T$ 
along the line orthogonal to $\rho$. Then we have to give also 
the set $\Gamma_-$, and we have \thetag{7.5}.  

Thus, we have the following main result: 

\proclaim{Basic Theorem 7.1} Table 4 below gives the complete
list (containing 66 cases numerated by $N$) of main (i.e.
even for even determinant) hyperbolic lattices $S$ of rank $3$ and
with square-free determinant which are reflective of hyperbolic 
type. In Table 4, for the lattice $S$
we give: 
invariants $(d,\,\eta)$ defining $S$ up to isomorphism;
the number $h$ of classes of central symmetries of the group 
$A(P(\M)_{\pr})$ of symmetries of the fundamental polygon $\M$ for  
$W(S)$;  
the matrix of $S$ for some standard bases; 
the generalized lattice Weyl vector $\rho$ and its square 
$\rho^2$;  
a subset $\Gamma_+\subset A(P(\M)_{\pr})$ (or 
two subsets $\Gamma_+\subset A(P(\M)_{\pr})$ and 
$\Gamma_-\subset A(P(\M)_{\pr})$, if  
$(d,\,\mu,\,h)=(65,\,0,\,0)$, $(74,\,1,\,0)$) of 
orthogonal vectors to consecutive sides of the $\M$ such that 
$P(\M)_{\pr}= A(P(\M)_{\pr})(\Gamma_+)$ 
(respectively 
($P(\M)_{\pr}= A(P(\M)_{\pr})(\Gamma_+\cup \Gamma_-))$; 
generators of the group $A(P(\M)_{\pr})$.  

As generators of $A(P(\M)_{\pr})$, we give central symmetries 
$C_1$ and $C_2$, if $h=2$; the sliding symmetry $S_1$ if 
$h=0$ and $(d,\,\mu)=(95,\,0)$, $(114,\,2)$, $(231,\,2)$; 
the translation $T$ if $h=0$ and $(d,\,\mu)=(65,\,0)$, $(74,\,1)$.  

If the lattice $S$ represents $0$ (i.e. there exists a non-zero
$x\in S$ with $x^2=0$), we give $S$ in the form
$S=U\oplus \langle -d \rangle $.
Then the fundamental polygon $\M$ contains an infinite vertex. 
In particular, Table 4 contains the complete list of reflective 
hyperbolic lattices of hyperbolic type which have the form 
$S=U\oplus\langle -d \rangle$ where $d$ is square-free. 
There are $8$ these cases corresponding to
$$ 
d=19,\,23,\,35,\,39,\,46,\,58,\,62,\,70. 
$$
If $S$ is given in different form, the lattice $S$
does not represent $0$ and $\M$ does not have infinite 
vertices. There are $58$ these cases.
\endproclaim

\demo{Idea of Proof}  
We use the same method as for the proof of Basic Theorem 2.3.2.1. 
If $S$ is hyperbolically reflective and has the rank three, then 
the invariant $h$ is equal to $0$ or $2$. Table 3 gives the 
list of all main hyperbolic lattices with $h=0$ or $h=1$ which 
have $d\le 100000$. It has 61 lattices with $h=0$. Using Vinberg's 
algorithm \cite{V2}, we check reflective type of them. Only five of them 
are hyperbolically reflective and are given in Table 4 (they are marked by 
"$hr$" in Table 3. All other are either elliptically reflective 
and are given in Table 1 (they are marked by $er$ in Table 3) or are 
not reflective (they are marked by $nr$ in Table 3). 
Using the formula for $h$ of Theorem 3.2.1, 
we get the list of all main lattices $S$ with $d\le 100000$ and 
$h=2$. It contains 259 lattices $S$ given in Table 6. Checking by the 
Vinberg's algorithm reflective type of these 259 lattices, we find 
the rest 61 lattices with $h=2$ of Table 4. 

Thus, Table 4 contains all main hyperbolic lattices with $d\le 100000$. 
Conjecturally the Tables 3 and 6 contain all main hyperbolic lattices 
with square-free determinant and $h\le 2$. To avoid this conjecture, in 
Sect. 8 we develop-e the ``method of narrow parts of polyhedra'' 
for polygons of restricted hyperbolic type in 2-dimensional 
hyperbolic spaces. It gives finiteness of the set of all 
hyperbolically reflective hyperbolic lattices of rank three. 
Applying this method to main hyperbolic lattices $S$ with 
square-free determinant, we get a finite list of their 
invariants $(d,\,\eta)$. 
We calculate the invariant $h$ for all of them, and we find that 
$d\le 100000$ if $h=0$ or $h=2$. 
It proves that the Tables 4 and 6 contain all main reflective 
hyperbolic lattices and finishes the proof of Basic Theorem 7.1. 
See details in Sect. 9.1. 
\enddemo 

From Basic Theorem 2.3.2.1 and Basic Theorem 7.1, we get (similarly to  
Theorem 2.3.3.1) classification of non-main hyperbolically reflective 
hyperbolic lattices of rank three and 
with square-free determinant.
 
\proclaim{Theorem 7.2} Table 5 below gives the complete list
(containing $21$ cases numerated by $N^\prime$) of non-main
(i. e. odd with even determinant) hyperbolic lattices $\widetilde{S}$
with square-free determinant and of rank three which are
hyperbolically reflective. In Table 7,
for each lattice $\widetilde{S}$
we give its invariants $(d,\text{odd},\eta)$, the number
$h$ of classes of central symmetries, matrix of $\widetilde{S}$ in
a standard bases, invariants $(d,\eta+\omega(d))$
of the main odd lattice $S$ of $\widetilde{S}$ (see Proposition
2.2.6).

If $S$ and $\widetilde{S}$ are equivariantly
equivalent (in Table 5 there are $11$ these cases),
we only give the lattice $S$. For this case,
calculation of $P(\M)_{\pr}$ and
$G(P(\M)_{\pr})$ for $\widetilde{S}$ follows from calculation
of similar sets for $S$ in Table  4 
(see Remark 2.3.3.2). In particular, all $\widetilde{S}$
representing $0$ (non-compact case) have this type, then
$\widetilde{S}=\langle 1 \rangle \oplus \langle -1 \rangle \oplus
\langle -d \rangle$ where
$$
d=38,\,46,\,70,\,78.
$$
If $S$ and $\widetilde{S}$ are not equivariantly equivalent
(there are 10 these cases), we (like for Table 4) give
$P(\M)_{\pr}$ and $G(P(\M)_{\pr})$ for $\widetilde{S}$ in
the standard bases.
\endproclaim

\demo{Proof} The proof is similar to proof of Theorem 2.3.3.1. 
See details in Sect. 9.1.2. 
\enddemo

From Theorems 7.1 and 7.2 we get (like in Sect. 2.3.4) 
classification of all elementary hyperbolic lattices of rank three 
which are hyperbolically reflective.

\proclaim{Theorem 7.3} A primitive elementary hyperbolic lattice
$F$ is hyperbolically reflective if and only if
$$
F\cong S^{\ast,m}(m),\ m\mid d
\tag{7.6}
$$
(see \thetag{2.2.5}) where $S$ is a lattice of the determinant $d=\det(S)$
(it is square-free) of Tables 4 or 5. Then $\det(F)=dm$.
In particular, for the fixed lattice $S$ the number
of lattices $F$ is equal to $2^t$ where $t$ is
the number of prime divisors of $d$.

Let $P(\M)_{\pr}$ be the set of orthogonal to $\M$ primitive roots of
$S$. Then the set $\widetilde{P}(\M)_{\pr}$ of orthogonal
to $\M$ primitive roots of $F$ is equal to
$$
\widetilde{P}(\M)_{\pr}=
\{\widetilde{\alpha}=\alpha/k_{\alpha,m}\ \mid \
\alpha \in P(\M)_{\pr}\}
\tag{7.7}
$$
where $k_{\alpha,m}$ is the greatest divisor of $m$ such that
$\alpha/k_{\alpha,m}\in S^\ast$ (equivalently, if
$(\alpha,S)=t\bz$, then $k_{\alpha,m}=\text{g.c.d}(t,m)$).
The Gram matrix of $\widetilde{P}(\M)_{\pr}$ for the lattice $F$
is equal to
$$
\left((\widetilde{\alpha},\widetilde{\beta})\right)=
\left({(\alpha,\beta)m\over k_{\alpha,m}k_{\beta,m}}\right),\ \
\widetilde{\alpha},\,\widetilde{\beta} \in \widetilde{P}(\M)_{\pr}.
\tag{7.8}
$$
Here, on the left hand side of \thetag{7.8}
we use the form of the
lattice $F$, and on the right hand side of \thetag{7.8}
we use the form of the lattice $S$.
\endproclaim

By Proposition 2.2.2, any hyperbolic lattice has 
an equivariant embedding (and even canonical embedding) to an elementary 
hyperbolic lattice (of the same rank). It follows that any 
reflective lattice has an equivariant embedding to a reflective elementary 
lattice. Thus all reflective lattices are equivariant sublattices 
of reflective elementary lattices. We have the following interesting 
statement which permits to describe many non-elementary hyperbolically 
reflective lattices. 

\proclaim{Theorem 7.4} Suppose that $S\subset F$ is an equivariant 
sublattice of a hyperbolically reflective elementary hyperbolic lattice 
$F$ of the rank 3.  
(All the lattices $F$ and their root systems are described in 
Theorems 7.1 --- 7.3.) Then $S$ 
is reflective (it is then automatically hyperbolically reflective),  
if and only if $W(S)=W(F)$ (equivalently, 
$S$ is $W(F)$-invariant). Equivalently, for any root 
$r\in F$, there exists $\lambda \in \bn$ such that $\lambda r$ is a root 
of $S$.   
\endproclaim    

\demo{Proof} See Sect. 9.3. 
\enddemo 

All other (different from Theorem 7.4) hyperbolically reflective 
hyperbolic lattices $S$ of the rank three have an equivariant embedding to 
reflective lattices $F$ classified in Theorem 2.3.4.1 (they 
are elliptically reflective). Classification of these lattices $S$ 
is more difficult.  

\newpage  

\centerline{\bf Table 4}

\centerline{The list of main hyperbolic lattices}
\centerline{of rank three and with square-free determinant }
\centerline{which are reflective of hyperbolic type}

\vskip20pt 

\vbox{\noindent
$N=1\ $ $d=17,\,\eta=1,\,h=2$: $S=\langle 51 \rangle \oplus A_2(1/3,1/3,-1/3)$;
\nobreak
\newline
$\rho= (3, -17, -17)$,\ $\rho^2=-119$; 
\newline 
$\Gamma_+=
\pmatrix{0}&{0}&{1}\cr{0}&{1}&{0}\cr
\endpmatrix,  
\hskip20pt 
G(\Gamma_+)=\pmatrix{-2}&{1}\cr{1}&{-2}\cr
\endpmatrix$;
\nobreak
\newline
$C_1=
\pmatrix{{31}\over{3}}&{{4}\over{3}}&{{2}\over{3}}\cr
{{-170}\over{3}}&{{-23}\over{3}}&{{-10}\over{3}}\cr
{{-136}\over{3}}&{{-16}\over{3}}&{{-11}\over{3}}\cr
\endpmatrix,
\hskip30pt
C_2=\pmatrix
{{65}\over{3}}&{{2}\over{3}}&{{10}\over{3}}\cr
{{-238}\over{3}}&{{-10}\over{3}}&{{-35}\over{3}}\cr
{{-374}\over{3}}&{{-11}\over{3}}&{{-58}\over{3}}\cr
\endpmatrix$.}

\vskip5pt

\vbox{\noindent
$N=2\ $ $d=19,\,\eta=1,\,h=2$: $S=U \oplus \langle -19 \rangle$;
\nobreak
\newline
$\rho= (57, 19, -11)$, $\rho^2=-133$;
\newline 
$\Gamma_+=
\pmatrix{9}&{1}&{-1}\cr{19}&{0}&{-1}\cr{0}&{0}&{1}\cr{-1}&{1}&{0}\cr
{3}&{3}&{-1}\cr\endpmatrix,  
\hskip20pt 
G(\Gamma_+)=
\pmatrix{-1}&{0}&{19}&{8}&{11}\cr{0}&{-19}&{19}&{19}&{38}\cr
{19}&{19}&{-19}&{0}&{19}\cr{8}&{19}&{0}&{-2}&{0}\cr{11}&{38}&{19}&{0}&{-1}\cr
\endpmatrix$;\nobreak
\newline
$C_1=
\pmatrix{38}&{169}&{494}\cr{9}&{38}&{114}\cr{-6}&{-26}&{-77}\cr
\endpmatrix,
\hskip30pt
C_2=
\pmatrix{19}&{50}&{190}\cr{8}&{19}&{76}\cr{-4}&{-10}&{-39}\cr
\endpmatrix$.}

\vskip5pt 

\vbox{\noindent
$N=3 \ $ $d=23 ,\,\eta=1 ,\,h=2$: $S=U \oplus \langle -23 \rangle$;
\nobreak
\newline
$\rho=(46, 23, -10)$,\ $\rho^2= -184$; 
\newline 
$\Gamma_+=
\pmatrix{15}&{3}&{-2}\cr{161}&{23}&{-18}\cr{11}&{1}&{-1}\cr
{23}&{0}&{-1}\cr{0}&{0}&{1}\cr{-1}&{1}&{0}\cr
\endpmatrix,  
\hskip20pt 
G(\Gamma_+)=
\pmatrix{-2}&{0}&{2}&{23}&{46}&{12}\cr{0}&{-46}&{0}&{115}&{414}&{138}\cr
{2}&{0}&{-1}&{0}&{23}&{10}\cr{23}&{115}&{0}&{-23}&{23}&{23}\cr
{46}&{414}&{23}&{23}&{-23}&{0}\cr{12}&{138}&{10}&{23}&{0}&{-2}\cr
\endpmatrix$;
\nobreak
\newline
$C_1=
\pmatrix{23}&{72}&{276}\cr{8}&{23}&{92}\cr{-4}&{-12}&{-47}\cr
\endpmatrix,
\hskip30pt
C_2=
\pmatrix{23}&{32}&{184}\cr{18}&{23}&{138}\cr{-6}&{-8}&{-47}\cr
\endpmatrix$.}

\vskip5pt

\vbox{\noindent
$N=4 \ $ $d=31 ,\,\eta=0 ,\,h=2 $: $S=\langle 31 \rangle \oplus 
\langle -1 \rangle \oplus \langle -1 \rangle$;
\nobreak
\newline
$\rho= (17,-31,-93) $,\ $\rho^2=-651$; 
\newline 
$\Gamma_+=
\pmatrix{1}&{-4}&{-4}\cr{0}&{-1}&{1}\cr{0}&{1}&{0}\cr
\endpmatrix,  
\hskip20pt 
G(\Gamma_+)=
\pmatrix{-1}&{0}&{4}\cr{0}&{-2}&{1}\cr{4}&{1}&{-1}\cr
\endpmatrix$;
\nobreak
\newline
$C_1=
\pmatrix{30}&{2}&{5}\cr{-62}&{-5}&{-10}\cr{-155}&{-10}&{-26}\cr
\endpmatrix,
\hskip30pt
C_2=
\pmatrix{123}&{2}&{22}\cr{-62}&{-2}&{-11}\cr{-682}&{-11}&{-122}\cr
\endpmatrix$.}

\vskip5pt 

\vbox{\noindent
$N=5 \ $ $d=35=5\cdot 7 ,\,\eta=1 ,\,h=2 $: 
$S=U\oplus \langle -35 \rangle$;
\nobreak
\newline
$\rho=(21, 7, -3) $,\ $\rho^2=-21 $; 
\newline 
$\Gamma_+=
\pmatrix{17}&{1}&{-1}\cr{35}&{0}&{-1}\cr{0}&{0}&{1}\cr{-1}&{1}&{0}\cr
{25}&{25}&{-6}\cr{23}&{19}&{-5}\cr
\endpmatrix,  
\hskip20pt 
G(\Gamma_+)=
\pmatrix{-1}&{0}&{35}&{16}&{240}&{171}\cr{0}&{-35}&{35}&{35}&{665}&{490}\cr
{35}&{35}&{-35}&{0}&{210}&{175}\cr{16}&{35}&{0}&{-2}&{0}&{4}\cr
{240}&{665}&{210}&{0}
&{-10}&{0}\cr{171}&{490}&{175}&{4}&{0}&{-1}\cr
\endpmatrix$;
\nobreak
\newline
$C_1=
\pmatrix{35}&{162}&{630}\cr{8}&{35}&{140}\cr{-4}&{-18}&{-71}\cr
\endpmatrix,
\hskip30pt
C_2=
\pmatrix{14}&{45}&{210}\cr{5}&{14}&{70}\cr{-2}&{-6}&{-29}\cr
\endpmatrix$.}

\vskip5pt 

\vbox{\noindent
$N=6 \ $ $d=37 ,\,\eta=1 ,\,h=2 $: 
$S=\langle 1 \rangle\oplus \langle -74 \rangle 
\oplus \langle -2 \rangle(0,1/2,1/2)$;
\nobreak
\newline
$\rho=(74, -17/2, -37/2)$,\ $\rho^2=-555$; 
\newline 
$\Gamma_+=
\pmatrix{43}&{-5}&{0}\cr{0}&{0}&{1}\cr{0}&{1}&{0}\cr{1}&{0}&{-1}\cr
\endpmatrix,  
\hskip20pt 
G(\Gamma_+)=
\pmatrix{-1}&{0}&{370}&{43}\cr{0}&{-2}&{0}&{2}\cr{370}&{0}&{-74}&{0}\cr
{43}&{2}&{0}&{-1}\cr
\endpmatrix$;
\nobreak
\newline
$C_1=
\pmatrix{168}&{1443}&{13}\cr
{{-39}\over{2}}&{{-335}\over{2}}&{{-3}\over{2}}\cr{{-13}
\over{2}}&{{-111}\over{2}}&{{-3}\over{2}}\cr
\endpmatrix,
\hskip30pt
C_2=
\pmatrix{24}&{185}&{15}\cr{{-5}\over{2}}&{{-39}\over{2}}&{{-3}\over{2}}\cr
{{-15}\over{2}}&{{-111}\over{2}}&{{-11}\over{2}}\cr
\endpmatrix$.}

\vskip5pt

\vbox{\noindent
$N=7 \ $ $d=39=3\cdot 13 ,\,\eta=1 ,\,h=2 $: 
$S= U \oplus \langle -39 \rangle$;
\nobreak
\newline
$\rho=(91,13,-8)$,\ $\rho^2=-130$; 
\newline 
$\Gamma_+=
\pmatrix{19}&{1}&{-1}\cr{39}&{0}&{-1}\cr{0}&{0}&{1}\cr{-1}&{1}&{0}\cr
{13}&{13}&{-3}\cr{6}&{3}&{-1}\cr{25}&{7}&{-3}\cr
\endpmatrix,  
\hskip10pt 
G(\Gamma_+)=
\pmatrix{-1}&{0}&{39}&{18}&{143}&{24}&{41}\cr
{0}&{-39}&{39}&{39}&{390}&{78}&{156}\cr
{39}&{39}&{-39}&{0}&{117}&{39}&{117}\cr{18}&{39}&{0}&{-2}&{0}&{3}&{18}\cr
{143}&{390}&{117}&{0}&{-13}&{0}&{65}\cr{24}&{78}&{39}&{3}&{0}&{-3}&{0}\cr
{41}&{156}&{117}&{18}&{65}&{0}&{-1}\cr
\endpmatrix$;
\nobreak
\newline
$C_1=
\pmatrix{351}&{3872}&{10296}\cr{32}&{351}&{936}\cr
{-24}&{-264}&{-703}\cr
\endpmatrix,
\hskip30pt
C_2=
\pmatrix{26}&{243}&{702}\cr{3}&{26}&{78}\cr{-2}&{-18}&{-53}\cr
\endpmatrix$.}

\vskip5pt

\vbox{\noindent
$N=8 \ $ $d=43 ,\,\eta=0 ,\,h=2 $: $S=\langle 43 \rangle\oplus 
\langle -1 \rangle \oplus  \langle -1 \rangle$;
\nobreak
\newline
$\rho=(47,-129,-301) $,\ $\rho^2= -12255 $; 
\newline 
$\Gamma_+=
\pmatrix{0}&{-1}&{1}\cr{0}&{1}&{0}\cr{59}&{0}&{-387}\cr{9}&{-2}&{-59}\cr
\endpmatrix,  
\hskip20pt 
G(\Gamma_+)=
\pmatrix{-2}&{1}&{387}&{57}\cr{1}&{-1}&{0}&{2}\cr{387}&{0}&{-86}&{0}\cr
{57}&{2}&{0}&{-2}\cr
\endpmatrix$;
\nobreak
\newline
$C_1=
\pmatrix{42}&{4}&{5}\cr{-172}&{-17}&{-20}\cr{-215}&{-20}&{-26}\cr
\endpmatrix,
\hskip30pt
C_2=
\pmatrix{171}&{2}&{26}\cr{-86}&{-2}&{-13}\cr{-1118}&{-13}&{-170}\cr
\endpmatrix$.}

\vskip5pt 

\vbox{\noindent
$N=9 \ $ $d=46=2\cdot 23 ,\,\eta=1 ,\,h=2 $: 
$S=U \oplus \langle -46 \rangle$;
\nobreak
\newline
$\rho= (46,23,-7) $,\ $\rho^2=-138$; 
\newline 
$\Gamma_+=
\pmatrix{11}&{2}&{-1}\cr{23}&{0}&{-1}\cr{0}&{0}&{1}\cr{-1}&{1}&{0}\cr
\endpmatrix,  
\hskip20pt 
G(\Gamma_+)=
\pmatrix{-2}&{0}&{46}&{9}\cr{0}&{-46}&{46}&{23}\cr{46}&{46}&{-46}&{0}\cr
{9}&{23}&{0}&{-2}\cr
\endpmatrix$;
\nobreak
\newline
$C_1=
\pmatrix{23}&{64}&{368}\cr{9}&{23}&{138}\cr{-3}&{-8}&{-47}\cr
\endpmatrix,
\hskip30pt
C_2=
\pmatrix{23}&{36}&{276}\cr{16}&{23}&{184}\cr{-4}&{-6}&{-47}\cr
\endpmatrix$.}

\vskip5pt

\vbox{\noindent
$N=10 \ $ $d=47,\,\eta=0,\,h=2$: 
$S=\langle 47\rangle\oplus \langle -1 \rangle 
\oplus 
\langle -1 \rangle$;
\nobreak
\newline
$\rho= (39,-141,-235) $,\ $\rho^2=-3619$; 
\newline 
$\Gamma_+=
\pmatrix{0}&{-1}&{1}\cr{0}&{1}&{0}\cr{1}&{0}&{-7}\cr{3}&{-5}&{-20}\cr
\endpmatrix,  
\hskip20pt 
G(\Gamma_+)=
\pmatrix{-2}&{1}&{7}&{15}\cr{1}&{-1}&{0}&{5}\cr{7}&{0}&{-2}&{1}\cr
{15}&{5}&{1}&{-2}\cr
\endpmatrix$;
\nobreak
\newline
$C_1=
\pmatrix{422}&{42}&{45}\cr{-1974}&{-197}&{-210}\cr{-2115}&{-210}&{-226}\cr
\endpmatrix,
\hskip30pt
C_2=
\pmatrix{46}&{3}&{6}\cr{-141}&{-10}&{-18}\cr{-282}&{-18}&{-37}\cr
\endpmatrix$.}

\vskip5pt

\vbox{\noindent
$N=11 \ $ $d=53 ,\,\eta=1 ,\,h=2 $: 
$S=\langle 1 \rangle\oplus \langle -106 \rangle \oplus 
\langle -2\rangle(0,1/2,1/2)$;
\nobreak
\newline
$\rho=(477,-41,-159)$,\ $\rho^2=-1219$; 
\newline 
$\Gamma_+=
\pmatrix{5}&{{-1}\over{2}}&{{-1}\over{2}}\cr{0}&{0}&{1}\cr{0}&{1}&{0}\cr
{1}&{0}&{-1}\cr{7}&{{-1}\over{2}}&{{-7}\over{2}}\cr
\endpmatrix,  
\hskip20pt 
G(\Gamma_+)=
\pmatrix{-2}&{1}&{53}&{4}&{5}\cr{1}&{-2}&{0}&{2}&{7}\cr
{53}&{0}&{-106}&{0}&{53}\cr{4}&{2}&{0}&{-1}&{0}\cr{5}&{7}&{53}&{0}&{-2}\cr
\endpmatrix$;
\nobreak
\newline
$C_1=
\pmatrix{647}&{5724}&{468}\cr{-54}&{-478}&{-39}\cr{-234}&{-2067}&{-170}\cr
\endpmatrix,
\hskip30pt
C_2=
\pmatrix{1155}&{10812}&{680}\cr{-102}&{-955}&{-60}\cr{-340}&{-3180}&{-201}\cr
\endpmatrix$.}

\vskip5pt

\vbox{\noindent
$N=12 \ $ $d=55=5\cdot 11 ,\,\eta=1 ,\,h=2 $: 
$S=\langle 11 \rangle\oplus \langle -10 \rangle \oplus 
\langle -2 \rangle(0,1/2,1/2)$;
\nobreak
\newline
$\rho= (25,-11,-55)$,\ $\rho^2=-385$; 
\newline 
$\Gamma_+=
\pmatrix{1}&{-1}&{-1}\cr{10}&{-11}&{0}\cr{0}&{0}&{1}\cr{0}&{1}&{0}\cr
{14}&{0}&{-33}\cr{9}&{-1}&{-21}\cr
\endpmatrix,  
\hskip20pt 
G(\Gamma_+)=
\pmatrix{-1}&{0}&{2}&{10}&{88}&{47}\cr{0}&{-110}&{0}&{110}&{1540}&{880}\cr
{2}&{0}&{-2}&{0}&{66}&{42}\cr{10}&{110}&{0}&{-10}&{0}&{10}\cr
{88}&{1540}&{66}&{0}&{-22}&{0}\cr{47}&{880}&{42}&{10}&{0}&{-1}\cr
\endpmatrix$;
\nobreak
\newline
$C_1=
\pmatrix{43}&{20}&{16}\cr{-22}&{-11}&{-8}\cr{-88}&{-40}&{-33}\cr
\endpmatrix,
\hskip30pt
C_2=
\pmatrix{87}&{20}&{36}\cr{-22}&{-6}&{-9}\cr{-198}&{-45}&{-82}\cr
\endpmatrix$.}

\vskip5pt

\vbox{\noindent
$N=13 \ $ $d=57=3\cdot 19,\,\eta=3 ,\,h=2 $: 
$S=\langle 1\rangle\oplus \langle -19 \rangle \oplus \langle -3 \rangle$;
\nobreak
\newline
$\rho= (95,-21,-19) $,\ $\rho^2=-437$; 
\newline 
$\Gamma_+=\pmatrix{13}&{-3}&{0}\cr{0}&{0}&{1}\cr{0}&{1}&{0}\cr{1}&{0}&{-1}\cr
\endpmatrix,  
\hskip20pt 
G(\Gamma_+)=\pmatrix 
{-2}&{0}&{57}&{13}\cr{0}&{-3}&{0}&{3}\cr{57}&{0}&{-19}&{0}\cr
{13}&{3}&{0}&{-2}\cr
\endpmatrix$;
\nobreak
\newline
$C_1=
\pmatrix{74}&{285}&{60}\cr{-15}&{-58}&{-12}\cr{-20}&{-76}&{-17}\cr
\endpmatrix,
\hskip30pt
C_2=
\pmatrix{1295}&{4788}&{1188}\cr{-252}&{-932}&{-231}\cr{-396}&{-1463}&{-364}\cr
\endpmatrix$.}

\vskip5pt 

\vbox{\noindent
$N=14 \ $ $d=58=2\cdot 29 ,\,\eta=1 ,\,h=2 $: 
$S=U \oplus \langle -58 \rangle$;
\nobreak
\newline
$\rho=(58,29,-8)$,\ $\rho^2=-348$; 
\newline 
$\Gamma_+=
\pmatrix{23}&{5}&{-2}\cr{14}&{2}&{-1}\cr{29}&{0}&{-1}\cr{0}&{0}&{1}\cr
{-1}&{1}&{0}\cr
\endpmatrix,  
\hskip20pt 
G(\Gamma_+)=
\pmatrix{-2}&{0}&{29}&{116}&{18}\cr{0}&{-2}&{0}&{58}&{12}\cr
{29}&{0}&{-58}&{58}&{29}\cr{116}&{58}&{58}&{-58}&{0}\cr
{18}&{12}&{29}&{0}&{-2}\cr
\endpmatrix$;
\nobreak
\newline
$C_1=
\pmatrix{29}&{100}&{580}\cr{9}&{29}&{174}\cr{-3}&{-10}&{-59}\cr
\endpmatrix,
\hskip30pt
C_2=
\pmatrix{29}&{36}&{348}\cr{25}&{29}&{290}\cr{-5}&{-6}&{-59}\cr
\endpmatrix$.}

\vskip5pt

\vbox{\noindent
$N=15 \ $ $d=59 ,\,\eta=0 ,\,h=2 $: $S=\langle 59 \rangle\oplus 
\langle -1 \rangle \oplus \langle -1 \rangle$;
\nobreak
\newline
$\rho=(24,-59,-177) $,\ $\rho^2=-826$; 
\newline 
$\Gamma_+=
\pmatrix{1}&{-5}&{-6}\cr{0}&{-1}&{1}\cr{0}&{1}&{0}\cr{23}&{0}&{-177}\cr
{3}&{-2}&{-23}\cr
\endpmatrix,  
\hskip20pt 
G(\Gamma_+)=
\pmatrix{-2}&{1}&{5}&{295}&{29}\cr
{1}&{-2}&{1}&{177}&{21}\cr{5}&{1}&{-1}&{0}&{2}\cr
{295}&{177}&{0}&{-118}&{0}\cr{29}&{21}&{2}&{0}&{-2}\cr
\endpmatrix$;
\nobreak
\newline
$C_1=
\pmatrix{117}&{6}&{14}\cr{-354}&{-19}&{-42}\cr{-826}&{-42}&{-99}\cr
\endpmatrix,
\hskip30pt
C_2=
\pmatrix{235}&{6}&{30}\cr{-354}&{-10}&{-45}\cr{-1770}&{-45}&{-226}\cr
\endpmatrix$.}

\vskip5pt 

\vbox{\noindent
$N=16 \ $ $d=62=2\cdot 31 ,\,\eta=1 ,\,h=2 $: 
$S=U \oplus \langle -62 \rangle$;
\nobreak
\newline
$\rho=(93,31,-10) $,\ $\rho^2=-434$; 
\newline 
$\Gamma_+=\pmatrix 
{15}&{2}&{-1}\cr{31}&{0}&{-1}\cr{0}&{0}&{1}\cr
{-1}&{1}&{0}\cr{6}&{5}&{-1}\cr
\endpmatrix,  
\hskip20pt 
G(\Gamma_+)=\pmatrix 
{-2}&{0}&{62}&{13}&{25}\cr{0}&{-62}&{62}&{31}&{93}\cr
{62}&{62}&{-62}&{0}&{62}\cr
{13}&{31}&{0}&{-2}&{1}\cr{25}&{93}&{62}&{1}&{-2}\cr
\endpmatrix$;
\nobreak
\newline
$C_1=
\pmatrix{124}&{625}&{3100}\cr{25}&{124}&{620}\cr{-10}&{-50}&{-249}\cr
\endpmatrix,
\hskip30pt
C_2=
\pmatrix{31}&{64}&{496}\cr{16}&{31}&{248}\cr{-4}&{-8}&{-63}\cr
\endpmatrix$.}

\vskip5pt 

\vbox{\noindent
$N=17 \ $ $d=65=5\cdot 13 ,\,\eta=0 ,\,h=0 $: 
$S=\langle 1 \rangle\oplus \langle -13 \rangle \oplus 
\langle -5 \rangle$;
\nobreak
\newline
$\rho=(15,-4,-4)$,\ $\rho^2=-63$; 
\newline 
$\Gamma_+=
\pmatrix{18}&{-5}&{0}\cr{0}&{0}&{1}\cr{0}&{1}&{0}\cr{2}&{0}&{-1}\cr
\endpmatrix,  
\hskip20pt 
G(\Gamma_+)=
\pmatrix{-1}&{0}&{65}&{36}\cr{0}&{-5}&{0}&{5}\cr{65}&{0}&{-13}&{0}\cr
{36}&{5}&{0}&{-1}\cr
\endpmatrix$;
\nobreak
\newline 
$\Gamma_-=
\pmatrix{19}&{-3}&{-7}\cr{260}&{-45}&{-91}\cr
{4}&{-1}&{-1}\cr{11}&{-3}&{-1}\cr
\endpmatrix,  
\hskip20pt 
G(\Gamma_-)=
\pmatrix{-1}&{0}&{2}&{57}\cr{0}&{-130}&{0}&{650}\cr{2}&{0}&{-2}&{0}\cr
{57}&{650}&{0}&{-1}\cr
\endpmatrix$;
\nobreak
\newline
$T=
\pmatrix{289}&{1040}&{40}\cr{-40}&{-144}&{-5}\cr{-112}&{-403}&{-16}\cr
\endpmatrix$.}

\vskip5pt

\vbox{\noindent
$N=18 \ $ $d=65=5\cdot 13 ,\,\eta=1 ,\,h=2 $: 
$S=\langle 2 \rangle\oplus \langle -13 \rangle \oplus 
\langle -10 \rangle(1/2,0,1/2)$;
\nobreak
\newline
$\rho=(91/2,-17,-13/2)$,\ $\rho^2=-39$; 
\newline 
$\Gamma_+=
\pmatrix{5}&{-2}&{0}\cr{0}&{0}&{1}\cr{0}&{1}&{0}\cr
{{1}\over{2}}&{0}&{{-1}\over{2}}\cr
\endpmatrix,  
\hskip20pt 
G(\Gamma_+)=
\pmatrix{-2}&{0}&{26}&{5}\cr{0}&{-10}&{0}&{5}\cr{26}&{0}&{-13}&{0}\cr
{5}&{5}&{0}&{-2}\cr
\endpmatrix$;
\nobreak
\newline
$C_1=
\pmatrix{255}&{624}&{160}\cr{-96}&{-235}&{-60}\cr{-32}&{-78}&{-21}\cr
\endpmatrix,
\hskip30pt
C_2=
\pmatrix{{603}\over{2}}&{715}&{{495}\over{2}}\cr
{-110}&{-261}&{-90}\cr{{-99}\over{2}}&{-117}&{{-83}\over{2}}\cr
\endpmatrix$.}

\vskip5pt

\vbox{\noindent
$N=19 \ $ $d=69=3\cdot 23 ,\,\eta=0 ,\,h=2 $: 
$S=\langle 46 \rangle\oplus \langle -6 \rangle \oplus 
\langle -1\rangle(1/2,1/2,0)$;
\nobreak
\newline
$\rho=(21/2, -23/2, -69) $,\ $\rho^2=-483$; 
\newline 
$\Gamma_+=
\pmatrix{{1}\over{2}}&{{-3}\over{2}}&{0}\cr{0}&{0}&{1}\cr{0}&{1}&{0}\cr
{{23}\over{2}}&{{-1}\over{2}}&{-78}\cr
\endpmatrix,  
\hskip20pt 
G(\Gamma_+)=
\pmatrix{-2}&{0}&{9}&{260}\cr{0}&{-1}&{0}&{78}\cr{9}&{0}&{-6}&{3}\cr
{260}&{78}&{3}&{-2}\cr
\endpmatrix$;
\nobreak
\newline
$C_1=
\pmatrix{22}&{3}&{3}\cr{-23}&{-4}&{-3}\cr{-138}&{-18}&{-19}\cr
\endpmatrix,
\hskip30pt
C_2=
\pmatrix{{205}\over{2}}&{{9}\over{2}}&{15}\cr
{{-69}\over{2}}&{{-5}\over{2}}&{-5}\cr{-690}&{-30}&{-101}\cr
\endpmatrix$.}

\vskip5pt 

\vbox{\noindent
$N=20 \ $ $d=70=2\cdot 5\cdot 7 ,\,\eta=0 ,\,h=2 $: 
$S=U \oplus \langle -70 \rangle$;
\nobreak
\newline
$\rho=(30,10,-3)$,\ $\rho^2=-30$; 
\newline 
$\Gamma_+=
\pmatrix{35}&{0}&{-1}\cr{0}&{0}&{1}\cr{-1}&{1}&{0}\cr{35}&{35}&{-6}\cr
\endpmatrix,  
\hskip20pt 
G(\Gamma_+)=\pmatrix 
{-70}&{70}&{35}&{805}\cr{70}&{-70}&{0}&{420}\cr{35}&{0}&{-2}&{0}\cr
{805}&{420}&{0}&{-70}\cr
\endpmatrix$;
\nobreak
\newline
$C_1=
\pmatrix{20}&{63}&{420}\cr{7}&{20}&{140}\cr{-2}&{-6}&{-41}\cr
\endpmatrix,
\hskip30pt
C_2=
\pmatrix{35}&{81}&{630}\cr{16}&{35}&{280}\cr{-4}&{-9}&{-71}\cr
\endpmatrix$.}

\vskip5pt 

\vbox{\noindent
$N=21 \ $ $d=74=2\cdot 37 ,\,\eta=1 ,\,h=0 $: 
$S= U \oplus \langle -74 \rangle$;
\nobreak
\newline
$\rho=(30,10,-3) $,\ $\rho^2=-66 $; 
\newline 
$\Gamma_+=
\pmatrix{37}&{0}&{-1}\cr{0}&{0}&{1}\cr{-1}&{1}&{0}\cr{6}&{6}&{-1}\cr
{259}&{185}&{-36}\cr
\endpmatrix,  
\hskip20pt 
G(\Gamma_+)=\pmatrix 
{-74}&{74}&{37}&{148}&{4181}\cr{74}&{-74}&{0}&{74}&{2664}\cr
{37}&{0}&{-2}&{0}&{74}\cr{148}&{74}&{0}&{-2}&{0}\cr
{4181}&{2664}&{74}&{0}&{-74}\cr
\endpmatrix$;
\nobreak
\newline 
$\Gamma_-=
\pmatrix{12}&{3}&{-1}\cr{9}&{4}&{-1}\cr
\endpmatrix,  
\hskip20pt 
G(\Gamma_-)=
\pmatrix{-2}&{1}\cr{1}&{-2}\cr
\endpmatrix$;
\nobreak
\newline 
$T=\pmatrix{37}&{225}&{1110}\cr{25}&{148}&{740}\cr{-5}&{-30}&{-149}\cr
\endpmatrix$.}

\vskip5pt

\vbox{\noindent
$N=22 \ $ $d=77=7\cdot 11 ,\,\eta=0 ,\,h=2 $: $S=\langle 14\rangle\oplus 
\langle -22 \rangle \oplus \langle -1 \rangle(1/2,1/2,0)$;
\nobreak
\newline
$\rho= (99/2, -49/2, -154) $,\ $\rho^2=-2618$; 
\newline 
$\Gamma_+=
\pmatrix{{1}\over{2}}&{{-1}\over{2}}&{0}\cr{0}&{0}&{1}\cr{0}&{1}&{0}\cr
{1}&{0}&{-4}\cr
\endpmatrix,  
\hskip20pt 
G(\Gamma_+)=
\pmatrix{-2}&{0}&{11}&{7}\cr{0}&{-1}&{0}&{4}\cr{11}&{0}&{-22}&{0}\cr
{7}&{4}&{0}&{-2}\cr
\endpmatrix$;
\nobreak
\newline
$C_1=
\pmatrix{{173}\over{2}}&{{165}\over{2}}&{15}\cr
{{-105}\over{2}}&{{-101}\over{2}}&{-9}\cr{-210}&{-198}&{-37}\cr
\endpmatrix,
\hskip30pt
C_2=
\pmatrix{62}&{33}&{15}\cr{-21}&{-12}&{-5}\cr{-210}&{-110}&{-51}\cr
\endpmatrix$.}

\vskip5pt

\vbox{\noindent
$N=23 \ $ $d=85=5\cdot 17 ,\,\eta=2 ,\,h=2 $: 
$S=\langle 1\rangle\oplus \langle -170 \rangle \oplus  
\langle -2\rangle(0,1/2,1/2)$;
\nobreak
\newline
$\rho=(34,-5/2,-17/2)$,\ $\rho^2=-51 $; 
\newline 
$\Gamma_+=
\pmatrix{13}&{-1}&{0}\cr{0}&{0}&{1}\cr{0}&{1}&{0}\cr{1}&{0}&{-1}\cr
\endpmatrix,  
\hskip20pt 
G(\Gamma_+)=
\pmatrix{-1}&{0}&{170}&{13}\cr{0}&{-2}&{0}&{2}\cr{170}&{0}&{-170}&{0}\cr
{13}&{2}&{0}&{-1}\cr
\endpmatrix$;
\nobreak
\newline
$C_1=
\pmatrix{48}&{595}&{21}\cr{{-7}\over{2}}&{{-87}\over{2}}&{{-3}\over{2}}\cr
{{-21}\over{2}}&{{-255}\over{2}}&{{-11}\over{2}}\cr
\endpmatrix,
\hskip30pt
C_2=
\pmatrix{89}&{1020}&{60}\cr{-6}&{-69}&{-4}\cr{-30}&{-340}&{-21}\cr
\endpmatrix$.}

\vskip5pt

\vbox{\noindent
$N=24 \ $ $d=87=3\cdot 29 ,\,\eta=0 ,\,h=2 $: 
$S=\langle 3\rangle\oplus \langle -29 \rangle 
\oplus  \langle -1 \rangle$;
\nobreak
\newline
$\rho=(116,-34,-87)$,\ $\rho^2=-725$; 
\newline 
$\Gamma_+=
\pmatrix{68}&{-21}&{-33}\cr{3}&{-1}&{0}\cr{0}&{0}&{1}\cr{0}&{1}&{0}\cr
{1}&{0}&{-3}\cr
\endpmatrix,  
\hskip20pt 
G(\Gamma_+)=
\pmatrix{-6}&{3}&{33}&{609}&{105}\cr{3}&{-2}&{0}&{29}&{9}\cr
{33}&{0}&{-1}&{0}&{3}\cr{609}&{29}&{0}&{-29}&{0}\cr{105}&{9}&{3}&{0}&{-6}\cr
\endpmatrix$;
\nobreak
\newline
$C_1=
\pmatrix{199}&{580}&{40}\cr{-60}&{-175}&{-12}\cr{-120}&{-348}&{-25}\cr
\endpmatrix,
\hskip30pt
C_2=
\pmatrix{362}&{957}&{110}\cr{-99}&{-262}&{-30}\cr{-330}&{-870}&{-101}\cr
\endpmatrix$.}

\vskip5pt 

\vbox{\noindent
$N=25 \ $ $d=91=7\cdot 13 ,\,\eta=0 ,\,h=2 $: $S=\langle 7 \rangle\oplus 
\langle -13 \rangle \oplus \langle -1 \rangle$;
\nobreak
\newline
$\rho= (13, -9, -13)$,\ $\rho^2=-39$; 
\newline 
$\Gamma_+=
\pmatrix{0}&{0}&{1}\cr{0}&{1}&{0}\cr{1}&{0}&{-3}\cr{39}&{-14}&{-91}\cr
{2}&{-1}&{-4}\cr
\endpmatrix,  
\hskip20pt 
G(\Gamma_+)=
\pmatrix{-1}&{0}&{3}&{91}&{4}\cr{0}&{-13}&{0}&{182}&{13}\cr
{3}&{0}&{-2}&{0}&{2}\cr
{91}&{182}&{0}&{-182}&{0}\cr{4}&{13}&{2}&{0}&{-1}\cr
\endpmatrix$;
\nobreak
\newline
$C_1=
\pmatrix{62}&{78}&{9}\cr{-42}&{-53}&{-6}\cr{-63}&{-78}&{-10}\cr
\endpmatrix,
\hskip30pt
C_2=
\pmatrix{120}&{143}&{22}\cr{-77}&{-92}&{-14}\cr{-154}&{-182}&{-29}\cr
\endpmatrix$.}

\vskip5pt

\vbox{\noindent
$N=26 \ $ $d=93=3\cdot 31 ,\,\eta=3 ,\,h=2 $: 
$S=\langle 1\rangle\oplus \langle -31 \rangle 
\oplus \langle -3 \rangle$;
\nobreak
\newline
$\rho= (155,-27,-31)$,\ $\rho^2=-1457 $; 
\newline 
$\Gamma_+=
\pmatrix{67}&{-12}&{-3}\cr{39}&{-7}&{-1}\cr{217}&{-39}&{0}\cr{0}&{0}&{1}\cr
{0}&{1}&{0}\cr{1}&{0}&{-1}\cr
\endpmatrix,  
\hskip20pt 
G(\Gamma_+)=
\pmatrix{-2}&{0}&{31}&{9}&{372}&{58}\cr{0}&{-1}&{0}&{3}&{217}&{36}\cr
{31}&{0}&{-62}&{0}&{1209}&{217}\cr{9}&{3}&{0}&{-3}&{0}&{3}\cr
{372}&{217}&{1209}&{0}&{-31}&{0}\cr{58}&{36}&{217}&{3}&{0}&{-2}\cr
\endpmatrix$;
\nobreak
\newline
$C_1=
\pmatrix{35}&{186}&{18}\cr{-6}&{-32}&{-3}\cr{-6}&{-31}&{-4}\cr
\endpmatrix,
\hskip30pt
C_2=\pmatrix{1199}&{5580}&{1140}\cr{-180}&{-838}&{-171}\cr
{-380}&{-1767}&{-362}\cr
\endpmatrix$.}

\vskip5pt

\vbox{\noindent
$N=27 \ $ $d=95=5\cdot 19 ,\,\eta=0 ,\,h=0 $: 
$S=\langle 19 \rangle\oplus \langle -5 \rangle \oplus 
\langle -1 \rangle$;
\nobreak
\newline
$\rho= (5,-6,-20)$,\ $\rho^2=-105$; 
\newline 
$\Gamma_+=
\pmatrix{1}&{-2}&{-1}\cr{0}&{0}&{1}\cr{0}&{1}&{0}\cr{13}&{0}&{-57}\cr
{9}&{-2}&{-39}\cr
\endpmatrix,  
\hskip20pt 
G(\Gamma_+)=
\pmatrix{-2}&{1}&{10}&{190}&{112}\cr{1}&{-1}&{0}&{57}&{39}\cr
{10}&{0}&{-5}&{0}&{10}\cr{190}&{57}&{0}&{-38}&{0}\cr{112}&{39}&{10}&{0}&{-2}\cr
\endpmatrix$;
\nobreak
\newline
$S_1=\pmatrix{115}&{10}&{26}\cr{-190}&{-16}&{-43}\cr{-266}&{-25}&{-60}\cr
\endpmatrix$.}

\vskip 5pt

\vbox{\noindent
$N=28 \ $ $d=114=2\cdot 3\cdot 19 ,\,\eta=1 ,\,h=2 $: 
$S=\langle 6 \rangle\oplus 
\langle -38\rangle \oplus \langle -2 \rangle$(0,1/2,1/2);
\nobreak
\newline
$\rho=(19, -7, -19)$,\ $\rho^2=-418$; 
\newline 
$\Gamma_+=
\pmatrix{190}&{-75}&{-38}\cr{5}&{-2}&{0}\cr{0}&{0}&{1}\cr{0}&{1}&{0}\cr
\endpmatrix,  
\hskip20pt 
G(\Gamma_+)=
\pmatrix{-38}&{0}&{76}&{2850}\cr{0}&{-2}&{0}&{76}\cr{76}&{0}&{-2}&{0}\cr
{2850}&{76}&{0}&{-38}\cr
\endpmatrix$;
\nobreak
\newline
$C_1=
\pmatrix{31}&{76}&{4}\cr{-12}&{{-59}\over{2}}&{{-3}\over{2}}\cr
{-12}&{{-57}\over{2}}&{{-5}\over{2}}\cr
\endpmatrix,
\hskip30pt
C_2=
\pmatrix{23}&{38}&{10}\cr{-6}&{{-21}\over{2}}&{{-5}\over{2}}\cr
{-30}&{{-95}\over{2}}&{{-27}\over{2}}\cr
\endpmatrix$.}

\vskip5pt

\vbox{\noindent
$N=29 \ $ $d=114=2\cdot 3\cdot 19 ,\,\eta=2 ,\,h=0 $: 
$S=U \oplus \langle -114 \rangle$;
\nobreak
\newline
$\rho=(95,19,-6)$,\ $\rho^2=-494$; 
\newline 
$\Gamma_+=
\pmatrix{321}&{30}&{-13}\cr{28}&{2}&{-1}\cr{57}&{0}&{-1}\cr
{0}&{0}&{1}\cr{-1}&{1}&{0}\cr{9}&{6}&{-1}\cr
\endpmatrix,  
\hskip10pt 
G(\Gamma_+)=
\pmatrix{-6}&{0}&{228}&{1482}&{291}&{714}\cr{0}&{-2}&{0}&{114}&{26}&{72}\cr
{228}&{0}&{-114}&{114}&{57}&{228}\cr{1482}&{114}&{114}&{-114}&{0}&{114}\cr
{291}&{26}&{57}&{0}&{-2}&{3}\cr{714}&{72}&{228}&{114}&{3}&{-6}\cr
\endpmatrix$;
\nobreak
\newline
$S_1=
\pmatrix{48}&{475}&{2280}\cr{19}&{192}&{912}\cr{-4}&{-40}&{-191}\cr
\endpmatrix$.}

\vskip5pt 

\vbox{\noindent
$N=30 \ $ $d=115=5\cdot 23 ,\,\eta=0 ,\,h=2 $: 
$S=\langle 23 \rangle\oplus 
\langle -5 \rangle \oplus \langle -1 \rangle$;
\nobreak
\newline
$\rho= (35, -69, -115)$,\ $\rho^2=-8855 $; 
\newline 
$\Gamma_+=
\pmatrix{0}&{0}&{1}\cr{0}&{1}&{0}\cr{1}&{0}&{-5}\cr{75}&{-46}&{-345}\cr
{4}&{-3}&{-18}\cr
\endpmatrix,  
\hskip20pt 
G(\Gamma_+)=
\pmatrix{-1}&{0}&{5}&{345}&{18}\cr{0}&{-5}&{0}&{230}&{15}\cr
{5}&{0}&{-2}&{0}&{2}\cr{345}&{230}&{0}&{-230}&{0}\cr{18}&{15}&{2}&{0}&{-1}\cr
\endpmatrix$;
\nobreak
\newline
$C_1=
\pmatrix{22}&{10}&{1}\cr{-46}&{-21}&{-2}\cr{-23}&{-10}&{-2}\cr
\endpmatrix,
\hskip30pt
C_2=
\pmatrix{22}&{5}&{4}\cr{-23}&{-6}&{-4}\cr{-92}&{-20}&{-17}\cr
\endpmatrix$.}

\vskip5pt

\vbox{\noindent
$N=31 \ $ $d=115=5\cdot 23 ,\,\eta=3 ,\,h=2 $: 
$S=\langle 115 \rangle\oplus 
\langle -1 \rangle \oplus \langle -1 \rangle$;
\nobreak
\newline
$\rho= (97, -115, -1035)$,\ $\rho^2=-2415 $; 
\newline 
$\Gamma_+=
\pmatrix{61}&{-115}&{-644}\cr{1}&{-5}&{-10}\cr{3}&{-23}&{-23}\cr{0}&{-1}&{1}\cr
{0}&{1}&{0}\cr{15}&{0}&{-161}\cr
\endpmatrix,  
\hskip5pt 
G(\Gamma_+)=
\pmatrix{-46}&{0}&{3588}&{529}&{115}&{1541}\cr{0}&{-10}&{0}&{5}&{5}&{115}\cr
{3588}&{0}&{-23}&{0}&{23}&{1472}\cr{529}&{5}&{0}&{-2}&{1}&{161}\cr
{115}&{5}&{23}&{1}&{-1}&{0}\cr{1541}&{115}&{1472}&{161}&{0}&{-46}\cr
\endpmatrix$;
\nobreak
\newline
$C_1=
\pmatrix{2874}&{40}&{265}\cr{-4600}&{-65}&{-424}\cr{-30475}&{-424}&{-2810}\cr
\endpmatrix,
\hskip30pt
C_2=
\pmatrix{1034}&{9}&{96}\cr{-1035}&{-10}&{-96}\cr{-11040}&{-96}&{-1025}\cr
\endpmatrix$.}

\vskip5pt 

\vbox{\noindent
$N=32 \ $ $d=119=7\cdot 17 ,\,\eta=0 ,\,h=2 $: 
$S=\langle 7 \rangle\oplus \langle -17 \rangle \oplus 
\langle -1 \rangle$;
\nobreak
\newline
$\rho=(13,-7,-21)$,\ $\rho^2=-91$; 
\newline 
$\Gamma_+=
\pmatrix{8}&{-5}&{-5}\cr{11}&{-7}&{-4}\cr{17}&{-11}&{0}\cr{0}&{0}&{1}\cr
{0}&{1}&{0}\cr{1}&{0}&{-3}\cr
\endpmatrix,  
\hskip20pt 
G(\Gamma_+)=
\pmatrix{-2}&{1}&{17}&{5}&{85}&{41}\cr{1}&{-2}&{0}&{4}&{119}&{65}\cr
{17}&{0}&{-34}&{0}&{187}&{119}\cr{5}&{4}&{0}&{-1}&{0}&{3}\cr
{85}&{119}&{187}&{0}&{-17}&{0}\cr{41}&{65}&{119}&{3}&{0}&{-2}\cr
\endpmatrix$;
\nobreak
\newline
$C_1=
\pmatrix{118}&{170}&{17}\cr{-70}&{-101}&{-10}\cr{-119}&{-170}&{-18}\cr
\endpmatrix,
\hskip30pt
C_2=
\pmatrix{27}&{34}&{6}\cr{-14}&{-18}&{-3}\cr{-42}&{-51}&{-10}\cr
\endpmatrix$.}

\vskip5pt 

\vbox{\noindent
$N=33 \ $ $d=123=3\cdot 41 ,\,\eta=0 ,\,h=2 $: $S=\langle 3 \rangle\oplus 
\langle -41 \rangle \oplus  \langle -1 \rangle$;
\nobreak
\newline
$\rho= (369,-98,-123)$,\ $\rho^2=-410$; 
\newline 
$\Gamma_+\hskip-2pt=\hskip-2pt 
\pmatrix{71}&{-19}&{-18}\cr{369}&{-99}&{-82}\cr{11}&{-3}&{0}\cr{0}&{0}&{1}\cr
{0}&{1}&{0}\cr{1}&{0}&{-3}\cr{4}&{-1}&{-3}\cr
\endpmatrix,  
\hskip2pt 
G(\Gamma_+)\hskip-2pt=\hskip-2pt 
\pmatrix{-2}&{0}&{6}&{18}&{779}&{159}&{19}\cr
{0}&{-82}&{0}&{82}&{4059}&{861}&{123}\cr
{6}&{0}&{-6}&{0}&{123}&{33}&{9}\cr{18}&{82}&{0}&{-1}&{0}&{3}&{3}\cr
{779}&{4059}&{123}&{0}&{-41}&{0}&{41}\cr{159}&{861}&{33}&{3}&{0}&{-6}&{3}\cr
{19}&{123}&{9}&{3}&{41}&{3}&{-2}\cr
\endpmatrix$;
\nobreak
\newline
$C_1=
\pmatrix{6074}&{22140}&{585}\cr
{-1620}&{-5905}&{-156}\cr{-1755}&{-6396}&{-170}\cr
\endpmatrix,
\hskip30pt
C_2=
\pmatrix{2311}&{8364}&{272}\cr{-612}&{-2215}&{-72}\cr{-816}&{-2952}&{-97}\cr
\endpmatrix$.}

\vskip5pt 

\vbox{\noindent
$N=34 \ $ $d=123=3\cdot 41 ,\,\eta=3 ,\,h=2 $: $S=\langle 123 \rangle\oplus 
\langle -1 \rangle \oplus \langle -1 \rangle$;
\nobreak
\newline
$\rho=(21,-123,-205)$,\ $\rho^2=-2911$; 
\newline 
$\Gamma_+=
\pmatrix{298}&{-2337}&{-2337}\cr{0}&{-1}&{1}\cr{0}&{1}&{0}\cr{11}&{0}&{-123}\cr
\endpmatrix,  
\hskip20pt 
G(\Gamma_+)=
\pmatrix{-246}&{0}&{2337}&{115743}\cr
{0}&{-2}&{1}&{123}\cr{2337}&{1}&{-1}&{0}\cr
{115743}&{123}&{0}&{-246}\cr
\endpmatrix$;
\nobreak
\newline
$C_1=
\pmatrix{40}&{2}&{3}\cr{-246}&{-13}&{-18}\cr{-369}&{-18}&{-28}\cr
\endpmatrix,
\hskip30pt
C_2=
\pmatrix{491}&{14}&{42}\cr{-1722}&{-50}&{-147}\cr{-5166}&{-147}&{-442}\cr
\endpmatrix$.}

\vskip5pt

\vbox{\noindent
$N=35 \ $ $d=129=3\cdot 43 ,\,\eta=0 ,\,h=2 $: $S=\langle 3 \rangle\oplus 
\langle -86 \rangle \oplus \langle -2 \rangle(0,1/2,1/2)$;
\nobreak
\newline
$\rho= (86,-27/2,-129/2)$,\ $\rho^2=-1806$; 
\newline 
$\Gamma_+=
\pmatrix{8}&{{-3}\over{2}}&{{-3}\over{2}}\cr{0}&{0}&{1}\cr{0}&{1}&{0}\cr
{2}&{0}&{-3}\cr
\endpmatrix,  
\hskip20pt 
G(\Gamma_+)=
\pmatrix{-6}&{3}&{129}&{39}\cr{3}&{-2}&{0}&{6}\cr{129}&{0}&{-86}&{0}\cr
{39}&{6}&{0}&{-6}\cr
\endpmatrix$;
\nobreak
\newline
$C_1=
\pmatrix{53}&{258}&{18}\cr{-9}&{-44}&{-3}\cr{-27}&{-129}&{-10}\cr
\endpmatrix,
\hskip30pt
C_2=
\pmatrix{47}&{172}&{28}\cr{-6}&{{-45}\over{2}}&{{-7}\over{2}}\cr
{-42}&{{-301}\over{2}}&{{-51}\over{2}}\cr
\endpmatrix$.}

\vskip5pt 

\vbox{\noindent
$N=36 \ $ $d=133=7\cdot 19 ,\,\eta=3 ,\,h=2 $: $S=\langle 1 \rangle\oplus 
\langle -19 \rangle \oplus \langle -7 \rangle$;
\nobreak
\newline
$\rho=(1197,-259,-247)$,\ $\rho^2= -268793$; 
\newline 
$\Gamma_+=
\pmatrix{13}&{-3}&{0}\cr{0}&{0}&{1}\cr{0}&{1}&{0}\cr{7}&{0}&{-3}\cr
{9}&{-1}&{-3}\cr{68}&{-9}&{-21}\cr
\endpmatrix,  
\hskip20pt 
G(\Gamma_+)=
\pmatrix{-2}&{0}&{57}&{91}&{60}&{371}\cr{0}&{-7}&{0}&{21}&{21}&{147}\cr
{57}&{0}&{-19}&{0}&{19}&{171}\cr{91}&{21}&{0}&{-14}&{0}&{35}\cr
{60}&{21}&{19}&{0}&{-1}&{0}\cr{371}&{147}&{171}&{35}&{0}&{-2}\cr
\endpmatrix$;
\nobreak
\newline
$C_1=
\pmatrix{483}&{2090}&{154}\cr{-110}&{-476}&{-35}\cr{-22}&{-95}&{-8}\cr
\endpmatrix,
\hskip30pt
C_2=
\pmatrix{48}&{133}&{98}\cr{-7}&{-20}&{-14}\cr{-14}&{-38}&{-29}\cr
\endpmatrix$.}

\vskip5pt 

\vbox{\noindent
$N=37 \ $ $d=138=2\cdot 3\cdot 23 ,\,\eta=0 ,\,h=2 $: 
$S=\langle 46\rangle\oplus 
\langle -6 \rangle \oplus \langle -2 \rangle(0,1/2,1/2)$;
\nobreak
\newline
$\rho=(117,-92,-552)$,\ $\rho^2=-30498 $; 
\newline 
$\Gamma_+=
\pmatrix{54}&{{-145}\over{2}}&{{-453}\over{2}}\cr{13}&{-18}&{-54}\cr
{0}&{{-1}\over{2}}&{{1}\over{2}}\cr{0}&{1}&{0}\cr
\endpmatrix,  
\hskip20pt 
G(\Gamma_+)=
\pmatrix{-6}&{0}&{9}&{435}\cr{0}&{-2}&{0}&{108}\cr{9}&{0}&{-2}&{3}\cr
{435}&{108}&{3}&{-6}\cr
\endpmatrix$;
\nobreak
\newline
$C_1=
\pmatrix{183}&{30}&{34}\cr{-230}&{{-77}\over{2}}&{{-85}\over{2}}\cr
{-782}&{{-255}\over{2}}&{{-291}\over{2}}\cr
\endpmatrix,
\hskip30pt
C_2=
\pmatrix{183}&{6}&{38}\cr{-46}&{{-5}\over{2}}&{{-19}\over{2}}\cr
{-874}&{{-57}\over{2}}&{{-363}\over{2}}\cr
\endpmatrix$.}

\vskip5pt 

\vbox{\noindent
$N=38 \ $ $d=143=11\cdot 13 ,\,\eta=3 ,\,h=2 $: 
$S=\langle 143 \rangle\oplus 
\langle -1 \rangle \oplus \langle -1 \rangle$;
\nobreak
\newline
$\rho=(19,-117,-195)$,\ $\rho^2=-91$; 
\newline 
$\Gamma_+=
\pmatrix{1}&{-8}&{-9}\cr{0}&{-1}&{1}\cr{0}&{1}&{0}\cr{1}&{0}&{-13}\cr
{7}&{-33}&{-77}\cr{4}&{-21}&{-43}\cr
\endpmatrix,  
\hskip20pt 
G(\Gamma_+)=
\pmatrix{-2}&{1}&{8}&{26}&{44}&{17}\cr{1}&{-2}&{1}&{13}&{44}&{22}\cr
{8}&{1}&{-1}&{0}&{33}&{21}\cr{26}&{13}&{0}&{-26}&{0}&{13}\cr
{44}&{44}&{33}&{0}&{-11}&{0}\cr{17}&{22}&{21}&{13}&{0}&{-2}\cr
\endpmatrix$;
\nobreak
\newline
$C_1=
\pmatrix{1871}&{84}&{132}\cr{-12012}&{-540}&{-847}\cr
{-18876}&{-847}&{-1332}\cr
\endpmatrix,
\hskip30pt
C_2=
\pmatrix{1286}&{54}&{93}\cr{-7722}&{-325}&{-558}\cr{-13299}&{-558}&{-962}\cr
\endpmatrix$.}

\vskip5pt

\vbox{\noindent
$N=39 \ $ $d=159=3\cdot 53 ,\,\eta=3 ,\,h=2 $: $S=\langle 159 \rangle\oplus 
\langle -1 \rangle \oplus \langle -1 \rangle$;
\nobreak
\newline
$\rho=(13,-53,-159) $,\ $\rho^2=-1219$; 
\newline 
$\Gamma_+=
\pmatrix{1}&{-9}&{-9}\cr{0}&{-1}&{1}\cr{0}&{1}&{0}\cr{21}&{0}&{-265}\cr
{5}&{-3}&{-63}\cr
\endpmatrix,  
\hskip20pt 
G(\Gamma_+)=
\pmatrix{-3}&{0}&{9}&{954}&{201}\cr{0}&{-2}&{1}&{265}&{60}\cr
{9}&{1}&{-1}&{0}&{3}\cr{954}&{265}&{0}&{-106}&{0}\cr
{201}&{60}&{3}&{0}&{-3}\cr
\endpmatrix$;
\nobreak
\newline
$C_1=
\pmatrix{52}&{1}&{4}\cr{-159}&{-4}&{-12}\cr{-636}&{-12}&{-49}\cr
\endpmatrix,
\hskip30pt
C_2=
\pmatrix{635}&{6}&{50}\cr{-954}&{-10}&{-75}\cr{-7950}&{-75}&{-626}\cr
\endpmatrix$.}

\vskip5pt 

\vbox{\noindent
$N=40 \ $ $d=165=3\cdot 5\cdot 11 ,\,\eta=2 ,\,h=2 $: 
$S=\langle 165 \rangle\oplus 
\langle -1 \rangle \oplus \langle -1 \rangle$;
\nobreak
\newline
$\rho=(7,-22,-88)$,\ $\rho^2=-143$; 
\newline 
$\Gamma_+ = 
\left(\smallmatrix{6}&{-30}&{-71}\cr
{23}&{-120}&{-270}\cr{18}&{-99}&{-209}\cr
{4}&{-25}&{-45}\cr{2}&{-15}&{-21}\cr
{6}&{-55}&{-55}\cr{0}&{-1}&{1}\cr{0}&{1}&{0}\cr
\endsmallmatrix\right),  
\hskip20pt 
G(\Gamma_+) = 
\left(\smallmatrix{-1}&{0}&{11}&{15}&{39}&{385}&{41}&{30}\cr
{0}&{-15}&{0}&{30}&{120}&{1320}&{150}&{120}\cr
{11}&{0}&{-22}&{0}&{66}&{880}&{110}&{99}\cr
{15}&{30}&{0}&{-10}&{0}&{110}&{20}&{25}\cr
{39}&{120}&{66}&{0}&{-6}&{0}&{6}&{15}\cr
{385}&{1320}&{880}&{110}&{0}&{-110}&{0}&{55}\cr
{41}&{150}&{110}&{20}&{6}&{0}&{-2}&{1}\cr
{30}&{120}&{99}&{25}&{15}&{55}&{1}&{-1}\cr
\endsmallmatrix\right)$;
\nobreak
\newline
$C_1=
\pmatrix{439}&{12}&{32}\cr{-1980}&{-55}&{-144}\cr{-5280}&{-144}&{-385}\cr
\endpmatrix,
\hskip30pt
C_2=
\pmatrix{131}&{2}&{10}\cr{-330}&{-6}&{-25}\cr{-1650}&{-25}&{-126}\cr
\endpmatrix$.}

\vskip5pt

\vbox{\noindent
$N=41 \ $ $d=174=2\cdot 3\cdot 29 ,\,\eta=2 ,\,h=2 $: 
$S=\langle 2\rangle\oplus \langle -58 \rangle \oplus 
\langle -6 \rangle(0,1/2,1/2)$;
\nobreak
\newline
$\rho= (58,-21/2,-29/2)$,\ $\rho^2=-928$; 
\newline 
$\Gamma_+=
\pmatrix{35}&{{-13}\over{2}}&{{-1}\over{2}}\cr{0}&{0}&{1}\cr{0}&{1}&{0}\cr
{3}&{0}&{-2}\cr{5}&{{-1}\over{2}}&{{-5}\over{2}}\cr
\endpmatrix,  
\hskip20pt 
G(\Gamma_+)=
\pmatrix{-2}&{3}&{377}&{204}&{154}\cr{3}&{-6}&{0}&{12}&{15}\cr
{377}&{0}&{-58}&{0}&{29}\cr{204}&{12}&{0}&{-6}&{0}\cr
{154}&{15}&{29}&{0}&{-2}\cr
\endpmatrix$;
\nobreak
\newline
$C_1=
\pmatrix{17}&{87}&{9}\cr{-3}&{{-31}\over{2}}&{{-3}\over{2}}\cr
{-3}&{{-29}\over{2}}&{{-5}\over{2}}\cr
\endpmatrix,
\hskip30pt
C_2=
\pmatrix{449}&{1740}&{540}\cr{-60}&{-233}&{-72}\cr{-180}&{-696}&{-217}\cr
\endpmatrix$.}

\vskip5pt

\vbox{\noindent
$N=42 \ $ $d=183=3\cdot 61 ,\,\eta=0 ,\,h=2 $: $S=\langle 3 \rangle\oplus 
\langle -61 \rangle \oplus \langle -1 \rangle$;
\nobreak
\newline
$\rho=(305,-66,-122)$,\ $\rho^2=-1525$; 
\newline 
$\Gamma_+=
\pmatrix{9}&{-2}&{-1}\cr{0}&{0}&{1}\cr{0}&{1}&{0}\cr
{1}&{0}&{-3}\cr{5}&{-1}&{-4}\cr
\endpmatrix,  
\hskip20pt 
G(\Gamma_+)=
\pmatrix{-2}&{1}&{122}&{24}&{9}\cr{1}&{-1}&{0}&{3}&{4}\cr
{122}&{0}&{-61}&{0}&{61}\cr{24}&{3}&{0}&{-6}&{3}\cr{9}&{4}&{61}&{3}&{-2}\cr
\endpmatrix$;
\nobreak
\newline
$C_1=
\pmatrix{391}&{1708}&{56}\cr{-84}&{-367}&{-12}\cr{-168}&{-732}&{-25}\cr
\endpmatrix,
\hskip30pt
C_2=
\pmatrix{8111}&{34892}&{1404}\cr{-1716}&{-7382}&{-297}\cr
{-4212}&{-18117}&{-730}\cr
\endpmatrix$.}

\vskip5pt 

\vbox{\noindent
$N=43 \ $ $d=195=3\cdot 5\cdot 13 ,\,\eta=0 ,\,h=2 $: 
$S=\langle 3 \rangle\oplus 
\langle -13 \rangle \oplus \langle -5 \rangle$;
\nobreak
\newline
$\rho=(65,-27,-26)$,\ $\rho^2=-182 $; 
\newline 
$\Gamma_+\hskip-2pt=\hskip-2pt
\pmatrix{13}&{-6}&{-3}\cr{2}&{-1}&{0}\cr{0}&{0}&{1}\cr{0}&{1}&{0}\cr
{1}&{0}&{-1}\cr{65}&{-20}&{-39}\cr{17}&{-6}&{-9}\cr
\endpmatrix,  
\hskip5pt 
G(\Gamma_+)\hskip-2pt=\hskip-2pt
\pmatrix{-6}&{0}&{15}&{78}&{24}&{390}&{60}\cr
{0}&{-1}&{0}&{13}&{6}&{130}&{24}\cr
{15}&{0}&{-5}&{0}&{5}&{195}&{45}\cr{78}&{13}&{0}&{-13}&{0}&{260}&{78}\cr
{24}&{6}&{5}&{0}&{-2}&{0}&{6}\cr{390}&{130}&{195}&{260}&{0}&{-130}&{0}\cr
{60}&{24}&{45}&{78}&{6}&{0}&{-6}\cr
\endpmatrix$;
\nobreak
\newline
$C_1=
\pmatrix{244}&{455}&{140}\cr{-105}&{-196}&{-60}\cr{-84}&{-156}&{-49}\cr
\endpmatrix,
\hskip30pt
C_2=
\pmatrix{506}&{845}&{390}\cr{-195}&{-326}&{-150}\cr{-234}&{-390}&{-181}\cr
\endpmatrix$.}

\vskip5pt

\vbox{\noindent
$N=44 \ $ $d=195=3\cdot 5\cdot 13 ,\,\eta=1 ,\,h=2 $: 
$S=\langle 5 \rangle\oplus 
\langle -13 \rangle \oplus \langle -3 \rangle$;
\nobreak
\newline
$\rho=(39,-21,-26)$,\ $\rho^2=-156$; 
\newline 
$\Gamma_+=
\left(\smallmatrix{207}&{-120}&{-95}\cr{17}&{-10}&{-7}\cr{15}&{-9}&{-5}\cr
{247}&{-150}&{-65}\cr{18}&{-11}&{-4}\cr{78}&{-48}&{-13}\cr{8}&{-5}&{0}\cr
{0}&{0}&{1}\cr{0}&{1}&{0}\cr{3}&{0}&{-5}\cr
\endsmallmatrix\right),  
\hskip20pt 
G(\Gamma_+)=
\left(\smallmatrix{-30}&{0}&{60}&{3120}&{330}&
{2145}&{480}&{285}&{1560}&{1680}\cr
{0}&{-2}&{0}&{130}&{16}&{117}&{30}&{21}&{130}&{150}\cr
{60}&{0}&{-3}&{0}&{3}&{39}&{15}&{15}&{117}&{150}\cr
{3120}&{130}&{0}&{-130}&{0}&{195}&{130}&{195}&{1950}&{2730}\cr
{330}&{16}&{3}&{0}&{-1}&{0}&{5}&{12}&{143}&{210}\cr
{2145}&{117}&{39}&{195}&{0}&{-39}&{0}&{39}&{624}&{975}\cr
{480}&{30}&{15}&{130}&{5}&{0}&{-5}&{0}&{65}&{120}\cr
{285}&{21}&{15}&{195}&{12}&{39}&{0}&{-3}&{0}&{15}\cr
{1560}&{130}&{117}&{1950}&{143}&{624}&{65}&{0}&{-13}&{0}\cr
{1680}&{150}&{150}&{2730}&{210}&{975}&{120}&{15}&{0}&{-30}\cr
\endsmallmatrix\right)$;
\nobreak
\newline
$C_1=
\pmatrix{161}&{234}&{54}\cr{-90}&{-131}&{-30}\cr{-90}&{-130}&{-31}\cr
\endpmatrix,
\hskip30pt
C_2=
\pmatrix{159}&{208}&{72}\cr{-80}&{-105}&{-36}\cr{-120}&{-156}&{-55}\cr
\endpmatrix$.}

\vskip5pt

\vbox{\noindent
$N=45 \ $ $d=195=3\cdot 5\cdot 13 ,\,\eta=7 ,\,h=2 $: 
$S=\langle 1 \rangle\oplus 
\langle -65 \rangle \oplus \langle -3 \rangle$;
\nobreak
\newline
$\rho=(780,-93,-130)$,\ $\rho^2=-4485$; 
\newline 
$\Gamma_+=
\left(\smallmatrix{123}&{-15}&{-13}\cr{65}&{-8}&{-5}\cr
{8}&{-1}&{0}\cr{0}&{0}&{1}\cr
{0}&{1}&{0}\cr{1}&{0}&{-1}\cr{39}&{-4}&{-13}\cr{27}&{-3}&{-7}\cr
\endsmallmatrix\right),  
\hskip5pt 
G(\Gamma_+)\hskip-2pt=\hskip-2pt
\left(\smallmatrix{-3}&{0}&{9}&{39}&{975}&{84}&{390}&{123}\cr
{0}&{-10}&{0}&{15}&{520}&{50}&{260}&{90}\cr
{9}&{0}&{-1}&{0}&{65}&{8}&{52}&{21}\cr
{39}&{15}&{0}&{-3}&{0}&{3}&{39}&{21}\cr
{975}&{520}&{65}&{0}&{-65}&{0}&{260}&{195}\cr
{84}&{50}&{8}&{3}&{0}&{-2}&{0}&{6}\cr
{390}&{260}&{52}&{39}&{260}&{0}&{-26}&{0}\cr
{123}&{90}&{21}&{21}&{195}&{6}&{0}&{-3}\cr
\endsmallmatrix\right)$;
\nobreak
\newline
$C_1=
\pmatrix{2177}&{17160}&{792}\cr{-264}&{-2081}&{-96}\cr{-264}&{-2080}&{-97}\cr
\endpmatrix,
\hskip30pt
C_2=
\pmatrix{288}&{2210}&{153}\cr{-34}&{-261}&{-18}\cr{-51}&{-390}&{-28}\cr
\endpmatrix$.}

\vskip5pt 

\vbox{\noindent
$N=46 \ $ $d=210=2\cdot 3\cdot 5\cdot 7 ,\,
\eta=1 ,\,h=2 $: $S=\langle 14 \rangle\oplus 
\langle -30 \rangle\oplus 
\langle -2 \rangle$(0,1/2,1/2);
\nobreak
\newline
$\rho=(15,-7,-35)$,\ $\rho^2=-770$; 
\newline 
$\Gamma_+=
\pmatrix{10}&{-7}&{0}\cr{0}&{0}&{1}\cr{0}&{1}&{0}\cr{3}&{0}&{-8}\cr
{80}&{-7}&{-210}\cr
\endpmatrix,  
\hskip20pt 
G(\Gamma_+)=
\pmatrix{-70}&{0}&{210}&{420}&{9730}\cr{0}&{-2}&{0}&{16}&{420}\cr
{210}&{0}&{-30}&{0}&{210}\cr{420}&{16}&{0}&{-2}&{0}\cr
{9730}&{420}&{210}&{0}&{-70}\cr
\endpmatrix$;
\nobreak
\newline
$C_1=
\pmatrix{13}&{15}&{3}\cr{-7}&{{-17}\over{2}}&{{-3}\over{2}}\cr
{-21}&{{-45}\over{2}}&{{-11}\over{2}}\cr
\endpmatrix,
\hskip30pt
C_2=
\pmatrix{41}&{15}&{15}\cr
{-7}&{{-7}\over{2}}&{{-5}\over{2}}\cr
{-105}&{{-75}\over{2}}&{{-77}\over{2}}\cr
\endpmatrix$.}

\vskip5pt

\vbox{\noindent
$N=47 \ $ $d=210=2\cdot 3\cdot 5\cdot 7 ,\,\eta=7 ,\,h=2 $: 
$S=U \oplus \langle -210 \rangle$;
\nobreak
\newline
$\rho=(182,14,-5)$,\ $\rho^2=-154$; 
\newline 
$\Gamma_+=
\left(\smallmatrix{425}&{20}&{-9}\cr{52}&{2}&{-1}\cr
{105}&{0}&{-1}\cr{0}&{0}&{1}\cr{-1}&{1}&{0}\cr{10}&{10}&{-1}\cr
{14}&{7}&{-1}\cr{20}&{5}&{-1}\cr
\endsmallmatrix\right),  
\hskip20pt 
G(\Gamma_+)=
\left(\smallmatrix{-10}&{0}&{210}&{1890}&{405}&{2560}&{1365}&{635}\cr
{0}&{-2}&{0}&{210}&{50}&{330}&{182}&{90}\cr
{210}&{0}&{-210}&{210}&{105}&{840}&{525}&{315}\cr
{1890}&{210}&{210}&{-210}&{0}&{210}&{210}&{210}\cr
{405}&{50}&{105}&{0}&{-2}&{0}&{7}&{15}\cr
{2560}&{330}&{840}&{210}&{0}&{-10}&{0}&{40}\cr
{1365}&{182}&{525}&{210}&{7}&{0}&{-14}&{0}\cr
{635}&{90}&{315}&{210}&{15}&{40}&{0}&{-10}\cr
\endsmallmatrix\right)$;
\nobreak
\newline
$C_1=
\pmatrix{84}&{1445}&{7140}\cr{5}&{84}&{420}\cr{-2}&{-34}&{-169}\cr
\endpmatrix,
\hskip30pt
C_2=
\pmatrix{35}&{432}&{2520}\cr{3}&{35}&{210}\cr{-1}&{-12}&{-71}\cr
\endpmatrix$.}

\vskip5pt

\vbox{\noindent
$N=48 \ $ $d=213=3\cdot 71 ,\,\eta=3 ,\,h=2 $: $S=\langle 1 \rangle\oplus 
\langle -71 \rangle \oplus \langle -3 \rangle$;
\nobreak
\newline
$\rho=(355,-41,-71)$,\ $\rho^2=-8449 $; 
\newline 
$\Gamma_+=
\left(\smallmatrix{59}&{-7}&{-1}\cr{497}&{-59}&{0}\cr
{0}&{0}&{1}\cr{0}&{1}&{0}\cr
{1}&{0}&{-1}\cr{5325}&{-516}&{-1775}\cr{257}&{-25}&{-85}\cr
\endsmallmatrix\right),  
\hskip20pt 
G(\Gamma_+)=
\left(\smallmatrix{-1}&{0}&{3}&{497}&{56}&{52398}&{2483}\cr
{0}&{-142}&{0}&{4189}&{497}&{485001}&{23004}\cr
{3}&{0}&{-3}&{0}&{3}&{5325}&{255}\cr
{497}&{4189}&{0}&{-71}&{0}&{36636}&{1775}\cr
{56}&{497}&{3}&{0}&{-2}&{0}&{2}\cr
{52398}&{485001}&{5325}&{36636}&{0}&{-426}&{0}\cr
{2483}&{23004}&{255}&{1775}&{2}&{0}&{-1}\cr
\endsmallmatrix\right)$;
\nobreak
\newline
$C_1=
\pmatrix{288}&{2414}&{51}\cr{-34}&{-285}&{-6}\cr{-17}&{-142}&{-4}\cr
\endpmatrix,
\hskip30pt
C_2=
\pmatrix{99}&{710}&{90}\cr{-10}&{-72}&{-9}\cr{-30}&{-213}&{-28}\cr
\endpmatrix$.}

\vskip5pt 

\vbox{\noindent
$N=49 \ $ $d=221=13\cdot 17 ,\,\eta=1 ,\,h=2 $: $S=\langle 2 \rangle\oplus 
\langle -26 \rangle \oplus \langle -17 \rangle(1/2,1/2,0)$;
\nobreak
\newline
$\rho=(91/2,-15/2,-13)$,\ $\rho^2=-195$; 
\newline 
$\Gamma_+\hskip-4pt=\hskip-4pt
\left(\smallmatrix{51}&{-11}&{-11}\cr{{77}\over{2}}&{{-17}\over{2}}&{-8}\cr
{{5967}\over{2}}&{{-1343}\over{2}}&{-598}\cr{106}&{-24}&{-21}\cr
{{273}\over{2}}&{{-63}\over{2}}&{-26}\cr{{3}\over{2}}&{{-1}\over{2}}&{0}\cr
{0}&{0}&{1}\cr{0}&{1}&{0}\cr{17}&{0}&{-6}\cr{15}&{-1}&{-5}\cr
\endsmallmatrix\right),   
G(\Gamma_+)\hskip-4pt=\hskip-4pt
\left(\smallmatrix{-1}&{0}&{442}&{21}&{52}&{10}&{187}&{286}&{612}&{309}\cr
{0}&{-2}&{0}&{2}&{13}&{5}&{136}&{221}&{493}&{254}\cr
{442}&{0}&{-442}&{0}&{221}&{221}&{10166}&{17459}&{40443}&{21216}\cr
{21}&{2}&{0}&{-1}&{0}&{6}&{357}&{624}&{1462}&{771}\cr
{52}&{13}&{221}&{0}&{-26}&{0}&{442}&{819}&{1989}&{1066}\cr
{10}&{5}&{221}&{6}&{0}&{-2}&{0}&{13}&{51}&{32}\cr
{187}&{136}&{10166}&{357}&{442}&{0}&{-17}&{0}&{102}&{85}\cr
{286}&{221}&{17459}&{624}&{819}&{13}&{0}&{-26}&{0}&{26}\cr
{612}&{493}&{40443}&{1462}&{1989}&{51}&{102}&{0}&{-34}&{0}\cr
{309}&{254}&{21216}&{771}&{1066}&{32}&{85}&{26}&{0}&{-1}\cr
\endsmallmatrix\right)$;
\nobreak
\newline
$C_1=
\pmatrix{{423}\over{2}}&{{1105}\over{2}}&{425}\cr
{{-85}\over{2}}&{{-223}\over{2}}&{-85}\cr{-50}&{-130}&{-101}\cr
\endpmatrix,
\hskip30pt
C_2=
\pmatrix{48}&{91}&{119}\cr{-7}&{-14}&{-17}\cr{-14}&{-26}&{-35}\cr
\endpmatrix$.}

\vskip5pt 

\vbox{\noindent
$N=50 \ $ $d=231=3\cdot 7\cdot 11 ,\,\eta=2 ,\,h=0 $: 
$S=\langle 3 \rangle\oplus \langle -22 \rangle \oplus 
\langle -14 \rangle(0,1/2,1/2)$;
\nobreak
\newline
$\rho=(21,-7,-6)$,\ $\rho^2=-259$; 
\newline 
$\Gamma_+=
\left(\smallmatrix{8}&{-3}&{0}\cr{0}&{0}&{1}\cr{0}&{1}&{0}\cr{2}&{0}&{-1}\cr
{7}&{-1}&{-3}\cr{308}&{{-105}\over{2}}&{{-253}\over{2}}\cr
{41}&{{-15}\over{2}}&{{-33}\over{2}}\cr
\endsmallmatrix\right),  
\hskip20pt 
G(\Gamma_+)=
\left(\smallmatrix{-6}&{0}&{66}&{48}&{102}&{3927}&{489}\cr
{0}&{-14}&{0}&{14}&{42}&{1771}&{231}\cr
{66}&{0}&{-22}&{0}&{22}&{1155}&{165}\cr{48}&{14}&{0}&{-2}&{0}&{77}&{15}\cr
{102}&{42}&{22}&{0}&{-1}&{0}&{3}\cr{3927}&{1771}&{1155}&{77}&{0}&{-77}&{0}\cr
{489}&{231}&{165}&{15}&{3}&{0}&{-6}\cr
\endsmallmatrix\right)$;
\nobreak
\newline
$S_1=
\pmatrix{112}&{297}&{49}\cr{{-45}\over{2}}&{{-119}\over{2}}&{{-21}\over{2}}\cr
{{-87}\over{2}}&{{-231}\over{2}}&{{-37}\over{2}}\cr
\endpmatrix$.}

\vskip5pt 

\vbox{\noindent
$N=51 \ $ $d=231=3\cdot 7\cdot 11 ,\,\eta=3 ,\,h=2 $: 
$S=\langle 11 \rangle\oplus 
\langle -7 \rangle \oplus \langle -3 \rangle$;
\nobreak
\newline
$\rho=(105,-99,-154)$,\ $\rho^2=-18480$; 
\newline 
$\Gamma_+\hskip-3pt =\hskip-3pt
\left(\smallmatrix{14}&{-17}&{-7}\cr{4}&{-5}&{-1}\cr{35}&{-44}&{0}\cr
{0}&{0}&{1}\cr{0}&{1}&{0}\cr{1}&{0}&{-2}\cr{30}&{-11}&{-55}\cr
{35}&{-15}&{-63}\cr{216}&{-99}&{-385}\cr{154}&{-73}&{-273}\cr
\endsmallmatrix\right),  
\hskip3pt 
G(\Gamma_+)\hskip-3pt=\hskip-3pt
\left(\smallmatrix
{-14}&{0}&{154}&{21}&{119}&{112}&{2156}&{2282}&{13398}&{9296}\cr
{0}&{-2}&{0}&{3}&{35}&{38}&{770}&{826}&{4884}&{3402}\cr
{154}&{0}&{-77}&{0}&{308}&{385}&{8162}&{8855}&{52668}&{36806}\cr
{21}&{3}&{0}&{-3}&{0}&{6}&{165}&{189}&{1155}&{819}\cr
{119}&{35}&{308}&{0}&{-7}&{0}&{77}&{105}&{693}&{511}\cr
{112}&{38}&{385}&{6}&{0}&{-1}&{0}&{7}&{66}&{56}\cr
{2156}&{770}&{8162}&{165}&{77}&{0}&{-22}&{0}&{132}&{154}\cr
{2282}&{826}&{8855}&{189}&{105}&{7}&{0}&{-7}&{0}&{28}\cr
{13398}&{4884}&{52668}&{1155}&{693}&{66}&{132}&{0}&{-66}&{0}\cr
{9296}&{3402}&{36806}&{819}&{511}&{56}&{154}&{28}&{0}&{-14}\cr
\endsmallmatrix\right)$;
\nobreak
\newline
$C_1=
\pmatrix{21}&{14}&{6}\cr{-22}&{-15}&{-6}\cr{-22}&{-14}&{-7}\cr
\endpmatrix,
\hskip30pt
C_2=
\pmatrix{351}&{112}&{168}\cr{-176}&{-57}&{-84}\cr{-616}&{-196}&{-295}\cr
\endpmatrix$.}

\vskip5pt 

\vbox{\noindent
$N=52 \ $ $d=231=3\cdot 7\cdot 11 ,\,\eta=4 ,\,h=2 $: 
$S=\langle 7 \rangle\oplus \langle -33 \rangle \oplus 
\langle -1 \rangle$;
\nobreak
\newline
$\rho=(16,-7,-14)$,\ $\rho^2=-21$; 
\newline 
$\Gamma_+\hskip-2pt=\hskip-2pt
\pmatrix{42}&{-19}&{-21}\cr{77}&{-35}&{-33}\cr{15}&{-7}&{0}\cr
{0}&{0}&{1}\cr{0}&{1}&{0}\cr{1}&{0}&{-3}\cr{3}&{-1}&{-6}\cr
\endpmatrix,  
\hskip5pt 
G(\Gamma_+)\hskip-2pt=\hskip-2pt
\pmatrix{-6}&{0}&{21}&{21}&{627}&{231}&{129}\cr
{0}&{-11}&{0}&{33}&{1155}&{440}&{264}\cr
{21}&{0}&{-42}&{0}&{231}&{105}&{84}\cr
{21}&{33}&{0}&{-1}&{0}&{3}&{6}\cr{627}&{1155}&{231}&{0}&{-33}&{0}&{33}\cr
{231}&{440}&{105}&{3}&{0}&{-2}&{3}\cr{129}&{264}&{84}&{6}&{33}&{3}&{-6}\cr
\endpmatrix$;
\nobreak
\newline
$C_1=
\pmatrix{377}&{792}&{36}\cr{-168}&{-353}&{-16}\cr{-252}&{-528}&{-25}\cr
\endpmatrix,
\hskip30pt
C_2=
\pmatrix{230}&{462}&{33}\cr{-98}&{-197}&{-14}\cr{-231}&{-462}&{-34}\cr
\endpmatrix$.}

\vskip5pt 

\vbox{\noindent
$N=53 \ $ $d=255=3\cdot 5\cdot 17 ,\,\eta=2 ,\,h=2 $: 
$S=\langle 1 \rangle\oplus \langle -17 \rangle \oplus 
\langle -15 \rangle$;
\nobreak
\newline
$\rho= (765,-135,-136)$,\ $\rho^2=-2040 $; 
\newline 
$\Gamma_+=
\left(\smallmatrix{265}&{-50}&{-43}\cr{323}&{-62}&{-51}\cr{75}&{-15}&{-11}\cr
{9}&{-2}&{-1}\cr{255}&{-60}&{-17}\cr{4}&{-1}&{0}\cr{0}&{0}&{1}\cr
{0}&{1}&{0}\cr{3}&{0}&{-1}\cr{65}&{-10}&{-13}\cr
\endsmallmatrix\right),  
\hskip20pt 
G(\Gamma_+)=
\left(\smallmatrix{-10}&{0}&{30}&{40}&{5610}&{210}&{645}&{850}&{150}&{340}\cr
{0}&{-34}&{0}&{34}&{6120}&{238}&{765}&{1054}&{204}&{510}\cr
{30}&{0}&{-15}&{0}&{1020}&{45}&{165}&{255}&{60}&{180}\cr
{40}&{34}&{0}&{-2}&{0}&{2}&{15}&{34}&{12}&{50}\cr
{5610}&{6120}&{1020}&{0}&{-510}&{0}&{255}&{1020}&{510}&{3060}\cr
{210}&{238}&{45}&{2}&{0}&{-1}&{0}&{17}&{12}&{90}\cr
{645}&{765}&{165}&{15}&{255}&{0}&{-15}&{0}&{15}&{195}\cr
{850}&{1054}&{255}&{34}&{1020}&{17}&{0}&{-17}&{0}&{170}\cr
{150}&{204}&{60}&{12}&{510}&{12}&{15}&{0}&{-6}&{0}\cr
{340}&{510}&{180}&{50}&{3060}&{90}&{195}&{170}&{0}&{-10}\cr
\endsmallmatrix\right)$;
\nobreak
\newline
$C_1=
\pmatrix{577}&{1734}&{1530}\cr{-102}&{-307}&{-270}\cr{-102}&{-306}&{-271}\cr
\endpmatrix,
\hskip30pt
C_2=
\pmatrix{4607}&{13056}&{12960}\cr{-768}&{-2177}&{-2160}\cr
{-864}&{-2448}&{-2431}\cr
\endpmatrix$.}

\vskip5pt  

\vbox{\noindent
$N=54 \ $ $d=255=3\cdot 5\cdot 17 ,\,\eta=3 ,\,h=2 $: $S=\langle 15 \rangle
\oplus \langle -17 \rangle \oplus \langle -1 \rangle$;
\nobreak
\newline
$\rho=(119,-95,-255)$,\ $\rho^2=-6035 $; 
\newline 
$\Gamma_+=
\left(\smallmatrix{156}&{-135}&{-235}\cr{31}&{-27}&{-45}\cr{51}&{-45}&{-68}\cr
{1}&{-1}&{0}\cr{0}&{0}&{1}\cr{0}&{1}&{0}\cr{1}&{0}&{-5}\cr
\endsmallmatrix\right),  
\hskip20pt 
G(\Gamma_+) =
\left(\smallmatrix{-10}&{0}&{85}&{45}&{235}&{2295}&{1165}\cr
{0}&{-3}&{0}&{6}&{45}&{459}&{240}\cr{85}&{0}&{-34}&{0}&{68}&{765}&{425}\cr
{45}&{6}&{0}&{-2}&{0}&{17}&{15}\cr{235}&{45}&{68}&{0}&{-1}&{0}&{5}\cr
{2295}&{459}&{765}&{17}&{0}&{-17}&{0}\cr
{1165}&{240}&{425}&{15}&{5}&{0}&{-10}\cr
\endsmallmatrix\right)$;
\nobreak
\newline
$C_1=
\pmatrix{734}&{714}&{77}\cr{-630}&{-613}&{-66}\cr{-1155}&{-1122}&{-122}\cr
\endpmatrix,
\hskip30pt
C_2=
\pmatrix{79}&{68}&{12}\cr{-60}&{-52}&{-9}\cr{-180}&{-153}&{-28}\cr
\endpmatrix$.}

\vskip5pt 

\vbox{\noindent
$N=55 \ $ $d=330=2\cdot 3 \cdot 5 \cdot 11 ,\,\eta=6 ,\,h=2 $: 
$S=\langle 22 \rangle\oplus \langle-10 \rangle \oplus 
\langle -6 \rangle(0,1/2,1/2)$;
\nobreak
\newline
$\rho=(63,-66,-88)$,\ $\rho^2=-2706$; 
\newline 
$\Gamma_+\hskip-2pt=\hskip-2pt
\pmatrix{15}&{-19}&{-15}\cr{1}&{{-3}\over{2}}&{{-1}\over{2}}\cr{0}&{0}&{1}\cr
{0}&{1}&{0}\cr{1}&{0}&{-2}\cr{3}&{{-3}\over{2}}&{{-11}\over{2}}\cr
{15}&{-11}&{-25}\cr
\endpmatrix,  
\hskip5pt 
G(\Gamma_+)\hskip-2pt=\hskip-2pt
\pmatrix{-10}&{0}&{90}&{190}&{150}&{210}&{610}\cr
{0}&{-2}&{3}&{15}&{16}&{27}&{90}\cr{90}&{3}&{-6}&{0}&{12}&{33}&{150}\cr
{190}&{15}&{0}&{-10}&{0}&{15}&{110}\cr{150}&{16}&{12}&{0}&{-2}&{0}&{30}\cr
{210}&{27}&{33}&{15}&{0}&{-6}&{0}\cr{610}&{90}&{150}&{110}&{30}&{0}&{-10}\cr
\endpmatrix$;
\nobreak
\newline
$C_1=
\pmatrix{197}&{105}&{63}\cr{-231}&{{-247}\over{2}}&{{-147}\over{2}}\cr
{-231}&{{-245}\over{2}}&{{-149}\over{2}}\cr
\endpmatrix,
\hskip30pt
C_2=
\pmatrix{109}&{45}&{45}\cr{-99}&{{-83}\over{2}}&{{-81}\over{2}}\cr
{-165}&{{-135}\over{2}}&{{-137}\over{2}}\cr
\endpmatrix$.}

\vskip5pt

\vbox{\noindent
$N=56 \ $ $d=390=2\cdot 3 \cdot 5 \cdot 13 ,\,
\eta=0 ,\,h=2 $: $S=\langle 130 \rangle\oplus 
\langle -6 \rangle \oplus \langle -2 \rangle(0,1/2,1/2)$;
\nobreak
\newline
$\rho= (17,-26,-130)$,\ $\rho^2=-286$; 
\newline 
$\Gamma_+=
\pmatrix{3}&{-7}&{-21}\cr{0}&{{-1}\over{2}}&{{1}\over{2}}\cr{0}&{1}&{0}\cr
{8}&{0}&{-65}\cr{3}&{-2}&{-24}\cr
\endpmatrix,  
\hskip20pt 
G(\Gamma_+)=
\pmatrix{-6}&{0}&{42}&{390}&{78}\cr{0}&{-2}&{3}&{65}&{18}\cr
{42}&{3}&{-6}&{0}&{12}\cr{390}&{65}&{0}&{-130}&{0}\cr{78}&{18}&{12}&{0}&{-6}\cr
\endpmatrix$;
\nobreak
\newline
$C_1=
\pmatrix{311}&{24}&{36}\cr{-520}&{-41}&{-60}\cr{-2340}&{-180}&{-271}\cr
\endpmatrix,
\hskip30pt
C_2=
\pmatrix{519}&{30}&{62}\cr{-650}&{{-77}\over{2}}&{{-155}\over{2}}\cr
{-4030}&{{-465}\over{2}}&{{-963}\over{2}}\cr
\endpmatrix$.}

\vskip5pt 

\vbox{\noindent
$N=57 \ $ $d=429=3\cdot 11\cdot 13,\,\eta=5 ,\,h=2 $: 
$S=\langle 1 \rangle\oplus \langle -143 \rangle \oplus 
\langle -3 \rangle$;
\nobreak
\newline
$\rho=(2431,-175,-715) $,\ $\rho^2=-3289 $; 
\newline 
$\Gamma_+=
\left(\smallmatrix{1352}&{-99}&{-377}\cr{528}&{-39}&{-143}\cr
{13}&{-1}&{-3}\cr{11}&{-1}&{0}\cr{0}&{0}&{1}\cr{0}&{1}&{0}\cr{1}&{0}&{-1}\cr
{117}&{-8}&{-39}\cr
\endsmallmatrix\right),  
\hskip20pt 
G(\Gamma_+)=
\left(\smallmatrix{-26}&{0}&{26}&{715}&{1131}&{14157}&{221}&{819}\cr
{0}&{-66}&{0}&{231}&{429}&{5577}&{99}&{429}\cr
{26}&{0}&{-1}&{0}&{9}&{143}&{4}&{26}\cr
{715}&{231}&{0}&{-22}&{0}&{143}&{11}&{143}\cr
{1131}&{429}&{9}&{0}&{-3}&{0}&{3}&{117}\cr
{14157}&{5577}&{143}&{143}&{0}&{-143}&{0}&{1144}\cr
{221}&{99}&{4}&{11}&{3}&{0}&{-2}&{0}\cr
{819}&{429}&{26}&{143}&{117}&{1144}&{0}&{-26}\cr
\endsmallmatrix\right)$;
\nobreak
\newline
$C_1=
\pmatrix{9074}&{94380}&{7755}\cr{-660}&{-6865}&{-564}\cr
{-2585}&{-26884}&{-2210}\cr
\endpmatrix,
\hskip30pt
C_2=
\pmatrix{4899}&{50050}&{4410}\cr{-350}&{-3576}&{-315}\cr
{-1470}&{-15015}&{-1324}\cr
\endpmatrix$.}

\vskip5pt 

\vbox{\noindent
$N=58 \ $ $d=462=2\cdot 3\cdot 7\cdot 11,\,\eta=0 ,\,h=2 $: 
$S=\langle 14 \rangle\oplus \langle -22 \rangle \oplus 
\langle -6\rangle(1/2,1/2,0)$;
\nobreak
\newline
$\rho=(33/2,-15/2,-22)$,\ $\rho^2=-330 $; 
\newline 
$\Gamma_+ = 
\left(\smallmatrix{{81}\over{2}}&{{-51}\over{2}}&{-38}\cr
{{1}\over{2}}&{{-1}\over{2}}&{0}\cr{0}&{0}&{1}\cr{0}&{1}&{0}\cr
{9}&{0}&{-14}\cr{66}&{-10}&{-99}\cr{{51}\over{2}}&{{-9}\over{2}}&{-38}\cr
\endsmallmatrix\right),  
\hskip20pt 
G(\Gamma_+) =
\left(\smallmatrix{-6}&{3}&{228}&{561}&{1911}&{9240}&{3270}\cr
{3}&{-2}&{0}&{11}&{63}&{352}&{129}\cr{228}&{0}&{-6}&{0}&{84}&{594}&{228}\cr
{561}&{11}&{0}&{-22}&{0}&{220}&{99}\cr{1911}&{63}&{84}&{0}&{-42}&{0}&{21}\cr
{9240}&{352}&{594}&{220}&{0}&{-22}&{0}\cr
{3270}&{129}&{228}&{99}&{21}&{0}&{-6}\cr
\endsmallmatrix\right)$;
\nobreak
\newline
$C_1=
\pmatrix{349}&{330}&{150}\cr{-210}&{-199}&{-90}\cr{-350}&{-330}&{-151}\cr
\endpmatrix,
\hskip30pt
C_2=
\pmatrix{{61}\over{2}}&{{33}\over{2}}&{18}\cr
{{-21}\over{2}}&{{-13}\over{2}}&{-6}\cr{-42}&{-22}&{-25}\cr
\endpmatrix$.}

\vskip5pt

\vbox{\noindent
$N=59 \ $ $d=462=2\cdot 3\cdot 7\cdot 11 ,\,\eta=6 ,\,h=2 $: 
$S=\langle 22 \rangle\oplus 
\langle -14 \rangle \oplus \langle -6 \rangle(1/2,1/2,0)$;
\nobreak
\newline
$\rho=(259/2,-121/2,-231)$,\ $\rho^2=-2464$; 
\newline 
$\Gamma_+=
\left(\smallmatrix{21}&{-13}&{-35}\cr
{{3}\over{2}}&{{-3}\over{2}}&{-2}\cr
{{7}\over{2}}&{{-11}\over{2}}&{0}\cr{0}&{0}&{1}\cr
{0}&{1}&{0}\cr{1}&{0}&{-2}\cr{210}&{-76}&{-385}\cr
\endsmallmatrix\right),  
\hskip20pt 
G(\Gamma_+)=
\left(\smallmatrix{-14}&{0}&{616}&{210}&{182}&{42}&{2338}\cr
{0}&{-6}&{0}&{12}&{21}&{9}&{714}\cr
{616}&{0}&{-154}&{0}&{77}&{77}&{10318}\cr
{210}&{12}&{0}&{-6}&{0}&{12}&{2310}\cr
{182}&{21}&{77}&{0}&{-14}&{0}&{1064}\cr
{42}&{9}&{77}&{12}&{0}&{-2}&{0}\cr
{2338}&{714}&{10318}&{2310}&{1064}&{0}&{-14}\cr
\endsmallmatrix\right)$;
\nobreak
\newline
$C_1=\pmatrix{351}&{112}&{168}\cr{-176}&{-57}&{-84}\cr
{-616}&{-196}&{-295}\cr
\endpmatrix 
\hskip30pt
C_2=
\pmatrix{{13473}\over{2}}&{{3185}\over{2}}&{3360}\cr
{{-5005}\over{2}}&{{-1185}\over{2}}&{-1248}\cr
{-12320}&{-2912}&{-6145}\cr\endpmatrix$.}

\vskip5pt 

\vbox{\noindent
$N=60 \ $ $d=555=3\cdot 5\cdot 37 ,\,\eta=3 ,\,h=2 $: 
$S=\langle 15 \rangle\oplus \langle -37 \rangle \oplus 
\langle -1 \rangle$;
\nobreak
\newline
$\rho=(37,-23,-37)$,\ $\rho^2=-407$; 
\newline 
$\Gamma_+=
\left(\smallmatrix{11}&{-7}&{-2}\cr{259}&{-165}&{0}\cr{0}&{0}&{1}\cr
{0}&{1}&{0}\cr{1}&{0}&{-5}\cr{37}&{-15}&{-111}\cr{2}&{-1}&{-5}\cr
\endsmallmatrix\right),  
\hskip20pt 
G(\Gamma_+)=
\left(\smallmatrix{-2}&{0}&{2}&{259}&{155}&{1998}&{61}\cr
{0}&{-1110}&{0}&{6105}&{3885}&{52170}&{1665}\cr
{2}&{0}&{-1}&{0}&{5}&{111}&{5}\cr{259}&{6105}&{0}&{-37}&{0}&{555}&{37}\cr
{155}&{3885}&{5}&{0}&{-10}&{0}&{5}\cr
{1998}&{52170}&{111}&{555}&{0}&{-111}&{0}\cr
{61}&{1665}&{5}&{37}&{5}&{0}&{-2}\cr
\endsmallmatrix\right)$;
\nobreak
\newline
$C_1=
\pmatrix{124}&{185}&{10}\cr{-75}&{-112}&{-6}\cr{-150}&{-222}&{-13}\cr
\endpmatrix,
\hskip30pt
C_2=
\pmatrix{866}&{1258}&{85}\cr{-510}&{-741}&{-50}\cr{-1275}&{-1850}&{-126}\cr
\endpmatrix$.}

\vskip5pt 

\vbox{\noindent
$N=61 \ $ $d=561=3\cdot 11\cdot 17 ,\,\eta=6 ,\,h=2 $: 
$S=\langle 1 \rangle\oplus \langle -51 \rangle \oplus 
\langle -11 \rangle$;
\nobreak
\newline
$\rho=(119,-15,-17)$,\ $\rho^2=-493$; 
\newline 
$\Gamma_+=
\left(\smallmatrix{81}&{-11}&{-6}\cr{561}&{-77}&{-34}\cr
{7}&{-1}&{0}\cr
{0}&{0}&{1}\cr{0}&{1}&{0}\cr{3}&{0}&{-1}\cr{12}&{-1}&{-3}\cr
\endsmallmatrix\right),  
\hskip20pt 
G(\Gamma_+)=
\left(\smallmatrix{-6}&{0}&{6}&{66}&{561}&{177}&{213}\cr
{0}&{-374}&{0}&{374}&{3927}&{1309}&{1683}\cr
{6}&{0}&{-2}&{0}&{51}&{21}&{33}\cr
{66}&{374}&{0}&{-11}&{0}&{11}&{33}\cr
{561}&{3927}&{51}&{0}&{-51}&{0}&{51}\cr
{177}&{1309}&{21}&{11}&{0}&{-2}&{3}\cr
{213}&{1683}&{33}&{33}&{51}&{3}&{-6}\cr
\endsmallmatrix\right)$;
\nobreak
\newline
$C_1=
\pmatrix{63}&{408}&{88}\cr{-8}&{-52}&{-11}\cr{-8}&{-51}&{-12}\cr
\endpmatrix,
\hskip30pt
C_2=
\pmatrix{296}&{1683}&{594}\cr{-33}&{-188}&{-66}\cr
{-54}&{-306}&{-109}\cr
\endpmatrix$.}

\vskip5pt 

\vbox{\noindent
$N=62 \ $ $d=609=3\cdot 7\cdot 29 ,\,\eta=3 ,\,h=2 $: 
$S=\langle 29\rangle\oplus \langle -7 \rangle \oplus 
\langle -3 \rangle$;
\nobreak
\newline
$\rho=(77,-87,-203)$,\ $\rho^2=-4669$; 
\newline 
$\Gamma_+=
\left(\smallmatrix{1}&{-2}&{-1}\cr{14}&{-29}&{0}\cr
{0}&{0}&{1}\cr{0}&{1}&{0}\cr
{9}&{0}&{-29}\cr{7}&{-4}&{-21}\cr
{191}&{-145}&{-551}\cr{40}&{-31}&{-115}\cr
\endsmallmatrix\right),  
\hskip20pt 
G(\Gamma_+)=
\left(\smallmatrix{-2}&{0}&{3}&{14}&{174}&{84}&{1856}&{381}\cr
{0}&{-203}&{0}&{203}&{3654}&{2030}&{48111}&{9947}\cr
{3}&{0}&{-3}&{0}&{87}&{63}&{1653}&{345}\cr
{14}&{203}&{0}&{-7}&{0}&{28}&{1015}&{217}\cr
{174}&{3654}&{87}&{0}&{-174}&{0}&{1914}&{435}\cr
{84}&{2030}&{63}&{28}&{0}&{-14}&{0}&{7}\cr
{1856}&{48111}&{1653}&{1015}&{1914}&{0}&{-29}&{0}\cr
{381}&{9947}&{345}&{217}&{435}&{7}&{0}&{-2}\cr
\endsmallmatrix\right)$;
\nobreak
\newline
$C_1=
\pmatrix{86}&{21}&{24}\cr{-87}&{-22}&{-24}\cr{-232}&{-56}&{-65}\cr
\endpmatrix,
\hskip30pt
C_2=
\pmatrix{1043}&{210}&{306}\cr{-870}&{-176}&{-255}\cr{-2958}&{-595}&{-868}\cr
\endpmatrix$.}

\vskip5pt

\vbox{\noindent
$N=63 \ $ $d=663=3\cdot 13\cdot 17 ,\,\eta=3 ,\,h=2 $: 
$S=\langle 6\rangle\oplus 
\langle -26\rangle \oplus \langle -17\rangle(1/2,1/2,0)$;
\nobreak
\newline
$\rho=(884,-357,-351)$,\ $\rho^2= -719355$; 
\newline 
$\Gamma_+\hskip-2pt=\hskip-2pt
\left(\smallmatrix{{1071}\over{2}}&{{-493}\over{2}}&{-91}\cr
{91}&{-42}&{-15}\cr{2}&{-1}&{0}\cr{0}&{0}&{1}\cr{0}&{1}&{0}\cr
{5}&{0}&{-3}\cr{{221}\over{2}}&{{-17}\over{2}}&{-65}\cr
{{7}\over{2}}&{{-1}\over{2}}&{-2}\cr{{85}\over{2}}&{{-17}\over{2}}&{-23}\cr
\endsmallmatrix\right),  
\hskip2pt 
G(\Gamma_+)\hskip-2pt=\hskip-2pt
\left(\smallmatrix{-34}&{0}&{17}&{1547}&{6409}&
{11424}&{200005}&{4947}&{46495}\cr
{0}&{-3}&{0}&{255}&{1092}&{1965}&{34476}&{855}&{8058}\cr
{17}&{0}&{-2}&{0}&{26}&{60}&{1105}&{29}&{289}\cr
{1547}&{255}&{0}&{-17}&{0}&{51}&{1105}&{34}&{391}\cr
{6409}&{1092}&{26}&{0}&{-26}&{0}&{221}&{13}&{221}\cr
{11424}&{1965}&{60}&{51}&{0}&{-3}&{0}&{3}&{102}\cr
{200005}&{34476}&{1105}&{1105}&{221}&{0}&{-442}&{0}&{884}\cr
{4947}&{855}&{29}&{34}&{13}&{3}&{0}&{-1}&{0}\cr
{46495}&{8058}&{289}&{391}&{221}&{102}&{884}&{0}&{-34}\cr
\endsmallmatrix\right)$;
\nobreak
\newline
$C_1=
\pmatrix{{361}\over{2}}&{{715}\over{2}}&{{187}\over{2}}\cr
{{-165}\over{2}}&{{-327}\over{2}}&{{-85}\over{2}}\cr{-33}&{-65}&{-18}\cr
\endpmatrix,
\hskip30pt
C_2=
\pmatrix{95}&{104}&{136}\cr{-24}&{-27}&{-34}\cr{-48}&{-52}&{-69}\cr
\endpmatrix$.}

\vskip5pt 

\vbox{\noindent
$N=64 \ $ $d=705=3\cdot 5\cdot 47 ,\,\eta=6 ,\,h=2 $: 
$S=\langle 1\rangle\oplus 
\langle -47\rangle \oplus \langle -15\rangle$;
\nobreak
\newline
$\rho= (705, -90, -94)$,\ $\rho^2=-16215 $; 
\newline 
$\Gamma_+\hskip-2pt=\hskip-2pt\left(\smallmatrix
{2961}&{-411}&{-235}\cr{35}&{-5}&{-2}\cr{47}&{-7}&{0}\cr{0}&{0}&{1}\cr
{0}&{1}&{0}\cr{3}&{0}&{-1}\cr{31}&{-3}&{-6}\cr{1269}&{-129}&{-235}\cr
\endsmallmatrix\right),  
\hskip5pt 
G(\Gamma_+)\hskip-2pt=\hskip-2pt
\left(\smallmatrix{-141}&{0}&{3948}&{3525}&{19317}&{5358}&{12690}&{437241}\cr
{0}&{-10}&{0}&{30}&{235}&{75}&{200}&{7050}\cr
{3948}&{0}&{-94}&{0}&{329}&{141}&{470}&{17202}\cr
{3525}&{30}&{0}&{-15}&{0}&{15}&{90}&{3525}\cr
{19317}&{235}&{329}&{0}&{-47}&{0}&{141}&{6063}\cr
{5358}&{75}&{141}&{15}&{0}&{-6}&{3}&{282}\cr
{12690}&{200}&{470}&{90}&{141}&{3}&{-2}&{0}\cr
{437241}&{7050}&{17202}&{3525}&{6063}&{282}&{0}&{-141}\cr
\endsmallmatrix\right)$;
\nobreak
\newline
$C_1=
\pmatrix{967}&{6204}&{1320}\cr{-132}&{-847}&{-180}\cr{-88}&{-564}&{-121}\cr
\endpmatrix,
\hskip30pt
C_2=
\pmatrix{63}&{376}&{120}\cr{-8}&{-48}&{-15}\cr{-8}&{-47}&{-16}\cr
\endpmatrix$.}

\vskip5pt

\vbox{\noindent
$N=65 \ $ $d=1365=3\cdot 5\cdot 7\cdot 13 ,\,\eta=9 ,\,h=2 $: 
$S=\langle 5 \rangle\oplus \langle -91\rangle \oplus 
\langle -3 \rangle$;
\nobreak
\newline
$\rho= (117,-25,-65)$,\ $\rho^2=-1105$; $\Gamma_+=$
\newline 
$
\left(\smallmatrix{1638}&{-375}&{-455}\cr{13}&{-3}&{-3}\cr{21}&{-5}&{0}\cr
{0}&{0}&{1}\cr{0}&{1}&{0}\cr{3}&{0}&{-5}\cr{13}&{-2}&{-13}\cr
{11}&{-2}&{-9}\cr{3003}&{-570}&{-2275}\cr
\endsmallmatrix\right),\hskip0pt    
G(\Gamma_+)\hskip-2pt = \hskip-2pt 
\left(\smallmatrix
{-2730}&{0}&{1365}&{1365}&{34125}&{17745}&{20475}&{9555}&{2037945}\cr
{0}&{-1}&{0}&{9}&{273}&{150}&{182}&{88}&{19110}\cr
{1365}&{0}&{-70}&{0}&{455}&{315}&{455}&{245}&{55965}\cr
{1365}&{9}&{0}&{-3}&{0}&{15}&{39}&{27}&{6825}\cr
{34125}&{273}&{455}&{0}&{-91}&{0}&{182}&{182}&{51870}\cr
{17745}&{150}&{315}&{15}&{0}&{-30}&{0}&{30}&{10920}\cr
{20475}&{182}&{455}&{39}&{182}&{0}&{-26}&{0}&{2730}\cr
{9555}&{88}&{245}&{27}&{182}&{30}&{0}&{-2}&{0}\cr
{2037945}&{19110}&{55965}&{6825}&{51870}&{10920}&{2730}&{0}&{-2730}\cr
\endsmallmatrix\right)$;
\nobreak
\newline
$C_1=
\pmatrix{139}&{546}&{42}\cr{-30}&{-118}&{-9}\cr{-70}&{-273}&{-22}\cr
\endpmatrix,
\hskip30pt
C_2=
\pmatrix{374}&{1365}&{150}\cr{-75}&{-274}&{-30}\cr{-250}&{-910}&{-101}\cr
\endpmatrix$.}

\vskip5pt 

\vbox{\noindent
$N=66 \ $ $d=1995=3\cdot 5\cdot 7\cdot 19 ,\,\eta=7 ,\,h=2 $: 
$S=\langle 35 \rangle\oplus 
\langle -57 \rangle \oplus \langle -1 \rangle$;
\nobreak
\newline
$\rho= (171,-105,-665)$,\ $\rho^2=-47215 $; $\Gamma_+=$
\newline 
$\left(\smallmatrix
{4}&{-3}&{-7}\cr{285}&{-217}&{-399}\cr{13}&{-10}&{-15}\cr
{9}&{-7}&{-7}\cr{19}&{-15}&{0}\cr{0}&{0}&{1}\cr{0}&{1}&{0}\cr
{1}&{0}&{-7}\cr{19}&{-8}&{-95}\cr{153}&{-70}&{-735}\cr
{47}&{-22}&{-223}\cr
\endsmallmatrix\right),  
\hskip2pt 
\left(\smallmatrix
{-2}&{0}&{5}&{14}&{95}&{7}&{171}&{91}&{627}&{4305}&{1257}\cr
{0}&{-399}&{0}&{399}&{3990}&{399}&{12369}&{7182}&{52668}&{367080}&{107730}\cr
{5}&{0}&{-10}&{0}&{95}&{15}&{570}&{350}&{2660}&{18690}&{5500}\cr
{14}&{399}&{0}&{-7}&{0}&{7}&{399}&{266}&{2128}&{15120}&{4466}\cr
{95}&{3990}&{95}&{0}&{-190}&{0}&{855}&{665}&{5795}&{41895}&{12445}\cr
{7}&{399}&{15}&{7}&{0}&{-1}&{0}&{7}&{95}&{735}&{223}\cr
{171}&{12369}&{570}&{399}&{855}&{0}&{-57}&{0}&{456}&{3990}&{1254}\cr
{91}&{7182}&{350}&{266}&{665}&{7}&{0}&{-14}&{0}&{210}&{84}\cr
{627}&{52668}&{2660}&{2128}&{5795}&{95}&{456}&{0}&{-38}&{0}&{38}\cr
{4305}&{367080}&{18690}&{15120}&{41895}&{735}&{3990}&{210}&{0}&{-210}&{0}\cr
{1257}&{107730}&{5500}&{4466}&{12445}&{223}&{1254}&{84}&{38}&{0}&{-2}\cr
\endsmallmatrix\right)$;
\nobreak
\newline
$C_1=
\pmatrix{104}&{114}&{9}\cr{-70}&{-77}&{-6}\cr{-315}&{-342}&{-28}\cr
\endpmatrix,
\hskip30pt
C_2=
\pmatrix{139}&{114}&{18}\cr{-70}&{-58}&{-9}\cr{-630}&{-513}&{-82}\cr
\endpmatrix$.}

\newpage

\centerline{\bf Table 5.}
\centerline{The list of odd hyperbolic lattices of rank three}
\centerline{and with even square-free determinant (i.e. non-main)}
\centerline{which are reflective of hyperbolic type.}

\vskip30pt

\vbox{\noindent
$N^\prime=1\ $ $d=30=2\cdot 13,\,\text{odd},\,\eta=0,\,h=2$:
$\langle 13 \rangle \oplus \langle -1 \rangle \oplus \langle -2 \rangle$
($=\widetilde{(13,1,h=1)}$). 
\nobreak
\newline
$\rho=(8,-26,-13)$,\ $\rho^2=-182$; 
\newline 
$\Gamma_+=
\pmatrix{18}&{-65}&{0}\cr{0}&{0}&{1}\cr
{0}&{1}&{0}\cr{5}&{0}&{-13}\cr
\endpmatrix,  
\hskip20pt 
G(\Gamma_+)=
\pmatrix{-13}&{0}&{65}&{1170}\cr{0}&{-2}&{0}&{26}\cr
{65}&{0}&{-1}&{0}\cr{1170}&{26}&{0}&{-13}\cr
\endpmatrix$;
\nobreak
\newline
$C_1=
\pmatrix{12}&{3}&{2}\cr{-39}&{-10}&{-6}\cr{-13}&{-3}&{-3}\cr
\endpmatrix,
\hskip30pt
C_2=\pmatrix
{25}&{4}&{8}\cr{-52}&{-9}&{-16}\cr{-52}&{-8}&{-17}\cr
\endpmatrix$.}

\vskip5pt 

\noindent
$N^\prime=2\ $
$d=38=2\cdot 19,\,\text{odd},\,\eta=0,\,h=2$:
$\langle 1 \rangle \oplus \langle -1 \rangle \oplus \langle -38 \rangle$.
It is equivariantly equivalent to
$d=19,\,\eta=1,\,h=2$: $U\oplus \langle -19 \rangle$.

\vskip5pt

\vbox{\noindent
$N^\prime=3\ $ $d=38=2\cdot 19,\,\text{odd},\,\eta=1 ,\,h=2$:
$\langle 38 \rangle \oplus \langle -1 \rangle \oplus \langle -1 \rangle$
($=\widetilde{(19,0,h=1)}$).  
\nobreak
\newline
$\rho=(24,-95,-114)$,\ $\rho^2=-133$; 
\newline 
$\Gamma_+=
\pmatrix{39}&{-170}&{-170}\cr{0}&{-1}&{1}\cr
{0}&{1}&{0}\cr{3}&{0}&{-19}\cr{1}&{-2}&{-6}\cr
\endpmatrix,  
\hskip20pt 
G(\Gamma_+)=
\pmatrix{-2}&{0}&{170}&{1216}&{122}\cr{0}&{-2}&{1}&{19}&{4}\cr
{170}&{1}&{-1}&{0}&{2}\cr{1216}&{19}&{0}&{-19}&{0}\cr{122}&{4}&{2}&{0}&{-2}\cr
\endpmatrix$;
\nobreak
\newline
$C_1=
\pmatrix{341}&{36}&{42}\cr{-1368}&{-145}&{-168}\cr
{-1596}&{-168}&{-197}\cr
\endpmatrix,
\hskip30pt
C_2=
\pmatrix{1899}&{180}&{250}\cr{-6840}&{-649}&{-900}\cr
{-9500}&{-900}&{-1251}\cr
\endpmatrix$.}

\vskip5pt 

\noindent
$N^\prime= 4\ $ $d=46=2\cdot 23,\,\text{odd},\,\eta=1,\,h=2$:
$\langle 1 \rangle\oplus \langle -1 \rangle \oplus \langle -46 \rangle$.
It is equivariantly equivalent to
$d=23,\,\eta=1,\,h=2$: $U\oplus \langle -23 \rangle $.

\vskip5pt

\noindent
$N^\prime =5\ $ $d=70=2\cdot 5\cdot 7 ,\,\text{odd},\,\eta=0,\,h=2$:
$\langle 1 \rangle\oplus \langle -1 \rangle \oplus \langle -70 \rangle  $.
It is equivariantly equivalent to
$d=35,\,\eta=1,\,h=2$: $U\oplus \langle -35\rangle $.

\vskip5pt 

\noindent
$N^\prime =6\ $ $d=78=2\cdot 3\cdot 13,\,\text{odd},\,\eta=2,\,h=2$:
$\langle 1\rangle\oplus \langle -1 \rangle \oplus \langle -78 \rangle  $.
It is equivariantly equivalent to
$d=39,\,\eta=1,\,h=2$: $U\oplus \langle -39 \rangle$.

\vskip5pt

\vbox{\noindent
$N^\prime=7\ $ $d=102=2\cdot 3\cdot 17,\,\text{odd},\,\eta=1,\,h=2$:
$\langle 6 \rangle \oplus \langle -34 \rangle \oplus \langle-2\rangle
(0,1/2,1/2)$ ($=\widetilde{(51,0,h=0)}$). 
\nobreak  
\newline 
$\rho=(17,-6,-17)$,\ $\rho^2=-68$; 
\newline 
$\Gamma_+=
\left(\smallmatrix{6}&{{-5}\over{2}}&{{-3}\over{2}}\cr
{68}&{{-57}\over{2}}&{{-17}\over{2}}\cr
{25}&{{-21}\over{2}}&{{-3}\over{2}}\cr
{119}&{-50}&{0}\cr{0}&{0}&{1}\cr{0}&{1}&{0}\cr
{1}&{0}&{-2}\cr{7}&{{-3}\over{2}}&{{-21}\over{2}}\cr
{14}&{{-7}\over{2}}&{{-39}\over{2}}\cr
\endsmallmatrix\right),  
\hskip20pt 
G(\Gamma_+)=
\left(\smallmatrix{-1}&{0}&{3}&{34}&{3}&{85}&{30}&{93}&{148}\cr
{0}&{-17}&{0}&{102}&{17}&{969}&{374}&{1224}&{1989}\cr
{3}&{0}&{-3}&{0}&{3}&{357}&{144}&{483}&{792}\cr
{34}&{102}&{0}&{-34}&{0}&{1700}&{714}&{2448}&{4046}\cr
{3}&{17}&{3}&{0}&{-2}&{0}&{4}&{21}&{39}\cr
{85}&{969}&{357}&{1700}&{0}&{-34}&{0}&{51}&{119}\cr
{30}&{374}&{144}&{714}&{4}&{0}&{-2}&{0}&{6}\cr
{93}&{1224}&{483}&{2448}&{21}&{51}&{0}&{-3}&{0}\cr
{148}&{1989}&{792}&{4046}&{39}&{119}&{6}&{0}&{-1}\cr
\endsmallmatrix\right)$;
\nobreak
\newline
$C_1=
\pmatrix{127}&{272}&{32}\cr{-48}&{-103}&{-12}\cr{-96}&{-204}&{-25}\cr
\endpmatrix,
\hskip30pt
C_2=
\pmatrix{53}&{102}&{18}\cr{-18}&{-35}&{-6}\cr{-54}&{-102}&{-19}\cr
\endpmatrix$.}

\hskip5pt 

\vbox{\noindent
$N^\prime=8\ $ $d=102=2\cdot 3\cdot 17 ,\,\text{odd},\,\eta=2,\,h=2 $:
$\langle 102 \rangle \oplus \langle -1 \rangle \oplus \langle -1 \rangle$
($=\widetilde{(51,3,h=1)}$). 
\nobreak
\newline
$\rho=(6,-34,-51)$,\ $\rho^2=-85$; 
\newline 
$\Gamma_+=
\pmatrix{7}&{-50}&{-50}\cr{0}&{-1}&{1}\cr{0}&{1}&{0}\cr
{5}&{0}&{-51}\cr{1}&{-2}&{-10}\cr
\endpmatrix,  
\hskip20pt 
G(\Gamma_+)=
\pmatrix{-2}&{0}&{50}&{1020}&{114}\cr{0}&{-2}&{1}&{51}&{8}\cr
{50}&{1}&{-1}&{0}&{2}\cr{1020}&{51}&{0}&{-51}&{0}\cr
{114}&{8}&{2}&{0}&{-2}\cr
\endpmatrix$;
\nobreak
\newline
$C_1=
\pmatrix{1087}&{72}&{80}\cr{-7344}&{-487}&{-540}\cr
{-8160}&{-540}&{-601}\cr
\endpmatrix,
\hskip30pt
C_2=
\pmatrix{101}&{6}&{8}\cr{-612}&{-37}&{-48}\cr{-816}&{-48}&{-65}\cr
\endpmatrix$.}

\vskip5pt 

\noindent
$N^\prime =9\ $ $d=110=2\cdot 5\cdot 11,\,\text{odd},\,\eta=2,\,h=2$:
$\langle 22 \rangle\oplus \langle -5 \rangle \oplus \langle  -1 \rangle  $.
It is equivariantly equivalent to
$d=55,\,\eta=1,\,h=2$: 
$\langle 11\rangle \oplus \langle -10 \rangle \oplus \langle -2 \rangle 
(0,1/2,1/2)$.

\vskip5pt

\vbox{\noindent
$N^\prime=10\ $ $d=110=2\cdot 5\cdot 11,\,\text{odd},\,\eta=3 ,\,h=2 $:
$\langle 22 \rangle \oplus \langle -10 \rangle \oplus \langle -2 \rangle
(0,1/2,1/2)$ ($=\widetilde{(55,0,h=1)}$).  
\nobreak
\newline
$\rho=(12,-11/2,-77/2)$,\ $\rho^2=-99$;  
\newline 
$\Gamma_+=
\left(\smallmatrix{2}&{-2}&{-5}\cr{5}&{-6}&{-10}\cr
{4}&{{-11}\over{2}}&{{-11}\over{2}}\cr
{1}&{{-3}\over{2}}&{{-1}\over{2}}\cr
{0}&{0}&{1}\cr{0}&{1}&{0}\cr{3}&{0}&{-10}\cr
{20}&{-3}&{-66}\cr
\endsmallmatrix\right),  
\hskip20pt 
G(\Gamma_+)=
\left(\smallmatrix{-2}&{0}&{11}&{9}&{10}&{20}&{32}&{160}\cr
{0}&{-10}&{0}&{10}&{20}&{60}&{130}&{700}\cr
{11}&{0}&{-11}&{0}&{11}&{55}&{154}&{869}\cr
{9}&{10}&{0}&{-1}&{1}&{15}&{56}&{329}\cr
{10}&{20}&{11}&{1}&{-2}&{0}&{20}&{132}\cr
{20}&{60}&{55}&{15}&{0}&{-10}&{0}&{30}\cr
{32}&{130}&{154}&{56}&{20}&{0}&{-2}&{0}\cr
{160}&{700}&{869}&{329}&{132}&{30}&{0}&{-2}\cr
\endsmallmatrix\right)$;
\nobreak
\newline
$C_1=
\pmatrix{109}&{30}&{30}\cr{-66}&{-19}&{-18}\cr
{-330}&{-90}&{-91}\cr
\endpmatrix,
\hskip30pt
C_2=
\pmatrix{175}&{20}&{52}\cr{-44}&{-6}&{-13}\cr
{-572}&{-65}&{-170}\cr
\endpmatrix$.}

\vskip5pt

\noindent
$N^\prime =11\ $ $d=130=2\cdot 5\cdot 13,\,\text{odd},\,\eta=3,\,h=0$:
$\langle 2 \rangle\oplus \langle -26 \rangle \oplus \langle -10 \rangle
(0,1/2,1/2) $. It is equivariantly equivalent to
$d=65,\,\eta=0,\,h=0$: $\langle 1 \rangle\oplus \langle -13 \rangle 
\oplus \langle -5 \rangle$.

\vskip5pt

\vbox{\noindent
$N^\prime=12\ $ $d=210=2\cdot 3\cdot 5\cdot 7,\,\text{odd},\,
\eta=5 ,\,h=2$:
$\langle 2 \rangle \oplus \langle -30 \rangle \oplus \langle -14 \rangle
(0,1/2,1/2)$
($=\widetilde{(105,6,h=0)}$). 
\nobreak
\newline
$\rho=(15,-2,-5)$,\ $\rho^2=-20$; 
\newline 
$\Gamma_+=
\left(\smallmatrix{7}&{{-3}\over{2}}&{{-3}\over{2}}\cr
{15}&{{-7}\over{2}}&{{-5}\over{2}}\cr{14}&{{-7}\over{2}}&{{-3}\over{2}}\cr
{15}&{-4}&{0}\cr{0}&{0}&{1}\cr{0}&{1}&{0}\cr{21}&{0}&{-8}\cr
{12}&{{-1}\over{2}}&{{-9}\over{2}}\cr{37}&{{-5}\over{2}}&{{-27}\over{2}}\cr
\endsmallmatrix\right),  
\hskip20pt 
G(\Gamma_+)=
\left(\smallmatrix{-1}&{0}&{7}&{30}&{21}&{45}&{126}&{51}&{122}\cr
{0}&{-5}&{0}&{30}&{35}&{105}&{350}&{150}&{375}\cr
{7}&{0}&{-7}&{0}&{21}&{105}&{420}&{189}&{490}\cr
{30}&{30}&{0}&{-30}&{0}&{120}&{630}&{300}&{810}\cr
{21}&{35}&{21}&{0}&{-14}&{0}&{112}&{63}&{189}\cr
{45}&{105}&{105}&{120}&{0}&{-30}&{0}&{15}&{75}\cr
{126}&{350}&{420}&{630}&{112}&{0}&{-14}&{0}&{42}\cr
{51}&{150}&{189}&{300}&{63}&{15}&{0}&{-3}&{0}\cr
{122}&{375}&{490}&{810}&{189}&{75}&{42}&{0}&{-1}\cr
\endsmallmatrix\right)$;
\nobreak
\newline
$C_1=
\pmatrix{167}&{420}&{336}\cr{-28}&{-71}&{-56}\cr{-48}&{-120}&{-97}\cr
\endpmatrix,
\hskip30pt
C_2=
\pmatrix{53}&{90}&{126}\cr{-6}&{-11}&{-14}\cr{-18}&{-30}&{-43}\cr
\endpmatrix$.}

\vskip5pt 

\vbox{\noindent
$N^\prime=13\ $ $d=330=2\cdot 3\cdot 5\cdot 11,\,\text{odd},\,\eta=1 ,\,h=0$:
$\langle 2 \rangle \oplus \langle -30 \rangle \oplus \langle -22 \rangle
(0,1/2,1/2)$
($=\widetilde{(165,6,h=0)}$). 
\nobreak
\newline
$\rho=(30,-11/2,-15/2)$,\ $\rho^2=-345$; 
\newline 
$\Gamma_+\hskip-2pt=\hskip-2pt 
\left(\smallmatrix
{100}&{{-47}\over{2}}&{{-25}\over{2}}\cr{231}&{-55}&{-27}\cr
{6}&{{-3}\over{2}}&{{-1}\over{2}}\cr{15}&{-4}&{0}\cr{0}&{0}&{1}\cr
{0}&{1}&{0}\cr{33}&{0}&{-10}\cr{15}&{{-1}\over{2}}&{{-9}\over{2}}\cr
{17}&{-1}&{-5}\cr{320}&{{-47}\over{2}}&{{-185}\over{2}}\cr
\endsmallmatrix\right),  
\hskip2pt 
G(\Gamma_+)\hskip-2pt=\hskip-2pt
\left(\smallmatrix
{-5}&{0}&{5}&{180}&{275}&{705}&{3850}&{1410}&{1320}&{21995}\cr
{0}&{-66}&{0}&{330}&{594}&{1650}&{9306}&{3432}&{3234}&{54120}\cr
{5}&{0}&{-1}&{0}&{11}&{45}&{286}&{108}&{104}&{1765}\cr
{180}&{330}&{0}&{-30}&{0}&{120}&{990}&{390}&{390}&{6780}\cr
{275}&{594}&{11}&{0}&{-22}&{0}&{220}&{99}&{110}&{2035}\cr
{705}&{1650}&{45}&{120}&{0}&{-30}&{0}&{15}&{30}&{705}\cr
{3850}&{9306}&{286}&{990}&{220}&{0}&{-22}&{0}&{22}&{770}\cr
{1410}&{3432}&{108}&{390}&{99}&{15}&{0}&{-3}&{0}&{90}\cr
{1320}&{3234}&{104}&{390}&{110}&{30}&{22}&{0}&{-2}&{0}\cr
{21995}&{54120}&{1765}&{6780}&{2035}&{705}&{770}&{90}&{0}&{-5}\cr
\endsmallmatrix\right)$;
\nobreak
\newline
$S=\pmatrix 
{285}&{990}&{418}\cr{-22}&{-76}&{-33}\cr{-82}&{-285}&{-120}\cr
\endpmatrix$.}

\vskip3pt 

\noindent
$N^\prime =14\ $ $d=330=2\cdot 3\cdot 5\cdot 11,\,\text{odd},\,\eta=5,\,h=2$:
$\langle 330 \rangle\oplus \langle -2 \rangle \oplus \langle -2 \rangle
(0,1/2,1/2)$.
It is equivariantly equivalent to
$d=165,\,\eta=2,\,h=2$: $\langle 165 \rangle\oplus \langle -1 \rangle 
\oplus \langle -1 \rangle $.

\vskip5pt

\vbox{\noindent
$N^\prime=15\ $ $d=330=2\cdot 3\cdot 5\cdot 11,\,\text{odd},\,\eta=7 ,\,h=2$:
$\langle 10 \rangle \oplus \langle -3 \rangle \oplus \langle -11\rangle$
($=\widetilde{(165,0,h=0)}$). 
\nobreak
\newline
$\rho=(33,-44,-24)$,\ $\rho^2=-1254$; 
\newline 
$\Gamma_+=\left(\smallmatrix
{231}&{-385}&{-90}\cr{1}&{-2}&{0}\cr{0}&{0}&{1}\cr{0}&{1}&{0}\cr
{1}&{0}&{-1}\cr{99}&{-55}&{-90}\cr
\endsmallmatrix\right),  
\hskip20pt 
G(\Gamma_+)=
\left(\smallmatrix{-165}&{0}&{990}&{1155}&{1320}&{76065}\cr
{0}&{-2}&{0}&{6}&{10}&{660}\cr{990}&{0}&{-11}&{0}&{11}&{990}\cr
{1155}&{6}&{0}&{-3}&{0}&{165}\cr{1320}&{10}&{11}&{0}&{-1}&{0}\cr
{76065}&{660}&{990}&{165}&{0}&{-165}\cr
\endsmallmatrix\right)$;
\nobreak
\newline
$C_1=
\pmatrix{39}&{18}&{22}\cr{-60}&{-28}&{-33}\cr{-20}&{-9}&{-12}\cr
\endpmatrix,
\hskip30pt
C_2=
\pmatrix{71}&{18}&{66}\cr{-60}&{-16}&{-55}\cr{-60}&{-15}&{-56}\cr
\endpmatrix$.}

\vskip5pt

\noindent
$N^\prime =16\ $ $d=390=2\cdot 3\cdot 5\cdot 13,\,\text{odd},\,\eta=6,\,h=2$:
$\langle 10 \rangle\oplus \langle -26 \rangle \oplus \langle -6 \rangle
(1/2,0,1/2)$.
It is equivariantly equivalent to
$d=195,\,\eta=1,\,h=2$: 
$\langle 5 \rangle\oplus \langle -13 \rangle \oplus \langle -3 \rangle $.

\vskip5pt

\noindent
$N^\prime =17\ $ $d=390=2\cdot 3\cdot 5\cdot 13 ,\,\text{odd},\,\eta=0,\,h=2$:
$\langle 2 \rangle\oplus \langle -130 \rangle \oplus \langle -6 \rangle
(1/2,0,1/2)$.
It is equivariantly equivalent to
$d=195,\,\eta=7,\,h=2$: 
$\langle 1 \rangle\oplus \langle -65 \rangle \oplus \langle -3 \rangle$.

\vskip5pt

\vbox{\noindent
$N^\prime=18\ $ $d=462=2\cdot 3\cdot 7\cdot 11,\,\text{odd},\,\eta=2 ,\,h=0$:
$\langle 3 \rangle \oplus \langle -7 \rangle \oplus \langle-22\rangle$
($=\widetilde{( 231,7,h=0)}$). 
\nobreak
\newline
$\rho=(22,-11,-7)$,\ $\rho^2=-473$; 
\newline 
$\Gamma_+=
\left(\smallmatrix{11}&{-7}&{-1}\cr{3}&{-2}&{0}\cr{0}&{0}&{1}\cr
{0}&{1}&{0}\cr{8}&{0}&{-3}\cr{385}&{-44}&{-140}\cr
{80}&{-10}&{-29}\cr{235}&{-31}&{-85}\cr
\endsmallmatrix\right),  
\hskip20pt 
G(\Gamma_+)=
\left(\smallmatrix{-2}&{1}&{22}&{49}&{198}&{7469}&{1512}&{4366}\cr
{1}&{-1}&{0}&{14}&{72}&{2849}&{580}&{1681}\cr
{22}&{0}&{-22}&{0}&{66}&{3080}&{638}&{1870}\cr
{49}&{14}&{0}&{-7}&{0}&{308}&{70}&{217}\cr
{198}&{72}&{66}&{0}&{-6}&{0}&{6}&{30}\cr
{7469}&{2849}&{3080}&{308}&{0}&{-77}&{0}&{77}\cr
{1512}&{580}&{638}&{70}&{6}&{0}&{-2}&{0}\cr
{4366}&{1681}&{1870}&{217}&{30}&{77}&{0}&{-2}\cr
\endsmallmatrix\right)$;
\nobreak
\newline
$S=\pmatrix{349}&{504}&{308}\cr{-48}&{-69}&{-44}\cr{-126}&{-182}&{-111}\cr
\endpmatrix$.}

\vskip5pt 

\noindent
$N^\prime =19\ $ $d=462=2\cdot 3\cdot 7\cdot 11,\,\text{odd},\,
\eta=6,\,h=2$:
$\langle 22 \rangle\oplus \langle -14 \rangle \oplus \langle -6 \rangle
(0,1/2,1/2)$.
It is equivariantly equivalent to
$d=231,\,\eta=3,\,h=2$: 
$\langle 11 \rangle\oplus \langle -7 \rangle \oplus \langle -3 \rangle
 $.

\vskip5pt

\noindent
$N^\prime =20\ $ $d=510=2\cdot 3\cdot 5\cdot 17,\,\text{odd},\,\eta=1,\,h=2$:
$\langle 2 \rangle\oplus \langle -34 \rangle \oplus \langle -30 \rangle
(1/2,0,1/2)$.
It is equivariantly equivalent to
$d=255,\,\eta=2,\,h=2$: $
\langle 1 \rangle\oplus \langle -17 \rangle \oplus \langle -15 \rangle $.

\vskip5pt

\vbox{\noindent
$N^\prime=21\ $ $d=510=2\cdot 3\cdot 5\cdot 17,\,\text{odd},\,\eta=5,\,h=0$:
$\langle 6 \rangle \oplus \langle -170 \rangle \oplus \langle -2 \rangle
(0,1/2,1/2)$
($=\widetilde{(255,6,h=0)}$). 
\nobreak
\newline
$\rho=(35,-6,-30)$,\ $\rho^2=-570$; 
\newline 
$\Gamma_+=
\left(\smallmatrix{19}&{{-7}\over{2}}&{{-13}\over{2}}\cr
{8}&{{-3}\over{2}}&{{-3}\over{2}}\cr{85}&{-16}&{0}\cr
{0}&{0}&{1}\cr{0}&{1}&{0}\cr{1}&{0}&{-2}\cr
{5}&{{-1}\over{2}}&{{-15}\over{2}}\cr{119}&{-15}&{-153}\cr
{19}&{{-5}\over{2}}&{{-47}\over{2}}\cr
\endsmallmatrix\right),  
\hskip20pt 
G(\Gamma_+)=
\left(\smallmatrix{-1}&{0}&{170}&{13}&{595}&{88}&{175}&{2652}&{373}\cr
{0}&{-3}&{0}&{3}&{255}&{42}&{90}&{1428}&{204}\cr
{170}&{0}&{-170}&{0}&{2720}&{510}&{1190}&{19890}&{2890}\cr
{13}&{3}&{0}&{-2}&{0}&{4}&{15}&{306}&{47}\cr
{595}&{255}&{2720}&{0}&{-170}&{0}&{85}&{2550}&{425}\cr
{88}&{42}&{510}&{4}&{0}&{-2}&{0}&{102}&{20}\cr
{175}&{90}&{1190}&{15}&{85}&{0}&{-5}&{0}&{5}\cr
{2652}&{1428}&{19890}&{306}&{2550}&{102}&{0}&{-102}&{0}\cr
{373}&{204}&{2890}&{47}&{425}&{20}&{5}&{0}&{-1}\cr
\endsmallmatrix\right)$;
\nobreak
\newline
$S=\pmatrix{263}&{1360}&{36}\cr{-36}&{-186}&{-5}\cr
{-312}&{-1615}&{-42}\cr
\endpmatrix$.}

\newpage 

\centerline{\bf Table 6}
\centerline{The list of main hyperbolic lattices of the rank $3$}
\centerline{with square-free determinant $d\le 100000$ and $h=2$}

\vskip20pt

{\settabs 12 \columns
\+$n=1 ,\,d=17 ,\,\eta=1 ,\,h=2$  &&&&&&&     
$\langle 51 \rangle\oplus A_2(1/3,1/3,-1/3)$ &&&&&& $hr$\cr
\+$n=2 ,\,d=19 ,\,\eta=1 ,\,h=2$  &&&&&&&      
$U\oplus\langle-19\rangle$  &&&&&& $hr$\cr
\+$n=3 ,\,d=23 ,\,\eta=1 ,\,h=2$  &&&&&&&      
$U\oplus \langle -23 \rangle$ &&&&&& $hr$\cr
\+$n=4 ,\,d=31 ,\,\eta=0 ,\,h=2$  &&&&&&&      
$\langle 31\rangle \oplus \langle -1 \rangle \oplus \langle-1\rangle$ 
&&&&&& $hr$\cr
\+$n=5 ,\,d=35=5\cdot 7 ,\,\eta=1 ,\,h=2$  &&&&&&&     
 $U\oplus \langle -35 \rangle$ &&&&&& $hr$\cr
\+$n=6 ,\,d=37 ,\,\eta=0 ,\,h=2$  &&&&&&&      
$U\oplus \langle-37\rangle$ &&&&&& $nr$\cr
\+$n=7 ,\,d=37 ,\,\eta=1 ,\,h=2$  &&&&&&&      
$\langle 1 \rangle\oplus \langle -74 \rangle\oplus 
\langle -2 \rangle(0,1/2,1/2)$ &&&&&& $hr$\cr
\+$n=8 ,\,d=39=3\cdot 13,\,\eta=1 ,\,h=2$  &&&&&&&      
$U\oplus \langle -39 \rangle$ &&&&&& $hr$\cr
\+$n=9 ,\,d=43 ,\,\eta=0 ,\,h=2$  &&&&&&&      
$\langle 43 \rangle\oplus \langle-1 \rangle\oplus 
\langle -1 \rangle$ &&&&&& $hr$\cr
\+$n=10 ,\,d=46=2\cdot 23 ,\,\eta=1 ,\,h=2$  &&&&&&&     
$U\oplus \langle -46 \rangle$ &&&&&& $hr$\cr
\+$n=11 ,\,d=47 ,\,\eta=0 ,\,h=2$  &&&&&&&     
$\langle 47\rangle\oplus \langle -1 \rangle\oplus 
\langle -1 \rangle$ &&&&&& $hr$\cr
\+$n=12 ,\,d=51=3\cdot 17 ,\,\eta=2 ,\,h=2$  &&&&&&&     
$U\oplus \langle -51 \rangle$  &&&&&& $nr$\cr
\+$n=13 ,\,d=53 ,\,\eta=0 ,\,h=2$  &&&&&&&    
$U\oplus \langle - 53 \rangle$ &&&&&& $nr$\cr
\+$n=14 ,\,d=53 ,\,\eta=1 ,\,h=2$  &&&&&&&     
$\langle 1 \rangle\oplus \langle -106 \rangle\oplus 
\langle -2 \rangle(0,1/2,1/2)$ &&&&&& $hr$\cr
\+$n=15 ,\,d=55=5\cdot 11 ,\,\eta=1 ,\,h=2$  &&&&&&&     
$\langle 11 \rangle\oplus \langle -10 \rangle\oplus 
\langle -2 \rangle(0,1/2,1/2)$ &&&&&& $hr$\cr
\+$n=16 ,\,d=57=3\cdot 19 ,\,\eta=3 ,\,h=2$  &&&&&&&    
$\langle 1 \rangle\oplus \langle -19 \rangle\oplus 
\langle -3 \rangle$ &&&&&& $hr$\cr
\+$n=17 ,\,d=58=2\cdot 29 ,\,\eta=1 ,\,h=2$  &&&&&&&     
$U\oplus \langle -58 \rangle$ &&&&&& $hr$\cr
\+$n=18 ,\,d=59 ,\,\eta=0 ,\,h=2$  &&&&&&&     
$\langle 59 \rangle\oplus \langle -1 \rangle\oplus 
\langle -1 \rangle$ &&&&&& $hr$\cr
\+$n=19 ,\,d=62=2\cdot 31 ,\,\eta=1 ,\,h=2$  &&&&&&&    
$U\oplus \langle -62 \rangle$ &&&&&& $hr$\cr
\+$n=20 ,\,d=65=5\cdot 13 ,\,\eta=1 ,\,h=2$  &&&&&&&    
$\langle 2 \rangle\oplus \langle -13 \rangle 
\oplus \langle -10 \rangle(1/2,0,1/2)$ &&&&&& $hr$\cr
\+$n=21 ,\,d=69=3\cdot 23 ,\,\eta=0 ,\,h=2$  &&&&&&&     
$\langle 46 \rangle\oplus \langle -6 \rangle\oplus 
\langle -1\rangle(1/2,1/2,0)$ &&&&&& $hr$\cr
\+$n=22 ,\,d=69=3\cdot 23 ,\,\eta=2 ,\,h=2$  &&&&&&&     
$U\oplus \langle -69 \rangle$ &&&&&& $nr$\cr
\+$n=23 ,\,d=70=2\cdot 5\cdot 7 ,\,\eta=0 ,\,h=2$  &&&&&&&     
$U\oplus \langle -70 \rangle$ &&&&&& $hr$\cr
\+$n=24 ,\,d=77=7\cdot 11 ,\,\eta=0 ,\,h=2$  &&&&&&&     
$\langle 14 \rangle\oplus \langle -22 \rangle\oplus 
\langle -1 \rangle(1/2,1/2,0)$ &&&&&& $hr$\cr
\+$n=25 ,\,d=79 ,\,\eta=0 ,\,h=2$  &&&&&&&     
$\langle 79 \rangle \oplus \langle -1 \rangle\oplus 
\langle -1 \rangle$ &&&&&& $nr$\cr
\+$n=26 ,\,d=85=5\cdot 17 ,\,\eta=2 ,\,h=2$  &&&&&&&     
$\langle 1 \rangle \oplus \langle -170 \rangle \oplus 
\langle -2 \rangle(0,1/2,1/2)$ &&&&&& $hr$\cr
\+$n=27 ,\,d=85=5\cdot 17 ,\,\eta=3 ,\,h=2$  &&&&&&&     
$U\oplus\langle -85 \rangle$ &&&&&& $nr$\cr
\+$n=28 ,\,d=87=3\cdot 29 ,\,\eta=0 ,\,h=2$  &&&&&&&     
$\langle 3 \rangle\oplus\langle -29 \rangle
\oplus\langle -1 \rangle$ &&&&&& $hr$\cr
\+$n=29 ,\,d=91=7\cdot 13 ,\,\eta=0 ,\,h=2$  &&&&&&&     
$\langle 7 \rangle\oplus\langle -13 \rangle\oplus
\langle -1 \rangle$ &&&&&& $hr$\cr
\+$n=30 ,\,d=91=7\cdot 13 ,\,\eta=1 ,\,h=2$  &&&&&&&    
$\langle 1\rangle\oplus\langle -13 \rangle\oplus
\langle -7 \rangle$ &&&&&& $nr$\cr
\+$n=31 ,\,d=93=3\cdot 31,\,\eta=3 ,\,h=2$  &&&&&&&     
$\langle 1 \rangle\oplus\langle -31 \rangle\oplus
\langle -3 \rangle$ &&&&&& $hr$\cr
\+$n=32 ,\,d=106=2\cdot 53 ,\,\eta=1 ,\,h=2$  &&&&&&&    
$U\oplus\langle -106 \rangle$ &&&&&& $nr$\cr
\+$n=33 ,\,d=107 ,\,\eta=0 ,\,h=2$  &&&&&&&    
$\langle 107 \rangle\oplus\langle -1 \rangle\oplus
\langle -1 \rangle$ &&&&&& $nr$\cr
\+$n=34 ,\,d=113 ,\,\eta=0 ,\,h=2$  &&&&&&&    
$U\oplus\langle -113 \rangle$&&&&&& $nr$\cr
\+$n=35 ,\,d=114=2\cdot 3\cdot 19 ,\,\eta=1 ,\,h=2$  &&&&&&&    
$\langle 6 \rangle\oplus\langle -38 \rangle\oplus
\langle -2 \rangle(0,1/2,1/2)$ &&&&&& $hr$\cr
\+$n=36 ,\,d=115=5\cdot 23 ,\,\eta=0 ,\,h=2$  &&&&&&&    
$\langle 23 \rangle\oplus\langle -5 \rangle\oplus
\langle -1 \rangle$ &&&&&& $hr$\cr
\+$n=37 ,\,d=115=5\cdot 23 ,\,\eta=3 ,\,h=2$  &&&&&&&   
$\langle 115 \rangle\oplus\langle -1 \rangle\oplus
\langle -1 \rangle$ &&&&&& $hr$\cr
\+$n=38 ,\,d=119=7\cdot 17 ,\,\eta=0 ,\,h=2$  &&&&&&&    
$\langle 7 \rangle\oplus\langle -17 \rangle\oplus
\langle -1 \rangle$ &&&&&& $hr$\cr
\+$n=39 ,\,d=122=2\cdot 61 ,\,\eta=1 ,\,h=2$  &&&&&&&   
$U\oplus\langle -122 \rangle$  &&&&&& $nr$\cr
\+$n=40 ,\,d=123=3\cdot 41 ,\,\eta=0 ,\,h=2$  &&&&&&&    
$\langle 3 \rangle\oplus\langle -41 \rangle\oplus
\langle -1 \rangle$ &&&&&& $hr$\cr
\+$n=41 ,\,d=123=3\cdot 41 ,\,\eta=1 ,\,h=2$  &&&&&&&    
$\langle 41 \rangle\oplus\langle -3 \rangle\oplus
\langle -1 \rangle$ &&&&&& $nr$\cr
\+$n=42 ,\,d=123=3\cdot 41 ,\,\eta=3 ,\,h=2$  &&&&&&&    
$\langle 123 \rangle\oplus\langle -1 \rangle\oplus
\langle -1 \rangle$ &&&&&& $hr$\cr
\+$n=43 ,\,d=129=3\cdot 43 ,\,\eta=0 ,\,h=2$  &&&&&&&    
$\langle 3 \rangle\oplus\langle -86 \rangle\oplus
\langle -2 \rangle(0,1/2,1/2)$ &&&&&& $hr$\cr
\+$n=44 ,\,d=129=3\cdot 43 ,\,\eta=2 ,\,h=2$  &&&&&&&   
$\langle 3 \rangle\oplus\langle -43 \rangle\oplus
\langle -1 \rangle$ &&&&&& $nr$\cr
\+$n=45 ,\,d=133=7\cdot 19 ,\,\eta=1 ,\,h=2$  &&&&&&&    
$\langle 19\rangle\oplus\langle -7 \rangle\oplus
\langle -1 \rangle$ &&&&&& $nr$\cr
\+$n=46 ,\,d=133=7\cdot 19 ,\,\eta=3 ,\,h=2$  &&&&&&&   
$\langle 1 \rangle \oplus\langle -19 \rangle\oplus
\langle -7\rangle$ &&&&&& $hr$\cr
\+$n=47 ,\,d=138=2\cdot 3\cdot 23 ,\,\eta=0 ,\,h=2$  &&&&&&&  
$\langle 46 \rangle\oplus\langle -6 \rangle\oplus
\langle -2 \rangle(0,1/2,1/2)$ &&&&&& $hr$\cr
\+$n=48 ,\,d=141=3\cdot 47 ,\,\eta=0 ,\,h=2$  &&&&&&&  
$\langle 94 \rangle\oplus\langle -6 \rangle\oplus
\langle-1\rangle(1/2,1/2,0)$ &&&&&& $nr$\cr
\+$n=49 ,\,d=143=11\cdot 13 ,\,\eta=3 ,\,h=2$  &&&&&&&  
$\langle 143 \rangle\oplus\langle -1 \rangle
\oplus\langle -1 \rangle$ &&&&&& $hr$\cr
\+$n=50 ,\,d=145=5\cdot 29 ,\,\eta=3 ,\,h=2$  &&&&&&&  
$\langle 1 \rangle\oplus\langle -290 \rangle
\oplus\langle -2 \rangle(0,1/2,1/2)$ &&&&&& $nr$\cr
\+$n=51 ,\,d=155=5\cdot 31 ,\,\eta=0 ,\,h=2$  &&&&&&&   
$\langle 31 \rangle\oplus\langle -5 \rangle\oplus
\langle -1 \rangle$ &&&&&& $nr$\cr
\+$n=52 ,\,d=159=3\cdot 53 ,\,\eta=0 ,\,h=2$  &&&&&&&  
$\langle 3 \rangle\oplus\langle -53 \rangle\oplus
\langle -1 \rangle$ &&&&&& $nr$\cr
\+$n=53 ,\,d=159=3\cdot 53 ,\,\eta=3 ,\,h=2$  &&&&&&&  
$\langle 159 \rangle\oplus\langle -1 \rangle\oplus
\langle -1 \rangle$ &&&&&& $hr$\cr
\+$n=54 ,\,d=165=3\cdot 5\cdot 11 ,\,\eta=2 ,\,h=2$  &&&&&&&  
$\langle 165 \rangle\oplus\langle -1 \rangle\oplus
\langle -1 \rangle$ &&&&&& $hr$\cr
\+$n=55 ,\,d=174=2\cdot 3\cdot 29 ,\,\eta=2 ,\,h=2$  &&&&&&&  
$\langle 2 \rangle\oplus\langle -58 \rangle\oplus
\langle -6 \rangle(0,1/2,1/2)$ &&&&&& $hr$\cr
\+$n=56 ,\,d=177=3\cdot 59 ,\,\eta=1 ,\,h=2$  &&&&&&&  
$\langle 59 \rangle\oplus\langle -3 \rangle
\oplus\langle -1 \rangle$ &&&&&& $nr$\cr
\+$n=57 ,\,d=182=2\cdot 7\cdot 13 ,\,\eta=0 ,\,h=2$  &&&&&&&  
$U\oplus\langle -182 \rangle$ &&&&&& $nr$\cr
\+$n=58 ,\,d=183=3\cdot 61 ,\,\eta=0 ,\,h=2$  &&&&&&&  
$\langle 3 \rangle\oplus\langle -61 \rangle
\oplus\langle -1 \rangle$ &&&&&& $hr$\cr
\+$n=59 ,\,d=185=5\cdot 37 ,\,\eta=0 ,\,h=2$  &&&&&&&   
$\langle 1 \rangle\oplus\langle -37 \rangle
\oplus\langle -5 \rangle$ &&&&&& $nr$\cr
\+$n=60 ,\,d=191 ,\,\eta=0 ,\,h=2$  &&&&&&&   
$\langle 191 \rangle\oplus\langle -1 \rangle
\oplus\langle -1 \rangle$ &&&&&& $nr$\cr
\+$n=61 ,\,d=195=3\cdot 5\cdot 13 ,\,\eta=0 ,\,h=2$  &&&&&&&   
$\langle 3 \rangle\oplus\langle -13 \rangle
\oplus\langle -5 \rangle$ &&&&&& $hr$\cr
\+$n=62 ,\,d=195=3\cdot 5\cdot 13 ,\,\eta=1 ,\,h=2$  &&&&&&&  
$\langle 5 \rangle\oplus\langle -13 \rangle\oplus
\langle -3 \rangle$ &&&&&& $hr$\cr
\+$n=63 ,\,d=195=3\cdot 5\cdot 13 ,\,\eta=7 ,\,h=2$  &&&&&&&   
$\langle 1 \rangle\oplus\langle -65 \rangle\oplus
\langle -3 \rangle$ &&&&&& $hr$\cr
\+$n=64 ,\,d=205=5\cdot 41 ,\,\eta=0 ,\,h=2$  &&&&&&&  
$U\oplus \langle -205\rangle$ &&&&&& $nr$\cr
\+$n=65 ,\,d=210=2\cdot 3\cdot 5\cdot 7 ,\,\eta=1 ,\,h=2$  &&&&&&&
$\langle 14\rangle \oplus \langle -30\rangle \oplus
\langle -2\rangle (0,1/2,1/2)$ &&&&&& $hr$\cr
\+$n=66 ,\,d=210=2\cdot 3\cdot 5\cdot 7 ,\,\eta=7 ,\,h=2$  &&&&&&& 
$U\oplus\langle -210\rangle $ &&&&&& $hr$\cr
\+$n=67 ,\,d=213=3\cdot 71 ,\,\eta=3 ,\,h=2$  &&&&&&&
$\langle 1\rangle \oplus \langle  -71 \rangle  \oplus 
\langle  -3 \rangle $ &&&&&& $hr$\cr
\+$n=68 ,\,d=217=7\cdot 31 ,\,\eta=3 ,\,h=2$  &&&&&&&
$\langle  1 \rangle \oplus \langle  -31 \rangle \oplus
\langle -7\rangle $ &&&&&& $nr$\cr
\+$n=69 ,\,d=218=2\cdot 109 ,\,\eta=1 ,\,h=2$  &&&&&&&
$U\oplus \langle -218\rangle $ &&&&&& $nr$\cr
\+$n=70 ,\,d=221=13\cdot 17 ,\,\eta=0 ,\,h=2$  &&&&&&&
$U\oplus \langle -221\rangle $ &&&&&& $nr$\cr
\+$n=71 ,\,d=221=13\cdot 17 ,\,\eta=1 ,\,h=2$  &&&&&&&
$\langle 2\rangle \oplus \langle -26\rangle \oplus 
\langle -17\rangle (1/2,1/2,0)$ &&&&&& $hr$\cr
\+$n=72 ,\,d=222=2\cdot 3\cdot 37 ,\,\eta=2 ,\,h=2$  &&&&&&&
$U\oplus \langle -222\rangle $&&&&&& $nr$\cr
\+$n=73 ,\,d=230=2\cdot 5\cdot 23 ,\,\eta=0 ,\,h=2$  &&&&&&&
$U\oplus \langle  -230\rangle $ &&&&&& $nr$\cr
\+$n=74 ,\,d=231=3\cdot 7\cdot 11 ,\,\eta=3 ,\,h=2$  &&&&&&&
$\langle  11\rangle \oplus \langle -7\rangle  \oplus 
\langle -3\rangle $ &&&&&& $hr$\cr
\+$n=75 ,\,d=231=3\cdot 7\cdot 11 ,\,\eta=4 ,\,h=2$  &&&&&&&
$\langle  7\rangle \oplus \langle  -33\rangle  \oplus 
\langle -1\rangle $ &&&&&& $hr$\cr
\+$n=76 ,\,d=235=5\cdot 47 ,\,\eta=0 ,\,h=2$  &&&&&&&
$\langle 47\rangle \oplus \langle -5\rangle \oplus 
\langle  -1\rangle $ &&&&&& $nr$\cr
\+$n=77 ,\,d=239 ,\,\eta=0 ,\,h=2$  &&&&&&&
$\langle 239\rangle \oplus \langle -1\rangle \oplus 
\langle -1\rangle $ &&&&&& $nr$\cr
\+$n=78 ,\,d=255=3\cdot 5\cdot 17 ,\,\eta=2 ,\,h=2$  &&&&&&&
$\langle 1\rangle \oplus \langle -17\rangle \oplus 
\langle -15\rangle $ &&&&&& $hr$\cr
\+$n=79 ,\,d=255=3\cdot 5\cdot 17 ,\,\eta=3 ,\,h=2$  &&&&&&&
$\langle 15\rangle \oplus \langle -17\rangle \oplus 
\langle -1\rangle $ &&&&&& $hr$\cr
\+$n=80 ,\,d=258=2\cdot 3\cdot 43 ,\,\eta=2 ,\,h=2$  &&&&&&&
$U\oplus \langle -258\rangle $ &&&&&& $nr$\cr
\+$n=81 ,\,d=273=3\cdot 7\cdot 13 ,\,\eta=2 ,\,h=2$  &&&&&&&
$\langle 3\rangle \oplus \langle -13\rangle 
\oplus \langle -7\rangle $ &&&&&& $nr$\cr
\+$n=82 ,\,d=282=2\cdot 3\cdot 47 ,\,\eta=0 ,\,h=2$  &&&&&&&
$\langle 94\rangle \oplus \langle -6\rangle \oplus 
\langle -2\rangle (0,1/2,1/2)$ &&&&&& $nr$\cr
\+$n=83 ,\,d=287=7\cdot 41 ,\,\eta=3 ,\,h=2$  &&&&&&&
$\langle 287\rangle \oplus \langle -1\rangle \oplus 
\langle -1\rangle $ &&&&&& $nr$\cr
\+$n=84 ,\,d=290=2\cdot 5\cdot 29 ,\,\eta=0 ,\,h=2$  &&&&&&&
$\langle 870\rangle \oplus A_2(1/3,1/3,-1/3)$ &&&&&& $nr$\cr
\+$n=85 ,\,d=290=2\cdot 5\cdot 29 ,\,\eta=3 ,\,h=2$  &&&&&&&
$U\oplus \langle  -290\rangle $ &&&&&& $nr$\cr
\+$n=86 ,\,d=301=7\cdot 43 ,\,\eta=1 ,\,h=2$  &&&&&&&
$U\oplus \langle  -301\rangle $ &&&&&& $nr$\cr
\+$n=87 ,\,d=301=7\cdot 43 ,\,\eta=3 ,\,h=2$  &&&&&&&
$\langle  1\rangle \oplus \langle  -43\rangle  \oplus 
\langle  -7\rangle $ &&&&&& $nr$\cr
\+$n=88 ,\,d=303=3\cdot 101 ,\,\eta=3 ,\,h=2$  &&&&&&&
$\langle  303)\oplus \langle  -1\rangle \oplus 
\langle  -1\rangle $ &&&&&& $nr$\cr
\+$n=89 ,\,d=309=3\cdot 103 ,\,\eta=0 ,\,h=2$  &&&&&&&
$\langle  2\rangle \oplus \langle  -618\rangle  
\oplus \langle  -1\rangle (1/2,1/2,0)$ &&&&&& $nr$\cr
\+$n=90 ,\,d=314=2\cdot 157 ,\,\eta=1 ,\,h=2$  &&&&&&&
$U\oplus \langle  -314\rangle $ &&&&&& $nr$\cr
\+$n=91 ,\,d=319=11\cdot 29 ,\,\eta=3 ,\,h=2$  &&&&&&&
$\langle  319\rangle \oplus \langle  -1\rangle \oplus 
\langle  -1\rangle $ &&&&&& $nr$\cr
\+$n=92 ,\,d=323=17\cdot 19 ,\,\eta=0 ,\,h=2$  &&&&&&&
$\langle  19\rangle \oplus \langle  -17\rangle \oplus 
\langle  -1\rangle $ &&&&&& $nr$\cr
\+$n=93 ,\,d=329=7\cdot 47 ,\,\eta=3 ,\,h=2$  &&&&&&&
$\langle  1\rangle \oplus \langle  -47\rangle  \oplus 
\langle  -7\rangle $ &&&&&& $nr$\cr
\+$n=94 ,\,d=330=2\cdot 3\cdot 5\cdot 11 ,\,\eta=6 ,\,h=2$  &&&&&&&
$\langle  22\rangle \oplus \langle  -10\rangle \oplus 
\langle  -6\rangle (0,1/2,1/2)$ &&&&&& $hr$\cr
\+$n=95 ,\,d=345=3\cdot 5\cdot 23 ,\,\eta=2 ,\,h=2$  &&&&&&&
$\langle  3\rangle \oplus \langle  -115\rangle \oplus 
\langle  -1\rangle $ &&&&&& $nr$\cr
\+$n=96 ,\,d=345=3\cdot 5\cdot 23 ,\,\eta=4 ,\,h=2$  &&&&&&&
$\langle  1\rangle \oplus \langle  -69\rangle \oplus 
\langle  -5\rangle $ &&&&&& $nr$\cr
\+$n=97 ,\,d=345=3\cdot 5\cdot 23 ,\,\eta=5 ,\,h=2$  &&&&&&&
$\langle  5\rangle \oplus \langle  -23\rangle \oplus 
\langle  -3\rangle $ &&&&&& $nr$\cr
\+$n=98 ,\,d=357=3\cdot 7\cdot 17 ,\,\eta=6 ,\,h=2$  &&&&&&&
$\langle  1\rangle \oplus \langle  -51\rangle \oplus 
\langle  -7\rangle $ &&&&&& $nr$\cr
\+$n=99 ,\,d=366=2\cdot 3\cdot 61 ,\,\eta=2 ,\,h=2$  &&&&&&&
$U\oplus \langle  -366\rangle $ &&&&&& $nr$\cr
\+$n=100 ,\,d=374=2\cdot 11\cdot 17 ,\,\eta=0 ,\,h=2$  &&&&&&&
$\langle  2\rangle \oplus \langle  -34\rangle \oplus
 \langle  -22\rangle (0,1/2,1/2)$ &&&&&& $nr$\cr
\+$n=101 ,\,d=390=2\cdot 3\cdot 5\cdot 13 ,\,\eta=0 ,\,h=2$  &&&&&&&
$\langle  130\rangle \oplus \langle  -6\rangle \oplus 
\langle  -2\rangle (0,1/2,1/2)$ &&&&&& $hr$\cr
\+$n=102 ,\,d=402=2\cdot 3\cdot 67 ,\,\eta=2 ,\,h=2$  &&&&&&&
$U\oplus \langle  -402\rangle $ &&&&&& $nr$\cr
\+$n=103 ,\,d=406=2\cdot 7\cdot 29 ,\,\eta=3 ,\,h=2$  &&&&&&&
$U\oplus \langle  -406)$ &&&&&& $nr$\cr
\+$n=104 ,\,d=410=2\cdot 5\cdot 41 ,\,\eta=2 ,\,h=2$  &&&&&&&
$\langle  1230\rangle \oplus A_2(1/3,1/3,-1/3)$ &&&&&& $nr$\cr
\+$n=105 ,\,d=426=2\cdot 3\cdot 71 ,\,\eta=0 ,\,h=2$  &&&&&&&
$\langle  142\rangle \oplus \langle  -6\rangle \oplus 
\langle  -2\rangle (0,1/2,1/2)$ &&&&&& $nr$\cr
\+$n=106 ,\,d=429=3\cdot 11\cdot 13 ,\,\eta=0 ,\,h=2$  &&&&&&&
$\langle  13\rangle \oplus \langle  -22\rangle \oplus 
\langle  -6\rangle (0,1/2,1/2)$ &&&&&& $nr$\cr
\+$n=107 ,\,d=429=3\cdot 11\cdot 13 ,\,\eta=5 ,\,h=2$  &&&&&&&
$\langle  1\rangle \oplus \langle  -143\rangle \oplus 
\langle  -3\rangle $ &&&&&& $hr$\cr
\+$n=108 ,\,d=429=3\cdot 11\cdot 13 ,\,\eta=6 ,\,h=2$  &&&&&&&
$\langle  2\rangle \oplus \langle  -26\rangle 
\oplus \langle  -33\rangle (1/2,1/2,0)$ &&&&&& $nr$\cr
\+$n=109 ,\,d=431= ,\,\eta=0 ,\,h=2$  &&&&&&&
$\langle  431\rangle \oplus \langle  -1\rangle \oplus 
\langle  -1\rangle $ &&&&&& $nr$\cr
\+$n=110 ,\,d=434=2\cdot 7\cdot 31 ,\,\eta=1 ,\,h=2$  &&&&&&&
$\langle  62\rangle \oplus \langle  -14\rangle \oplus 
\langle  -2\rangle (0,1/2,1/2)$ &&&&&& $nr$\cr
\+$n=111 ,\,d=438=2\cdot 3\cdot 73 ,\,\eta=0 ,\,h=2$  &&&&&&&
$U\oplus \langle  -438\rangle $ &&&&&& $nr$\cr
\+$n=112 ,\,d=455=5\cdot 7\cdot 13 ,\,\eta=3 ,\,h=2$  &&&&&&&
$\langle  35\rangle \oplus \langle  -13\rangle \oplus 
\langle  -1\rangle $ &&&&&& $nr$\cr
\+$n=113 ,\,d=455=5\cdot 7\cdot 13 ,\,\eta=6 ,\,h=2$  &&&&&&&
$\langle  91\rangle \oplus \langle  -5\rangle \oplus 
\langle  -1\rangle $ &&&&&& $nr$\cr
\+$n=114 ,\,d=462=2\cdot 3\cdot 7\cdot 11 ,\,\eta=0 ,\,h=2$  &&&&&&&
$\langle  14\rangle \oplus \langle  -22\rangle \oplus 
\langle  -6\rangle (1/2,1/2,0)$ &&&&&& $hr$\cr
\+$n=115 ,\,d=462=2\cdot 3\cdot 7\cdot 11 ,\,\eta=6 ,\,h=2$  &&&&&&&
$\langle  22\rangle \oplus \langle  -14\rangle \oplus 
\langle  -6\rangle (1/2,1/2,0)$ &&&&&& $hr$\cr
\+$n=116 ,\,d=465=3\cdot 5\cdot 31 ,\,\eta=1 ,\,h=2$  &&&&&&&
$\langle  1\rangle \oplus \langle  -93\rangle \oplus 
\langle  -5\rangle $&&&&&& $nr$\cr
\+$n=117 ,\,d=465=3\cdot 5\cdot 31 ,\,\eta=6 ,\,h=2$  &&&&&&&
$\langle  1\rangle \oplus \langle  -31\rangle \oplus 
\langle  -15\rangle $&&&&&& $nr$\cr
\+$n=118 ,\,d=469=7\cdot 67 ,\,\eta=1 ,\,h=2$  &&&&&&&
$U\oplus \langle  -469\rangle $ &&&&&& $nr$\cr
\+$n=119 ,\,d=469=7\cdot 67 ,\,\eta=3 ,\,h=2$  &&&&&&&
$\langle  1\rangle \oplus \langle  -67\rangle \oplus 
\langle  -7\rangle $ &&&&&& $nr$\cr
\+$n=120 ,\,d=470=2\cdot 5\cdot 47 ,\,\eta=0 ,\,h=2$  &&&&&&&
$U\oplus \langle  -470\rangle $ &&&&&& $nr$\cr
\+$n=121 ,\,d=471=3\cdot 157 ,\,\eta=0 ,\,h=2$  &&&&&&&
$\langle  3\rangle \oplus \langle  -157\rangle \oplus 
\langle  -1\rangle $ &&&&&& $nr$\cr
\+$n=122 ,\,d=474=2\cdot 3\cdot 79 ,\,\eta=0 ,\,h=2$  &&&&&&&
$U\oplus \langle  -474\rangle $ &&&&&& $nr$\cr
\+$n=123 ,\,d=483=3\cdot 7\cdot 23 ,\,\eta=2 ,\,h=2$  &&&&&&&
$\langle  3\rangle \oplus \langle -161\rangle \oplus 
\langle  -1\rangle $ &&&&&& $nr$\cr
\+$n=124 ,\,d=515=5\cdot 103 ,\,\eta=0 ,\,h=2$  &&&&&&&
$\langle  103\rangle \oplus \langle  -5\rangle \oplus 
\langle  -1\rangle $ &&&&&& $nr$\cr
\+$n=125 ,\,d=518=2\cdot 7\cdot 37 ,\,\eta=3 ,\,h=2$  &&&&&&&
$U\oplus \langle  -518\rangle $ &&&&&& $nr$\cr
\+$n=126 ,\,d=527=17\cdot 31 ,\,\eta=3 ,\,h=2$  &&&&&&&
$\langle  527\rangle \oplus \langle  -1\rangle \oplus 
\langle  -1\rangle $ &&&&&& $nr$\cr
\+$n=127 ,\,d=554=2\cdot 277 ,\,\eta=1 ,\,h=2$  &&&&&&&
$U\oplus \langle  -554\rangle $ &&&&&& $nr$\cr
\+$n=128 ,\,d=555=3\cdot 5\cdot 37 ,\,\eta=3 ,\,h=2$  &&&&&&&
$\langle  15\rangle \oplus \langle  -37\rangle \oplus 
\langle  -1\rangle $ &&&&&& $hr$\cr
\+$n=129 ,\,d=555=3\cdot 5\cdot 37 ,\,\eta=5 ,\,h=2$  &&&&&&&
$\langle  555\rangle \oplus \langle  -1\rangle \oplus 
\langle  -1\rangle $ &&&&&& $nr$\cr
\+$n=130 ,\,d=555=3\cdot 5\cdot 37 ,\,\eta=6 ,\,h=2$  &&&&&&&
$\langle  3\rangle \oplus \langle  -185\rangle \oplus 
\langle  -1\rangle $ &&&&&& $nr$\cr
\+$n=131 ,\,d=561=3\cdot 11\cdot 17 ,\,\eta=6 ,\,h=2$  &&&&&&&
$\langle  1\rangle \oplus \langle  -51\rangle \oplus 
\langle  -11\rangle $ &&&&&& $hr$\cr
\+$n=132 ,\,d=582=2\cdot 3\cdot 97 ,\,\eta=0 ,\,h=2$  &&&&&&&
$U\oplus \langle  -582\rangle $ &&&&&& $nr$\cr
\+$n=133 ,\,d=595=5\cdot 7\cdot 17 ,\,\eta=3 ,\,h=2$  &&&&&&&
$\langle  35\rangle \oplus \langle  -17\rangle \oplus 
\langle  -1\rangle $ &&&&&& $nr$\cr
\+$n=134 ,\,d=609=3\cdot 7\cdot 29 ,\,\eta=3 ,\,h=2$  &&&&&&&
$\langle  29\rangle \oplus \langle  -7\rangle \oplus 
\langle  -3\rangle $ &&&&&& $hr$\cr
\+$n=135 ,\,d=609=3\cdot 7\cdot 29 ,\,\eta=6 ,\,h=2$  &&&&&&&
$\langle  1\rangle \oplus \langle  -87\rangle \oplus 
\langle  -7\rangle $ &&&&&& $nr$\cr
\+$n=136 ,\,d=618=2\cdot 3\cdot 103 ,\,\eta=0 ,\,h=2$  &&&&&&&
$U\oplus \langle  -618\rangle $ &&&&&& $nr$\cr
\+$n=137 ,\,d=623=7\cdot 89 ,\,\eta=3 ,\,h=2$  &&&&&&&
$\langle  623\rangle \oplus \langle  -1\rangle \oplus 
\langle  -1\rangle $ &&&&&& $nr$\cr
\+$n=138 ,\,d=627=3\cdot 11\cdot 19 ,\,\eta=1 ,\,h=2$  &&&&&&&
$\langle  11\rangle \oplus \langle -57\rangle \oplus 
\langle  -1\rangle $ &&&&&& $nr$\cr
\+$n=139 ,\,d=627=3\cdot 11\cdot 19 ,\,\eta=3 ,\,h=2$  &&&&&&&
$\langle  1\rangle \oplus \langle  -57\rangle \oplus 
\langle  -11\rangle $ &&&&&& $nr$\cr
\+$n=140 ,\,d=663=3\cdot 13\cdot 17 ,\,\eta=3 ,\,h=2$  &&&&&&&
$\langle  6\rangle \oplus \langle  -26\rangle \oplus 
\langle  -17\rangle (1/2,1/2,0)$ &&&&&& $hr$\cr
\+$n=141 ,\,d=665=5\cdot 7\cdot 19 ,\,\eta=0 ,\,h=2$  &&&&&&&
$\langle  2\rangle \oplus \langle  -226\rangle \oplus 
\langle  -5\rangle (1/2,1/2,0)$ &&&&&& $nr$\cr
\+$n=142 ,\,d=690=2\cdot 3\cdot 5\cdot 23 ,\,\eta=1 ,\,h=2$  &&&&&&&
$\langle  46\rangle \oplus \langle  -30\rangle \oplus 
\langle  -2\rangle (0,1/2,1/2)$ &&&&&& $nr$\cr
\+$n=143 ,\,d=690=2\cdot 3\cdot 5\cdot 23 ,\,\eta=2 ,\,h=2$  &&&&&&&
$U\oplus \langle  -690\rangle $ &&&&&& $nr$\cr
\+$n=144 ,\,d=690=2\cdot 3\cdot 5\cdot 23 ,\,\eta=7 ,\,h=2$  &&&&&&&
$\langle  6\rangle \oplus \langle  -46\rangle \oplus 
\langle  -10\rangle (0,1/2,1/2)$ &&&&&& $nr$\cr
\+$n=145 ,\,d=705=3\cdot 5\cdot 47 ,\,\eta=6 ,\,h=2$  &&&&&&&
$\langle  1\rangle \oplus \langle  -47\rangle \oplus 
\langle  -15\rangle $ &&&&&& $hr$\cr
\+$n=146 ,\,d=715=5\cdot 11\cdot 13 ,\,\eta=0 ,\,h=2$  &&&&&&&
$\langle  11\rangle \oplus \langle  -13\rangle \oplus 
\langle  -5\rangle $ &&&&&& $nr$\cr
\+$n=147 ,\,d=715=5\cdot 11\cdot 13 ,\,\eta=6 ,\,h=2$  &&&&&&&
$\langle  143\rangle \oplus \langle  -5\rangle \oplus 
\langle  -1\rangle $ &&&&&& $nr$\cr
\+$n=148 ,\,d=741=3\cdot 13\cdot 19 ,\,\eta=0 ,\,h=2$  &&&&&&&
$\langle  2\rangle \oplus \langle  -1482\rangle \oplus 
\langle  -1\rangle (1/2,1/2,0)$ &&&&&& $nr$\cr
\+$n=149 ,\,d=741=3\cdot 13\cdot 19 ,\,\eta=6 ,\,h=2$  &&&&&&&
$\langle  1\rangle \oplus \langle  -114\rangle \oplus 
\langle  -26\rangle (0,1/2,1/2)$ &&&&&& $nr$\cr
\+$n=150 ,\,d=770=2\cdot 5\cdot 7\cdot 11 ,\,\eta=2 ,\,h=2$  &&&&&&&
$\langle  770\rangle \oplus \langle  -2\rangle \oplus 
\langle  -2\rangle (1/2,1/2,0)$ &&&&&& $nr$\cr
\+$n=151 ,\,d=770=2\cdot 5\cdot 7\cdot 11 ,\,\eta=7 ,\,h=2$  &&&&&&&
$\langle  10\rangle \oplus \langle -154\rangle \oplus 
\langle  -2\rangle (1/2,0,1/2)$ &&&&&& $nr$\cr
\+$n=152 ,\,d=791=7\cdot 113 ,\,\eta=3 ,\,h=2$  &&&&&&&
$\langle  2373\rangle \oplus A_2(1/3,1/3,-1/3)$ &&&&&& $nr$\cr
\+$n=153 ,\,d=794=2\cdot 397 ,\,\eta=1 ,\,h=2$  &&&&&&&
$U\oplus \langle  -794\rangle $ &&&&&& $nr$\cr
\+$n=154 ,\,d=798=2\cdot 3\cdot 7\cdot 19 ,\,\eta=0 ,\,h=2$  &&&&&&&
$\langle  2\rangle \oplus \langle  -38\rangle \oplus 
\langle  -42\rangle (0,1/2,1/2)$ &&&&&& $nr$\cr
\+$n=155 ,\,d=798=2\cdot 3\cdot 7\cdot 19 ,\,\eta=6 ,\,h=2$  &&&&&&&
$U\oplus \langle  -798\rangle $ &&&&&& $nr$\cr
\+$n=156 ,\,d=806=2\cdot 13\cdot 31 ,\,\eta=0 ,\,h=2$  &&&&&&&
$U\oplus \langle  -806\rangle $ &&&&&& $nr$\cr
\+$n=157 ,\,d=834=2\cdot 3\cdot 139 ,\,\eta=2 ,\,h=2$  &&&&&&&
$U\oplus \langle  -834\rangle $ &&&&&& $nr$\cr
\+$n=158 ,\,d=858=2\cdot 3\cdot 11\cdot 13 ,\,\eta=0 ,\,h=2$  &&&&&&&
$\langle  858\rangle \oplus \langle  -2\rangle \oplus 
\langle  -2\rangle (1/2,1/2,0)$ &&&&&& $nr$\cr
\+$n=159 ,\,d=890=2\cdot 5\cdot 89 ,\,\eta=2 ,\,h=2$  &&&&&&&
$\langle  2670\rangle \oplus A_2(1/3,1/3,-1/3)$ &&&&&& $nr$\cr
\+$n=160 ,\,d=897=3\cdot 13\cdot 23 ,\,\eta=6 ,\,h=2$  &&&&&&&
$\langle  1\rangle \oplus \langle  -598\rangle \oplus 
\langle  -6\rangle (0,1/2,1/2)$ &&&&&& $nr$\cr
\+$n=161 ,\,d=903=3\cdot 7\cdot 43 ,\,\eta=2 ,\,h=2$  &&&&&&&
$\langle  3\rangle \oplus \langle  -301\rangle \oplus 
\langle  -1\rangle $ &&&&&& $nr$\cr
\+$n=162 ,\,d=906=2\cdot 3\cdot 151 ,\,\eta=0 ,\,h=2$  &&&&&&&
$U\oplus \langle  -906\rangle $ &&&&&& $nr$\cr
\+$n=163 ,\,d=915=3\cdot 5\cdot 61 ,\,\eta=5 ,\,h=2$  &&&&&&&
$\langle  2\rangle \oplus \langle  -1830\rangle \oplus 
\langle  -1\rangle (1/2,1/2,0)$ &&&&&& $nr$\cr
\+$n=164 ,\,d=935=5\cdot 11\cdot 17 ,\,\eta=0 ,\,h=2$  &&&&&&&
$\langle  11\rangle \oplus \langle  -17\rangle \oplus 
\langle  -5\rangle $ &&&&&& $nr$\cr
\+$n=165 ,\,d=959=7\cdot 137 ,\,\eta=3 ,\,h=2$  &&&&&&&
$\langle  2877\rangle \oplus A_2(1/3,1/3,-1/3)$ &&&&&& $nr$\cr
\+$n=166 ,\,d=966=2\cdot 3\cdot 7\cdot 23 ,\,\eta=4 ,\,h=2$  &&&&&&&
$\langle  14\rangle \oplus \langle  -46\rangle \oplus 
\langle  -6\rangle (1/2,1/2,0)$ &&&&&& $nr$\cr
\+$n=167 ,\,d=986=2\cdot 17\cdot 29 ,\,\eta=1 ,\,h=2$  &&&&&&&
$U\oplus \langle  -986\rangle $ &&&&&& $nr$\cr
\+$n=168 ,\,d=987=3\cdot 7\cdot 47 ,\,\eta=7 ,\,h=2$  &&&&&&&
$\langle  2\rangle \oplus \langle  -282\rangle \oplus 
\langle  -7\rangle (1/2,1/2,0)$ &&&&&& $nr$\cr
\+$n=169 ,\,d=1001=7\cdot 11\cdot 13 ,\,\eta=5 ,\,h=2$  &&&&&&&
$\langle  1\rangle \oplus \langle  -143\rangle \oplus 
\langle  -7\rangle $ &&&&&& $nr$\cr
\+$n=170 ,\,d=1010=2\cdot 5\cdot 101 ,\,\eta=0 ,\,h=2$  &&&&&&&
$\langle  3030\rangle \oplus A_2(1/3,1/3,-1/3)$ &&&&&& $nr$\cr
\+$n=171 ,\,d=1015=5\cdot 7\cdot 29 ,\,\eta=6 ,\,h=2$  &&&&&&&
$\langle  5\rangle \oplus \langle -406\rangle \oplus 
\langle  -2\rangle (0,1/2,1/2)$ &&&&&& $nr$\cr
\+$n=172 ,\,d=1023=3\cdot 11\cdot 31 ,\,\eta=1 ,\,h=2$  &&&&&&&
$\langle  11\rangle \oplus \langle  -93\rangle \oplus 
\langle  -1\rangle $ &&&&&& $nr$\cr
\+$n=173 ,\,d=1066=2\cdot 13\cdot 41 ,\,\eta=2 ,\,h=2$  &&&&&&&
$U\oplus \langle  -1066\rangle $ &&&&&& $nr$\cr
\+$n=174 ,\,d=1085=5\cdot 7\cdot 31 ,\,\eta=6 ,\,h=2$  &&&&&&&
$\langle  1\rangle \oplus \langle  -155\rangle \oplus 
\langle -7\rangle $ &&&&&& $nr$\cr
\+$n=175 ,\,d=1095=3\cdot 5\cdot 73 ,\,\eta=3 ,\,h=2$  &&&&&&&
$\langle  15\rangle \oplus \langle  -73\rangle \oplus 
\langle  -1\rangle $ &&&&&& $nr$\cr
\+$n=176 ,\,d=1110=2\cdot 3\cdot 5\cdot 37 ,\,\eta=0 ,\,h=2$  &&&&&&&
$\langle  2\rangle \oplus \langle  -370\rangle \oplus 
\langle  -6\rangle (0,1/2,1/2)$ &&&&&& $nr$\cr
\+$n=177 ,\,d=1122=2\cdot 3\cdot 11\cdot 17 ,\,\eta=1 ,\,h=2$  &&&&&&&
$\langle  6\rangle \oplus \langle  -34\rangle \oplus 
\langle  -22\rangle (0,1/2,1/2)$ &&&&&& $nr$\cr
\+$n=178 ,\,d=1130=2\cdot 5\cdot 113 ,\,\eta=2 ,\,h=2$  &&&&&&&
$U\oplus \langle  -1130\rangle $ &&&&&& $nr$\cr
\+$n=179 ,\,d=1155=3\cdot 5\cdot 7\cdot 11 ,\,\eta=10 ,\,h=2$  &&&&&&&
$\langle  1\rangle \oplus \langle  -35\rangle \oplus 
\langle  -33\rangle $ &&&&&& $nr$\cr
\+$n=180 ,\,d=1173=3\cdot 17\cdot 23 ,\,\eta=6 ,\,h=2$  &&&&&&&
$\langle  1\rangle \oplus \langle  -51\rangle \oplus 
\langle  -23\rangle $ &&&&&& $nr$\cr
\+$n=181 ,\,d=1194=2\cdot 3\cdot 199 ,\,\eta=0 ,\,h=2$  &&&&&&&
$U\oplus \langle  -1194\rangle $ &&&&&& $nr$\cr
\+$n=182 ,\,d=1218=2\cdot 3\cdot 7\cdot 29 ,\,\eta=2 ,\,h=2$  &&&&&&&
$\langle  174\rangle \oplus \langle  -14\rangle \oplus 
\langle  -2\rangle (0,1/2,1/2)$ &&&&&& $nr$\cr
\+$n=183 ,\,d=1235=5\cdot 13\cdot 19 ,\,\eta=6 ,\,h=2$  &&&&&&&
$\langle  247\rangle \oplus \langle  -5\rangle \oplus 
\langle  -1\rangle $ &&&&&& $nr$\cr
\+$n=184 ,\,d=1245=3\cdot 5\cdot 83 ,\,\eta=6 ,\,h=2$  &&&&&&&
$\langle  1\rangle \oplus \langle  -83\rangle \oplus 
\langle  -15\rangle $ &&&&&& $nr$\cr
\+$n=185 ,\,d=1254=2\cdot 3\cdot 11\cdot 19 ,\,\eta=2 ,\,h=2$  &&&&&&&
$\langle  2\rangle \oplus \langle  -418\rangle \oplus 
\langle  -6\rangle (1/2,1/2,0)$ &&&&&& $nr$\cr
\+$n=186 ,\,d=1271=31\cdot 41 ,\,\eta=3 ,\,h=2$  &&&&&&&
$\langle  3813\rangle \oplus A_2(1/3,1/3,-1/3)$ &&&&&& $nr$\cr
\+$n=187 ,\,d=1295=5\cdot 7\cdot 37 ,\,\eta=3 ,\,h=2$  &&&&&&&
$\langle  35\rangle \oplus \langle  -37\rangle \oplus 
\langle  -1\rangle $ &&&&&& $nr$\cr
\+$n=188 ,\,d=1302=2\cdot 3\cdot 7\cdot 31 ,\,\eta=2 ,\,h=2$  &&&&&&&
$U\oplus \langle  -1302\rangle $ &&&&&& $nr$\cr
\+$n=189 ,\,d=1302=2\cdot 3\cdot 7\cdot 31 ,\,\eta=4 ,\,h=2$  &&&&&&&
$\langle  2\rangle \oplus \langle  -62\rangle \oplus 
\langle  -42)(1/2,0,1/2)$ &&&&&& $nr$\cr
\+$n=190 ,\,d=1338=2\cdot 3\cdot 223 ,\,\eta=0 ,\,h=2$  &&&&&&&
$U\oplus \langle  -1338\rangle $ &&&&&& $nr$\cr
\+$n=191 ,\,d=1365=3\cdot 5\cdot 7\cdot 13 ,\,\eta=9 ,\,h=2$  &&&&&&&
$\langle  5\rangle \oplus \langle  -91\rangle \oplus 
\langle  -3\rangle $ &&&&&& $hr$\cr
\+$n=192 ,\,d=1365=3\cdot 5\cdot 7\cdot 13 ,\,\eta=10 ,\,h=2$  &&&&&&&
$\langle  1\rangle \oplus \langle  -91\rangle \oplus 
\langle  -15\rangle $ &&&&&& $nr$\cr
\+$n=193 ,\,d=1365=3\cdot 5\cdot 7\cdot 13 ,\,\eta=12 ,\,h=2$  &&&&&&&
$\langle  1\rangle \oplus \langle  -195\rangle \oplus 
\langle  -7\rangle $ &&&&&& $nr$\cr
\+$n=194 ,\,d=1365=3\cdot 5\cdot 7\cdot 13 ,\,\eta=13 ,\,h=2$  &&&&&&&
$\langle  1\rangle \oplus \langle  -273\rangle \oplus 
\langle  -5\rangle $ &&&&&& $nr$\cr
\+$n=195 ,\,d=1370=2\cdot 5\cdot 137 ,\,\eta=2 ,\,h=2$  &&&&&&&
$U\oplus \langle  -1370\rangle $ &&&&&& $nr$\cr
\+$n=196 ,\,d=1446=2\cdot 3\cdot 241 ,\,\eta=0 ,\,h=2$  &&&&&&&
$U\oplus \langle  -1446\rangle $ &&&&&& $nr$\cr
\+$n=197 ,\,d=1463=7\cdot 11\cdot 19 ,\,\eta=1 ,\,h=2$  &&&&&&&
$\langle  19\rangle \oplus \langle  -77\rangle \oplus 
\langle  -1\rangle $ &&&&&& $nr$\cr
\+$n=198 ,\,d=1479=3\cdot 17\cdot 29 ,\,\eta=6 ,\,h=2$  &&&&&&&
$\langle  3\rangle \oplus \langle  -493\rangle \oplus 
\langle  -1\rangle $ &&&&&& $nr$\cr
\+$n=199 ,\,d=1482=2\cdot 3\cdot 13\cdot 19 ,\,\eta=0 ,\,h=2$  &&&&&&&
$U\oplus \langle  -1482\rangle $ &&&&&& $nr$\cr
\+$n=200 ,\,d=1482=2\cdot 3\cdot 13\cdot 19 ,\,\eta=6 ,\,h=2$  &&&&&&&
$\langle  38\rangle \oplus \langle  -78\rangle \oplus 
\langle  -2\rangle (0,1/2,1/2)$ &&&&&& $nr$\cr
\+$n=201 ,\,d=1495=5\cdot 13\cdot 23 ,\,\eta=5 ,\,h=2$  &&&&&&&
$\langle  115\rangle \oplus \langle  -13\rangle \oplus 
\langle  -1\rangle $ &&&&&& $nr$\cr
\+$n=202 ,\,d=1526=2\cdot 7\cdot 109 ,\,\eta=3 ,\,h=2$  &&&&&&&
$U\oplus \langle  -1526\rangle $ &&&&&& $nr$\cr
\+$n=203 ,\,d=1551=3\cdot 11\cdot 47 ,\,\eta=4 ,\,h=2$  &&&&&&&
$\langle  11\rangle \oplus \langle -141\rangle \oplus 
\langle  -1\rangle $ &&&&&& $nr$\cr
\+$n=204 ,\,d=1554=2\cdot 3\cdot 7\cdot 37 ,\,\eta=2 ,\,h=2$  &&&&&&&
&&&&&& $nr$\cr
\+$n=205 ,\,d=1554=2\cdot 3\cdot 7\cdot 37 ,\,\eta=4 ,\,h=2$  &&&&&&&
$U\oplus \langle  -1554\rangle $ &&&&&& $nr$\cr
\+$n=206 ,\,d=1586=2\cdot 13\cdot 61 ,\,\eta=0 ,\,h=2$  &&&&&&&
&&&&&& $nr$\cr
\+$n=207 ,\,d=1586=2\cdot 13\cdot 61 ,\,\eta=3 ,\,h=2$  &&&&&&&
$U\oplus \langle  -1586\rangle $ &&&&&& $nr$\cr
\+$n=208 ,\,d=1605=3\cdot 5\cdot 107 ,\,\eta=6 ,\,h=2$  &&&&&&&
$\langle  1\rangle \oplus \langle  -107\rangle \oplus 
\langle  -15\rangle $ &&&&&& $nr$\cr
\+$n=209 ,\,d=1626=2\cdot 3\cdot 271 ,\,\eta=0 ,\,h=2$  &&&&&&&
$U\oplus \langle  -1626\rangle $ &&&&&& $nr$\cr
\+$n=210 ,\,d=1635=3\cdot 5\cdot 109 ,\,\eta=5 ,\,h=2$  &&&&&&&
$\langle  2\rangle \oplus \langle  -3270\rangle \oplus 
\langle  -1\rangle (1/2,1/2,0)$ &&&&&& $nr$\cr
\+$n=211 ,\,d=1677=3\cdot 13\cdot 43 ,\,\eta=6 ,\,h=2$  &&&&&&&
$\langle  2\rangle \oplus \langle  -3354\rangle 
\oplus \langle  -1\rangle (1/2,1/2,0)$ &&&&&& $nr$\cr
\+$n=212 ,\,d=1751=17\cdot 103 ,\,\eta=3 ,\,h=2$  &&&&&&&
$\langle  5253\rangle \oplus A_2(1/3,1/3,-1/3)$ &&&&&& $nr$\cr
\+$n=213 ,\,d=1771=7\cdot 11\cdot 23 ,\,\eta=7 ,\,h=2$  &&&&&&&
$\langle  2\rangle \oplus \langle  -3542\rangle \oplus 
\langle  -1\rangle (1/2,1/2,0)$ &&&&&& $nr$\cr
\+$n=214 ,\,d=1785=3\cdot 5\cdot 7\cdot 17 ,\,\eta=10 ,\,h=2$  &&&&&&&
$\langle  1\rangle \oplus \langle  -119\rangle \oplus 
\langle  -15\rangle $ &&&&&& $nr$\cr
\+$n=215 ,\,d=1794=2\cdot 3\cdot 13\cdot 23 ,\,\eta=4 ,\,h=2$  
&&&&&&&  &&&&&& $nr$\cr
\+$n=216 ,\,d=1806=2\cdot 3\cdot 7\cdot 43 ,\,\eta=6 ,\,h=2$  &&&&&&&
$\langle  86\rangle \oplus \langle  -14\rangle \oplus 
\langle  -6\rangle (1/2,1/2,0)$ &&&&&& $nr$\cr
\+$n=217 ,\,d=1898=2\cdot 13\cdot 73 ,\,\eta=2 ,\,h=2$  &&&&&&&
$U\oplus \langle  -1898\rangle $ &&&&&& $nr$\cr
\+$n=218 ,\,d=1986=2\cdot 3\cdot 331 ,\,\eta=2 ,\,h=2$  &&&&&&&
$U\oplus \langle  -1986\rangle $ &&&&&& $nr$\cr
\+$n=219 ,\,d=1995=3\cdot 5\cdot 7\cdot 19 ,\,\eta=7 ,\,h=2$  &&&&&&&
$\langle  35\rangle \oplus \langle  -57\rangle \oplus 
\langle  -1\rangle $ &&&&&& $hr$\cr
\+$n=220 ,\,d=1995=3\cdot 5\cdot 7\cdot 19 ,\,\eta=14 ,\,h=2$  &&&&&&&
$\langle  3\rangle \oplus \langle  -665\rangle \oplus 
\langle  -1\rangle$ &&&&&& $nr$\cr
\+$n=221 ,\,d=2015=5\cdot 13\cdot 31 ,\,\eta=6 ,\,h=2$  &&&&&&&
$\langle  403\rangle \oplus \langle  -5\rangle \oplus 
\langle  -1\rangle $ &&&&&& $nr$\cr
\+$n=222 ,\,d=2046=2\cdot 3\cdot 11\cdot 31 ,\,\eta=0 ,\,h=2$  &&&&&&&
$\langle  62\rangle \oplus \langle  -22\rangle \oplus 
\langle  -6\rangle (1/2,1/2,0)$ &&&&&& $nr$\cr
\+$n=223 ,\,d=2090=2\cdot 5\cdot 11\cdot 19 ,\,\eta=0 ,\,h=2$  &&&&&&&
$\langle  6270\rangle \oplus A_2(1/3,1/3,-1/3)$ &&&&&& $nr$\cr
\+$n=224 ,\,d=2145=3\cdot 5\cdot 11\cdot 13 ,\eta=14,h=2$  &&&&&&&
$\langle  1\rangle \oplus \langle -65\rangle \oplus 
\langle -33\rangle $ &&&&&& $nr$\cr
\+$n=225 ,\,d=2226=2\cdot 3\cdot 7\cdot 53 ,\,\eta=2 ,\,h=2$  &&&&&&&
$\langle  318\rangle \oplus \langle  -14\rangle 
\oplus \langle -2\rangle (0,1/2,1/2)$ &&&&&& $nr$\cr
\+$n=226 ,\,d=2415=3\cdot 5\cdot 7\cdot 23 ,\,\eta=1 ,\,h=2$  &&&&&&&
$\langle  7\rangle \oplus \langle  -345\rangle \oplus 
\langle  -1\rangle $ &&&&&& $nr$\cr
\+$n=227 ,\,d=2415=3\cdot 5\cdot 7\cdot 23 ,\,\eta=4 ,\,h=2$  &&&&&&&
$\langle  3\rangle \oplus \langle -161\rangle \oplus 
\langle  -5\rangle $ &&&&&& $nr$\cr
\+$n=228 ,\,d=2415=3\cdot 5\cdot 7\cdot 23 ,\,\eta=8 ,\,h=2$  &&&&&&&
$\langle  3\rangle \oplus \langle  -805\rangle \oplus 
\langle  -1\rangle $ &&&&&& $nr$\cr
\+$n=229 ,\,d=2454=2\cdot 3\cdot 409 ,\,\eta=0 ,\,h=2$  &&&&&&&
$U\oplus \langle  -2454\rangle $ &&&&&& $nr$\cr
\+$n=230 ,\,d=2562=2\cdot 3\cdot 7\cdot 61 ,\,\eta=2 ,\,h=2$  &&&&&&&
$U\oplus \langle  -2562\rangle $ &&&&&& $nr$\cr
\+$n=231 ,\,d=2570=2\cdot 5\cdot 257 ,\,\eta=2 ,\,h=2$  &&&&&&&
$U\oplus \langle  -2570\rangle $ &&&&&& $nr$\cr
\+$n=232 ,\,d=2639=7\cdot 13\cdot 29 ,\,\eta=5 ,\,h=2$  &&&&&&&
$\langle  13\rangle \oplus \langle  -406\rangle \oplus 
\langle  -2\rangle (0,1/2,1/2)$ &&&&&& $nr$\cr
\+$n=233 ,\,d=2730=2\cdot 3\cdot 5\cdot 7\cdot 13,\eta=0,h=2$  &&&&&&&
\hskip10pt $\langle  14\rangle \hskip-2pt\oplus\hskip-2pt 
\langle  -130\rangle \hskip-2pt\oplus\hskip-2pt 
\langle  -6\rangle (1/2,0,1/2)$ &&&&&& $nr$\cr
\+$n=234 ,\,d=2730=2\cdot 3\cdot 5\cdot 7\cdot 13,\eta=10,h=2$  &&&&&&&
\hskip10pt $\langle  14\rangle\hskip-2pt \oplus\hskip-2pt 
\langle  -78\rangle \hskip-2pt \oplus 
\hskip-2pt \langle  -10\rangle (0,1/2,1/2)$ &&&&&& $nr$\cr
\+$n=235 ,\,d=2805=3\cdot 5\cdot 11\cdot 17 ,\,\eta=3 ,\,h=2$  &&&&&&&
$\langle  1\rangle \oplus \langle  -935\rangle \oplus 
\langle  -3\rangle $ &&&&&& $nr$\cr
\+$n=236 ,\,d=2829=3\cdot 23\cdot 41 ,\,\eta=6 ,\,h=2$  &&&&&&&
$\langle  1\rangle \oplus \langle  -123\rangle \oplus 
\langle  -23\rangle $ &&&&&& $nr$\cr
\+$n=237 ,\,d=2886=2\cdot 3\cdot 13\cdot 37 ,\,\eta=0 ,\,h=2$  &&&&&&&
$U\oplus \langle  -2886\rangle $ &&&&&& $nr$\cr
\+$n=238 ,\,d=3003=3\cdot 7\cdot 11\cdot 13,\eta=11,h=2$  &&&&&&&
$\langle  11\rangle \oplus \langle  -273\rangle \oplus 
\langle  -1\rangle $ &&&&&& $nr$\cr
\+$n=239 ,\,d=3045=3\cdot 5\cdot 7\cdot 29 ,\,\eta=10 ,\,h=2$  &&&&&&&
$\langle  1\rangle \oplus \langle  -87\rangle \oplus 
\langle  -35\rangle $ &&&&&& $nr$\cr
\+$n=240 ,\,d=3066=2\cdot 3\cdot 7\cdot 73 ,\,\eta=6 ,\,h=2$  &&&&&&&
$U\oplus \langle  -3066\rangle $ &&&&&& $nr$\cr
\+$n=241 ,\,d=3135=3\cdot 5\cdot 11\cdot 19 ,\,\eta=2 ,\,h=2$  &&&&&&&
$\langle  11\rangle \oplus \langle  -285\rangle \oplus 
\langle  -1\rangle $ &&&&&& $nr$\cr
\+$n=242 ,\,d=3135=3\cdot 5\cdot 11\cdot 19 ,\,\eta=8 ,\,h=2$  &&&&&&&
$\langle  3\rangle \oplus \langle  -1045\rangle \oplus 
\langle  -1\rangle $ &&&&&& $nr$\cr
\+$n=243 ,\,d=3315=3\cdot 5\cdot 13\cdot 17 ,\,\eta=0 ,\,h=2$  &&&&&&&
$\langle  3\rangle \oplus \langle  -221\rangle \oplus 
\langle  -5\rangle $ &&&&&& $nr$\cr
\+$n=244 ,\,d=3318=2\cdot 3\cdot 7\cdot 79 ,\,\eta=2 ,\,h=2$  &&&&&&&
$\langle  2\rangle \oplus \langle  -474\rangle 
\oplus \langle  -14\rangle (0,1/2,1/2)$ &&&&&& $nr$\cr
\+$n=245 ,\,d=3354=2\cdot 3\cdot 13\cdot 43 ,\,\eta=0 ,\,h=2$  &&&&&&&
&&&&&& $nr$\cr
\+$n=246 ,\,d=3354=2\cdot 3\cdot 13\cdot 43 ,\,\eta=6 ,\,h=2$  &&&&&&&
$U\oplus \langle  -3354\rangle $ &&&&&& $nr$\cr
\+$n=247 ,\,d=3410=2\cdot 5\cdot 11\cdot 31 ,\,\eta=4 ,\,h=2$  &&&&&&&
&&&&&& $nr$\cr
\+$n=248 ,\,d=4026=2\cdot 3\cdot 11\cdot 61 ,\,\eta=0 ,\,h=2$  &&&&&&&
$\langle  1342\rangle \oplus \langle  -6\rangle 
\oplus \langle  -2\rangle (0,1/2,1/2)$ &&&&&& $nr$\cr
\+$n=249 ,\,d=4074=2\cdot 3\cdot 7\cdot 97 ,\,\eta=6 ,\,h=2$  &&&&&&&
$U\oplus \langle  -4074\rangle $ &&&&&& $nr$\cr
\+$n=250 ,\,d=4290=2\cdot 3\cdot 5\cdot 11\cdot 13 ,\eta=1,h=2$  &&&&&&&
\hskip10pt $\langle  286\rangle  \hskip-2pt\oplus  \hskip-2pt
\langle  -30\rangle \hskip-2pt \oplus \hskip-2pt 
\langle  -2\rangle (0,1/2,1/2)$ &&&&&& $nr$\cr
\+$n=251 ,\,d=4326=2\cdot 3\cdot 7\cdot 103 ,\,\eta=2 ,\,h=2$  &&&&&&&
$U\oplus \langle  -4326\rangle $ &&&&&& $nr$\cr
\+$n=252 ,\,d=4902=2\cdot 3\cdot 19\cdot 43 ,\,\eta=4 ,\,h=2$  &&&&&&&
$\langle  2\rangle \oplus \langle  -1634\rangle \oplus 
\langle  -6\rangle (0,1/2,1/2)$ &&&&&& $nr$\cr
\+$n=253 ,\,d=4991=7\cdot 23\cdot 31 ,\,\eta=7 ,\,h=2$  &&&&&&&
$\langle  4991\rangle \oplus \langle  -1\rangle \oplus 
\langle  -1\rangle $ &&&&&& $nr$\cr
\+$n=254 ,\,d=5226=2\cdot 3\cdot 13\cdot 67 ,\,\eta=0 ,\,h=2$  &&&&&&&
$U\oplus \langle  -5226\rangle $ &&&&&& $nr$\cr
\+$n=255 ,\,d=5334=2\cdot 3\cdot 7\cdot 127 ,\,\eta=2 ,\,h=2$  &&&&&&&
$\langle  2\rangle \oplus \langle  -762\rangle \oplus 
\langle  -14\rangle (0,1/2,1/2)$ &&&&&& $nr$\cr
\+$n=256 ,\,d=6006=2\cdot 3\cdot 7\cdot 11\cdot 13 ,\eta=2,h=2$  &&&&&&&
\hskip10pt$\langle  286\rangle 
\hskip-2pt\oplus\hskip-2pt \langle-14\rangle \hskip-2pt\oplus\hskip-2pt 
\langle  -6\rangle (1/2,1/2,0)$ &&&&&& $nr$\cr
\+$n=257 ,\,d=7590=2\cdot 3\cdot 5\cdot 11\cdot 23,\eta=8 ,h=2$  &&&&&&&
\hskip10pt $\langle  230\rangle \hskip-2pt\oplus \hskip-2pt\langle-22\rangle 
\hskip-2pt\oplus\hskip-2pt 
\langle  -6\rangle (1/2,1/2,0)$ &&&&&& $nr$\cr
\+$n=258 ,\,d=10374=2\cdot 3\cdot 7\cdot 13\cdot 19 ,\eta=2,h=2$  &&&&&&&
\hskip10pt $\langle  494\rangle \hskip-2pt \oplus\hskip-2pt 
\langle  -14\rangle \hskip-2pt\oplus \hskip-2pt
\langle  -6\rangle (1/2,1/2,0)$ &&&&&& $nr$\cr
\+$n=259 ,\,d=29526=2\cdot 3\cdot 7\cdot 19\cdot 37 ,\,\eta=2 ,\,h=2$  &&&&&&&
\hskip20pt $U\oplus \langle  -29526\rangle $ &&&&&& $nr$\cr}

\newpage

\head
8. Narrow parts of restricted hyperbolic convex polygons on the
hyperbolic plane. Application to hyperbolically reflective
hyperbolic lattices
\endhead

Here we consider narrow parts of restricted hyperbolic convex polygons 
on the hyperbolic plane and apply them to hyperbolically reflective 
hyperbolic lattices of rank three. Like in \cite{N4}, \cite{N5},
\cite{N9}, these results can be generalized on restricted hyperbolic 
convex polyhedra in hyperbolic spaces of arbitrary dimension and on 
hyperbolically reflective hyperbolic lattices of rank $\ge 3$.

\subhead
8.1. Some formulae using cross-ratio
\endsubhead

We consider a hyperbolic plane and cross-ratio for points at
infinity. The infinity is a non-singular rational curve of 
degree 2. We normalize the curvature $\kappa=-1$. 

\proclaim{Lemma 8.1.1} Let $(AB)$ and $(CD)$ are two
lines in a hyperbolic plane with terminals $A$, $B$, $C$,
$D$ at infinity, and $O$ a point of the hyperbolic plane which
does not belong to each line $(AB)$ and $(CD)$ and orientations of the
triangles $OAB$ and $OCD$ coincide.
Let $\delta_1$ and $\delta_2$ are orthogonal vectors with
square $-2$ to lines $(AB)$ and $(CD)$ respectively such that
$O$ is contained in both half-planes $\Ha_{\delta_1}^+$ and
$\Ha_{\delta_2}^+$.

Then
$$
(\delta_1,\,\delta_2)=4[A,D,C,B]-2.
$$

As a corollary, we get:

1) If lines $(AB)$ and $(CD)$ do not intersect each other and
$\Ha_{\delta_1}^+\cap\Ha_{\delta_2}^+$ is a band
between lines $(AB)$ and $(CD)$, then
$$
2 \ch~{\rho }=
(\delta_1,\delta_2)=4[A,D,C,B]-2
\tag{8.1.1}
$$
where $\rho$ is the distance between lines $(AB)$ and $(CD)$.

2) If $\Ha_{\delta_1}^+\cap \Ha_{\delta_2}^+$ is
an angle $\alpha$, then
$$
2 \cos{\alpha}=
(\delta_1,\delta_2)=4[A,D,C,B]-2.
\tag{8.1.2}
$$
\endproclaim

\demo{Proof} One can prove this by direct calculation
(also see the proof of Lemma 4.1.2).
\enddemo

We have
$[A,D,C,B]=[A,D,B,C]^{-1}=(1-[A,B,D,C])^{-1}=$
\newline
$(1-[A,B,C,D]^{-1})^{-1}=[A,B,C,D]/([A,B,C,D]-1)$.
Thus, we have
$$
[A,B,C,D]={\ch~\rho+1\over \ch~\rho-1}>1
\tag{8.1.3}
$$
for the case 1) of Lemma 8.1.1,
and
$$
[A,B,C,D]={\cos{\alpha}+1\over \cos{\alpha}-1}<0
\tag{8.1.4}
$$
for the case 2) of Lemma 8.1.1.

\smallpagebreak

Further we fix a line $l=(MN)$ on the hyperbolic plane with
terminals $M$ and $N$ at infinity. This line will be called
{\it axise}. We shall consider pairs of points $A$ and $B$
at ifinity. In general, we assume that they are different from
$M$ and $N$ (i. e. they do not
belong to $l$) but one can also consider cases when one of them is
equal to $M$ or $N$ taking an appropriate limit. The axis $l$ 
divides the infinity in two arcs containing in two different 
half-planes  bounded by $l$.

First suppose that $A$ and $B$ belong to one half-plane
bounded by $l$. Additionally suppose that
$M$, $N$, $A$ and $B$ are four consecutive points at infinity
(here we consider the infinity as a circle).
Considering a point $O$ between
the lines $l$ and $(AB)$ and using Lemma 8.1.1 and
\thetag{8.1.3}, we get
$$
[M,N,A,B]={\ch~\rho(l,(AB))+1\over \ch~\rho(l,(AB))-1}\ge 1.
\tag{8.1.5}
$$
It follows
$$
0<[M,N,B,A]=[M,N,A,B]^{-1}={\ch~\rho(l,(AB))-1\over \ch~\rho(l,(AB))+1}\le 1.
\tag{8.1.6}
$$
By definition of the cross-ratio, 
$$
[M,N,A,B]\cdot [M,N,B,C]=[M,N,A,C].
\tag{8.1.7}
$$
Further, for any two points $A$, $B$ at infinity which
are contained in one half-plane bounded by the axis $l$, we  
introduce the {\it ``angle''}
$$
\theta(A,B)=\ln{[M,N,A,B]}\in \br.
\tag{8.1.8}
$$
The ``angle'' $\theta(A,B)$ does behave like an angle: by
\thetag{8.1.6} and \thetag{8.1.7},
$$
\theta(B,A)=-\theta(A,B),\ \ \theta(A,C)=\theta(A,B)+\theta(B,C)
\tag{8.1.9}
$$
where we suppose that all points $A$, $B$ and $C$ belong to
one half-plane bounded by the axis $l$.

The axis $l$ divides the infinity on two arcs and defines the orthogonal
projection $X\mapsto X^\prime$ with respect to $l$
of one  arc onto another: by definition, the lines $l$ and
$(XX')$ are perpendicular. By \thetag{8.1.4} and \thetag{8.1.7},
we have
$$
[M,N,X,X']=-1,\ \text{and}\ [M,N,A,B']=-[M,N,A,B].
\tag{8.1.10}
$$
For points
$A$ and $B$ at infinity which belong to different half-planes
bounded by the axis $l=(MN)$, we then define the ``angle'' as follows:
$$
\theta(A,\,B)=:\theta(A,\,B^\prime)
\tag{8.1.11}
$$
where $\theta(A,\,B^\prime)$ was defined in \thetag{8.1.8}.
Remark that $\theta(A,\,B^\prime)=\theta(A^\prime,\,B)$.
Thus, we defined the angle $\theta(A,\,B)$ for all points $A$ and $B$
at infinity. For all of them one has
\thetag{8.1.9}.

For points $A$ and $B$ at infinity
we define the invariant
$$
\epsilon(A,B)=\pm 1
\tag{8.1.12}
$$
where {$\epsilon(A,B)=1$ if $A$ and $B$ belong to one half-plane
bounded by $l$, and $\epsilon(A,B)=-1$ otherwise.} Obviously,
$$
\epsilon(A,B)=\epsilon(B,A)\ \text{and}\ 
\epsilon(A,B)\epsilon(B,C)=\epsilon(A,C).
\tag{8.1.13}
$$
From definition above, it follows
$$
[M,N,A,B]=\epsilon(A,B)\exp{\theta(A,B)}
\tag{8.1.14}
$$
for any two points $A$ and $B$ at infinity.

Now assume that $A$, $B$, $C$ and $D$ are
points at infinity.
We choose an affine coordinate at infinity
such that $M=\infty$, $N=0$, $A=1$. Then $B=u$, $C=uv$, $D=uvw$
where
$[M,N,A,B]=u=\epsilon(A,B)\exp{\theta(A,B)}$,
$[M,N,B,C]=v=\epsilon(B,C)\exp{\theta(B,C)}$,
$[M,N,C,D]=w=\epsilon(C,D)\exp{\theta(C,D)}$.
It follows
$$
[A,D,C,B]={uv-1\over uv-uvw}:{u-1\over u-uvw}=
{(uv-1)(vw-1)\over (u-1)(w-1)v}=
$$
$$
\left[{\exp{\left({\theta(A,B)+\theta(B,C)\over 2}\right)}-
\epsilon(A,C)\exp{\left(-{\theta(A,B)+\theta(B,C)\over 2}\right)}\over 2}
\right.\times
$$
$$
\left.\left.{\exp{\left({\theta(B,C)+\theta(C,D)\over 2}\right)}-
\epsilon(B,D)\exp{\left(-{\theta(B,C)+\theta(C,D)\over 2}\right)}\over 2}
\right]\right/
$$
$$
\left[{\exp{\left({\theta(A,B)\over 2}\right)}\hskip-2pt-\hskip-2pt
\epsilon(A,B)\exp{\left(-{\theta(A,B)\over 2}\right)}\over 2}\times
{\exp{\left({\theta(B,C)\over 2}\right)}\hskip-2pt-\hskip-2pt
\epsilon(B,C)\exp{\left(-{\theta(B,C)\over 2}\right)}\over 2}\right].
\tag{8.1.15}
$$
We remark that functions used in this expression
are hyperbolic sinus or cosinus. E. g.
$$
{\exp{\left({\theta(A,B)\over 2}\right)}\hskip-2pt-\hskip-2pt
\epsilon(A,B)\exp{\left(-{\theta(A,B)\over 2}\right)}\over 2}=
\cases
\sh~{\theta(A,B)\over 2} &\text{if $\epsilon(A,B)=1$}\\
\ch~{\theta(A,B)\over 2} &\text{if $\epsilon(A,B)=-1$}
\endcases.
\tag{8.1.16}
$$

Using Lemma 8.1.1 and \thetag{8.1.15}, we get the following
analog of Lemmas 4.1.2 and 4.2.4 (for Lemmas 4.1.2 and 4.2.4,
a finite or infinite point was used instead of the axis $l$).

\proclaim{Lemma 8.1.2}  Let $l=(MN)$ be a line (called axis ) of
a hyperbolic plane with terminals $M$ and $N$ at infinity.
Let $(AB)$ and $(CD)$ are two
lines of the hyperbolic plane with terminals $A$, $B$, $C$,
$D$ at infinity, and there exists a point $O$ of the hyperbolic plane
which does not belong to each line $(AB)$ and $(CD)$ and orientations
of the triangles $OAB$ and $OCD$ coincide.
Let $\delta_1$ and $\delta_2$ are orthogonal vectors with
square $-2$ to lines $(AB)$ and $(CD)$ respectively such that
$O$ is contained in both half-planes $\Ha_{\delta_1}^+$ and
$\Ha_{\delta_2}^+$.

Then
$$
(\delta_1,\,\delta_2)=4[A,D,C,B]-2
$$
where
$$
[A,D,C,B]=
$$
$$
\left[{\exp{\left({\theta(A,B)+\theta(B,C)\over 2}\right)}-
\epsilon(A,C)\exp{\left(-{\theta(A,B)+\theta(B,C)\over 2}\right)}\over 2}
\right.\times
$$
$$
\left.\left.{\exp{\left({\theta(B,C)+\theta(C,D)\over 2}\right)}-
\epsilon(B,D)\exp{\left(-{\theta(B,C)+\theta(C,D)\over 2}\right)}\over 2}
\right]\right/
$$
$$
\left[{\exp{\left({\theta(A,B)\over 2}\right)}\hskip-2pt-\hskip-2pt
\epsilon(A,B)\exp{\left(-{\theta(A,B)\over 2}\right)}\over 2}\times
{\exp{\left({\theta(B,C)\over 2}\right)}\hskip-2pt-\hskip-2pt
\epsilon(B,C)\exp{\left(-{\theta(B,C)\over 2}\right)}\over 2}\right].
\tag{8.1.17}
$$
Here we use invariants $\theta$ and $\epsilon$ introduced in
\thetag{8.1.8}, \thetag{8.1.11} and \thetag{8.1.12} 
(also see \thetag{8.1.5}). 
\endproclaim

Below we give two the most important for us particular
cases of Lemma 8.1.2.

\proclaim{Lemma 8.1.3}  Let $l=(MN)$ be a line (called axis) of a
hyperbolic plane with terminals $M$ and $N$ at infinity.
Let $(AB)$ and $(CD)$ are two lines of the
plane which belong to one half-plane bounded by the axis
$l$ and $A$, $B$, $C$, $D$ are their terminals at infinity.  We
suppose that $M$, $N$, $A$, $B$ are four consecutive points
at infininity and $M$, $N$, $C$, $D$ are four consecutive points
at infinity either.
We consider angles
$\theta_1=\theta(A,B)$, $\theta_2=\theta(C,D)$ and
$\theta_{12}=\theta(B,C)$ defined by \thetag{8.1.5} and
\thetag{8.1.8}. Let $\delta_1$ and $\delta_2$ are orthogonal vectors with
square $-2$ to lines $(AB)$ and $(CD)$ respectively such that
$l$ is contained in both half-planes $\Ha_{\delta_1}^+$ and
$\Ha_{\delta_2}^+$. Then
$$
(\delta_1,\,\delta_2)=
4{\sh~{\theta_1+\theta_{12}\over 2}\sh~{\theta_2+\theta_{12}\over 2}
\over
\sh~{\theta_1\over 2}\sh~{\theta_2\over 2}}-2.
\tag{8.1.18}
$$

As a corollary we get:

If lines $(AB)$ and $(CD)$ do not intersect and
$\Ha_{\delta_1}^+\cap \Ha_{\delta_2}^+$ is the band bounded by
these lines, we have
$$
(\delta_1,\,\delta_2)\ =\ 2\ch~\rho\ =\
4{\sh~{\theta_1+\theta_{12}\over 2}\sh~{\theta_2+\theta_{12}\over 2}
\over
\sh~{\theta_1\over 2}\sh~{\theta_2\over 2}}-2
\tag{8.1.19}
$$
where $\rho$ is the distance between lines $(AB)$ and $(CD)$.

If lines $(AB)$ and $(CD)$ define an angle
$\alpha=\Ha_{\delta_1}^+\cap \Ha_{\delta_2}^+$, we have
$$
(\delta_1,\,\delta_2)\ =\ 2\cos{\alpha}=\
4{\sh~{\theta_1-\theta_{21}\over 2}\sh~{\theta_2-\theta_{21}\over 2}
\over
\sh~{\theta_1\over 2}\sh~{\theta_2\over 2}}-2
\tag{8.1.20}
$$
where $\theta_{21}=-\theta_{12}=\theta(C,B)$.
\endproclaim

\proclaim{Lemma 8.1.4}  Let $l=(MN)$ be a line on a hyperbolic
plane (called axis) with terminals $M$ and $N$ at infinity.
Let $(AB)$ and $(CD)$ are two lines of the
plane which belong to different half-planes bounded by the axis
$l$, and $A$, $B$, $C$, $D$ are their terminals at infinity.
We suppose that $M$, $N$, $A$, $B$ are four consecutive points at
infinity and $M$, $N$, $D$, $C$ are four consecutive points at
infinity either.  We consider angles
$\theta_1=\theta(A,B)$, $\theta_2=\theta(D,C)$ defined by
\thetag{8.1.5} and \thetag{8.1.8}, and $\theta_{12}=\theta(B,C)$
defined by \thetag{8.1.11}.
Let $\delta_1$ and $\delta_2$ are orthogonal
vectors with square $-2$ to lines $(AB)$ and $(CD)$ respectively
such that $l$ is contained in both half-planes $\Ha_{\delta_1}^+$ and
$\Ha_{\delta_2}^+$. Then
$$
(\delta_1,\,\delta_2)=2\ch~\rho=
4{\ch~{\theta_1+\theta_{12}\over 2}
\ch~{\theta_2-\theta_{12}\over 2}
\over
\sh~{\theta_1\over 2}\sh~{\theta_2\over 2}}-2.
\tag{8.1.21}
$$
Here $\rho$ is the distance between lines $(AB)$ and $(CD)$.
\endproclaim

\subhead
8.2. Narrow parts of restricted hyperbolic convex poligons in a
hyperbolic plane: general results
\endsubhead

We use notation of Sect. 1 for convex locally finite
polyhdera in hyperbolic spaces. Remind that $P(\M)$ denote the set of
orthogonal vectors to faces of $\M$ of highest dimension and
directed outwards of $\M$. We normalize $\delta\in P(\M)$
by the condition $\delta^2=-2$.

We use definition from \cite{N9}
of hyperbolic polyhedra in hyperbolic spaces.
Let $\T\subset \La$ be a hyperbolic subspace of a hyperbolic space
$\La$ such that $1\le \dim \T\le \dim \La -1$. For $K\subset \T$,
we denote by $C_K$ the orthogonal cylinder with the base $K$.
Thus, $C_K$ is union of all lines orthogonal to $\T$ and containing
a point of $K$.

\definition{Definition 8.2.1}
A convex polyhedron $\M$ in $\La$ is called
{\it hyperbolic (relative to a subspace $\T\subset \La$, 
$1\le \dim \T\le \dim \La -1$, which is called axis)}
if it is non-degenerate (i. e. it is not contained 
in a proper subspace of $\La$), locally finite, and the conditions 
(1) and (2) below hold:

(1) $\M$ is finite at infinite points, i. e. the set
$\{ \delta\in P(\M)\ \mid (\delta,\,c)=0\}$ is finite
for any infinite point $\br_{++}c\in \La_{\infty}$, $c^2=0$. 

(2) For any compact elliptic polyhedron $\N\subset \T$, the polyhedron
$C_\N\cap \M$ is elliptic if it is non-degenerate. (We remind that
a polyhedron is called elliptic if it is a convex envelope of a finite
set of points and is non-degenerate.)

A hyperbolic polyhedron is called {\it restricted} if additionally

(3) The set of angles and distances with the axis $\T$
defined by all half-spaces $\Ha_\delta^+$, $\delta\in P(\M)$,
is finite.
\enddefinition

Further we shall consider only restricted hyperbolic convex
polygons $\M$ in a hyperbolic plane. Then axis is a line $l$.

If a polygon $\M$ is elliptic, then $\M$ is restricted hyperbolic
with respect to any axis $l$. If $\M$ is hyperbolic and is finite
(i. e. the set $P(\M)$ is finite), then $\M$ is elliptic.
Narrow parts of elliptic and parabolic polygons have been considered in
Sect. 4. Therefore, further we can suppose that $\M$ is infinite
(i. e. it has infinite set of sides, equivalently the set
$P(\M)$ is infinite).

The axis $l$ divides $\M$ in two half-s. We denote by $A_i$
vertices of one half of $\M$, and by $B_j$ vertices of another half of $\M$.
We prove the following basic result about narrow parts of restricted
hyperbolic convex polygons in a hyperbolic plane.

\proclaim{Theorem 8.2.2 (about narrow parts of types (II) and (III)
of restricted hyperbolic convex polygons)}
We fix a constant $\eta >0$.

Any infinite restricted hyperbolic convex polygon $\M$ has a narrow
part of one of types
(AII), (AIII), (BII$_1$), (BII$_2$), (BIII) defined below:

(AII) There exist five consecutive vertices
$A_0$, $A_1$, $A_2$, $A_3$ and $A_4$ of one half of $\M$
such that for orthogonal
vectors $\delta_1$, $\delta_2$, $\delta_3$ and $\delta_4$
to lines $(A_0A_1)$, $(A_1A_2)$, $(A_2A_3)$ and $(A_3A_4)$
respectively directed outwards of $\M$ and with
$\delta_1^2=\delta_2^2=\delta_3^2=\delta_4^2=-2$, one has
$(\delta_i,\,\delta_{i+1})=2\cos{\alpha_i}$, $i=1,\,2,\,3$,
where $\alpha_i=\angle~A_{i-1}A_iA_{i+1}$, and
$$
2<(\delta_1,\,\delta_3)\le 8\ch~{\eta\over 2}+6,
\tag{8.2.1}
$$
$$
2<(\delta_1,\,\delta_4)\le 4(2\ch~{\eta\over 2}+1)^2-2.
\tag{8.2.2}
$$
Moreover, the set $\{\delta_1,\,\delta_3,\,\delta_4\}$
generates the 3-dimensional hyperbolic vector space and
has a connected Gram graph.
Additionally, open half-planes $\Ha_{\delta_i}^+$,
$i=1,\,2,\,3,\,4$, contain the axis $l$.

(AIII) There exist six consecutive vertices
$A_0$, $A_1$, $A_2$, $A_3$, $A_4$ and $A_5$ of one half of $\M$
such that for orthogonal vectors $\delta_1$, $\delta_2$, $\delta_3$,
$\delta_4$ and $\delta_5$
to lines $(A_0A_1)$, $(A_1A_2)$, $(A_2A_3)$, $(A_3A_4)$ and
$(A_4A_5)$ respectively directed
outwards of $\M$ and with
$\delta_1^2=\delta_2^2=\delta_3^2=\delta_4^2=\delta_5^2=-2$
one has $(\delta_i,\delta_{i+1})=2\cos{\alpha_i}$, $i=1,\,2,\,3,\,4$,
where $\alpha_i=\angle~A_{i-1}A_iA_{i+1}$, and
$$
2<(\delta_1,\delta_3)\le 8\ch~{\eta\over 2}+6,
\tag{8.2.3}
$$
$$
2<(\delta_3,\delta_5)\le 8\ch~{\eta\over 2}+6,
\tag{8.2.4}
$$
and
$$
2<(\delta_1,\delta_5)\le
32(\ch^3{\eta \over 2}+\ch^2{\eta \over 2})-2.
\tag{8.2.5}
$$
Moreover, the
set $\{\delta_1,\,\delta_3,\,\delta_5\}$ generates the 3-dimensional
hyperbolic vector space and has a connected Gram graph.
Additionally, open half-planes $\Ha_{\delta_i}^+$,
$i=1,\,2,\,3,\,4,$ $5$, contain the axis $l$.

(BII$_1$) There exist four consecutive vertices
$A_0$, $A_1$, $A_2$, $A_3$ of one half of $\M$
and two consecutive vertices $B_0$, $B_1$ of another half of $\M$
such that for orthogonal
vectors $\delta_1$, $\delta_2$, $\delta_3$ and $\delta_4$
to lines $(A_0A_1)$, $(A_1A_2)$, $(A_2A_3)$ and $(B_0B_1)$
respectively directed outwards of $\M$ and with
$\delta_1^2=\delta_2^2=\delta_3^2=\delta_4^2=-2$, one has
$(\delta_i,\,\delta_{i+1})=2\cos{\alpha_i}$, $i=1,\,2$,
where $\alpha_i=\angle~A_{i-1}A_iA_{i+1}$, and
$$
2<(\delta_1,\,\delta_4)\le
{8\over \ch~{\eta\over 2}-1}+6,
\tag{8.2.6}
$$
$$
2<(\delta_2,\,\delta_4)\le
{8\over \ch~{\eta\over 2}-1}+6,
\tag{8.2.7}
$$
$$
2<(\delta_3,\,\delta_4)\le
{8\over \ch~{\eta\over 2}-1}+6.
\tag{8.2.8}
$$
Moreover, the
set $\{\delta_1,\,\delta_2,\,\delta_4\}$ generates the 3-dimensional
hyperbolic vector space and has a connected Gram graph. The same
is valid for $\{\delta_2,\,\delta_3,\,\delta_4\}$.
Additionally, open half-planes $\Ha_{\delta_i}^+$,
$i=1,\,2,\,3$, and closed half-plane $\Ha_{\delta_4}^+$
contain the axis $l$; both lines orthogonal to the axis $l$ and
containing $A_1$ and $A_2$ intersect the line $(B_0B_1)$.

(BII$_2$) There exist three consecutive vertices
$A_0$, $A_1$, $A_2$ of one half of $\M$
and three consecutive vertices $B_0$, $B_1$, $B_2$ of
another half of $\M$ such that for orthogonal
vectors $\delta_1$, $\delta_2$, $\delta_3$ and $\delta_4$
to lines $(A_0A_1)$, $(A_1A_2)$, $(B_0B_1)$ and $(B_1B_2)$
respectively directed outwards of $\M$ and with
$\delta_1^2=\delta_2^2=\delta_3^2=\delta_4^2=-2$, one has
$(\delta_1,\,\delta_2)=2\cos{\alpha}$ where $\alpha=\angle~A_0A_1A_2$,
$(\delta_3,\,\delta_4)=2\cos{\beta}$ where
$\beta=\angle~B_0B_1B_2$, and
$$
2<(\delta_1,\,\delta_3)\le
{8\over \ch~{\eta\over 2}-1}+6,
\tag{8.2.9}
$$
$$
2<(\delta_2,\,\delta_3)\le
{8\over \ch~{\eta\over 2}-1}+6,
\tag{8.2.10}
$$
$$
2<(\delta_2,\,\delta_4)\le
{8\over \ch~{\eta\over 2}-1}+6. 
\tag{8.2.11}
$$
Moreover, the
set $\{\delta_1,\,\delta_2,\,\delta_3\}$ generates the 3-dimensional
hyperbolic vector space and has a connected Gram graph. The same
is valid for $\{\delta_2,\,\delta_3,\,\delta_4\}$.
Additionally, open half-planes $\Ha_{\delta_i}^+$,
$i=1,\,2,\,3,\,4$, contain the axis $l$; the line
orthogonal to the axis $l$ and containing $A_1$ intersects the line 
$(B_0B_1)$, and the line orthogonal to $l$ and containing $B_1$
intersects the line $(A_1A_2)$.

(BIII) There exist four consecutive vertices
$A_0$, $A_1$, $A_2$, $A_3$ of one half of $\M$,
and two consecutive vertices $B_0$, $B_1$ of
another half of $\M$ such that for orthogonal
vectors $\delta_1$, $\delta_2$, $\delta_3$ and $\delta_4$
to lines $(A_0A_1)$, $(A_1A_2)$, $(A_2A_3)$ and $(B_0B_1)$
respectively directed outwards of $\M$ and with
$\delta_1^2=\delta_2^2=\delta_3^2=\delta_4^2=-2$, one has
$(\delta_i,\,\delta_{i+1})=2\cos{\alpha_i}$ where
$\alpha_i=\angle~A_{i-1}A_iA_{i+1}$, $i=1,\,2$, and
$$
2<(\delta_1,\,\delta_3)\le
8\ch~{\eta\over 2}+6, 
\tag{8.2.12}
$$
$$
2<(\delta_1,\,\delta_4)\le
8\ch~{\eta\over 2}+10+{8\over \ch~{\eta\over 2}-1}, 
\tag{8.2.13}
$$
$$
2<(\delta_3,\,\delta_4)\le
8\ch~{\eta\over 2}+10+{8\over \ch~{\eta\over 2}-1}.
\tag{8.2.14}
$$
Moreover, the
set $\{\delta_1,\,\delta_2,\,\delta_4\}$ generates the
3-dimensional hyperbolic vector space and has a connected Gram graph.
The same is valid for the set $\{\delta_2,\,\delta_3,\,\delta_4\}$.
Additionally, open half-planes $\Ha_{\delta_i}^+$,
$i=1,\,2,\,3,\,4$, contain the axis $l$; there exists a line
orthogonal to the axis $l$ and intersecting the lines
$(A_1A_2)$ and $(B_0B_1)$.
\endproclaim

\demo{Proof} Compare with \cite{N4, Appendix}, \cite{N13}
and the proof of Theorem 4.1.5.

Further we shall normalize elements $\delta \in P(\M)$
by the condition $\delta^2=-2$.

Let us denote by $E\subset P(\M)$ the set of all elements
$e \in P(\M)$ such that the line $\Ha_e$ intersects the axis $l$
(in finite or infinite points). By condition (3), the set
$E$ is finite. Otherwise
there exist two different elements $e_1,\,e_2\in E$ such that
angles of the half-planes $\Ha_{e_1}^+$ and $\Ha_{e_2}^+$ with $l$
are equal. It then follows that one of two half-planes
$\Ha_{e_1}^+$, $\Ha_{e_2}^+$ contains another. This contradicts the
definition of $P(\M)$. Let
$$
\M_0=\bigcap_{e\in E}{\Ha_e^+}.
\tag{8.2.15}
$$
The $\M_0$ is a finite convex polyhedron containing $\M$.
If $l\cap \M_0$ is a finite or empty interval of $l$, then
orthogonal projection (here and in what follows it is
the orthogonal projection relative to the axis $l$)
of $\M_0$ into $l$ is
compact. By (2), then $\M$ is finite. We get a
contradiction. Thus $l\cap \M_0$ is either the axis
$l=(NM)$ or a half-axis $(KM)\subset l$ where $M$ and $N$
are terminals of $l$ at infinity and $K$ is a finite point of $l$.

Suppose that $\delta\in P(\M)-E$. Then either the open half-plane
$\Ha_\delta^+$ contains the axis $l$ or the closed half-plane
$\Ha_\delta^+$ does not intersect the axis $l$. For the last case
we also obtain that $\M$ is finite (by (2)).
It then follows that for any $\delta\in P(\M)-E$ the open
half-plane $\Ha_\delta^+$ contains
the axis $l$. Remark that then the line $\Ha_\delta$ has a compact
orthogonal projection into $l$ and
defines (by (2)) a side $(C_iC_{i+1})$ of $\M$ where
$C_i$ and $C_{i+1}$ are vertices of $\M$.

By the considerations above, there exists an infinite
chain $A_i$, $i\to +\infty$ of consecutive vertices of
$\M$ with the following properties. All lines
$(A_{i-1}A_{i})$ do not intersect the axis $l$ and the open
half-plane of $\La$ bounded by $(A_{i-1}A_{i})$ contains $l$.
Thus all $A_i$ belong to one half of $\M$. The vertices
$A_i$, $i\to +\infty$ tend to the terminal $M$ of $l$ at infinity.

Around the terminal $M$ of $l$, another half of $\M$ is
either another similar infinite chain
$B_j$, $j\to +\infty$ and $B_j\to M$,
or a side $(B_0M)$ of $\M$, or the axis $l$. Further we shall
consider only the first case of the infinite chain $B_j$ leaving
more simple last two cases to a reader (for these cases, one should
just argue with one chain $A_i$ and replace the chain
$B_j$ by $(B_0M)$ or $l$ respectively if necessary).

We shall denote by $e_i\in P(\M)$, $i\in I$, the orthogonal vector to
the side $(A_{i-1}A_i)$ of $\M$. Respectively $f_j\in P(\M)$,
$j\in J$, denote the orthogonal vector to the side $(B_{j-1}B_j)$ of
$\M$. We denote by $A_{i1}$ and $A_{i2}$ terminals of the
line $(A_{i-1}A_i)$ at infinity such that $A_{i1}$, $A_{i-1}$,
$A_i$, $A_{i2}$ are consecutive points of the line
$(A_{i-1}A_i)$. Similarly, we denote by $B_{j1}$ and
$B_{j2}$ terminals of the
line $(B_{j-1}B_j)$ at infinity such that $B_{j1}$, $B_{j-1}$,
$B_j$, $B_{j2}$ are consecutive points of the line
$(B_{j-1}B_j)$. We denote $\theta^1_i=
\theta(A_{i1}, A_{i2})$, $i\in I$,
and $\theta^2_j=\theta(B_{j1},\,B_{j2})$, $j\in J$. 
See \thetag{8.1.8}. 

In considerations below, without lost of generality,
we suppose that the sets $I=\bz$ and $J=\bz$, i. e. chains $A_i$
and $B_j$ are infinite in both directions. One can always suppose this
considering parts of these chains in a sufficiently small neighbourhood
of the infinite point $M\in l$ of $\M$.

We denote
$$
\theta=\min{(\,\min_{i\in I}{(\,\theta^1_{i-1}+\theta^1_i\,)},\
\min_{j\in J}{(\,\theta^2_{j-1}+\theta^2_j\,)}\,)}.
\tag{8.2.16}
$$
We can assume (changing numerations if necessary) that
$$
\theta=\theta^1_2+\theta^1_3\  \text{and}\
\theta^1_2\le \theta^1_3.
\tag{8.2.17}
$$
By definition of $\theta$, we have
$$
\theta^1_1\ge\theta^1_3\ge \theta^1_2\ \text{and}\
\theta^1_1\ge \theta/2,\
\theta^1_2\le\theta/2,\
\theta^1_3\ge \theta/2.
\tag{8.2.18}
$$
By \thetag{8.1.19},
$$
(e_1,\,e_3)\le
4{\sh~{\theta^1_1+\theta^1_2\over 2}\sh~{\theta^1_3+\theta^1_2\over 2}
\over
\sh~{\theta^1_1\over 2}\sh~{\theta^1_3\over 2}}-2.
\tag{8.2.19}
$$
Using \thetag{8.2.17} and \thetag{8.2.18}, we then get
$$
\split
&{\sh~{\theta^1_1+\theta^1_2\over 2}\sh~{\theta^1_3+\theta^1_2\over 2}
\over
\sh~{\theta^1_1\over 2}\sh~{\theta^1_3\over 2}}\ \le\
{(\sh~{\theta^1_1\over 2}\ch~{\theta^1_2\over 2}+
\ch~{\theta^1_1\over 2}\sh~{\theta^1_2\over 2})\sh~{\theta\over 2}
\over
\sh~{\theta^1_1\over 2}\sh~{\theta\over 4}}\ \le\ \\
&(\ch~{\theta^1_2\over 2}+
\sh~{\theta^1_2\over 2}\times
{\ch~{\theta^1_1\over 2}\over \sh~{\theta^1_1\over 2}})\times
{\sh~{\theta\over 2}\over \sh~{\theta\over 4}}
\le
(\ch~{\theta\over 4}+
\sh~{\theta\over 4}\times
{\ch~{\theta\over 4}\over \sh~{\theta\over 4}})\times
{\sh~{\theta\over 2}\over \sh~{\theta\over 4}}=\\
&4\ch^2{\theta\over 4}=2\ch~{\theta\over 2}+2.
\endsplit
$$
It follows
$$
(e_1,\,e_3)\le 8\ch~{\theta\over 2}+6.
\tag{8.2.20}
$$
We also have the inequality $2<(e_1,\,e_3)$ because otherwise
the lines $(A_0A_1)$ and $(A_2A_3)$ intersect each other and together
with the line $(A_1A_3)$ define a triangle containing the axis $l$
inside, which is impossible. Similarly, in considerations below,
one can prove inequalities of the type $2<(.\,,\,.)$ of the theorem,
and further we omit these considerations.

Further we consider two cases:

{\bf Case II:} Assume that
$$
\theta^1_4\ge \theta/2.
\tag{8.2.21}
$$
Then, by \thetag{8.1.19},
$$
(e_1,\,e_4)\le
4{\sh~{\theta^1_1+\theta\over 2}\sh~{\theta+\theta^1_4\over 2}
\over
\sh~{\theta^1_1\over 2}\sh~{\theta^1_4\over 2}}-2.
\tag{8.2.22}
$$
Using \thetag{8.2.17} and \thetag{8.2.21}, we then get like above
$$
{\sh~{\theta^1_1+\theta\over 2}\sh~{\theta+\theta^1_4\over 2}
\over
\sh~{\theta^1_1\over 2}\sh~{\theta^1_4\over 2}}\le
{\sh~{3\theta\over 4}\sh~{3\theta\over 4}\over
\sh~{\theta\over 4}\sh~{\theta\over 4}}=(2\ch~{\theta\over 2}+1)^2.
$$
It follows
$$
(e_1,\,e_4)\le 4(2\ch~{\theta\over 2}+1)^2-2.
\tag{8.2.23}
$$
Like above, one can prove that $2<(e_1,\,e_4)$. Elements
$e_1$, $e_3$ and $e_4$ generate the hyperbolic vector space
defining the hyperbolic plane. Otherwise, the lines
$(A_0A_1)$, $(A_2A_3)$ and $(A_3A_4)$ either have a common point or
are orthogonal to a line which is impossible. The Gram graph of
these three elements is connected because $2<(e_1,\,e_3)$ and
$2<(e_1,\,e_4)$. Similarly one can prove all statements about
connectideness of Gram graphs of the theorem, and further we
omit these considerations.

{\bf Case (III):} Assume that
$$
\theta^1_4\le \theta/2.
\tag{8.2.24}
$$
By definition of $\theta$, then $\theta^1_5\ge\theta/2$.
Since $\theta^1_3\ge \theta/2$, $\theta^1_4\le \theta/2$
and $\theta^1_5\ge \theta/2$, like for \thetag{8.2.20}, we get
$$
(e_3,\,e_5)\le
4{\sh~{\theta^1_3+\theta^1_4\over 2}\sh~{\theta^1_5+\theta^1_4\over 2}
\over
\sh~{\theta^1_3\over 2}\sh~{\theta^1_5\over 2}}-2\le
8\ch~{\theta\over 2}+6.
\tag{8.2.25}
$$
We have $\theta^1_2+\theta^1_3+\theta^1_4=\theta+\theta^1_4\le
3\theta/2$. It then follows like above, that
$$
(e_1,\,e_5)\le
4{\sh~{\theta^1_1+3\theta/2\over 2}\sh~{\theta^1_5+3\theta/2\over 2}
\over
\sh~{\theta^1_1\over 2}\sh~{\theta^1_5\over 2}}-2.
\tag{8.2.26}
$$
Here
$$
{\sh~{\theta^1_1+3\theta/2\over 2}\sh~{\theta^1_5+3\theta/2\over 2}
\over
\sh~{\theta^1_1\over 2}\sh~{\theta^1_5\over 2}}\le
{\sh^2\theta\over \sh^2{\theta\over 4}}=
8(\ch^3{\theta\over 2}+\ch^2{\theta\over 2}).
$$
Thus,
$$
(e_1,\,e_5)\le
32(\ch^3{\theta\over 2}+\ch^2{\theta\over 2})-2.
\tag{8.2.27}
$$

{\bf Cases (A):} Suppose that
$$
\theta\le \eta.
\tag{8.2.28}
$$
Then for the case (II) above we get the case (AII) of
the theorem, and for the case (III) above we get the case (AIII)
of the theorem.

{\bf Cases (B)}: Now assume that
$$
\theta\ge \eta.
\tag{8.2.29}
$$
Equivalently, then
$$
\theta^1_{i-1}+\theta^1_i\ge \eta\ \forall i\in I\ \
\text{and}\ \
\theta^2_{j-1}+\theta^2_j\ge \eta\ \ \forall j\in J.
\tag{8.2.30}
$$

We consider two cases:

{\bf Case (BII):} Under the condition \thetag{8.2.29} (equivalently,
\thetag{8.2.30}), we additionally suppose that
$$
\theta^1_i\ge \eta/2\ \forall\ i\in I\ \ \text{and}\ \
\theta^2_j\ge \eta/2\ \forall j\in J.
\tag{8.2.31}
$$

We then consider orthogonal projections of all vertices
$A_i$, $i\in I$, and $B_j$, $j\in J$, into the axis $l$. The set
of these projections gives a locally finite set of points in $l$.
We choose two neighbouring points of this set and consider
their preimages. We then get
(after changing numerations if necessary)
one of two possibilities (BII$_1$) or (BII$_2$)
which we describe below.

{\bf Case (BII$_1$):} There exist four vertices $A_0$, $A_1$, $A_2$
and $A_3$ and two vertices $B_0$ and $B_1$
such that orthogonal projections of $A_1$ and $A_2$
into $l$ belong to the orthogonal projection of the line $(B_0B_1)$
into $l$. Equivalently, two lines containing $A_1$ and $A_2$ and
orthogonal to $l$ intersect the line $(B_0B_1)$.
We want to show that then all conditions of the case (BII$_1$)
of the theorem are satisfied.

Using \thetag{8.1.21}, we get
$$
(e_1,\,f_1)=
4{\ch~{\theta^1_1+\theta_{12}\over 2}
\ch~{\theta^2_1-\theta_{12}\over 2}
\over
\sh~{\theta^1_1\over 2}\sh~{\theta^2_1\over 2}}-2.
\tag{8.2.32}
$$
where $-\theta^1_1\le \theta_{12}=\theta(A_{12},\,B_{12})\le \theta^2_1$
because there exists a line orthogonal to the axis $l$ and intersecting
both lines $(A_0A_1)$ and $(B_0B_1)$.
It follows
$$
{\ch~{\theta^1_1+\theta_{12}\over 2}
\ch~{\theta^2_1-\theta_{12}\over 2}
\over
\sh~{\theta^1_1\over 2}\sh~{\theta^2_1\over 2}}=
{{1\over 2}\left[\ch~({\theta^1_1+\theta^2_1\over 2})+
\ch~({\theta^1_1-\theta^2_1\over 2}+\theta_{12})\right]
\over
\sh~{\theta^1_1\over 2}\sh~{\theta^2_1\over 2}}\le
$$
$$
{{1\over 2}\left[\ch~({\theta^1_1+\theta^2_1\over 2})+
\max(\ch~({-\theta^1_1-\theta^2_1\over 2}),\
\ch~({\theta^1_1+\theta^2_1\over 2})\right]
\over
\sh~{\theta^1_1\over 2}\sh~{\theta^2_1\over 2}}=
{\ch~{\theta^1_1+\theta^2_1\over 2}
\over
\sh~{\theta^1_1\over 2}\sh~{\theta^2_1\over 2}}\le
$$
$$
{\ch~{\eta\over 2}
\over
\sh~{\eta\over 4}\sh~{\eta\over 4}}=
{2\ch~{\eta\over 2}\over \ch~{\eta\over 2}-1}=
2+{2\over \ch~{\eta\over 2}-1}
$$
because $\theta^1_1\ge \eta/2$ and $\theta^1_2\ge \eta/2$. It follows
$$
(e_1,\,f_1)\le {8\over \ch~{\eta\over 2}-1}+6.
\tag{8.2.33}
$$
Similarly one can prove that
$$
(e_2,\,f_1)\le {8\over \ch~{\eta\over 2}-1}+6\ \text{and}\
(e_3,\,f_1)\le {8\over \ch~{\eta\over 2}-1}+6.
\tag{8.2.34}
$$
It proves that we do get the case (BII$_1$) of the theorem.

{\bf Case (BII$_2$):} There exist three vertices $A_0$, $A_1$, $A_2$
and three vertices $B_0$, $B_1$, $B_2$
such that orthogonal projection of $A_1$ into $l$ belongs
to the orthogonal projection of the line $(B_0B_1)$
into $l$, and
orthogonal projection of $B_1$ into $l$ belongs
to the orthogonal projection of the line $(A_1A_2)$
into $l$.
We want to show that then all conditions of the case (BII$_2$)
of the theorem are valid. Really, like for the case (BII$_1$) above,
one can prove that
$$
(e_1,\,f_1),\ (e_2,\,f_1)
 \le {8\over \ch~{\eta\over 2}-1}+6
\tag{8.2.35}
$$
because there exists a line orthogonal to $l$ (it contains $A_1$)
which intersects the line $(B_0B_1)$ and intersects lines
$(A_0A_1)$ and $(A_1A_2)$. Similarly,
$$
(e_2,\,f_1),\ (e_2,\,f_2)
\le {8\over \ch~{\eta\over 2}-1}+6.
\tag{8.2.36}
$$

We now consider the opposite (to the case (BII)) case:

{\bf Case (BIII):} Under the condition \thetag{8.2.29}
(equivalently, \thetag{8.2.30}),
we assume that there exists  either $i\in I$ such that
$\theta^1_i\le \eta/2$ or $j\in J$ such that $\theta^2_j\le \eta/2$.
We then define similar to the invariant $\theta$ above, an invariant
$$
\theta_\eta=\max(
\max_{\{i\in I\,|\,\theta^1_i\le \eta/2\}}{\theta^1_i},\
\max_{\{j\in J\,|\,\theta^2_j\le \eta/2\}}{\theta^2_j}) .
\tag{8.2.37}
$$
Obviously, $\theta_\eta\le \eta/2$.
For any $i\in I$, if $\theta^1_i>\theta_\eta$ 
then $\theta^1_i>\eta/2$. For any $j\in J$,  
if $\theta^2_j>\theta_\eta$ then $\theta^2_j>\eta/2$.

For the case (BIII), after changing numerations if necessary,
there exist vertices
$A_0$, $A_1$, $A_2$, $A_3$ such that
$\theta^1_2=\theta_\eta$. Since $\theta_\eta\le \eta/2$, one has
$\theta^1_1\ge \eta/2$ because $\theta^1_1+\theta^1_2\ge \theta\ge \eta$.
Similarly, $\theta^1_3\ge \eta/2$. Like in \thetag{8.2.20},
we then have
$$
(e_1,\,e_3)\le 8\,\ch~{\eta\over 2}+6.
\tag{8.2.38}
$$
Suppose that there exists a vertex $B_j$ such that its orthogonal
projection into $l$ is contained in the orthogonal projection of the line
$(A_1A_2)$ into $l$. Then either the line $(B_{j-1}B_j)$ has the
angle $\theta^2_j\ge \theta/2\ge \eta/2$ or the line $(B_jB_{j+1})$
has the angle $\theta^2_{j+1}\ge \theta/2\ge \eta/2$. Remark that
both these lines have a line (containing $B_j$) which is orthogonal to
$l$ and intersects $(A_1A_2)$. If the vertex $B_j$ does not exist,
then there exists a line $(B_{j-1}B_j)$ such that the orthogonal
projection of the line $(A_1A_2)$ into $l$ is contained inside
the interval defined by orthogonal projections of $[B_{j-1}B_j]$.
It follows that $\theta^2_j>\theta^1_2=\theta_\eta$. By definition
of $\theta_\eta$, then $\theta^2_j>\eta/2$. Obviously, there
exists a line orthogonal to $l$ which intersects both lines
$(A_1A_2)$ and $(B_{j-1}B_j)$.

Finally (after changing numeration of $B_j$ if necessary)
we get that there exist vertices $B_0$ and $B_1$ such that
$\theta^2_1\ge \eta/2$ and there exists a line which is
orthogonal to $l$ and intersects both
lines $(A_1A_2)$ and $(B_0B_1)$. Let us show that
vertices $A_0$, $A_1$, $A_2$, $A_3$ and $B_0$, $B_1$ satisfy
conditions of the case (BIII) of the theorem.

We have proved \thetag{8.2.12}. By \thetag{8.1.21}, 
$$
(e_1,\,f_1)=4{\ch~{\theta^1_1+\theta_{12}\over 2}
\ch~{\theta^2_1-\theta_{12}\over 2}
\over
\sh~{\theta^1_1\over 2}\sh~{\theta^2_1\over 2}}-2
$$
where
$-\theta^1_2\le \theta_{12}\le \theta^1_2+\theta^2_1$ 
because there exists a line which is orthogonal to $l$ 
and intersects both lines $(A_1A_2)$ and $(B_0B_1)$. 
We have
$$
{\ch~{\theta^1_1+\theta_{12}\over 2}
\ch~{\theta^2_1-\theta_{12}\over 2}
\over
\sh~{\theta^1_1\over 2}\sh~{\theta^2_1\over 2}}=
{{1\over 2}\left[\ch~({\theta^1_1+\theta^2_1\over 2})+
\ch~({\theta^1_1-\theta^2_1\over 2}+\theta_{12})\right]
\over
\sh~{\theta^1_1\over 2}\sh~{\theta^2_1\over 2}}\le
$$
$$
{{1\over 2}\left[\ch~({\theta^1_1+\theta^2_1\over 2})+
\max(\ch~({\theta^1_1-2\theta^1_2-\theta^2_1\over 2}),\
\ch~({\theta^1_1+2\theta^1_2+\theta^2_1\over 2}))\right]
\over
\sh~{\theta^1_1\over 2}\sh~{\theta^2_1\over 2}}=
$$
$$
{{1\over 2}\left[\ch~{\theta^1_1+\theta^2_1\over 2}+
\ch~({\theta^1_1+2\theta^1_2+\theta^2_1\over 2})\right]
\over
\sh~{\theta^1_1\over 2}\sh~{\theta^2_1\over 2}}\le
$$
$$
{{1\over 2}\left[\ch~{\eta \over 2}+
\ch~(\eta )\right]
\over
\sh^2{\eta\over 4}}={\ch~{\eta\over 2}+2\ch^2{\eta\over 2}-1
\over \ch~{\eta\over 2}-1}=2\ch~{\eta\over 2}+3+
{2\over \ch~{\eta\over 2}-1}
$$
because $\theta^1_1\ge \eta/2$, $\theta^2_1\ge \eta/2$
and $\theta^1_2\le \eta/2$.
It follows,
$$
(e_1,\,f_1)\le
8\ch~{\eta\over 2}+10+
{8\over \ch~{\eta\over 2}-1}.
\tag{8.2.40}
$$
Similarly, we have
$$
(e_3,\,f_1)\le
8\ch~{\eta\over 2}+10+
{8\over \ch~{\eta\over 2}-1}.
\tag{8.2.41}
$$

This finishes the proof of Theorem 8.2.2.
\enddemo

Applying Theorem 8.2.2 to $\eta=\beta$ where 
$\beta=0.983986...$ is the solution of the equation
$$
32(\ch^3{\eta \over 2}+\ch^2{\eta \over 2})-2=
8\ch~{\eta\over 2}+10+{8\over \ch~{\eta\over 2}-1},
\tag{8.2.42}
$$
we get the important in applications 

\proclaim{Theorem 8.2.3 (about narrow part of restricted hyperbolic 
convex polygons in a hyperbolic plane)} 
For any restricted hyperbolic convex polygon $\M$ in a hyperbolic plane 
there are $\delta_1,\,\delta_2,\,\delta_3\in P(\M)$ with 
$\delta_1^2=\delta_2^2=\delta_3^2=-2$ such that 
$\delta_1,\,\delta_2,\,\delta_3$ generate the hyperbolic vector space 
(defining the hyperbolic space), they have connected Gram graph 
and 
$$
(\delta_i,\,\delta_j)\le 
32(\ch^3{\beta \over 2}+\ch^2{\beta \over 2})-2
< 83.7706
\tag{8.2.43}
$$
for all $1\le i,j,\le 3$. Here $\beta=0.983986...$ is the solution of 
the equation \thetag{8.2.42} 
\endproclaim

We remind that the same statement is valid for elliptic and parabolic 
polygons in hyperbolic plane with the better constant $62$ 
(see Theorems 4.1.5 and 4.2.2) and even $30$ 
(see Theorems 4.1.8 and 4.2.3) instead of $83.7706$.   

Similar method permits to prove Theorem 8.2.3 for restricted hyperbolic 
polyhedra in hyperbolic spaces of arbitrary dimension $n\ge 2$  
(here $3$ should be replaced by $n+1$).  

\subhead
8.3. Narrow parts of restricted hyperbolic convex polygons in a
hyperbolic plane: refined results
\endsubhead

Here we improve Theorems 8.2.2 and 8.2.3 to get better estimates 
depending on the angles of the polygon $\M$.

We first solve the equation
$$
\sh~(u-x)\sh~(v-x)=a\,\sh~u\,\sh~v
\tag{8.3.1}
$$
where $0\le a={1+\cos{\alpha}\over 2}=\cos^2{\alpha\over 2} \le 1$
and $u,v\ge 0$.

We have
$$
(\sh~u\,\ch~x-\ch~u\,\sh~x)(\sh~v\,\ch~x-\ch~v\,\sh~x)=
a(\ch^2~x-\sh^2~x)\sh~u\,\sh~v.
$$
Equivalently 
$$
\sh~u\,\ch~x\,\sh~v\,\ch~x-\sh~u\,\ch~x\,\ch~v\,\sh~x-
\ch~u\,\sh~x\,\sh~v\,\ch~x+\ch~u\,\sh~x\,\ch~v\,\sh~x=
$$
$$
=a(\ch^2~x-\sh^2~x)\,\sh~u\,\sh~v.
$$
Equivalently 
$$
\split
&(1-a)\sh~u\,\sh~v\,\ch^2~x-(\sh~u\,\ch~v+\ch~u\,\sh~v)\ch~x\,\sh~x+\\
&+(\ch~u\,\ch~v+a\,\sh~u\,\sh~v)\,\sh^2~x=0.
\endsplit
\tag{8.3.2}
$$
It follows that the smallest $x=g(\alpha,u,v)=
\widetilde{g}(a,u,v)$ is defined as follows:
$$
\th~x=\th~g(\alpha,u,v)=
{\sh~(u+v)-
\sqrt{\left(a\,\ch~(u+v)+(1-a)\,\ch~(u-v)\right)^2-1}\over
(a+1)\,\ch~(u+v)+(1-a)\,\ch~(u-v)}
\tag{8.3.3}
$$
and
$$
x=g(\alpha,u,v)=\widetilde{g}(a,u,v)=
{1\over 2}\ln{1+\th~g(\alpha,u,v)\over 1-\th~g(\alpha,u,v)}=
$$
$$
\split
{1\over 2}\ln &\Bigg[(a+1)\,\ch~(u+v)+(1-a)\,\ch~(u-v)+\sh~(u+v)-\\
&-\sqrt{\left(a\,\ch~(u+v)+(1-a)\,\ch~(u-v)\right)^2-1}\Bigg]\Bigg/\\
&\Bigg[(a+1)\,\ch~(u+v)+(1-a)\,\ch~(u-v)-\sh~(u+v)+\\
&+\sqrt{\left(a\,\ch~(u+v)+(1-a)\,\ch~(u-v)\right)^2-1}\Bigg].
\endsplit
\tag{8.3.4}
$$
In particular, for $u=v$, we get
$$
\th~x=\th~g(\alpha,u,u)=
{\sh~2u-
\sqrt{\left(a\,\ch~2u+(1-a)\right)^2-1}\over
(a+1)\,\ch~2u+(1-a)}
\tag{8.3.5}
$$
and
$$
\split
&x=g(\alpha,u,u)=\widetilde{g}(a,u,u)=\\
&{1\over 2}\ln{{(a+1)\,\ch~2u+(1-a)+\sh~2u-
\sqrt{\left(a\,\ch~2u+(1-a)\right)^2-1}\over
(a+1)\,\ch~2u+(1-a)-\sh~2u+
\sqrt{\left(a\,\ch~2u+(1-a)\right)^2-1}}}.\\
\endsplit
\tag{8.3.6}
$$

We introduce functions which use the function
$g(\alpha,u,v)=\widetilde{g}(a,u,v)$. We set
$$
f_{h1}(\alpha_1,\alpha_2,\theta)=\widetilde{f}_{h1}(a_1,a_2,\theta)=
{\sh^2({\theta\over 2}-g(\alpha_1,{\theta\over 4},{\theta\over 4})-
g(\alpha_2,{\theta\over 4},{\theta\over 4}))
\over
\sh^2({\theta\over 4})};
\tag{8.3.7}
$$
$$
\split
&f_{h2}(\alpha_1,\alpha_2,\alpha_3,\theta,t)=
\widetilde{f}_{h2}(a_1,a_2,a_3,\theta,t)=\\
&\sh~{\Big({3\theta\over 4}+{t\over 4}-
g(\alpha_1,{\theta+t\over 4},{\theta-t\over 4})-
g(\alpha_2,{\theta-t\over 4},{\theta+t\over 4})-
g(\alpha_3,{\theta+t\over 4},{\theta\over 4})\Big)\times}\\
&\sh~\Big({3\theta\over 4}-
g(\alpha_1,{\theta+t\over 4},{\theta-t\over 4})-
g(\alpha_2,{\theta-t\over 4},{\theta+t\over 4})-
g(\alpha_3,{\theta+t\over 4},{\theta\over 4})\Big)
\Bigg/ \\
&\Big(\sh~{\theta+t\over 4}\sh~{\theta\over 4}\Big)\ ;
\endsplit
\tag{8.3.8}
$$
$$
\split
&f_{h3}(\alpha_1,\alpha_2,\alpha_3,\alpha_4,\theta,z,w)=
\widetilde{f}_{h3}(a_1,a_2,a_3,a_4,\theta,z,w)=\\
&\sh~{\Big(\theta\hskip-2pt+\hskip-2pt
{z\hskip-2pt-\hskip-2ptw\over 4}-}\\
&{-g(\alpha_1,\hskip-2pt
{\theta \hskip-2pt+\hskip-2pt z\over 4},\hskip-2pt
{\theta\hskip-2pt-\hskip-2pt z\over 4})\hskip-2pt-\hskip-2pt
g(\alpha_2,\hskip-2pt{\theta\hskip-2pt-\hskip-2pt
z\over 4},\hskip-2pt
{\theta\hskip-2pt+\hskip-2pt z\over 4})\hskip-2pt-\hskip-2pt  
g(\alpha_3,\hskip-2pt {\theta\hskip-2pt+\hskip-2pt z\over 4},
\hskip-2pt
{\theta\hskip-2pt-\hskip-2pt w\over 4})\hskip-2pt-\hskip-2pt
g(\alpha_4,\hskip-2pt{\theta\hskip-2pt -\hskip-2pt w\over 4},\hskip-2pt
{\theta\hskip-2pt +\hskip-2pt w\over 4})
\Big)}\times \\
&\sh~{\hskip-2pt\Big(\theta\hskip-2pt-\hskip-2pt
g(\alpha_1,\hskip-2pt {\theta \hskip-3pt+\hskip-3pt z\over 4},\hskip-2pt
{\theta\hskip-3pt-\hskip-3pt z\over 4})\hskip-2pt-\hskip-2pt
g(\alpha_2,\hskip-1pt{\theta\hskip-3pt-\hskip-3pt
z\over 4},\hskip-1pt {\theta\hskip-3pt+
\hskip-3pt z\over 4})\hskip-2pt-\hskip-2pt  g(\alpha_3,\hskip-2pt
{\theta\hskip-3pt+\hskip-3pt z\over 4},\hskip-2pt
{\theta\hskip-3pt-\hskip-3pt w\over 4})\hskip-2pt-\hskip-2pt
g(\alpha_4,\hskip-2pt
{\theta\hskip-3pt -\hskip-3pt w\over 4},\hskip-2pt
{\theta\hskip-2pt +\hskip-2pt w\over 4})
\Big)}\Bigg/\\
&\Big(\sh~{\theta+z\over 4}\sh~{\theta+w\over 4}\Big)\ .
\endsplit
\tag{8.3.9}
$$

We mention the following properties of the
function $g(\alpha,u,v)$. We have
$g(\alpha,u,v)$ $=g(\alpha,v,u)$ and $g(\alpha,u,v)\le u,\,v$.
For fixed $\alpha$ and $u$, the functions $g(\alpha,u,v)$ and
$v-g(\alpha,u,v)$ increase if $v$ increases. It follows that
the number $g(\alpha,u,+\infty)$ is defined and
$$
g(\alpha,u,+\infty)=u-\text{arc\,sh}~(a\,\sh~u).
\tag{8.3.10}
$$
These properties follow from geometric interpretation below
of the equation \thetag{8.3.1}, see \thetag{8.3.31}.

We also denote by
$$
\text{Ach}~(x)=
\cases
0, &\text{if $x<1$},\\
$\text{arc\,ch}~(x)$, &\text{if $x\ge 1$}.
\endcases
\tag{8.3.11}
$$

We prove the following refined result (Theorem 8.3.1 below)
where all estimates of Theorem 8.2.2 will be improved depending on
angles of $\M$. In Theorem 8.3.1, we did exact calculations of
maxima using computer only for angles $\alpha_i$ of the form
${\pi\over n}$ where $n=2,\,3,\,4,\,5,\,6$ or $+\infty$. Only 
this result we need for reflective lattices. 
But we believe that the same results are valid for arbitrary
$0\le \alpha_i<\pi$ (corresponding considerations should be
similar to the proof of Theorem 4.1.5).

\proclaim{Theorem 8.3.1 (about narrow parts of types (II) and (III)
of restricted hyperbolic convex polygons)}
We fix a constant $\eta$ such that $3\ge \eta >0$ (to get good
estimates, one can take the constant
$\eta=\eta_0=2\text{arc}\,\ch~{3\over 2}=$
\newline
$=1.92484730023841378999103565...$;
equivalently, $\ch~{\eta_0\over 2}={3\over 2}$).

Any infinite restricted hyperbolic convex polygon $\M$ has a narrow
part of one of types
(AII), (AIII), (BII$_1$), (BII$_2$), (BIII) defined below:

(AII) There exist five consecutive vertices
$A_0$, $A_1$, $A_2$, $A_3$ and $A_4$ of one half of $\M$
such that for orthogonal
vectors $\delta_1$, $\delta_2$, $\delta_3$ and $\delta_4$
to lines $(A_0A_1)$, $(A_1A_2)$, $(A_2A_3)$ and $(A_3A_4)$
respectively directed outwards of $\M$ and with
$\delta_1^2=\delta_2^2=\delta_3^2=\delta_4^2=-2$, one has
$(\delta_i,\,\delta_{i+1})=2\cos{\alpha_i}$, $i=1,\,2,\,3$,
where $\alpha_i=\angle~A_{i-1}A_iA_{i+1}$, and
$$
2<(\delta_1,\,\delta_3)
\le 4\,\max_{0\le \theta\le \eta}{f_{h1}(\alpha_1,\alpha_2,\theta)}-2=
4\,f_{h1}(\alpha_1,\alpha_2,\eta)-2,
\tag{8.3.12}
$$
$$
2<(\delta_1,\,\delta_4)\le
4\,\max_{0\le \theta\le \eta}{\ \max_{0\le t\le \theta}
{f_{h2}(\alpha_1,\alpha_2,\alpha_3,\theta,t)}}-2=
$$
$$
4\,\max_{0\le \theta\le \eta}
{\big(f_{h2}(\alpha_1,\alpha_2,\alpha_3,\theta,0),
f_{h2}(\alpha_1,\alpha_2,\alpha_3,\theta,\theta)\big)}-2=
$$
$$
4\,\max{\big(f_{h2}(\alpha_1,\alpha_2,\alpha_3,\eta,0),\,
f_{h2}(\alpha_1,\alpha_2,\alpha_3,\eta,\eta)\big)}-2.
\tag{8.3.13}
$$
Moreover, the set $\{\delta_1,\,\delta_3,\,\delta_4\}$
generates the 3-dimensional hyperbolic vector space and
has a connected Gram graph.
Additionally, open half-planes $\Ha_{\delta_i}^+$,
$i=1,\,2,\,3,\,4$, contain the axis $l$.

(AIII) There exist six consecutive vertices
$A_0$, $A_1$, $A_2$, $A_3$, $A_4$ and $A_5$ of one half of $\M$
such that for orthogonal vectors $\delta_1$, $\delta_2$, $\delta_3$,
$\delta_4$ and $\delta_5$
to lines $(A_0A_1)$, $(A_1A_2)$, $(A_2A_3)$, $(A_3A_4)$ and
$(A_4A_5)$ respectively directed
outwards of $\M$ and with
$\delta_1^2=\delta_2^2=\delta_3^2=\delta_4^2=\delta_5^2=-2$
one has $(\delta_i,\delta_{i+1})=2\cos{\alpha_i}$, $i=1,\,2,\,3,\,4$,
where $\alpha_i=\angle~A_{i-1}A_iA_{i+1}$, and
$$
2<(\delta_1,\delta_3)
\le 4\,\max_{0\le \theta \le \eta}{f_{h1}(\alpha_1,\alpha_2,\theta)}-2=
4\,f_{h1}(\alpha_1,\alpha_2,\eta)-2,
\tag{8.3.14}
$$
$$
2<(\delta_3,\,\delta_5)
\le 4\,\max_{0\le \theta \le \eta}{f_{h1}(\alpha_1,\alpha_2,\theta )}-2=
4\,f_{h1}(\alpha_1,\alpha_2,\eta)-2,
\tag{8.3.15}
$$
and
$$
2<(\delta_1,\,\delta_5)\le
4\,\max_{0\le \theta\le \eta}{\ \max_{0\le w\le z\le \theta}
{f_{h3}(\alpha_1,\alpha_2,\alpha_3,\alpha_4,\theta,z,w)}}-2=
$$
$$
\split
&4\max_{0\le \theta\le \eta}
{\big(f_{h3}(\alpha_1,\hskip-1pt  \alpha_2,\hskip-1pt \alpha_3,
\hskip-1pt \alpha_4,\hskip-1pt\theta,\hskip-1pt0,\hskip-1pt0),\,
f_{h3}(\alpha_1,\hskip-1pt \alpha_2,\hskip-1pt 
\alpha_3,\hskip-1pt \alpha_4,\hskip-1pt\theta,\hskip-1pt
\theta,\hskip-1pt 0),\
f_{h3}(\alpha_1,\hskip-1pt\alpha_2,\hskip-1pt\alpha_3,\hskip-1pt\alpha_4,
\hskip-1pt\theta,\hskip-1pt \theta,\hskip-1pt \theta)\big)}\\
&-2=
\endsplit
$$
$$
\split
&4\max
{\big(f_{h3}(\alpha_1,\hskip-1pt  \alpha_2,\hskip-1pt \alpha_3,
\hskip-1pt \alpha_4,\hskip-1pt\eta,\hskip-1pt0,\hskip-1pt0),\,
f_{h3}(\alpha_1,\hskip-1pt \alpha_2,\hskip-1pt 
\alpha_3,\hskip-1pt \alpha_4,\hskip-1pt\eta,\hskip-1pt\eta,\hskip-1pt0),\
f_{h3}(\alpha_1,\hskip-1pt\alpha_2,\hskip-1pt\alpha_3,\hskip-1pt\alpha_4,
\hskip-1pt\eta,\hskip-1pt \eta,\hskip-1pt \eta)\big)}\\
&-2.
\endsplit
\tag{8.3.16}
$$
Moreover, the
set $\{\delta_1,\,\delta_3,\,\delta_5\}$ generates the 3-dimensional
hyperbolic vector space and has a connected Gram graph.
Additionally, open half-planes $\Ha_{\delta_i}^+$,
$i=1,\,2,\,3,\,4$, $5$, contain the axis $l$.

(BII$_1$) There exist four consecutive vertices
$A_0$, $A_1$, $A_2$, $A_3$ of one half of $\M$
and two consecutive vertices $B_0$, $B_1$ of another half of $\M$
such that for orthogonal
vectors $\delta_1$, $\delta_2$, $\delta_3$ and $\delta_4$
to lines $(A_0A_1)$, $(A_1A_2)$, $(A_2A_3)$ and $(B_0B_1)$
respectively directed outwards of $\M$ and with
$\delta_1^2=\delta_2^2=\delta_3^2=\delta_4^2=-2$, one has
$(\delta_i,\,\delta_{i+1})=2\cos{\alpha_i}$, $i=1,\,2$,
where $\alpha_i=\angle~A_{i-1}A_iA_{i+1}$, and
$$
2<(\delta_1,\,\delta_4)\le
{8\over \ch~{\eta\over 2}-1}+6,
\tag{8.3.17}
$$
$$
2<(\delta_2,\,\delta_4)\le
{8\over \ch~{\eta\over 2}-1}+6,
\tag{8.3.18}
$$
$$
2<(\delta_3,\,\delta_4)\le
{8\over \ch~{\eta\over 2}-1}+6,
\tag{8.3.19}
$$
$$
\split
&\text{Ach}~\left({(\ch~{\eta\over 2}\hskip-2pt-\hskip-2pt 1)
((\delta_1,\,\delta_4)
\hskip-1pt+\hskip-1pt2)\over 4}\hskip-1pt-\hskip-1pt
\ch~{\eta\over 2}\right)\hskip-1pt+\hskip-1pt
\text{Ach}~\left({(\ch~{\eta\over 2}\hskip-2pt-\hskip-2pt1)
((\delta_2,\,\delta_4)
\hskip-1pt+\hskip-1pt 2)\over 4}\hskip-1pt-\hskip-1pt
\ch~{\eta\over 2}\right)\ \le\\
&{\eta\over 2}+2g(\alpha_1,\,{\eta\over 4},+\infty),
\endsplit
\tag{8.3.20}
$$
$$
\split
&\text{Ach}~\left({(\ch~{\eta\over 2}\hskip-2pt-\hskip-2pt 1)
((\delta_2,\,\delta_4)
\hskip-1pt+\hskip-1pt2)\over 4}\hskip-1pt-\hskip-1pt
\ch~{\eta\over 2}\right)\hskip-1pt+\hskip-1pt
\text{Ach}~\left({(\ch~{\eta\over 2}\hskip-2pt-\hskip-2pt1)
((\delta_3,\,\delta_4)
\hskip-1pt+\hskip-1pt 2)\over 4}\hskip-1pt-\hskip-1pt
\ch~{\eta\over 2}\right)\ \le\\
&{\eta\over 2}+2g(\alpha_2,\,{\eta\over 4},+\infty).
\endsplit
\tag{8.3.21}
$$
Moreover, the
set $\{\delta_1,\,\delta_2,\,\delta_4\}$ generates the 3-dimensional
hyperbolic vector space and has a connected Gram graph. The same
is valid for $\{\delta_2,\,\delta_3,\,\delta_4\}$.
Additionally, open half-planes $\Ha_{\delta_i}^+$,
$i=1,\,2,\,3$, and closed half-plane $\Ha_{\delta_4}^+$
contain the axis $l$; both lines orthogonal to the axis $l$ and
containing $A_1$ and $A_2$ intersect the line $(B_0B_1)$.

(BII$_2$) There exist three consecutive vertices
$A_0$, $A_1$, $A_2$ of one half of $\M$
and three consecutive vertices $B_0$, $B_1$, $B_2$ of
another half of $\M$ such that for orthogonal
vectors $\delta_1$, $\delta_2$, $\delta_3$ and $\delta_4$
to lines $(A_0A_1)$, $(A_1A_2)$, $(B_0B_1)$ and $(B_1B_2)$
respectively directed outwards of $\M$ and with
$\delta_1^2=\delta_2^2=\delta_3^2=\delta_4^2=-2$, one has
$(\delta_1,\,\delta_2)=2\cos{\alpha}$ where $\alpha=\angle~A_0A_1A_2$,
$(\delta_3,\,\delta_4)=2\cos{\beta}$ where
$\beta=\angle~B_0B_1B_2$, and
$$
2<(\delta_1,\,\delta_3)\le
{8\over \ch~{\eta\over 2}-1}+6,
\tag{8.3.22}
$$
$$
2<(\delta_2,\,\delta_3)\le
{8\over \ch~{\eta\over 2}-1}+6,
\tag{8.3.23}
$$
$$
2<(\delta_2,\,\delta_4)\le
{8\over \ch~{\eta\over 2}-1}+6,
\tag{8.3.24}
$$
$$
\split
&\text{Ach}~\left({(\ch~{\eta\over 2}\hskip-2pt-\hskip-2pt 1)
((\delta_1,\,\delta_3)
\hskip-1pt+\hskip-1pt2)\over 4}\hskip-1pt-\hskip-1pt
\ch~{\eta\over 2}\right)\hskip-1pt+\hskip-1pt
\text{Ach}~\left({(\ch~{\eta\over 2}\hskip-2pt-\hskip-2pt1)
((\delta_2,\,\delta_3)
\hskip-1pt+\hskip-1pt 2)\over 4}\hskip-1pt-\hskip-1pt
\ch~{\eta\over 2}\right)\ \le\\
&{\eta\over 2}+2g(\alpha,\,{\eta\over 4},+\infty),
\endsplit
\tag{8.3.25}
$$
$$
\split
&\text{Ach}~\left({(\ch~{\eta\over 2}\hskip-2pt-\hskip-2pt 1)
((\delta_3,\,\delta_2)
\hskip-1pt+\hskip-1pt2)\over 4}\hskip-1pt-\hskip-1pt
\ch~{\eta\over 2}\right)\hskip-1pt+\hskip-1pt
\text{Ach}~\left({(\ch~{\eta\over 2}\hskip-2pt-\hskip-2pt1)
((\delta_4,\,\delta_2)
\hskip-1pt+\hskip-1pt 2)\over 4}\hskip-1pt-\hskip-1pt
\ch~{\eta\over 2}\right)\ \le\\
&{\eta\over 2}+2g(\beta,\,{\eta\over 4},+\infty).
\endsplit
\tag{8.3.26}
$$
Moreover, the
set $\{\delta_1,\,\delta_2,\,\delta_3\}$ generates the 3-dimensional
hyperbolic vector space and has a connected Gram graph. The same
is valid for $\{\delta_2,\,\delta_3,\,\delta_4\}$.
Additionally, open half-planes $\Ha_{\delta_i}^+$,
$i=1,\,2,\,3,\,4$, contain the axis $l$; the line
orthogonal to the axis $l$ and containing $A_1$ intersects
$(B_0B_1)$, and the line orthogonal to $l$ and containing $B_1$
intersects $(A_1A_2)$.

(BIII) There exist four consecutive vertices
$A_0$, $A_1$, $A_2$, $A_3$ of one half of $\M$,
and two consecutive vertices $B_0$, $B_1$ of
another half of $\M$ such that for orthogonal
vectors $\delta_1$, $\delta_2$, $\delta_3$ and $\delta_4$
to lines $(A_0A_1)$, $(A_1A_2)$, $(A_2A_3)$ and $(B_0B_1)$
respectively directed outwards of $\M$ and with
$\delta_1^2=\delta_2^2=\delta_3^2=\delta_4^2=-2$, one has
$(\delta_i,\,\delta_{i+1})=2\cos{\alpha_i}$ where
$\alpha_i=\angle~A_{i-1}A_iA_{i+1}$, $i=1,\,2$, and
$$
2<(\delta_1,\,\delta_3)
\le 4\,\max_{0\le x\le \eta}{f_{h1}(\alpha_1,\alpha_2,x)}-2=
4\,f_{h1}(\alpha_1,\alpha_2,\eta)-2,
\tag{8.3.27}
$$
$$
2<(\delta_1,\,\delta_4)\le
4{\ch~{\eta\over 2}+
\ch~(\eta - 2\,g(\alpha_1,{\eta\over 4},{\eta\over 4}))
\over
\ch~{\eta\over 2}-1}-2,
\tag{8.3.28}
$$
$$
2<(\delta_3,\,\delta_4)\le
4{\ch~{\eta\over 2}+
\ch~(\eta - 2\,g(\alpha_2,{\eta\over 4},{\eta\over 4}))
\over
\ch~{\eta\over 2}-1}-2,
\tag{8.3.29}
$$
$$
\split
&\text{Ach}~\left({(\ch~{\eta\over 2}\hskip-1pt-1\hskip-1pt)
((\delta_1,\,\delta_4)\hskip-1pt+\hskip-1pt 2)\over 4}
\hskip-1pt-\hskip-1pt\ch~{\eta\over 2}\right)+
\text{Ach}~\left({(\ch~{\eta\over 2}\hskip-1pt-\hskip-1pt1)
((\delta_3,\,\delta_4)\hskip-1pt+\hskip-1pt2)\over 4}
\hskip-1pt-\hskip-1pt\ch~{\eta\over 2}\right)\le\\
&{\eta\over 2}+\max_{0\le u\le {\eta\over 2}}
{\left(u+2\,|g(\alpha_1,{\eta\over 4},{u\over 2})-
g(\alpha_2,{\eta\over 4},{u\over 2})|\right)}.
\endsplit
\tag{8.3.30}
$$
Moreover, the
set $\{\delta_1,\,\delta_2,\,\delta_4\}$ generates the
3-dimensional hyperbolic vector space and has a connected Gram graph.
The same is valid for the set $\{\delta_2,\,\delta_3,\,\delta_4\}$.
Additionally, open half-planes $\Ha_{\delta_i}^+$,
$i=1,\,2,\,3,\,4$, contain the axis $l$; there exists a line
orthogonal to the axis $l$ and intersecting the lines
$(A_1A_2)$ and $(B_0B_1)$.
\endproclaim

\demo{Proof} Compare the proof of
Theorem 4.1.5. The general line of the proof is the same as in
Theorem 8.2.2, but all estimates should be improved. Thus
we follow the proof of Theorem 8.2.2 and notations of the proof.

We introduce additional angles $\theta^1_{(i+1)i}=
\theta(A_{(i+1)1},\,A_{i2})$ and
$\theta^2_{(i+1)i}=\theta(B_{(i+1)1},$ $B_{i2})$.
We denote $\alpha_i=\angle A_{i-1}A_iA_{i+1}$, $a_i=\cos^2{\alpha_i\over 2}$.
Similarly, we denote $\beta_j=\angle B_{j-1}B_jB_{j+1}$,
$b_j=\cos^2{\beta_j\over 2}$. By \thetag{8.1.20}, we get
$$
\theta^1_{(i+1)i}=2g(\alpha_i,\,{\theta^1_i\over 2},\,
{\theta^1_{i+1}\over 2}),\ \
\theta^2_{(j+1)j}=2g(\beta_j,\,{\theta^2_j\over 2},\,
{\theta^2_{j+1}\over 2}),
\tag{8.3.31}
$$

Let us consider the proof of the estimate \thetag{8.2.20}.
For lines $m_1=(A_0A_1)$ and $m_3=(A_2A_3)$ with orthogonal
vectors $e_1$, $e_3$ we have (using  \thetag{8.1.19})
$$
(e_1,\,e_3)=
4{\sh~{\theta^1_1+\theta^1_2-\theta^1_{21}-\theta^1_{32}\over 2}
\sh~{\theta^1_3+\theta^1_2-\theta^1_{21}-\theta^1_{32}\over 2}
\over
\sh~{\theta^1_1\over 2}\sh~{\theta^1_3\over 2}}-2=
$$
$$
\split
&4\,{\sh\left({\theta^1_1\over 2}\hskip-2pt+\hskip-2pt
{\theta^1_2\over 2}\hskip-2pt-\hskip-2pt
g(\alpha_1,{\theta^1_1\over 2},{\theta^1_2\over 2})\hskip-2pt-\hskip-2pt
g(\alpha_2,{\theta^1_2\over 2},{\theta^1_3\over 2})\right)
\hskip-2pt\sh\left({\theta^1_3\over 2}\hskip-2pt+\hskip-2pt
{\theta^1_2\over 2}\hskip-2pt-\hskip-2pt
g(\alpha_1,{\theta^1_1\over 2},{\theta^1_2\over 2})\hskip-2pt-\hskip-2pt
g(\alpha_2,{\theta^1_2\over 2},{\theta^1_3\over 2})\right)
\over
\sh{\theta^1_1\over 2}\,\sh~{\theta^1_3\over 2}}\\
&-2.
\endsplit
\tag{8.3.32}
$$
Here $\theta^1_1\ge \theta/2$, $\theta^1_2\le \theta/2$ and
$\theta^1_3\ge \theta/2$. In particular,
$\theta^1_1\ge \theta^1_2$ and $\theta^1_3\ge \theta^1_2$.
Then we shall prove that
$$
(e_1,\,e_3)\le
4\,f_{h1}(\alpha_1,\alpha_2,2\theta^1_2)-2
\le 4\,\max_{0\le x\le \theta}{f_{h1}(\alpha_1,\alpha_2,x)}-2.
\tag{8.3.33}
$$
Equivalently, \thetag{8.3.32} increases if one puts
$\theta^1_1:=\theta^1_2$ and $\theta^1_3:=\theta^1_2$.

We argue
like in proof of Theorem 4.1.4. We know that
$\theta^1_1\ge \theta/2\ge \theta^1_2$ and
$\theta^1_3\ge \theta/2\ge \theta^1_2/2$.
We consider lines $m_i=(A_{i-1}A_i)=(A_{i1}A_{i2})$.
We know that lines $m_1$, $m_3$ do not intersect. Equivalently,
$(e_1,\,e_3)>2$.

Since $\theta^1_1\ge \theta^1_2$, there exists a
line $m_1'=(A_{11}'A_{12}')$ with terminals
$A_{11}'$ and $A_{12}'$ at infinity such that the line
$m_1'$ is contained in $\Ha_{-e_1}^+$, points $A_{11}$, $A_{11}'$
and $A_{12}$, $A_{12}'$ are contained in the same half-planes bounded
by $m_2$, the line $m_1'$ has the same angle
$\alpha_1=\angle A_{11}'A_1'A_{22}$ (as $m_1$) with the line $m_2$
(we denote by $A_1'$ their intersection point),
and the angle $\theta(A_{11}',A_{12}')$ of $m_1'$ is equal to
$\theta^1_2$. We
denote by $e_1'$ the orthogonal vector to $m_1'$
directed outwards of $\M$ and with $(e_1')^2=-2$, and we denote
$\theta_{21}'=\theta(A_{21},A_{12}')$.

Similarly,
there exists a line $m_3'=(A_{31}'A_{32}')$ with terminals
$A_{31}'$ and $A_{32}'$ at infinity such that the line
$m_3'$ is contained in $\Ha_{-e_3}^+$, points $A_{31}$, $A_{31}'$
and $A_{32}$, $A_{32}'$ are contained in the same half-planes bounded
by $m_2$, the line $m_3'$ has the same angle
$\alpha_2=\angle A_{21}A_2'A_{32}'$ (as $m_3$)
with the line $m_2$ (we denote by $A_2'$ their point of intersection),
and the angle $\theta(A_{31}',A_{32}')$ of $m_3'$ is equal to
$\theta^1_2$. We denote by $e_3'$ the orthogonal vector to $m_3'$
directed outwards of $\M$ and with $(e_3')^2=-2$, and we denote
$\theta_{32}'=\theta(A_{31}',A_{22})$.

Since the lines $m_1$ and $m_3$ do not intersect, by our construction,
any interval with terminals at $m_1'$ and $m_3'$ intersects
both lines $m_1$ and $m_3$. It follows that distance between lines
$m_1'$ and $m_3'$ is greater than distance between lines $m_1$
and $m_3$. It follows $(e_1,e_3)\le (e_1',e_3')$ and
by \thetag{8.1.19} and \thetag{8.1.20}, we have
$$
(e_1,\,e_3)\le (e_1',\,e_3')=4\,f_{h1}(\alpha_1,\alpha_2,2\theta^1_2)-2.
$$
Since $2\theta^1_2\le \theta$, we get \thetag{8.3.33}.

Thus, instead of \thetag{8.2.20}, we get \thetag{8.3.33}.
Following the line of the proof of Theorem 8.2.2,
using \thetag{8.3.33}, we get estimates
\thetag{8.3.12}, \thetag{8.3.14},
\thetag{8.3.15}, \thetag{8.3.27}.

Now we consider the proof of \thetag{8.2.23}.
Using \thetag{8.1.19} and \thetag{8.1.20}, we have
$$
(e_1,\,e_4)=
4\,{\sh~{\theta^1_1+\theta-\theta^1_{21}-\theta^1_{32}-\theta^1_{43}
\over 2}
\sh~{\theta^1_4+\theta-\theta^1_{21}-\theta^1_{32}-\theta^1_{43}
\over 2}
\over
\sh~{\theta^1_1\over 2}\sh~{\theta^1_4\over 2}}-2=
$$
$$
\split
4\,&\sh~{\Big({\theta^1_1\over 2}+{\theta\over 2}-
g(\alpha_1,{\theta^1_1\over 2},{\theta^1_2\over 2})-
g(\alpha_2,{\theta^1_2\over 2},{\theta^1_3\over 2})-
g(\alpha_3,{\theta^1_3\over 2},{\theta^1_4\over 2})\Big)}\times\\
&\sh~{\Big({\theta^1_4\over 2}+{\theta\over 2}-
g(\alpha_1,{\theta^1_1\over 2},{\theta^1_2\over 2})-
g(\alpha_2,{\theta^1_2\over 2},{\theta^1_3\over 2})-
g(\alpha_3,{\theta^1_3\over 2},{\theta^1_4\over 2})\Big)}\Bigg/\\
&\Big(\sh~{\theta^1_1\over 2}\sh~{\theta^1_4\over 2}\Big)-2\
\endsplit
\tag{8.3.34}
$$
where $\theta^1_2+\theta^1_3=\theta$, $\theta^1_2\le \theta^1_3$,
$\theta^1_1\ge \theta^1_3$ and $\theta^1_4\ge \theta/2$. Like above,
one can prove that \thetag{8.3.34} increases if one puts
$\theta^1_1:=\theta^1_3$ and $\theta^1_4=:\theta/2$. Thus, we get
$$
\split
&\ \ \ (e_1,\,e_4)\le\\
4\,&\sh~{\Big({\theta^1_3\over 2}+{\theta\over 2}-
g(\alpha_1,{\theta^1_3\over 2},{\theta^1_2\over 2})-
g(\alpha_2,{\theta^1_2\over 2},{\theta^1_3\over 2})-
g(\alpha_3,{\theta^1_3\over 2},{\theta\over 4})\Big)}\times\\
&\sh~{\Big({3\theta\over 4}-
g(\alpha_1,{\theta^1_3\over 2},{\theta^1_2\over 2})-
g(\alpha_2,{\theta^1_2\over 2},{\theta^1_3\over 2})-
g(\alpha_3,{\theta^1_3\over 2},{\theta\over 4})\Big)}\Bigg/
\Big(\sh~{\theta^1_3\over 2}\sh~{\theta\over 4}\Big)-2\
\endsplit
\tag{8.3.35}
$$
We set $t=\theta^1_3-\theta^1_2$ where $0\le t\le \theta$.
Then $\theta^1_2=(\theta-t)/2$, $\theta^1_3=(\theta+t)/2$.
We get
$$
\split
&(e_1,\,e_4)\le\\
&4\,\sh~{\Big({3\theta\over 4}+{t\over 4}-
g(\alpha_1,{\theta+t\over 4},{\theta-t\over 4})-
g(\alpha_2,{\theta-t\over 4},{\theta+t\over 4})-
g(\alpha_3,{\theta+t\over 4},{\theta\over 4})\Big)\times}\\
&\sh~\Big({3\theta\over 4}-
g(\alpha_1,{\theta+t\over 4},{\theta-t\over 4})-
g(\alpha_2,{\theta-t\over 4},{\theta+t\over 4})-
g(\alpha_3,{\theta+t\over 4},{\theta\over 4})\Big)
\Bigg/ \\
&\Big(\sh~{\theta+t\over 4}\sh~{\theta\over 4}\Big)\ -\ 2\ .
\endsplit
\tag{8.3.36}
$$
Thus, instead of \thetag{8.2.23}, we get (see \thetag{8.3.8})
$$
(e_1,\,e_4)\ \le\ 4\,\max_{0\le t\le \theta}
{f_{h2}(\alpha_1,\alpha_2,\alpha_3,\theta,t)}-2.
\tag{8.3.37}
$$

Following the proof of Theorem 8.2.2, using \thetag{8.3.37},
we get \thetag{8.3.13}.

Now we consider the proof of \thetag{8.2.27}.
Using \thetag{8.1.19} and \thetag{8.1.20}, we have
$$
(e_1,\,e_5)=
4\,{\sh~{\theta^1_1+\theta+\theta^1_4-
\theta^1_{21}-\theta^1_{32}-\theta^1_{43}-\theta^1_{54}
\over 2}
\sh~{\theta^1_5+\theta+\theta^1_4-
\theta^1_{21}-\theta^1_{32}-\theta^1_{43}-\theta^1_{54}
\over 2}
\over
\sh~{\theta^1_1\over 2}\sh~{\theta^1_5\over 2}}-2=
$$
$$
\split
4\,&\sh~{\Big({\theta^1_1\over 2}\hskip-2pt+\hskip-2pt
{\theta\over 2}\hskip-2pt+\hskip-2pt{\theta^1_4\over 2}\hskip-2pt-\hskip-2pt
g(\alpha_1,{\theta^1_1\over 2},{\theta^1_2\over 2})\hskip-1pt-\hskip-1pt
g(\alpha_2,{\theta^1_2\over 2},{\theta^1_3\over 2})\hskip-1pt-\hskip-1pt
g(\alpha_3,{\theta^1_3\over 2},{\theta^1_4\over 2})\hskip-1pt-\hskip-1pt
g(\alpha_4,{\theta^1_4\over 2},{\theta^1_5\over 2})
\Big)}\times\\
&\sh~{\Big({\theta^1_5\over 2}\hskip-2pt+\hskip-2pt
{\theta\over 2}\hskip-1pt+\hskip-1pt{\theta^1_4\over 2}\hskip-2pt-\hskip-2pt
g(\alpha_1,{\theta^1_1\over 2},{\theta^1_2\over 2})\hskip-1pt-\hskip-1pt
g(\alpha_2,{\theta^1_2\over 2},{\theta^1_3\over 2})\hskip-1pt-
g(\alpha_3,{\theta^1_3\over 2},{\theta^1_4\over 2})\hskip-1pt-\hskip-1pt
g(\alpha_4,{\theta^1_4\over 2},{\theta^1_5\over 2})
\Big)}\Bigg/\\
&\Big(\sh~{\theta^1_1\over 2}\sh~{\theta^1_5\over 2}\Big)-2\
\endsplit
\tag{8.3.38}
$$
where $\theta^1_2+\theta^1_3=\theta$, $\theta^1_2\le \theta^1_3$,
$\theta^1_1\ge \theta^1_3$ and $\theta^1_4\le \theta/2$. Then
$\theta^1_5\ge \theta/2$ and $\theta^1_4+\theta^1_5\ge \theta$.
Like above, \thetag{8.3.38} increases if one puts
$\theta^1_1:=\theta^1_3$ and $\theta^1_5:=\theta-\theta^1_4$.
Thus, we get
$$
\split
&\ \ \ (e_1,\,e_5)\le\\
4\,&\sh~{\Big({\theta^1_3\over 2}\hskip-2pt+\hskip-2pt
{\theta\over 2}\hskip-2pt+\hskip-2pt{\theta^1_4\over 2}
\hskip-2pt-\hskip-2pt
g(\alpha_1,\hskip-1pt{\theta^1_3\over 2},
\hskip-1pt{\theta^1_2\over 2})\hskip-2pt-\hskip-2pt
g(\alpha_2,\hskip-1pt {\theta^1_2\over 2},\hskip-1pt
{\theta^1_3\over 2})\hskip-1pt-\hskip-1pt g(\alpha_3,\hskip-1pt
{\theta^1_3\over 2},\hskip-1pt
{\theta^1_4\over 2})\hskip-1pt-\hskip-1pt
g(\alpha_4,\hskip-1pt
{\theta^1_4\over 2},\hskip-1pt
{\theta\hskip-2pt-\hskip-2pt\theta^1_4\over 2})
\Big)}\times\\
\times &\sh~{\Big(\theta-
g(\alpha_1,{\theta^1_3\over 2},{\theta^1_2\over 2})-
g(\alpha_2,{\theta^1_2\over 2},{\theta^1_3\over 2})-
g(\alpha_3,{\theta^1_3\over 2},{\theta^1_4\over 2})-
g(\alpha_4,{\theta^1_4\over 2},{\theta-\theta^1_4\over 2})
\Big)}\Bigg/\\
&\Big(\sh~{\theta^1_3\over 2}\sh~{\theta-\theta^1_4\over 2}\Big)\ -2\ .
\endsplit
\tag{8.3.39}
$$
We set $z=\theta^1_3-\theta^1_2$ and $w=\theta^1_5-\theta^1_4$ where
$\theta\ge z \ge w\ge 0$ (this follows from definition). Then
$\theta^1_1=\theta^1_3={\theta +z\over 2}$, $\theta^1_2={\theta-z\over 2}$,
$\theta^1_4={\theta-w\over 2}$, $\theta^1_5={\theta+w\over 2}$.
From \thetag{8.3.39}, we then get
$$
\split
&\ \ \ (e_1,\,e_5)\le 4\,\sh~{\Big(\theta\hskip-2pt+\hskip-2pt
{z\hskip-2pt-\hskip-2ptw\over 4}-}\\
&{-g(\alpha_1,\hskip-2pt
{\theta \hskip-2pt+\hskip-2pt z\over 4},\hskip-2pt
{\theta\hskip-2pt-\hskip-2pt z\over 4})\hskip-2pt-\hskip-2pt
g(\alpha_2,\hskip-2pt{\theta\hskip-2pt-\hskip-2pt
z\over 4},\hskip-2pt
{\theta\hskip-2pt+\hskip-2pt z\over 4})\hskip-2pt-\hskip-2pt  
g(\alpha_3,\hskip-2pt {\theta\hskip-2pt+\hskip-2pt z\over 4},
\hskip-2pt
{\theta\hskip-2pt-\hskip-2pt w\over 4})\hskip-2pt-\hskip-2pt
g(\alpha_4,\hskip-2pt{\theta\hskip-2pt -\hskip-2pt w\over 4},\hskip-2pt
{\theta\hskip-2pt +\hskip-2pt w\over 4})
\Big)}\times \\
&\sh~{\hskip-2pt\Big(\theta\hskip-2pt-\hskip-2pt
g(\alpha_1,\hskip-2pt {\theta \hskip-3pt+\hskip-3pt z\over 4},\hskip-2pt
{\theta\hskip-3pt-\hskip-3pt z\over 4})\hskip-2pt-\hskip-2pt
g(\alpha_2,\hskip-1pt{\theta\hskip-3pt-\hskip-3pt
z\over 4},\hskip-1pt {\theta\hskip-3pt+\hskip-3pt z\over 4})\hskip-2pt-
\hskip-2pt  g(\alpha_3,\hskip-2pt
{\theta\hskip-3pt+\hskip-3pt z\over 4},\hskip-2pt
{\theta\hskip-3pt-\hskip-3pt w\over 4})\hskip-2pt-\hskip-2pt
g(\alpha_4,\hskip-2pt
{\theta\hskip-3pt -\hskip-3pt w\over 4},\hskip-2pt
{\theta\hskip-2pt +\hskip-2pt w\over 4})
\Big)}\Bigg/\\
&\Big(\sh~{\theta+z\over 4}\sh~{\theta+w\over 4}\Big)-2\ .
\endsplit
\tag{8.3.40}
$$
It follows, that instead of \thetag{8.2.27} we have
$$
(e_1,\,e_5)\le 4\,\max_{0\le w\le z\le \theta}
{f_{h3}(\alpha_1,\alpha_2,\alpha_3,\alpha_4,\theta,z,w)}-2.
\tag{8.3.41}
$$
Following the proof of Theorem 8.2.2, using \thetag{8.3.41},
we get \thetag{8.3.16}.

Now we consider the proof of \thetag{8.2.33} and \thetag{8.2.34}.
Using \thetag{8.1.21}, we get
$$
(e_1,\,f_1)=
4{\ch~{\theta^1_1+\theta_{12}\over 2}
\ch~{\theta^2_1-\theta_{12}\over 2}
\over
\sh~{\theta^1_1\over 2}\sh~{\theta^2_1\over 2}}-2
\tag{8.3.42}
$$
and similar formular for $(e_2,\,f_1)$.
Here $\theta^1_1\ge {\eta\over 2}$, $\theta^1_2\ge {\eta\over 2}$
and $\theta^2_1\ge {\eta\over 2}$. There exists a line $(AB)$ which
is orthogonal to the axis, contains the vertex $A_1$ and intersects
the line $(B_0B_1)$. Here we denote
by $A$ and $B$ terminals of the line at infinity and we suppose that
$A$ belongs to the half-plane of the axis, which contains vertices
$A_i$. Then $B$ belongs to the half-plane of the axis containing
vertices $B_j$. We denote by $u_1=\theta(A_{21},A)$,
$u_2=\theta(A,A_{12})$, then
$\theta^1_{21}=u_1+u_2=2g(\alpha_1,{\theta^1_1\over 2},
{\theta^1_2\over 2})$ where $u_1\ge 0$ and $u_2\ge 0$. We denote by
$v_1=\theta(B_{11},B)$, $v_2=\theta(B,B_{12})$. Then
$\theta^2_1=v_1+v_2$ where $v_1\ge 0$ and $v_2\ge 0$.
We have $\theta_{12}=v_2-u_2$ in \thetag{8.3.42}.

We want to prove \thetag{8.3.20} (where $\delta_4=f_1$).
Like for the proof of \thetag{8.3.33},
we can suppose one of two cases (1) or (2) below:

Case (1). We can assume that
$\theta^1_1=\theta^1_2=\theta^2_1={\eta\over 2}$ (instead of
inequalities
$\theta^1_1,\,\theta^1_2,\,\theta^2_1\ge {\eta\over 2}$). Then
$u_1=u_2=\theta^1_{21}/2$, and,
by \thetag{8.3.42}, we have
$$
{(e_1,\,f_1)+2\over 4}\le
{\ch~{\eta\over 2}\over \ch~{\eta\over 2}-1}+
{\ch~{(v_2-u_2)}\over \ch~{\eta\over 2}-1}.
$$
Similarly, we have
$$
{(e_2,\,f_1)+2\over 4}\le
{\ch~{\eta\over 2}\over \ch~{\eta\over 2}-1}+
{\ch~{(v_1-u_1)}\over \ch~{\eta\over 2}-1}.
$$
We have
$|v_1-u_1|+|v_2-u_2|=|v_1-u_1|+|{\eta\over 2}-v_1-u_1|
\le {\eta\over 2}$ because
$u_1=u_2\le {\eta\over 2}$, $v_1+v_2={\eta\over 2}$ and
$v_1,\,v_2\ge 0$.
It follows
$$
\text{Ach}~\left({(\ch~{\eta\over 2}\hskip-2pt-\hskip-2pt 1)
((e_1,\,f_1)\hskip-1pt+\hskip-1pt2)\over 4}\hskip-1pt-\hskip-1pt
\ch~{\eta\over 2}\right)\hskip-1pt+\hskip-1pt
\text{Ach}~\left({(\ch~{\eta\over 2}\hskip-2pt-\hskip-2pt1)
((e_2,\,f_1)\hskip-1pt+\hskip-1pt 2)\over 4}\hskip-1pt-\hskip-1pt
\ch~{\eta\over 2}\right)\hskip-1pt\le \hskip-1pt{\eta\over 2}.
$$
It follows \thetag{8.3.20}.

Case (2). We can assume that
$\theta^1_1=\theta^2_1={\eta\over 2}$, $\theta^1_2\ge {\eta\over 2}$,
$v_2=0$ (or similar case
$\theta^1_2=\theta^2_1={\eta\over 2}$,
$\theta^1_1\ge {\eta\over 2}$ and $v_1=0$). Then
$v_1={\eta\over 2}$, and by \thetag{8.3.42}, one has
$$
{(e_1,\,f_1)+2\over 4}\le
{\ch~{\eta\over 2}\over \ch~{\eta\over 2}-1}+
{\ch~{u_2}\over \ch~{\eta\over 2}-1}
$$
and
$$
{(e_2,\,f_1)+2\over 4}\le
{\ch~{\theta^1_2+\eta/2\over 2}\over
2\ch~{\theta^1_2\over 2}\ch~{\eta\over 4}}+
{\ch~({\theta^1_2\over 2}+{\eta\over 4}-u_1)\over
2\ch~{\theta^1_2\over 2}\ch~{\eta\over 4}}\le
$$
$$
{\ch~{\eta\over 2}\over \ch~{\eta\over 2}-1}+
{\ch~({\eta\over 4}+|{\eta\over 4}-u_1|)\over \ch~{\eta\over 2}-1}.
$$
If $u_1\ge {\eta\over 4}$, we have
$u_2+{\eta\over 4}+|{\eta\over 4}-u_1|
= u_2+u_1=\theta^1_{21}\le {\eta\over 2}$. If $u_1\le {\eta\over 4}$,
we have
$u_2+{\eta\over 4}+|{\eta\over 4}-u_1|={\eta\over 2}+(u_2-u_1)\le
{\eta\over 2}+\theta^1_{21}={\eta\over 2}+
2\max_{x\ge {\eta\over 2}}g(\alpha_1,{\eta\over 4},{x\over 2})=
{\eta\over 2}+2g(\alpha_1,{\eta\over 4},+\infty)$. It follows
\thetag{8.3.20}. Similarly we get \thetag{8.3.21},
\thetag{8.3.25} and \thetag{8.3.26}. All other inequalities of
cases (BII$_1$) and (BII$_2$) of Theorem 8.3.1 are the same as
in Theorem 8.2.2.

Now we consider the proof of \thetag{8.2.40} and \thetag{8.2.41}.
Using \thetag{8.1.21}, we have
$$
(e_1,\,f_1)=
4{\ch~{\theta^1_1+\theta_{12}\over 2}
\ch~{\theta^2_1-\theta_{12}\over 2}
\over
\sh~{\theta^1_1\over 2}\sh~{\theta^2_1\over 2}}-2
\tag{8.3.43}
$$
where $\theta^1_1\ge {\eta\over 2}$, $\theta^1_2\le {\eta\over 2}$,
$\theta^2_1\ge {\eta\over 2}$ and there exists a line $(AB)$
which is orthogonal to the axis and intersects both lines
$(A_1A_2)$ and $(B_0B_1)$. (One has similar formular for
$(e_3,\,f_1)$, we shall consider it below). Here we denote by
$A$ and $B$ terminals of the line $(AB)$ at infinity, and
$A$ belongs to the half-plane of the axis containing vertices
$A_i$. Then $B$ belongs to the half-plane of the axis which contains
vertices $B_j$. We denote by $u_1=\theta(A_{21},A)$ and
$u_2=\theta(A,A_{22})$. Then
$\theta^1_2=u=u_1+u_2$ where
$u_1,\,u_2 \ge 0$. Similarly we introduce
$v_1=\theta(B_{11},B)$,
$v_2=\theta(B,B_{12})$. Then $\theta^2_1=v_1+v_2$ where
$v_1,\,v_2\ge 0$.

We want to prove \thetag{8.3.28}, \thetag{8.3.29} and \thetag{8.3.30}
(where $\delta_4=f_1$). Like for the proof of \thetag{8.3.33},
we can suppose that $\theta^1_1=\theta^1_3=\theta^2_1={\eta\over 2}$
(instead of inequalities
$\theta^1_1,\,\theta^1_3,\,\theta^2_1\ \ge\ {\eta\over 2}$).
Using \thetag{8.3.43}, we have
$$
{(e_1,\,f_1)+2\over 4}=
{\ch~{\theta^1_1+\theta_{12}\over 2}
\ch~{\theta^2_1-\theta_{12}\over 2}
\over
\sh~{\theta^1_1\over 2}\sh~{\theta^2_1\over 2}}\le
$$
$$
{\ch~{\eta/2+\theta_{12}\over 2}
\ch~{\eta/2-\theta_{12}\over 2}
\over
\sh~{\eta\over 4}\sh~{\eta\over 4}}=
{\ch~{\eta\over 2}\over \ch~{\eta\over 2}-1}+
{\ch~\theta_{12}\over \ch~{\eta\over 2}-1}
\tag{8.3.44}
$$
where $\theta_{12}=u_1-2g(\alpha_1,{\eta\over 4},{u\over 2})+v_2$.
It follows
$$
\split
&\text{Ach}~\left({(\ch~{\eta\over 2}-1)((e_1,\,f_1)+2)\over 4}-
\ch~{\eta\over 2}\right)\le
|v_2+u_1-2\,g(\alpha_1,{\eta\over 4},{u\over 2})|\le \\
&{\eta\over 2}+\max_{0\le u\le \eta/2}
{\left(u-2\,g(\alpha_1,{\eta\over 4},{u\over 2})\right)}=
\eta-2\,g(\alpha_1,{\eta\over 4},{\eta\over 4}).
\endsplit
$$
To get the last inequality, we use that
$0\le v_2\le {\eta\over 2}$, $0\le u_1 \le u\le {\eta\over 2}$ and
$0\le \theta_{21}^1=
2\,g(\alpha_1,{\eta\over 4},{u\over 2})\le u\le {\eta\over 2}$.
Similarly we have
$$
\split
&\text{Ach}~\left({(\ch~{\eta\over 2}-1)((e_3,\,f_1)+2)\over 4}-
\ch~{\eta\over 2}\right)\le
|u_2+v_1-2g(\alpha_2,{\eta\over 4},{u\over 2})|\le \\
&\le {\eta\over 2}+\max_{0\le u\le \eta/2}
{\left(u-2\,g(\alpha_2,{\eta\over 4},{u\over 2})\right)}=
\eta-2\,g(\alpha_2,{\eta\over 4},{\eta\over 4}).
\endsplit
$$
This proves \thetag{8.3.28} and \thetag{8.3.29}.
We have
$$
|u_1+v_2-2\,g(\alpha_1,{\eta\over 4},{u\over 2})|+
|u_2+v_1-2\,g(\alpha_2,{\eta\over 4},{u\over 2})|\le
$$
$$
\max_{0\le u\le {\eta\over 2}}
{\ \max_{0\le x\le {\eta\over 2}+u}
{\left(|x-2g(\alpha_1,{\eta\over 4},{u\over 2})|+
|{\eta\over 2}+u-x-2\,g(\alpha_2,{\eta\over 4},{u\over 2})|\right)}}=
$$
$$
\max_{0\le u\le {\eta\over 2}}
{\left({\eta\over 2}+u+
2\,|g(\alpha_1,{\eta\over 4},{u\over 2})-
g(\alpha_2,{\eta\over 4},{u\over 2})|\right)}.
$$
Here we use obviouse inequalities:
$\theta^1_{21}=2\,g(\alpha_1,{\eta\over 4},{u\over 2})\le u\le \eta/2$,
$\theta^1_{32}=2\,g(\alpha_2,{\eta\over 4},$ ${u\over 2}) \le u\le \eta/2$.
Thus, we get
$$
\text{Ach}~\left({(\ch~{\eta\over 2}\hskip-1pt-\hskip-1pt 1)
((e_1,\,f_1)\hskip-1pt+\hskip-1pt2)\over 4}\hskip-1pt-\hskip-1pt
\ch~{\eta\over 2}\right)+
\text{Ach}~\left({(\ch~{\eta\over 2}\hskip-1pt-\hskip-1pt1)
((e_3,\,f_1)\hskip-1pt+\hskip-1pt 2)\over 4}\hskip-1pt-\hskip-1pt
\ch~{\eta\over 2}\right)\le
$$
$$
{\eta\over 2}+\max_{0\le u\le {\eta\over 2}}
{\left(u+2\,|g(\alpha_1,{\eta\over 4},{u\over 2})-
g(\alpha_2,{\eta\over 4},{u\over 2})|\right)}.
$$
This proves \thetag{8.3.30}.

All other statements of Theorem 8.3.1 are the same as similar statements
of Theorem 8.2.2.

It finishes the proof.
\enddemo

From Theorem 8.3.1 we get an important for application statement 
which is analogous to Theorems 4.1.8 and 4.2.3. 

\proclaim{Theorem 8.3.2 (about narrow parts of types (I),
(II) and (III) of restricted hyperbolic convex polygons)} 
We fix a constant $0<\eta\le 3$.

Any infinite restricted hyperbolic convex polygon $\M$ has a narrow
part of one of types
(AI), (AII), (AIII), (BI), (BII$_1$),(BII$_2$), (BIII) 
defined below:

(AI) There exist four consecutive vertices
$A_0$, $A_1$, $A_2$, $A_3$ of one half of $\M$
such that for orthogonal
vectors $\delta_1$, $\delta_2$ and $\delta_3$
to lines $(A_0A_1)$, $(A_1A_2)$ and $(A_2A_3)$
respectively directed outwards of $\M$ and with
$\delta_1^2=\delta_2^2=\delta_3^2=-2$, one has
$(\delta_i,\,\delta_{i+1})=2\cos{\alpha_i}$, where  
$\alpha_i=\angle~A_{i-1}A_iA_{i+1}$, $i=1,\,2$, 
$\alpha_2\not={\pi\over 2}$, and
$$
2<(\delta_1,\,\delta_3) \le 4\,f_{h1}(\alpha_1,\alpha_2,\eta)-2.
\tag{8.3.45}
$$  
The set $\{\delta_1,\,\delta_2,\,\delta_3\}$
generates the 3-dimensional hyperbolic vector space and
has a connected Gram graph.
Additionally, open half-planes $\Ha_{\delta_i}^+$,
$i=1,\,2,\,3$, contain the axis $l$.

(AII) There exist five consecutive vertices
$A_0$, $A_1$, $A_2$, $A_3$ and $A_4$ of one half of $\M$
such that for orthogonal
vectors $\delta_1$, $\delta_2$, $\delta_3$ and $\delta_4$
to lines $(A_0A_1)$, $(A_1A_2)$, $(A_2A_3)$ and $(A_3A_4)$
respectively directed outwards of $\M$ and with
$\delta_1^2=\delta_2^2=\delta_3^2=\delta_4^2=-2$, one has
$(\delta_i,\,\delta_{i+1})=2\cos{\alpha_i}$ 
where $\alpha_i=\angle~A_{i-1}A_iA_{i+1}$, $i=1,\,2,\,3$, 
$\alpha_1=\alpha_2={\pi\over 2}$, and
$$
2<(\delta_1,\,\delta_3)
\le 4\,f_{h1}({\pi\over 2},{\pi\over 2},\eta)-2,
\tag{8.3.46}
$$
$$
2<(\delta_1,\,\delta_4)\le 
4\,\max{\big(f_{h2}({\pi\over 2}, {\pi\over 2},\alpha_3,\eta,\,0),\,
f_{h2}({\pi\over 2},{\pi\over 2},\alpha_3,\eta, \eta)\big)}-2.
\tag{8.3.47}
$$
The set $\{\delta_1,\,\delta_3,\,\delta_4\}$
generates the 3-dimensional hyperbolic vector space and
has a connected Gram graph. Additionally, open half-planes 
$\Ha_{\delta_i}^+$, $i=1,\,2,\,3,\,4$, contain the axis $l$.

(AIII) There exist six consecutive vertices
$A_0$, $A_1$, $A_2$, $A_3$, $A_4$ and $A_5$ of one half of $\M$
such that for orthogonal vectors $\delta_1$, $\delta_2$, $\delta_3$,
$\delta_4$ and $\delta_5$
to lines $(A_0A_1)$, $(A_1A_2)$, $(A_2A_3)$, $(A_3A_4)$ and
$(A_4A_5)$ respectively directed
outwards of $\M$ and with
$\delta_1^2=\delta_2^2=\delta_3^2=\delta_4^2=\delta_5^2=-2$
one has $(\delta_i,\,\delta_{i+1})=2\cos{\alpha_i}$, 
where $\alpha_i=\angle~A_{i-1}A_iA_{i+1}$, $i=1,\,2,\,3,\,4$, 
$\alpha_1=\alpha_2=\alpha_3=\alpha_4={\pi\over 2}$, and
$$
2<(\delta_1,\,\delta_3)
\le 4\,f_{h1}({\pi\over 2},\,{\pi\over 2},\,\eta)-2,
\tag{8.3.48}
$$
$$
2<(\delta_3,\,\delta_5)
\le 4\,f_{h1}({\pi\over 2},\,{\pi\over 2},\,\eta)-2,
\tag{8.3.49}
$$
and
$$
\split
&2<(\delta_1,\,\delta_5)\le\\
&4\,\max
{\big(f_{h3}({\pi\over 2},\hskip-1pt  {\pi\over 2},\hskip-1pt
{\pi\over 2},\hskip-1pt {\pi\over 2},
\hskip-1pt\eta,\hskip-1pt0,\hskip-1pt0),\,
f_{h3}({\pi\over 2},\hskip-1pt {\pi\over 2},\hskip-1pt
{\pi\over 2},
\hskip-1pt {\pi\over 2},\hskip-1pt\eta,\hskip-1pt\eta,\hskip-1pt0),\
f_{h3}({\pi\over 2},
\hskip-1pt {\pi\over 2},\hskip-1pt {\pi\over 2},\hskip-1pt
{\pi\over 2},
\hskip-1pt\eta,\hskip-1pt \eta,\hskip-1pt \eta)\big)}-2. 
\endsplit
\tag{8.3.50}
$$
The set $\{\delta_1,\,\delta_3,\,\delta_5\}$ generates the 3-dimensional
hyperbolic vector space and has a connected Gram graph.
Additionally, open half-planes $\Ha_{\delta_i}^+$,
$i=1,\,2,\,3,\,4$, $5$, contain the axis $l$.

(BI) There exist four consecutive vertices
$A_0$, $A_1$, $A_2$ of one half of $\M$
and two consecutive vertices $B_0$, $B_1$ of another half of $\M$
such that for orthogonal
vectors $\delta_1$, $\delta_2$ and  $\delta_3$
to lines $(A_0A_1)$, $(A_1A_2)$ and $(B_0B_1)$ 
respectively directed outwards of $\M$ and with
$\delta_1^2=\delta_2^2=\delta_3^2=-2$, one has
$(\delta_1,\,\delta_2)=2\cos{\alpha}$,
where $\alpha=\angle~A_0A_1A_2\not= {\pi\over 2}$, and
$$
2<(\delta_1,\,\delta_3)\le {8\over \ch~{\eta\over 2}-1}+6,
\tag{8.3.51}
$$
$$
2<(\delta_2,\,\delta_3)\le {8\over \ch~{\eta\over 2}-1}+6,
\tag{8.3.52}
$$
$$
\split
&\text{Ach}~\left({(\ch~{\eta\over 2}\hskip-2pt-\hskip-2pt 1)
((\delta_1,\,\delta_3)
\hskip-1pt+\hskip-1pt2)\over 4}\hskip-1pt-\hskip-1pt
\ch~{\eta\over 2}\right)\hskip-1pt+\hskip-1pt
\text{Ach}~\left({(\ch~{\eta\over 2}\hskip-2pt-\hskip-2pt1)
((\delta_2,\,\delta_3)
\hskip-1pt+\hskip-1pt 2)\over 4}\hskip-1pt-\hskip-1pt
\ch~{\eta\over 2}\right)\ \le\\
&{\eta\over 2}+2g(\alpha,\,{\eta\over 4},+\infty).
\endsplit
\tag{8.3.53}
$$
The set $\{\delta_1,\,\delta_2,\,\delta_3\}$ generates the 3-dimensional
hyperbolic vector space and has a connected Gram graph.
Additionally, open half-planes $\Ha_{\delta_i}^+$,
$i=1,\,2$, and closed half-plane $\Ha_{\delta_3}^+$
contain the axis $l$; the line orthogonal to the axis $l$ and
containing $A_1$ intersects the line $(B_0B_1)$.

(BII$_1$) There exist four consecutive vertices
$A_0$, $A_1$, $A_2$, $A_3$ of one half of $\M$
and two consecutive vertices $B_0$, $B_1$ of another half of $\M$
such that for orthogonal
vectors $\delta_1$, $\delta_2$, $\delta_3$ and $\delta_4$
to lines $(A_0A_1)$, $(A_1A_2)$, $(A_2A_3)$ and $(B_0B_1)$
respectively directed outwards of $\M$ and with
$\delta_1^2=\delta_2^2=\delta_3^2=\delta_4^2=-2$, 
one has $(\delta_i,\,\delta_{i+1})=2\cos{\alpha_i}$, 
where $\alpha_i=\angle~A_{i-1}A_iA_{i+1}$, $i=1,\,2$, 
$\alpha_1=\alpha_2={\pi\over 2}$, and 
$$
2<(\delta_1,\,\delta_4)\le {8\over \ch~{\eta\over 2}-1}+6,
\tag{8.3.54}
$$
$$
2<(\delta_2,\,\delta_4)\le {8\over \ch~{\eta\over 2}-1}+6,
\tag{8.3.55}
$$
$$
2<(\delta_3,\,\delta_4)\le {8\over \ch~{\eta\over 2}-1}+6,
\tag{8.3.56}
$$
$$
\split
&\text{Ach}~\left({(\ch~{\eta\over 2}\hskip-2pt-\hskip-2pt 1)
((\delta_1,\,\delta_4)
\hskip-1pt+\hskip-1pt2)\over 4}\hskip-1pt-\hskip-1pt
\ch~{\eta\over 2}\right)\hskip-1pt+\hskip-1pt
\text{Ach}~\left({(\ch~{\eta\over 2}\hskip-2pt-\hskip-2pt1)
((\delta_2,\,\delta_4)
\hskip-1pt+\hskip-1pt 2)\over 4}\hskip-1pt-\hskip-1pt
\ch~{\eta\over 2}\right)\ \le\\
&{\eta\over 2}+2g({\pi\over 2},\,{\eta\over 4},+\infty),
\endsplit
\tag{8.3.57}
$$
$$
\split
&\text{Ach}~\left({(\ch~{\eta\over 2}\hskip-2pt-\hskip-2pt 1)
((\delta_2,\,\delta_4)
\hskip-1pt+\hskip-1pt2)\over 4}\hskip-1pt-\hskip-1pt
\ch~{\eta\over 2}\right)\hskip-1pt+\hskip-1pt
\text{Ach}~\left({(\ch~{\eta\over 2}\hskip-2pt-\hskip-2pt1)
((\delta_3,\,\delta_4)
\hskip-1pt+\hskip-1pt 2)\over 4}\hskip-1pt-\hskip-1pt
\ch~{\eta\over 2}\right)\ \le\\
&{\eta\over 2}+2g({\pi\over 2},\,{\eta\over 4},+\infty).
\endsplit
\tag{8.3.58}
$$
The set $\{\delta_1,\,\delta_2,\,\delta_4\}$ generates the 3-dimensional
hyperbolic vector space and has a connected Gram graph. The same
is valid for $\{\delta_2,\,\delta_3,\,\delta_4\}$.
Additionally, open half-planes $\Ha_{\delta_i}^+$,
$i=1,\,2,\,3$, and closed half-plane $\Ha_{\delta_4}^+$
contain the axis $l$; both lines orthogonal to the axis $l$ and
containing $A_1$ and $A_2$ intersect the line $(B_0B_1)$.

(BII$_2$) There exist three consecutive vertices
$A_0$, $A_1$, $A_2$ of one half of $\M$
and three consecutive vertices $B_0$, $B_1$, $B_2$ of
another half of $\M$ such that for orthogonal
vectors $\delta_1$, $\delta_2$, $\delta_3$ and $\delta_4$
to lines $(A_0A_1)$, $(A_1A_2)$, $(B_0B_1)$ and $(B_1B_2)$
respectively directed outwards of $\M$ and with
$\delta_1^2=\delta_2^2=\delta_3^2=\delta_4^2=-2$,
one has $(\delta_1,\,\delta_2)=2\cos{\alpha}$,  
$(\delta_3,\,\delta_4)=2\cos{\beta}$, 
where $\alpha=\angle~A_{0}A_1A_2$,  $\beta=\angle~B_{0}B_1B_2$,
$\alpha =\beta ={\pi\over 2}$, and 
$$
2<(\delta_1,\,\delta_3)\le {8\over \ch~{\eta\over 2}-1}+6,
\tag{8.3.59}
$$
$$
2<(\delta_2,\,\delta_3)\le {8\over \ch~{\eta\over 2}-1}+6,
\tag{8.3.60}
$$
$$
2<(\delta_2,\,\delta_4)\le {8\over \ch~{\eta\over 2}-1}+6,
\tag{8.3.61}
$$
$$
\split
&\text{Ach}~\left({(\ch~{\eta\over 2}\hskip-2pt-\hskip-2pt 1)
((\delta_1,\,\delta_3)
\hskip-1pt+\hskip-1pt2)\over 4}\hskip-1pt-\hskip-1pt
\ch~{\eta\over 2}\right)\hskip-1pt+\hskip-1pt
\text{Ach}~\left({(\ch~{\eta\over 2}\hskip-2pt-\hskip-2pt1)
((\delta_2,\,\delta_3)
\hskip-1pt+\hskip-1pt 2)\over 4}\hskip-1pt-\hskip-1pt
\ch~{\eta\over 2}\right)\ \le\\
&{\eta\over 2}+2g({\pi\over 2},\,{\eta\over 4},+\infty),
\endsplit
\tag{8.3.62}
$$
$$
\split
&\text{Ach}~\left({(\ch~{\eta\over 2}\hskip-2pt-\hskip-2pt 1)
((\delta_2,\,\delta_3)
\hskip-1pt+\hskip-1pt2)\over 4}\hskip-1pt-\hskip-1pt
\ch~{\eta\over 2}\right)\hskip-1pt+\hskip-1pt
\text{Ach}~\left({(\ch~{\eta\over 2}\hskip-2pt-\hskip-2pt1)
((\delta_2,\,\delta_4)
\hskip-1pt+\hskip-1pt 2)\over 4}\hskip-1pt-\hskip-1pt
\ch~{\eta\over 2}\right)\ \le\\
&{\eta\over 2}+2g({\pi\over 2},\,{\eta\over 4},+\infty).
\endsplit
\tag{8.3.63}
$$
The
set $\{\delta_1,\,\delta_2,\,\delta_3\}$ generates the 3-dimensional
hyperbolic vector space and has a connected Gram graph. The same
is valid for $\{\delta_2,\,\delta_3,\,\delta_4\}$.
Additionally, open half-planes $\Ha_{\delta_i}^+$,
$i=1,\,2,\,3,\,4$, contain the axis $l$; the line
orthogonal to the axis $l$ and containing $A_1$ intersects
$(B_0B_1)$, and the line orthogonal to $l$ and containing $B_1$
intersects $(A_1A_2)$.

(BIII) There exist four consecutive vertices
$A_0$, $A_1$, $A_2$, $A_3$ of one half of $\M$,
and two consecutive vertices $B_0$, $B_1$ of
another half of $\M$ such that for orthogonal
vectors $\delta_1$, $\delta_2$, $\delta_3$ and $\delta_4$
to lines $(A_0A_1)$, $(A_1A_2)$, $(A_2A_3)$ and $(B_0B_1)$
respectively directed outwards of $\M$ and with
$\delta_1^2=\delta_2^2=\delta_3^2=\delta_4^2=-2$, one has
$(\delta_i,\,\delta_{i+1})=2\cos{\alpha_i}$ where
$\alpha_i=\angle~A_{i-1}A_iA_{i+1}$, $i=1,\,2$, 
$\alpha_1=\alpha_2={\pi\over 2}$ and
$$
2<(\delta_1,\,\delta_3)
\le
4\,f_{h1}({\pi\over 2}, {\pi\over 2},\eta)-2,
\tag{8.3.64}
$$
$$
2<(\delta_1,\,\delta_4)\le
4{\ch~{\eta\over 2}+
\ch~(\eta - 2\,g({\pi\over 2},{\eta\over 4},{\eta\over 4}))
\over
\ch~{\eta\over 2}-1}-2,
\tag{8.3.65}
$$
$$
2<(\delta_3,\,\delta_4)\le
4{\ch~{\eta\over 2}+
\ch~(\eta - 2\,g({\pi\over 2},{\eta\over 4},{\eta\over 4}))
\over
\ch~{\eta\over 2}-1}-2, 
\tag{8.3.66}
$$
$$
\text{Ach}\,\left({(\ch~{\eta\over 2}\hskip-1pt-1\hskip-1pt)
((\delta_1,\,\delta_4)\hskip-1pt+\hskip-1pt 2)\over 4}
\hskip-1pt-\hskip-1pt\ch~{\eta\over 2}\right)\hskip-1pt+\hskip-1pt
\text{Ach}\,\left({(\ch~{\eta\over 2}\hskip-1pt-\hskip-1pt1)
((\delta_3,\,\delta_4)\hskip-1pt+\hskip-1pt2)\over 4}
\hskip-1pt-\hskip-1pt\ch~{\eta\over 2}\right)\le \eta.
\tag{8.3.67}
$$
The set $\{\delta_1,\,\delta_2,\,\delta_4\}$ generates the
3-dimensional hyperbolic vector space and has a connected Gram graph.
The same is valid for the set $\{\delta_2,\,\delta_3,\,\delta_4\}$.
Additionally, open half-planes $\Ha_{\delta_i}^+$,
$i=1,\,2,\,3,\,4$, contain the axis $l$; there exists a line
orthogonal to the axis $l$ and intersecting the lines
$(A_1A_2)$ and $(B_0B_1)$.
\endproclaim

\demo{Proof} The polygon $\M$ satisfies to 
one of cases (AII), (AIII), (BII$_1$), (BII$_2$), (BIII) of 
Theorem 8.3.1. 

Suppose that $\M$ satisfies the case (AII) of Theorem 8.3.1 
and one of angles $\alpha_1$, $\alpha_2$ is not equal to 
$\pi/2$. Suppose $\alpha_2\not=\pi/2$. 
Then $\delta_1,\,\delta_2,\,\delta_3$ give the case (AI) of 
Theorem 8.3.2. If $\alpha_1=\alpha_2=\pi/2$, we get the case 
(AII) of Theorem 8.3.2. 

Suppose that $\M$ satisfies the case (AIII) of Theorem 8.3.1. 
Like above, part of elements $\delta_1, \dots ,\delta_5$ gives 
the case (AI) of Theorem 8.3.2, 
if at least one of angles $\alpha_i\not=\pi/2$, $i=1,\dots ,4$. 
Otherwise, we get the case (A3) of Theorem 8.3.2. 

For cases  (BII$_1$), (BII$_2$), (BIII) of Theorem 8.3.1 
the proof is the same. 

For special values of angles of Theorem 8.3.2, we can calculate 
maxima in Theorem 8.3.1 with the help of computer.

This finishes the proof of Theorem 8.3.2.
\enddemo 

Applying Theorem 8.3.2 to $\eta=\beta=1,4134...$ where $\beta$ 
is the root of the equation
$$
\split 
&\max
{\big(f_{h3}({\pi\over 2},\hskip-1pt  {\pi\over 2},\hskip-1pt
{\pi\over 2},\hskip-1pt {\pi\over 2},
\hskip-1pt\eta,\hskip-1pt0,\hskip-1pt0),\,
f_{h3}({\pi\over 2},\hskip-1pt {\pi\over 2},\hskip-1pt
{\pi\over 2},
\hskip-1pt {\pi\over 2},\hskip-1pt\eta,\hskip-1pt\eta,\hskip-1pt0),\
f_{h3}({\pi\over 2},
\hskip-1pt {\pi\over 2},\hskip-1pt {\pi\over 2},\hskip-1pt
{\pi\over 2},
\hskip-1pt\eta,\hskip-1pt \eta,\hskip-1pt \eta)\big)}=\\
&={\ch~{\eta\over 2}+
\ch~(\eta - 2\,g({\pi\over 2},{\eta\over 4},{\eta\over 4}))
\over
\ch~{\eta\over 2}-1},
\endsplit 
\tag{8.3.68}
$$ 
we improve Theorem 8.2.2:

\proclaim{Theorem 8.3.3 (about narrow part of restricted hyperbolic 
convex polygons in a hyperbolic plane)} 
For any restricted hyperbolic convex polygon $\M$ in a hyperbolic plane 
there are $\delta_1,\,\delta_2,\,\delta_3\in P(\M)$ with 
$\delta_1^2=\delta_2^2=\delta_3^2=-2$ such that 
$\delta_1,\,\delta_2,\,\delta_3$ generate the hyperbolic vector space 
(defining the hyperbolic space), they have connected Gram graph 
and 
$$
(\delta_i,\,\delta_j)\le 
4\,{\ch~{\beta\over 2}+
\ch~(\beta - 2\,g({\pi\over 2},{\beta\over 4},{\beta\over 4}))
\over
\ch~{\beta\over 2}-1} < 45,4629...\ 
\tag{8.3.69}
$$ 
for all $1\le i,j,\le 3$. Here $\beta=1,4134...$ is the root of 
the equation \thetag{8.3.68} 
\endproclaim

By computer experiments (see Sect. 8.5), we found that for classification of 
hyperbolically reflective hyperbolic lattices of the rank three 
it is good to apply 
Theorem 8.3.1 to $\eta=
\eta_0=2\,\text{arc}\,\ch~{3\over 2}$.  
Thus, finally we get the result which we  
later will apply to classification of hyperbolically reflective 
hyperbolic lattices of the rank three.

\proclaim{Theorem 8.3.4 (about narrow parts of types (I),
(II) and (III) of restricted hyperbolic convex polygons)}

We set
$\eta_0=2\text{arc}\,\ch~{3\over 2}=1.9248473002384137899910...$;
equivalently, $\ch~{\eta_0\over 2}={3\over 2}$).

Any infinite restricted hyperbolic convex polygon $\M$ has a narrow
part of one of types
(AI), (AII), (AIII), (BI), (BII$_1$),(BII$_2$), (BIII) 
defined below:

(AI) There exist four consecutive vertices
$A_0$, $A_1$, $A_2$, $A_3$ of one half of $\M$
such that for orthogonal
vectors $\delta_1$, $\delta_2$ and $\delta_3$
to lines $(A_0A_1)$, $(A_1A_2)$ and $(A_2A_3)$
respectively directed outwards of $\M$ and with
$\delta_1^2=\delta_2^2=\delta_3^2=-2$, one has
$(\delta_i,\,\delta_{i+1})=2\cos{\alpha_i}$, where   
$\alpha_i=\angle~A_{i-1}A_iA_{i+1}$, $i=1,\,2$, 
$\alpha_2\not={\pi\over 2}$, and
$$
2<(\delta_1,\,\delta_3)\le
r_{h1}(\alpha_1,\,\alpha_2):=4\,f_{h1}(\alpha_1,\alpha_2,\eta_0)-2 
\tag{8.3.70}
$$
where in particular:

\vskip5pt 

{\settabs 2\columns
\+$r_{h1}({\pi\over 2},\,{\pi\over 2})=7$ \ \ \ &
$r_{h1}({\pi\over 2},\,{\pi\over 3})= 9.412375826...$\cr
\+$r_{h1}({\pi\over 2},\,{\pi\over 4})= 10.37390342...$ &
$r_{h1}({\pi\over 2},\,{\pi\over 6})= 11.10113930...$\cr
\+$r_{h1}({\pi\over 2},\,0)=            11.70820393...$ &
$r_{h1}({\pi\over 3},\,{\pi\over 3})= 12.25;$ \cr
\+$r_{h1}({\pi\over 3},\,{\pi\over 4})= 13.37793021...$ &
$r_{h1}({\pi\over 3},\,{\pi\over 6})= 14.23012096...$\cr
\+$r_{h1}({\pi\over 3},\,0)=            14.94097150...$ &
$r_{h1}({\pi\over 4},\,{\pi\over 4})= 14.57106781...$\cr
\+$r_{h1}({\pi\over 4},\,{\pi\over 6})= 15.47225159...$&
$r_{h1}({\pi\over 4},\,0)=            16.22381115...$\cr
\+$r_{h1}({\pi\over 6},\,{\pi\over 6})= 16.41025403...$ &
$r_{h1}({\pi\over 6},\,0)=            17.19241152...$\cr
\+$r_{h1}(0,\,0)=18.$ & &\cr}

\vskip5pt 

\noindent
The set $\{\delta_1,\,\delta_2,\,\delta_3\}$
generates the 3-dimensional hyperbolic vector space and
has a connected Gram graph.
Additionally, open half-planes $\Ha_{\delta_i}^+$,
$i=1,\,2,\,3$, contain the axis $l$.

(AII) There exist five consecutive vertices
$A_0$, $A_1$, $A_2$, $A_3$ and $A_4$ of one half of $\M$
such that for orthogonal
vectors $\delta_1$, $\delta_2$, $\delta_3$ and $\delta_4$
to lines $(A_0A_1)$, $(A_1A_2)$, $(A_2A_3)$ and $(A_3A_4)$
respectively directed outwards of $\M$ and with
$\delta_1^2=\delta_2^2=\delta_3^2=\delta_4^2=-2$, one has
$(\delta_i,\,\delta_{i+1})=2\cos{\alpha_i}$ 
where $\alpha_i=\angle~A_{i-1}A_iA_{i+1}$, $i=1,\,2,\,3$, 
$\alpha_1=\alpha_2={\pi\over 2}$, and
$$
2<(\delta_1,\,\delta_3)
\le r_{h1}({\pi\over 2},\,{\pi\over 2})=7,
\tag{8.3.71}
$$
$$
2<(\delta_1,\,\delta_4)\le r_{h2}(\alpha_3):=
4\,\max{\big(f_{h2}({\pi\over 2}, {\pi\over 2},\alpha_3,\eta_0,0),\,
f_{h2}({\pi\over 2},{\pi\over 2},\alpha_3,\eta_0,\eta_0)\big)}-2 
\tag{8.3.72}
$$
where in particular: 
$$
\split
&r_{h2}({\pi\over 2})=31.15549442...\ \ \
r_{h2}({\pi\over 3})= 38.68043607...\ \ \
r_{h2}({\pi\over 4})=41.73090517...\\
&r_{h2}({\pi\over 6})=44.05297726...\ \ \
r_{h2}(0)=46.
\endsplit
\tag{8.3.73}
$$ 
The set $\{\delta_1,\,\delta_3,\,\delta_4\}$
generates the 3-dimensional hyperbolic vector space and
has a connected Gram graph. Additionally, open half-planes 
$\Ha_{\delta_i}^+$, $i=1,\,2,\,3,\,4$, contain the axis $l$.

(AIII) There exist six consecutive vertices
$A_0$, $A_1$, $A_2$, $A_3$, $A_4$ and $A_5$ of one half of $\M$
such that for orthogonal vectors $\delta_1$, $\delta_2$, $\delta_3$,
$\delta_4$ and $\delta_5$
to lines $(A_0A_1)$, $(A_1A_2)$, $(A_2A_3)$, $(A_3A_4)$ and
$(A_4A_5)$ respectively directed
outwards of $\M$ and with
$\delta_1^2=\delta_2^2=\delta_3^2=\delta_4^2=\delta_5^2=-2$
one has $(\delta_i,\,\delta_{i+1})=2\cos{\alpha_i}$, 
where $\alpha_i=\angle~A_{i-1}A_iA_{i+1}$, $i=1,\,2,\,3,\,4$, 
$\alpha_1=\alpha_2=\alpha_3=\alpha_4={\pi\over 2}$, and
$$
2<(\delta_1,\delta_3)
\le r_{h1}({\pi\over 2},\,{\pi\over 2})=7,
\tag{8.3.74}
$$
$$
2<(\delta_3,\,\delta_5)
\le r_{h1}({\pi\over 2},\,{\pi\over 2})=7,
\tag{8.3.75}
$$
and
$$
\split
&2<(\delta_1,\,\delta_5)\le r_{h3}:=\\
&4\max
{\big(f_{h3}({\pi\over 2},\hskip-1pt  {\pi\over 2},\hskip-1pt
{\pi\over 2},\hskip-1pt {\pi\over 2},
\hskip-1pt\eta_0,\hskip-1pt0,\hskip-1pt0),\,
f_{h3}({\pi\over 2},\hskip-1pt {\pi\over 2},\hskip-1pt
{\pi\over 2},
\hskip-1pt {\pi\over 2},\hskip-1pt\eta_0,\hskip-1pt\eta_0,\hskip-1pt0),\
f_{h3}({\pi\over 2},
\hskip-1pt {\pi\over 2},\hskip-1pt {\pi\over 2},\hskip-1pt
{\pi\over 2},
\hskip-1pt\eta_0,\hskip-1pt \eta_0,\hskip-1pt \eta_0)\big)}\\
&-2=68.1815011826...
\endsplit
\tag{8.3.76}
$$
The set $\{\delta_1,\,\delta_3,\,\delta_5\}$ generates the 3-dimensional
hyperbolic vector space and has a connected Gram graph.
Additionally, open half-planes $\Ha_{\delta_i}^+$,
$i=1,\,2,\,3,\,4$, $5$, contain the axis $l$.

(BI) There exist four consecutive vertices
$A_0$, $A_1$, $A_2$ of one half of $\M$
and two consecutive vertices $B_0$, $B_1$ of another half of $\M$
such that for orthogonal
vectors $\delta_1$, $\delta_2$ and  $\delta_3$
to lines $(A_0A_1)$, $(A_1A_2)$ and $(B_0B_1)$ 
respectively directed outwards of $\M$ and with
$\delta_1^2=\delta_2^2=\delta_3^2=-2$, one has
$(\delta_1,\,\delta_2)=2\cos{\alpha}$,
where $\alpha=\angle~A_0A_1A_2\not= {\pi\over 2}$, and
$$
2<(\delta_1,\,\delta_3)\le 22,
\tag{8.3.77}
$$
$$
2<(\delta_2,\,\delta_3)\le 22,
\tag{8.3.78}
$$
$$
\split
&\text{Ach}~\left({(\delta_1,\,\delta_3)\over 8}-{5\over 4}\right)+
\text{Ach}~\left({(\delta_2,\,\delta_3)\over 8}-{5\over 4}\right)\le\\
&r_{h4}(\alpha):={\eta_0\over 2}+2g(\alpha ,\,{\eta_0\over 4},+\infty)
=\eta_0-2\text{arc}\,\sh~({1+\cos \alpha\over 2}\sh~{\eta_0\over 4})
\endsplit
\tag{8.3.79}
$$
where in particular

\vskip5pt 

\noindent
{\settabs 3\columns
\+$r_{h4}({\pi\over 2})=1.429914377...$ &
$r_{h4}({\pi\over 3})=  1.191398091...$ &
$r_{h4}({\pi\over 4})=1.095286061...$ \cr 
\+$r_{h4}({\pi\over 6})=1.022736500...$ &
$r_{h4}(0)={\eta_0\over 2}=0.962423650...\ .$ && \cr}

\vskip5pt 

\noindent
The set $\{\delta_1,\,\delta_2,\,\delta_3\}$
generates the 3-dimensional hyperbolic vector space and has
a connected Gram graph.
Additionally, open half-planes $\Ha_{\delta_i}^+$,
$i=1,\,2$, and closed half-plane $\Ha_{\delta_3}^+$
contain the axis $l$; the line orthogonal to the axis $l$ and
containing $A_1$ intersects the line $(B_0B_1)$. 

(BII$_1$) There exist four consecutive vertices
$A_0$, $A_1$, $A_2$, $A_3$ of one half of $\M$
and two consecutive vertices $B_0$, $B_1$ of another half of $\M$
such that for orthogonal
vectors $\delta_1$, $\delta_2$, $\delta_3$ and $\delta_4$
to lines $(A_0A_1)$, $(A_1A_2)$, $(A_2A_3)$ and $(B_0B_1)$
respectively directed outwards of $\M$ and with
$\delta_1^2=\delta_2^2=\delta_3^2=\delta_4^2=-2$, 
one has $(\delta_i,\,\delta_{i+1})=2\cos{\alpha_i}$, 
where $\alpha_i=\angle~A_{i-1}A_iA_{i+1}$, $i=1,\,2$, 
$\alpha_1=\alpha_2={\pi\over 2}$, and 
$$
2<(\delta_1,\,\delta_4)\le 22,
\tag{8.3.80}
$$
$$
2<(\delta_2,\,\delta_4)\le 22,
\tag{8.3.81}
$$
$$
2<(\delta_3,\,\delta_4)\le 22,
\tag{8.3.82}
$$
$$
\text{Ach}~\left({(\delta_1,\,\delta_4)\over 8}-{5\over 4}\right)+
\text{Ach}~\left({(\delta_2,\,\delta_4)\over 8}-{5\over 4}\right)\le 
rh_4({\pi\over 2})=1.429914377...\ ,
\tag{8.3.83}
$$
$$
\text{Ach}~\left({(\delta_2,\,\delta_4)\over 8}-{5\over 4}\right)+
\text{Ach}~\left({(\delta_3,\,\delta_4)\over 8}-{5\over 4}\right)\le 
r_{h4}({\pi\over 2})=1.429914377...\ .
\tag{8.3.84}
$$
The set $\{\delta_1,\,\delta_2,\,\delta_4\}$ generates the 3-dimensional
hyperbolic vector space and has a connected Gram graph. The same
is valid for $\{\delta_2,\,\delta_3,\,\delta_4\}$.
Additionally, open half-planes $\Ha_{\delta_i}^+$,
$i=1,\,2,\,3$, and closed half-plane $\Ha_{\delta_4}^+$
contain the axis $l$; both lines orthogonal to the axis $l$ and
containing $A_1$ and $A_2$ intersect the line $(B_0B_1)$.

(BII$_2$) There exist three consecutive vertices
$A_0$, $A_1$, $A_2$ of one half of $\M$
and three consecutive vertices $B_0$, $B_1$, $B_2$ of
another half of $\M$ such that for orthogonal
vectors $\delta_1$, $\delta_2$, $\delta_3$ and $\delta_4$
to lines $(A_0A_1)$, $(A_1A_2)$, $(B_0B_1)$ and $(B_1B_2)$
respectively directed outwards of $\M$ and with
$\delta_1^2=\delta_2^2=\delta_3^2=\delta_4^2=-2$,
one has $(\delta_1,\,\delta_2)=2\cos{\alpha}$,  
$(\delta_3,\,\delta_4)=2\cos{\beta}$, 
where $\alpha=\angle~A_{0}A_1A_2$,  $\beta=\angle~B_{0}B_1B_2$,
$\alpha =\beta ={\pi\over 2}$, and 
$$
2<(\delta_1,\,\delta_3)\le 22,
\tag{8.3.85}
$$
$$
2<(\delta_2,\,\delta_3)\le 22,
\tag{8.3.86}
$$
$$
2<(\delta_2,\,\delta_4)\le 22,
\tag{8.3.87}
$$
$$
\text{Ach}~\left({(\delta_1,\,\delta_3)\over 8}-{5\over 4}\right)+
\text{Ach}~\left({(\delta_2,\,\delta_3)\over 8}-{5\over 4}\right)\le 
r_{h4}({\pi\over 2})=1.429914377...\ ,
\tag{8.3.88}
$$
$$
\text{Ach}~\left({(\delta_3,\,\delta_2)\over 8}-{5\over 4}\right)+
\text{Ach}~\left({(\delta_4,\,\delta_2)\over 8}-{5\over 4}\right)\le 
r_{h4}({\pi\over 2})=1.429914377...\ .
\tag{8.3.89}
$$
The
set $\{\delta_1,\,\delta_2,\,\delta_3\}$ generates the 3-dimensional
hyperbolic vector space and has a connected Gram graph. The same
is valid for $\{\delta_2,\,\delta_3,\,\delta_4\}$.
Additionally, open half-planes $\Ha_{\delta_i}^+$,
$i=1,\,2,\,3,\,4$, contain the axis $l$; the line
orthogonal to the axis $l$ and containing $A_1$ intersects
$(B_0B_1)$, and the line orthogonal to $l$ and containing $B_1$
intersects $(A_1A_2)$.

(BIII) There exist four consecutive vertices
$A_0$, $A_1$, $A_2$, $A_3$ of one half of $\M$,
and two consecutive vertices $B_0$, $B_1$ of
another half of $\M$ such that for orthogonal
vectors $\delta_1$, $\delta_2$, $\delta_3$ and $\delta_4$
to lines $(A_0A_1)$, $(A_1A_2)$, $(A_2A_3)$ and $(B_0B_1)$
respectively directed outwards of $\M$ and with
$\delta_1^2=\delta_2^2=\delta_3^2=\delta_4^2=-2$, one has
$(\delta_i,\,\delta_{i+1})=2\cos{\alpha_i}$ where
$\alpha_i=\angle~A_{i-1}A_iA_{i+1}$, $i=1,\,2$, 
$\alpha_1=\alpha_2={\pi\over 2}$ and
$$
2<(\delta_1,\,\delta_3)
\le r_{h1}({\pi\over 2},\,{\pi\over 2})=7,
\tag{8.3.90}
$$
$$
2<(\delta_1,\,\delta_4)\le
4{\ch~{\eta_0\over 2}+
\ch~(\eta_0 - 2\,g({\pi\over 2},{\eta_0\over 4},{\eta_0\over 4}))
\over
\ch~{\eta_0\over 2}-1}-2=31.70820393...\ ,
\tag{8.3.91}
$$
$$
2<(\delta_3,\,\delta_4)\le
4{\ch~{\eta_0\over 2}+
\ch~(\eta_0 - 2\,g({\pi\over 2},{\eta_0\over 4},{\eta_0\over 4}))
\over
\ch~{\eta_0\over 2}-1}-2=31.70820393...\ ,
\tag{8.3.92}
$$
$$
\text{Ach}~\left({(\delta_1,\,\delta_4)\over 8}-{5\over 4}\right)+
\text{Ach}~\left({(\delta_3,\,\delta_4)\over 8}-{5\over 4}\right)\le
\eta_0=1.9248473002...\ .
\tag{8.3.93}
$$
The set $\{\delta_1,\,\delta_2,\,\delta_4\}$ generates the
3-dimensional hyperbolic vector space and has a connected Gram graph.
The same is valid for the set $\{\delta_2,\,\delta_3,\,\delta_4\}$.
Additionally, open half-planes $\Ha_{\delta_i}^+$,
$i=1,\,2,\,3,\,4$, contain the axis $l$; there exists a line
orthogonal to the axis $l$ and intersecting the lines
$(A_1A_2)$ and $(B_0B_1)$.
\endproclaim

\subhead
8.4. Some general applications 
\endsubhead

\subsubhead
8.4.1. Reflective hyperbolic lattices of hyperbolic type over  
purely real algebraic number fields 
\endsubsubhead
Let $\bk$ be a purely real algebraic number field of finite degree 
$[\bk:\bq]$ and $\bo$ its ring of integers. A {\it lattice} $S$ over $\bk$ 
is a projective module over $\bo$ of a finite rank equipped with a 
non-degenerate symmetric bilinear form with values in $\bo$. A lattice $S$ 
is called hyperbolic if there exists an embedding $\sigma^{(+)}:\bk\to \br$ 
(it is called the {\it geometric embedding)} such that the real symmetric 
bilinear form $S\otimes \br$ is hyperbolic (i.e. it has exactly one positive 
square), and for all other embeddings $\sigma\not=\sigma^{(+)}$ the 
real form $S\otimes \br$ is negative definite. Using the geometric embedding, 
one can define the hyperbolic space $\La(S)$, the automorphism group 
$O^+(S)$, the group $W(S)\subset O^+(S)$ 
generated by reflections of the lattice $S$.  
Here an automorphism of $O^{(+)}(S)$ is called  
{\it reflection} if it acts as reflection with respect to a hyperplane of 
$\La(S)$. For a hyperbolic lattice $S$, the groups $O^+(S)$ 
is discrete in $\La(S)$ and has a fundamental domain of finite volume. 
Let $\M\subset \La(S)$ 
be a fundamental polyhedron of $W(S)$ and 
$A(\M)=\{ \phi\in O^+(S)\ |\ \phi(\M)=\M\}$ its group of symmetries. 
Then $O^+(S)=W\rtimes A(\M)$ is the semi-direct product. The lattice $S$ is 
called {\it reflective of hyperbolic type (or hyperbolically reflective)} 
if there exists a subgroup of finite index $A^\prime \subset A(\M)$ 
which keeps invariant a proper subspace $\T\subset \La(S)$  
(here ``proper'' means that $0<\dim \T <\dim \La(S)$). Using 
standard discrete groups arguments, one can see that the polyhedron 
$\M$ is then restricted hyperbolic with respect to $\T$.  
From Theorem 8.2.3 (or 8.3.3) 
we get (like for usual reflective lattices, see \cite{N4}, \cite{N5, 
Appendix})

\proclaim{Theorem 8.4.1.1} For the fixed $N=[\bk:\bq]$ and 
$\rk S=\dim S\otimes \br=3$, the set 
of reflective hyperbolic lattices $S$ of hyperbolic type is finite up 
to multiplication of the form of $S$ by elements of $\bk$.    
\endproclaim

We remark that the generalization of Theorem 8.2.3 for 
arbitrary  dimension of hyperbolic space (using similar method) 
permits to prove 
Theorem 8.4.1.1 for arbitrary $\rk S\ge 3$. 
This gives the partial answer to the problem posed in 
\cite{N8, Sect. 3} about description of generalized arithmetic 
reflection groups (see also \cite{N9}, \cite{N13}).   
For the fixed degree $[\bk:\bq]$, finiteness 
of the set of fields $\bk$ such that 
there exists a reflective hyperbolic lattice $S$  
of hyperbolic type over $\bk$ and of rank $\rk S \ge 3$ 
was proved in \cite{N13}.   

\subsubhead
8.4.2. Algebraic surfaces with almost finite polyhedral Mori cone
\endsubsubhead
We refer to \cite{N14} about this subject. 

Let $X$ be a non-singular projective algebraic
surface over algebraically closed field. 
Let $\rho(X)=\rk \NS(X) \ge 3$ be the Picard number of $X$. We 
denote by $\Exc(X)$ the set of all exceptional (i. e. irreducible 
and with negative self-intersection) curves of $X$. 

We say that $X$ has {\it almost finite
polyhedral Mori cone} $\overline{NE}(X)$ if (1) and (2) below
hold: 

(1) There exist finite maxima: 
$$
\delta_E(X)=\max_{C\in \Exc(X)}{-C^2}\ \ \text{and}\ \
p_E(X)=\max_{C\in \Exc(X)}{p_a(C)}
$$
where $p_a(C)$ is the arithmetic genus of the curve $C$. 

(2) There exists a non-zero 
$r \in \NS(X)$ such that any extremal ray of $\overline{\NE}(X)$
is either generated by an exceptional curve $E\in \Exc(X)$ or
by $c\in \NS(X)\otimes \br$ such that $c^2=0$ and $c\cdot r=0$.
We remark that if the canonical class 
$K\not\equiv 0$, one can always put $r=K$. 

By elementary considerations, one can see that the dual 
cone $\text{NEF}(X)$ to the Mori cone 
$\overline{NE}(X)$ defines the polyhedron  
$\M=\text{NEF}(X)/\br_{++}\subset \La(X)$ of restricted 
hyperbolic type with respect to the hyperplane orthogonal to $r$, 
if $r^2<0$. If $r^2>0$ then $\M$ will be 
elliptic, and if $r^2=0$ then $\M$ will be parabolic. Like in 
\cite{N14} for elliptic case (when $r^2>0$), from  
Theorems 4.2.3 and 8.2.3 we get   

\proclaim{Theorem 8.4.2.1} For $\rho=3$, 
there are constants $N(\rho,\,\delta_E)$ and 
$N'(\rho,\,\delta_E,\,p_E)$ 
depending only on $(\rho,\,\delta_E)$ and $(\rho,\,\delta_E,\,p_E)$
respectively such that for any $X$ with almost finite polyhedral
Mori cone and $\rho(X)=\rho$, $\delta_E(X)=\delta_E$ and
$p_E(X)=p_E$, there exists an ample effective divisor $h$ such that
$h^2\le N(\rho,\,\delta_E)$, and if the ground field is $\bc$,
there exists a very ample divisor $h'$ such that
${h'}^2\le N'(\rho,\,\delta_E,\,p_E)$.
\endproclaim

The generalization of Theorem 8.2.3 to arbitrary dimension 
(using the same method) permits 
to prove Theorem 8.4.2.1 for arbitrary $\rho\ge 3$. 

\subhead
8.5. Description of narrow parts of fundamental polygons $\M$ of
reflection subgroups $W\subset W(S)$ of hyperbolic type where 
$\rk S=3$. Application to reflective lattices of 
hyperbolic type 
\endsubhead

Here we want to make similar calculations for hyperbolic type 
as we did for elliptic and parabolic 
types in Sect. 4.3. We use the same notations. 

Let $S$ be a primitive hyperbolic lattice of $\rk S=3$ and
$W\subset W(S)$ its reflection subgroup  of hyperbolic type
with a fundamental polygon $\M$ (see Sects. 1.3, 1.4).
Remind that it means that $W$ and $P(\M)_{\pr}$ have 
restricted arithmetic type and $P(\M)_{\pr}$ has 
a generalized lattice Weyl vector $\rho$ with 
$\rho^2<0$ (and does not have a generalized lattice Weyl vector 
$\rho$ with $\rho^2\ge 0$. It is easy to see that
the polygon $\M$ is then restricted hyperbolic. We remind 
that the lattice $S$ having a reflection subgroups $W\subset W(S)$ 
of hyperbolic type is called reflective of hyperbolic type.
Here we apply Theorem 8.3.4 to describe narrow parts of $\M$. Using
this description, we shall get a finite list of lattices such that
any hyperbolically reflective lattice $S$ belongs to
the list. Let $a(S^\ast/S)$ be the exponent of the discriminant 
group $S^\ast/S$: i. e. $a(S^\ast/S)$ is the 
least natural $a$ such that $aS^\ast/S=0$. We shall get an
estimate of  $a(S^\ast/S)$.

In Sect. 4.3 we gave description of narrow parts of fundamental 
polygons $\M$ for elliptic and parabolic types according to 
Theorems 4.1.8 and 4.2.3. Here we similarly want to 
describe narrow parts for hyperbolic type according to Theorem 
8.3.4. 
We use the same notations as in Sect. 4.3. 
We remind that we describe all possible matrices 
$B=((r_i,\,r_j))_{pr}$, $i=1,\dots ,k$ where $r_i\in S$ are orthogonal 
roots to $\M$ and $\delta_i=2r_i/\sqrt{-2r_i^2}$, $i=1,\dots ,k$ 
define the narrow part of $\M$ described in Theorem 8.3.4. For the 
matrices $B$ we calculate the invariant $a(B)$ which is equal to 
$a(K^\ast/K)$ for the lattice $K$ defined by $B$, the invariant 
$a_1(B)$ which is equal to to the product of all odd prime divisors of $a(B)$, 
and the invariant $a_2(B)$ which is equal to the greatest odd prime 
divisor of $a(B)$.   

\subsubhead
8.5.1. Matrices $B$ of narrow parts of the type (AI1)
\endsubsubhead
This is a particular case of the type (AI) in Theorem 8.3.4 when 
additionally $\alpha_1\not=\pi/2$. 
This case is similar to Sect. 4.3.1. We should only 
change the condition \thetag{4.3.1.2} by the condition
$$
\split
&1\le \alpha_{12}=\alpha_{21}\le 4,\ \
\alpha_{12}\le \alpha_{23}=\alpha_{32}\le 4,\\
&4< \alpha_{13}=\alpha_{31}\le 
r_{h1}\big(\arc \cos{{\sqrt{\alpha_{12}}\over 2}},\ 
\arc \cos{{\sqrt{\alpha_{23}}\over 2}}\big)^2.
\endsplit
\tag{8.5.1.1}
$$ 
In Appendix 2, we give the Program 13: funda11.gen which enumerates all 
matrices $B$ of this case. It finds the number $nAI1$ of these 
matrices and 
$$
aAI1=\max_{B}{a(B)},\ aAI1_1=\max_{B}{a_1(B)},\ aAI1_2=\max_B{a_2(B)}.
\tag{8.5.1.2}
$$
Calculation using this program gives
$$
nAI1=332,\ aAI1=5832,\ aAI1_1=759,\ aAI1_2=181.
\tag{8.5.1.3}
$$

\subsubhead
8.5.2. Matrices $B$ of narrow parts of the type (AI0)
\endsubsubhead
It is a particular case of Type (AI) of Theorem 8.3.4 when
additionally the angle $\alpha_1=\pi/2$. This case is similar to 
Sect. 4.3.2. We should only change the condition \thetag{4.3.2.2} 
by the condition 
$$
1\le \alpha_{23}=\alpha_{32}\le 4,\ \
4< \alpha_{13}=\alpha_{31}\le
r_{h1}\big({\pi\over 2},\ 
\arc \cos{{\sqrt{\alpha_{23}}\over 2}}\big)^2. 
\tag{8.5.2.1}
$$
 In Appendix 2, we give the Program 14: funda10.gen which enumerates all 
matrices $B$ of this case. It finds the number $nAI0$ of these 
matrices and 
$$
aAI0=\max_{B}{a(B)},\ aAI0_1=\max_{B}{a_1(B)},\ aAI0_2=\max_B{a_2(B)}.
\tag{8.5.2.2}
$$
Calculation using this program gives
$$
nAI0=4788,\ aAI0=154568,\ aAI0_1=21063,\ aAI0_2=139.
\tag{8.5.2.3}
$$

\subsubhead
8.5.3. Matrices $B$ of narrow parts of the type (AII1)
\endsubsubhead
 It is a particular case of Type (AII) of Theorem 8.3.4 when
additionally the angle $\alpha_3\not=\pi/2$. 
This case is similar to 
Sect. 4.3.3. We should only change the condition \thetag{4.3.3.2} 
by the condition 
$$
\split
&1\le \alpha_{34}=\alpha_{43}\le 4,\ \
4< \alpha_{13}=\alpha_{31}\le r_{h1}(\pi/2,\pi/2)^2=7^2,\\
&4\le \alpha_{14}=\alpha_{41}\le 
r_{h2}\big(\arc \cos{{\sqrt{\alpha_{34}}\over 2}}\big)^2,\\
&r_{h1}\big({\pi\over 2},\ 
\arc \cos{{\sqrt{\alpha_{34}}\over 2}}\big)^2
<\alpha_{24}=\alpha_{42}. 
\endsplit
\tag{8.5.3.1}
$$
In Appendix 2, we give the Program 15: funda21.gen which enumerates all 
matrices $B$ of this case. It finds the number $nAII1$ of these 
matrices and 
$$
aAII1=\max_{B}{a(B)},\ aAII1_1=\max_{B}{a_1(B)},\ aAII1_2=\max_B{a_2(B)}.
\tag{8.5.3.2}
$$
Calculation using this program gives
$$
nAII1=33488,\ aAII1=220116,\ aAII1_1=55029,\ aAII1_2=2113.
\tag{8.5.3.3}
$$

\subsubhead
8.5.4. Matrices $B$ of narrow parts of the type (AII0)
\endsubsubhead
 It is a particular case of Type (AII) of Theorem 8.3.4 when
additionally the angle $\alpha_3=\pi/2$. This case is similar to 
Sect. 4.3.4. We should only change the condition \thetag{4.3.4.2} 
by the condition 
$$
\split
&4< \alpha_{13}=\alpha_{31}\le 7^2,\ \
4< \alpha_{14}=\alpha_{41}\le r_{h2}({\pi\over 2})^2=(31.15549442...)^2\,,\\
&\alpha_{13}\le \alpha_{24}=\alpha_{42}.
\endsplit 
\tag{8.5.4.1}
$$
In Appendix 2, we give the Program 16: funda20.gen which enumerates all 
matrices $B$ of this case. It finds the number $nAII0$ of these 
matrices and 
$$
aAII0=\max_{B}{a(B)},\ aAII0_1=\max_{B}{a_1(B)},\ aAII0_2=\max_B{a_2(B)}.
\tag{8.5.4.2}
$$
Calculation using this program gives
$$
nAII0=2321854,\ aAII0=11430900,\ aAII0_1=2834349,\ aAII0_2=977.
\tag{8.5.4.3}
$$

\subsubhead
8.5.5. Matrices $B$ of narrow parts of the type (AIII)
\endsubsubhead
This case is similar to 
Sect. 4.3.5. We should only change the condition \thetag{4.3.5.2} 
by the condition 
$$
\split
&4< \alpha_{13}=\alpha_{31}\le 7^2,\ \ 
\alpha_{31}\le \alpha_{35}=\alpha_{53}\le 7^2,\\
&4< \alpha_{15}=\alpha_{51}\le (r_{h3})^2=(68.1815011826...)^2,\\
&r_{h2}({\pi\over 2})^2=(31.15549442...)^2\hskip-2pt <\hskip-2pt 
\alpha_{14}=\alpha_{41},\ 
r_{h2}({\pi\over 2})^2=(31.15549442...)^2\hskip-2pt 
<\hskip-2pt \alpha_{25}=\alpha_{52},
\endsplit
\tag{8.5.5.1}
$$
In Appendix 2, we give Program 17: funda3.gen which enumerates all 
matrices $B$ of this case. It finds the number $nAIII$ of these 
matrices and 
$$
aAIII=\max_{B}{a(B)},\ aAIII_1=\max_{B}{a_1(B)},\ aAIII_2=\max_B{a_2(B)}.
\tag{8.5.5.2}
$$
Calculation using this program gives
$$
nAIII=802291,\ aAIII=3150000,\ aAIII_1=219453,\ aAIII_2=4561.
\tag{8.5.5.3}
$$

\subsubhead
8.5.6. Matrices $B$ of narrow parts of the type (BI)
\endsubsubhead
This case is similar to 
Sect. 4.3.1. We should only change the condition \thetag{4.3.1.1} 
by the condition 
$$
\split
&1 \le \alpha_{12}=\alpha_{21}\le 4,\\
&4< \alpha_{13}=\alpha_{31}\le 22^2,\ \ 
\alpha_{13} \le \alpha_{23}=\alpha_{32}\le 22^2,\\
&\text{Ach}~\left({\sqrt{\alpha_{13}}\over 8}-{5\over 4}\right)+
\text{Ach}~\left({\sqrt{\alpha_{23}}\over 8}-{5\over 4}\right)\le
r_{h4}(\arc \cos{{\sqrt{\alpha_{12}}\over 2}}).
\endsplit
\tag{8.5.6.1}
$$
In Appendix 2, we give Program 18: fundb1.gen which enumerates all 
matrices $B$ of this case. It finds the number $nBI$ of these 
matrices and 
$$
aBI=\max_{B}{a(B)},\ aBI_1=\max_{B}{a_1(B)},\ aBI_2=\max_B{a_2(B)}.
\tag{8.5.6.2}
$$
Calculation using this program gives
$$
nBI=24329,\ aBI=1159692,\ aBI_1=289923,\ aBI_2=1429.
\tag{8.5.6.3}
$$

\subsubhead
8.5.7. Matrices $B$ of narrow parts of the type (BII$_1$)
\endsubsubhead
This case is similar to 
Sect. 4.3.3. We should only change the condition \thetag{4.3.3.2}  
by the condition 
$$
\split
&4\,<\,\alpha_{14}=\alpha_{41}\le 22^2,\ 
4\,<\,\alpha_{24}=\alpha_{42}\le 22^2,\  
\alpha_{14}\le \alpha_{34}=\alpha_{43}\le 22^2,\\
&\text{Ach}~\left({\sqrt{\alpha_{14}}\over 8}-{5\over 4}\right)+
\text{Ach}~\left({\sqrt{\alpha_{24}}\over 8}-{5\over 4}\right)\le
r_{h4}\left({\pi\over 2}\right)=1.429914377...\ ,\\
&\text{Ach}~\left({\sqrt{\alpha_{24}}\over 8}-{5\over 4}\right)+
\text{Ach}~\left({\sqrt{\alpha_{34}}\over 8}-{5\over 4}\right)\le
r_{h4}\left({\pi\over 2}\right)=1.429914377...\ ,\\
&4<\alpha_{13}=\alpha_{31}.
\endsplit
\tag{8.5.7.1}
$$
The natural numbers $\alpha_{ij}$ satisfy conditions \thetag{4.3.3.3} and 
\thetag{4.3.3.4}. We denote 
$$
u=\sqrt{\alpha_{13}\alpha_{14}\alpha_{34}}.
\tag{8.5.7.2} 
$$
We can rewrite the conditions \thetag{4.3.3.3} and \thetag{4.3.3.4} 
in the form 
$$
\split 
&(\alpha_{24}-4)u^2-4\alpha_{14}\alpha_{34}u+4\alpha_{14}\alpha_{34}
(4-\alpha_{14}-\alpha_{34}-\alpha_{24})=0\\
&u\in \bn\ \ \text{and}\ \  \alpha_{13}={u^2\over 
\alpha_{14}\alpha_{34}}\in \bn.
\endsplit 
\tag{8.5.7.3} 
$$
It follows,
$$
u={2\alpha_{14}\alpha_{34}+2\sqrt{D}\over \alpha_{24}-4}\ 
\tag{8.5.7.4}
$$
where 
$$
D=\alpha_{14}\alpha_{34}(\alpha_{14}\alpha_{34}+
(\alpha_{24}-4)(\alpha_{14}+\alpha_{24}+
\alpha_{34}-4))
\tag{8.5.7.5}
$$
is a perfect square. 
All other considerations are the same as in Sect. 4.3.3. 
In Appendix 2, we give the Program 19: fundb21.gen which enumerates all 
matrices $B$ of this case. It finds the number $nBII1$ of these 
matrices and 
$$
aBII1=\max_{B}{a(B)},\ aBII1_1=\max_{B}{a_1(B)},\ 
aBII1_2=\max_B{a_2(B)}.
\tag{8.5.7.6}
$$
Calculation using this program gives
$$
nBII1=1291199,\ aBII1=16144800,\ aBII1_1=1248555,\ aBII1_2=487.
\tag{8.5.7.7}
$$

\subsubhead
8.5.8. Matrices $B$ of narrow parts of the type (BII$_2$)
\endsubsubhead
For this case
$k=4$ and $\delta_1,\,\delta_2,\,\delta_3,\,\delta_4$
generate the $3$-dimensional hyperbolic vector space and any three of them
give a bases of the space. We get that the
matrix $\Cal A$ is a symmetric matrix
$$
{\Cal A}=\pmatrix
4          & 0         &\alpha_{13}&\alpha_{14}\\
0          & 4         &\alpha_{23}& \alpha_{24}\\
\alpha_{31}&\alpha_{32}& 4         & 0\\
\alpha_{41}&\alpha_{42}& 0         &  4
\endpmatrix
\tag{8.5.8.1}
$$
with integral non-negative coefficients where
$$
\split
&4\,<\,\alpha_{13}=\alpha_{31}\le 22^2,\ 
4\,<\,\alpha_{23}=\alpha_{32}\le 22^2,\  
\alpha_{13}\le \alpha_{24}=\alpha_{42}\le 22^2,\\
&\text{Ach}~\left({\sqrt{\alpha_{13}}\over 8}-{5\over 4}\right)+
\text{Ach}~\left({\sqrt{\alpha_{23}}\over 8}-{5\over 4}\right)\le
r_{h4}\left({\pi\over 2}\right)=1.429914377...\ ,\\
&\text{Ach}~\left({\sqrt{\alpha_{23}}\over 8}-{5\over 4}\right)+
\text{Ach}~\left({\sqrt{\alpha_{24}}\over 8}-{5\over 4}\right)\le
r_{h4}\left({\pi\over 2}\right)=1.429914377...\ ,\\
&4<\alpha_{14}=\alpha_{41}.
\endsplit
\tag{8.5.8.2}
$$
Here we use the conditions (BII$_2$) of Theorem 8.3.4 and add some
inequalities to avoid symmetric cases. We also have
$$
det(\Gamma )=
\alpha_{14}(\alpha_{23}-4)-
2\sqrt{\alpha_{14}\alpha_{13}\alpha_{23}\alpha_{24}}+
(\alpha_{13}-4)(\alpha_{24}-4)-4\alpha_{23}=0.
\tag{8.5.8.3}
$$
where
$$
\alpha_{14}\alpha_{13}\alpha_{23}\alpha_{24}\ \text{is a perfect square}.
\tag{8.5.8.4}
$$
We denote 
$$
u=\sqrt{\alpha_{14}\alpha_{13}\alpha_{23}\alpha_{24}}. 
\tag{8.5.8.5}
$$
The conditions \thetag{8.5.8.3} and \thetag{8.5.8.4} are equivalent to 
$$
\split
&u^2\,(\alpha_{23}-4)-u\,2\alpha_{13}\alpha_{23}\alpha_{24}+
\alpha_{13}\alpha_{23}\alpha_{24}\big((\alpha_{13}-4)(\alpha_{24}-4)-
4\alpha_{23}\big)=0,\\
&u\in \bn\ \ \text{and}\ \ \alpha_{14}={u^2\over 
\alpha_{13}\alpha_{23}\alpha_{24}}\in \bn.
\endsplit
\tag{8.5.8.6} 
$$
It follows
$$
u=
{\alpha_{13}\alpha_{23}\alpha_{24}\pm 
\sqrt{D}\over \alpha_{23}-4}
\tag{8.5.8.7}
$$
where 
$$
D=\alpha_{13}\alpha_{23}\alpha_{24}
\Big(\alpha_{13}\alpha_{23}\alpha_{24}-(\alpha_{23}-4)
\big((\alpha_{13}-4)(\alpha_{24}-4)-4\alpha_{23}\big)\Big)
\tag{8.5.8.8}
$$ 
is a perfect square. 
We can easily enumerate the finite set of all possible
matrices $\Cal A$ satisfying these conditions. For
each $\Cal A$ we find all symmetrizable generalized Cartan
matrices
$$
A=
\pmatrix
-2      &  0      &a_{13}   & a_{14}\\
0       & -2      &a_{23}   & a_{24}\\
a_{31}  &  a_{32} &   -2    & 0     \\
a_{41}  &  a_{42} &   0     & -2
\endpmatrix
\tag{8.5.8.9}
$$
using relations
$$
\split
&a_{13}a_{31}=\alpha_{13},\ a_{14}a_{41}=\alpha_{14},\
a_{23}a_{32}=\alpha_{23},\ a_{24}a_{42}=\alpha_{24},\\ 
&a_{13}a_{32}a_{24}a_{41}=
\sqrt{\alpha_{13}\alpha_{23}\alpha_{24}\alpha_{14}}.
\endsplit
\tag{8.5.8.10}
$$
For the diagonal matrix
$$
\Lambda=\text{diag}(a_{13}a_{14}a_{32},\ a_{31}a_{14}a_{23},\
a_{31}a_{14}a_{32},\ a_{13}a_{41}a_{32}), 
\tag{8.5.8.11}
$$
the matrix $A\Lambda$ is symmetric. Thus, finally we get that
$$
B=(A\Lambda)_{\pr}.
\tag{8.5.8.12}
$$
In Appendix 2, we give the Program 20: fundb22.gen which uses this algorithm
to enumerate all the matrices $B$. For each of them it calculates
the invariants $a(B)$, $a_1(B)$, $a_2(B)$ and finds the number
$nBII2$ of all the matrices $B$, and the numbers
$$
aBII2=\max_{B}{a(B)},\ aBII2_1=\max_{B}{a_1(B)},\ aBII2_2=\max_B{a_2(B)}.
\tag{8.5.8.13}
$$
Calculation using this program gives
$$
nBII2=929616,\ aBII2=3873744,\ aBII2_1=255255,\ aBII2_2=487.
\tag{8.5.8.14}
$$

\subsubhead
8.5.9. Matrices $B$ of narrow parts of the type (BIII)
\endsubsubhead
This case is similar to 
Sect. 4.3.3. We should only change the condition \thetag{4.3.3.2} 
by the condition 
$$
\split
&4< \alpha_{13}=\alpha_{31}\le r_{h1}({\pi\over 2},\,{\pi\over 2})^2=7^2,\\
&4< \alpha_{14}=\alpha_{41}\le (31.70820393...)^2\,,\ \ 
\alpha_{14}\le \alpha_{34}=\alpha_{43}\le  (31.70820393...)^2\\
&\text{Ach}~\left({\sqrt{\alpha_{14}}\over 8}-{5\over 4}\right)+
\text{Ach}~\left({\sqrt{\alpha_{34}}\over 8}-{5\over 4}\right)\le
\eta_0=1.9248473002...\ .
\endsplit 
\tag{8.5.9.1}
$$
In Appendix 2, we give the Program 21: fundb3.gen which uses this algorithm
to enumerate all the matrices $B$. For each of them it calculates
the invariants $a(B)$, $a_1(B)$, $a_2(B)$ and finds the number
$nBIII$ of all the matrices $B$, and the numbers
$$
aBIII=\max_{B}{a(B)},\ aBIII_1=\max_{B}{a_1(B)},\ aBIII_2=\max_B{a_2(B)}.
\tag{8.5.9.2}
$$
Calculation using this program gives
$$
nBIII=1878630,\ aBIII=7013160,\ aBIII_1=1627563,\ aBIII_2=2851.
\tag{8.5.9.3}
$$

\subsubhead
8.5.10. The global estimate of invariants of primitive reflective
hyperbolic lattices of the rank $3$ having hyperbolic type
\endsubsubhead Like in Sect. 4.3.6, we get some global estimates. 
For a lattice $L$ we denote by $a(L^\ast/L)$ the exponent of
the discriminant group $L^\ast/L$, we denote by $a_1(L^\ast/L)$
the product of all different odd prime divisors of $a(L^\ast/L)$, and
we denote by $a_2(L^\ast/L)$ the greatest prime divisor of
$a(L^\ast/L)$.

Using Proposition 4.3.1 and calculations \thetag{8.5.1.3},
\thetag{8.5.2.3}, \thetag{8.5.3.3}, \linebreak \thetag{8.5.4.3},
\thetag{8.5.5.3},\thetag{8.5.6.3}, \thetag{8.5.7.7}, 
\thetag{8.5.8.14}, \thetag{8.5.9.3}, we get

\proclaim{Theorem 8.5.10.1} For any primitive reflective hyperbolic
lattice $S$ of $\rk S=3$ having hyperbolic type we
have estimates:
$$
a(S^\ast/S)\le 8(16144800)^2,\ \
a_1(S^\ast/S)\le 2834349,\ \
a_2(S^\ast/S)\le 4561.
$$
\endproclaim

Since $\det (S)\le a(S^\ast/S)^2$ and number of lattices with the
fixed rank and determinant is finite (e. g. see \cite{C}),
Theorem 8.5.10.1 gives a finite list of lattices which contains
all the reflective lattices $S$.

The estimates of Theorem 8.5.10.1 are very preliminary,
and we shall significantly improve them below.

\head
9. Classification of hyperbolically reflective elementary 
hyperbolic lattices of the rank 3: Proofs 
\endhead

We are now ready to prove classification results of Sect. 7. 
The proofs are very similar to Sect. 5 but their details are more 
complicated and require much more calculations.

\subhead
9.1. Proof of Basic Theorem 7.1
\endsubhead 

Below (see the proof of Theorem 9.1.1) we will show that if a lattice $S$ 
of the rank 3 is hyperbolically reflective, then its invariant $h$ is 
equal to $0$ or $2$. Table 3 contains all main 
hyperbolic lattices with square-free determinant $d$ which 
have $h\le 1$ and $d\le 100000$. Using 
Theorem 3.2.1 and Program 3: refh3 (see Appendix), we can similarly 
find all lattices of this type with $h=2$ and $d\le 100000$. 
(One should only change $h\le 1$ by $h=2$ in two places of Program 3.) 
The result of these calculations is given in Table 6. It contains 
259 lattices. The maximal $d=29526$. Conjecturally, Tables 3 and 6 
contain all main hyperbolic lattices with 
square-free determinant $d$ (without the inequality $d\le 100000$) and 
$h\le 2$. To avoid this conjecture, we apply results of Sect. 8.   

We prove analog of Theorem 5.1.1. 

\proclaim{Theorem 9.1.1} Any hyperbolically reflective main 
hyperbolic lattice of the rank 3 and with square-free determinant $d$  
belongs either to the list of lattices with $h=0$ of Table 3  or to the 
list of lattices of Table 6 (all of them have $h=2$). We remind that 
the tables 3 and 6 contain 
all main hyperbolic lattices $S$ of the rank 3 with square-free 
$d=\det (S)\le 100000$ and $h=hnr(S)\le 2$.    
\endproclaim

\demo{Proof} Assume that a lattice $S$ of the rank 3 is hyperbolically 
reflective, and $\M\subset \La(S)$ is a fundamental polygon 
for the reflection group 
$W(S)$ and $A(\M)$ is its group of symmetries. 
Let us show that then either $h=0$ or $2$. Remind that $h$ is the 
number of conjugacy classes of central symmetries from $A(\M)$.  
Suppose that $\rho_1$ and $\rho_2$ are two primitive 
generalized lattice Weyl vectors for $A(\M)$ 
such that $\rho_1\not= \pm\rho_2$. Then 
a subgroup $G\subset A(\M)$ acts trivially on the sublattice 
$[\rho_1,\,\rho_2]$ of the rank 2 generated by $\rho_1$ and $\rho_2$, and 
any non-zero element of $[\rho_1,\,\rho_2]$ is a generalized 
lattice Weyl vector. All of them should have negative square since 
$S$ is hyperbolically reflective.  
It then follows that $[\rho_1,\,\rho_2]$ is negative definite. 
Its orthogonal complement is a one-dimensional hyperbolic 
lattice generated by an element $h$ with $h^2>0$. The $h$ is then 
also invariant with respect to a subgroup of finite index of $A(\M)$ and 
is a generalized lattice Weyl vector. It is impossible since $S$ is 
hyperbolically reflective. Thus $\rho_1=\pm \rho_2$. It follows that 
a primitive generalized lattice Weyl vector $\rho$ of $S$ is 
unique up to replacing by  $-\rho$. If $g\in A(\M)$, 
then $g(\rho)$ is also a generalized lattice Weyl vector of $S$. Thus, 
$g(\rho)=\pm \rho$. Let $l$ be the line of $\La(S)$ which is orthogonal 
to $\rho$. The group $A(\M)$ preserves this line and acts as a discrete 
group of motions of the line. In particular, the center of any central 
symmetry of $A(\M)$ belongs to the line $l$. The group $A(\M)$ does not 
have reflections and is infinite. It follows that its action on the line $l$ 
has only trivial kernel. An element of $A(\M)$ is a central symmetry if 
and only if it has a fixed point on the line $l$ (it is a central symmetry of 
the line $l$). Any discrete group of motions of a line either has no 
central symmetries or is generated by two central symmetries which give 
two different conjugacy classes of central symmetries of the group. 
Thus, the same is valid for $A(\M)$, and the invariant $h$ is equal to 
$0$ or $2$. 

Now, to prove Theorem 9.1.1, we should show that 
$d\le 100000$ for main hyperbolically reflective hyperbolic lattices 
with square-free determinant $d$. 

Let $D\subset l$ be a fundamental interval for the action of 
$A(\M)$ in $l$ and $C_D$ is its orthogonal to $l$ cylinder. 
The $C_D\cap \M$ is a fundamental domain for the group $O^+(S)$. 
By the theory of arithmetic 
groups or by the theory of automorphism groups of integral quadratic forms 
(e. g. see \cite{R}), the domain $C_D\cap \M$ is a finite polygon.    
It then follows that the locally finite polygon 
$\M$ is restricted hyperbolic relative to the line $l$ which is its axis. 

Now we apply to $S$ and $\M$ results of Sect. 8.5. Depending on the 
type of the narrow part of $\M$, these results give estimates on 
$d$. By Proposition 4.3.1, $d\le 2a_1$ where 
$a_1=\max_B{a_1(B)}$ where $B$ runs through all matrices $B$ 
describing possible narrow parts  
of the corresponding type of $\M$. They are given in Sect. 8.5. 
Unfortunately, for all types of narrow parts 
(except for the types (AI1) and (A10)), 
the invariant $a_1>50000$, and we don't get the desirable estimate 
($\le 100000$) for $d$. 

In Sect. 5.1 we had the same problem   
for the type (II0). To overcome this problem, we had considered  
matrices $B$ only for main hyperbolic lattices with square-free $d$. 
Using the matrices $B$, we  
had calculated the invariants $(d,\,\eta, h)$ of $S$ and 
found that if $h\le 1$, than $d\le 100000$.  
Actually, for the type (II0) we found all triplets  
$(d,\,\eta, h)$ with $h\le 1$. There were 132 the triplets, and all of 
them were contained in Table 3. 

Here we have to do the same for all types (AI1), (AI0), (AII1), 
(AII0), (AIII), (BI), (BII$_1$), (BII$_2$), (BIII) of narrow parts of 
$\M$. 

The following elementary Lemma is very useful. Actually, we 
had used its statements (a) and (d) without proof (since they are  
very simple).  

\proclaim{Lemma 9.1.2} Let $S$ be a main hyperbolic lattice of 
the rank three and with 
square-free determinant $d$, and $p$ denote any odd prime $p|d$. 
We have: 

(a) $r^2$ is square-free for any primitive root $r\in S$. If 
odd prime $p|r^2$, then $p|d$.  

(b) Let $r_1,\,r_2$ are two primitive roots of $S$ such that 
$r_2\not=\lambda r_1$, $\lambda \in \bq$, and 
${4(r_1,\,r_2)^2\over r_1^2r_2^2}=0$, $2$ or $3$.       
Then $\text{g.c.d}(r_1^2,\,r_2^2)\le 2$. If 
${4(r_1,\,r_2)^2\over r_1^2r_2^2}=0$, then  
$S\otimes \bz_p\cong \langle r_1^2 \rangle\oplus 
\langle r_2^2 \rangle \oplus 
\langle d/(r_1^2r_2^2)\rangle$. 
If ${4(r_1,\,r_2)^2\over r_1^2r_2^2}=2$, then 
$S\otimes \bz_p\cong \langle -1 \rangle\oplus 
\langle -1 \rangle \oplus 
\langle d \rangle$.  
If ${4(r_1,\,r_2)^2\over r_1^2r_2^2}=3$, then $3|d$ and  
$S\otimes \bz_p\cong \langle -2 \rangle\oplus 
\langle -6 \rangle \oplus 
\langle d/3 \rangle$. 

(c) Let primitive roots $r_1,\,r_2\in P(\M)$ and     
${4(r_1,\,r_2)^2\over r_1^2r_2^2} = 1$. Then 
$\text{g.c.d}(r_1^2,\,r_2^2)$ $\le 2$. Moreover 
$d\equiv 2\mod 3$ and 
$S\otimes \bz_p\cong \langle -2 \rangle\oplus 
\langle -6 \rangle\oplus \langle 3d \rangle$. 

(d) Assume that $S$ represents $0$. Then 
$S\cong U\oplus \langle -d \rangle$. In particular, 
this is true if there are two primitive roots $r_1,\,r_2\in P(\M)$  
such that ${4(r_1,\,r_2)^2\over r_1^2r_2^2}=4$.  
\endproclaim

\demo{Proof} (a). Assume that $r$ is a primitive root of $S$ and $K$ is 
its orthogonal complement. Then either $S=\bz r\oplus K$ or 
$S=[r, K, (r+k)/2]$ where $k\in K$ is a primitive element.  
Since the determinant of $S$ is square-free, it follows  
that $p^2\nmid r^2$ for any odd prime $p$, and 
$2^2\nmid r^2$ if $S=\bz r\oplus K$. Suppose that 
$S=[r, K, (r+k)/2]$ where $k\in K$ is a primitive element. 
Suppose that $4|r^2$. Since $((r+k)/2)^2=(r^2+k^2)/4 \in \bz$, 
it follows $4|k^2$.   
Let $k,\,s$ be a bases in $K$.  
We have $(k,s)/2=((r+k)/2,\,s)$ in $\bz$ and $2|(k,s)$. It follows 
$4|det(K)=k^2s^2-(k,s)^2$ and $16|det(\bz r\oplus K)=
r^2det(K)$. Since $[S:\bz r\oplus K]=2$, we have  
$det(S)=det(\bz r\oplus K)/4$ and $4|det(S)$. We get the contradiction. 
Thus $2^2 \nmid r^2$, and $r^2$ is square-free for any 
primitive root $r\in S$.    
  
(b)  
$4(r_1,\,r_2)^2/(r_1^2r_2^2)\in \bz$ because $r_1$ and $r_2$ 
are roots and $r_1^2|2(r_1,r_2)$, $r_2^2|2(r_1,r_2)$. 

Suppose that $(r_1,\,r_2)=0$ and $\bz h$ is the orthogonal complement 
to $[r_1,r_2]$. Since $r_1$ and $r_2$ define 
reflections of $S$, we have 
$S\otimes \bz_p=\bz_pr_1\oplus \bz_pr_2\oplus \bz_p h$ for any odd prime 
$p$. It follows that 
$p$ can divide only one of three numbers: $r_1^2$, $r_2^2$, $h^2$. 
Since $r_1^2$ and $r_2^2$ are square-free, it follows  
$\text{g.c.d}(r_1^2,\,r_2^2)\le 2$.  

Suppose that ${4(r_1,\,r_2)^2\over r_1^2r_2^2}=2$. Then we have 
(up to changing numeration and replacing $r_1$ by $-r_1$) 
that $r_1^2=-k$, $r_2^2=-2k$ where 
$k=(r_1,\,r_2)>0$.  
Then $s_{r_1}(r_2)=r_2^\prime$ is a primitive root such that 
$(r_2^\prime, \,r_2)=0$ and $(r_2^\prime)^2=r_2^2=-2k$. 
By consideration above, $\text{g.c.d}(r_2^2,\, r_2^2)\le 2$, $k=1$ 
and $r_2^2=-2$.  It follows that $\text{g.c.d}(r_1^2,r_2^2)\le 2$. 
We have $r_1^2=-1$, $r_2^2=-2$, $(r_1,\,r_2)=1$. 
Let $r_1^\prime =s_{r_2}(r_1)$. Then $(r_1^\prime,\,r_1)=0$ and 
$S\otimes \bz_p\cong \langle -1 \rangle\oplus \langle -1 \rangle
\oplus \langle d \rangle$. 

Suppose that ${4(r_1,\,r_2)^2\over r_1^2r_2^2}=3$. Then we have 
(up to changing numeration and replacing $r_1$ by $-r_1$) 
$r_1^2=-2k$, $r_2^2=-6k$, $(r_1,\,r_2)=3k$ where $k\in \bn$. 
Let $r_1^\prime=s_{r_2}(r_1)$, $r_2^\prime =s_{r_1}(r_2)$.    
Then $(r_1^\prime,\,r_2^\prime)=0$, and, as above, we get  
$\text{g.c.d}(r_1^2,\,r_2^2)\le 2$. Actually, we get $k=1$,  
$r_1^2=-2$, $r_2^2=-6$, $(r_1,\,r_2)=3$. Moreover, we get $3|d$ 
and $S\otimes \bz_p\cong \langle -2 \rangle \oplus \langle -6 \rangle 
\oplus \langle d/3 \rangle$.  

(c) We have $r_1^2=r_2^2=-2k$ and $(r_1,\,r_2)=k \in \bn$ where $k$ is 
odd and square-free. Suppose 
that $v=ar_1+br_2\in S$ where $a,\,b\in \bq$. Since 
$r_1^2|2(r_1,\,v)$, $r_2^2|2(r_2,\,v)$ and $v^2\in \bz$, 
simple calculations show 
that either $v=0$ or $3|k$ and $v=(r_1-r_2)/3$ modulo $\bz r_1+\bz r_2$. 
It follows that either $r_1$ and $r_2$ generate a primitive 
sublattice $[r_1,\,r_2]$ of $S$, or 
$3|k$ and $[r_1,\,r_2,\,(r_1-r_2)/3]$ is a primitive sublattice of $S$. 
Let us consider the first case and suppose that $r_1,\,r_2,\,u$ is 
a bases of $S$. Since $-2k=r_1^2|2(r_1,\,u)$ and $-2k=r_2^2|2(r_2,\,u)$, 
we then get that $k^2|det(S)$. It follows that $k=1$ and 
$\text{g.c.d}(r_1^2,\,r_2^2)\le 2$. Suppose that 
$3|k$ and  $[r_1,\,r_2,\,(r_1-r_2)/3]$ is a primitive sublattice of $S$. 
Like above, we then get $k=3$. Then $r_3=(r_1-r_2)/3$ is a root of 
$S$ because $r_3^2=-2$. We have $(r_1,\,r_3)=-3$ and 
$(r_2,\,r_3)=3$. It is impossible because $r_1,\,r_2\in P(\M)$ are 
orthogonal to sides of $\M$ defining a vertex of $\M$ with the angle 
$\pi/3$. But the line orthogonal to the root 
$r_3$ contains the same vertex and divides the angle. 
Thus, we have proved that $k=1$ and $r_1^2=r_2^2=-2$,  
$(r_1,\,r_2)=1$. Moreover, the sublattice $[r_1,\,r_2]$ is  
primitive in $S$. Let $\bz h$ be the orthogonal complement to 
$[r_1,\,r_2]$. Suppose that $S=[r_1,\,r_2]\oplus \bz h$. Then 
$r_3=(r_1-r_2)$ is a root of $S$ such that $(r_3,\,r_1)=-3$ and 
$(r_3,\,r_2)=3$. Like above, we get a contradiction. Thus, 
the lattice $S$ is bigger than its sublattice of finite 
index $[r_1,\,r_2]\oplus \bz h$. There are rational (not all integral) 
$a,\,b,\,c$ such that $v=ar_1+br_2+ch\in S$. Simple calculations 
show that modulo $[r_1,\,r_2]\oplus \bz h$ the element 
$v=(r_1-r_2)/3 \pm h/3$ where $h^2=n\equiv 6 \mod 9$. It follows that 
$d=n/3\equiv 2 \mod 3$. Simple calculations show that  
$S\otimes \bz_p\cong \langle -2 \rangle\oplus \langle -6 \rangle
\oplus \langle d/3 \rangle$.     

(d) Suppose that $S$ represents zero. Then for any prime $p$, the lattice 
$S\otimes \bz_p$ represents zero, and  
$S\otimes \bz_p\cong U\otimes \bz_p\oplus \langle -d \rangle$ because 
$d$ is square-free. It follows that $S$ has the same genus as 
$U\oplus \langle -d \rangle$. 
By Proposition 2.2.4, $S\cong U\oplus \langle -d \rangle$. 

Suppose that there are two primitive roots $r_1,\,r_2\in P(\M)$  
such that ${4(r_1,\,r_2)^2\over (r_1^2r_2^2)}=4$. 
Then we have: either $r_1^2=r_2^2=-k$, $(r_1,\,r_2)=k>0$ and 
$(r_1+r_2)^2=0$ or
$r_1^2=-k$, $r_2^2=-4k$, $(r_1,\,r_2)=2k>0$ and $(2r_1+r_2)^2=0$. 
In both cases $S$ represents $0$. 

This finishes the proof of the lemma.  
\enddemo 

Now suppose that $S$ is a hyperbolically reflective main hyperbolic 
lattice of the rank 3 with square-free determinant $d$. Let 
$r_1,\dots, r_k \in P(\M)$ are primitive roots defining a narrow parts 
of $\M$ (from Sect. 8.5) and $\Gamma=((r_i,\,r_j))$, $1\le i,\,j\le k$. 
In Sect. 8.5 we described possible matrices $B=\Gamma_{\pr}$ for all 
hyperbolically reflective lattices $S$ of the rank 3. 
Here we add additional 
conditions when $S$ is a main hyperbolic lattice with square-free 
determinant, and describe all possible 
matrices $\Gamma=\lambda B$, $\lambda \in \bn$, for these lattices. 
We denote by $G_0$ a lattice with the matrix $B$ (it is defined 
by the matrix $B$ modulo its kernel). We denote 
$d_1=\det(G_0)$, and $d_2$ the product of all odd prime $p\vert d_1$ 
such that $\nu_p(d_1)\equiv 0 \mod2$.  

By considerations above and Lemma 9.1.2, 
we should add to conditions of Sect. 8.5  
additional conditions (1) --- (10) below: 

{\it 
(1) $b_{ii}$ is square-free;   
if odd prime $p|b_{ii}$, then $\nu_p(d_1)\equiv 1\mod 2$, $1\le i\le k$; 
if $2|(b_{11}\cdots b_{kk})$ and 
$\text{g.c.d}(2,b_{11},b_{22},...,b_{kk})=1$, then 
$\nu_2(d_1)\equiv 0\mod 2$.     

(2) If $\alpha_{ij}\le 3$, then 
$\text{g.c.d}(b_{ii},\,b_{jj})\le 2$, $1\le i<j\le k$. 

(3) $\lambda|2d_2$;  
$\text{g.c.d}(\lambda, \, b_{11}b_{22}\cdots b_{kk})=1)$;   
$\lambda=1$ or $2$ if there exists $1\le i<j\le k$ such that 
$\alpha_{ij}\le 3$.  

(4) The invariant $d=\det(S)$ is the product of all prime $p$ such that 
$\nu_p(d_1\lambda)\equiv 1 \mod 2$. 

(5) If $\alpha_{ij}=0$, then the invariant $\eta$ of $S$ is defined by 
$$
(-1)^{\eta_p}=
\cases
\left({\lambda b_{ii}/p\over p}\right) &\text{if $p\vert b_{ii}$,}\\
\left({\lambda b_{jj}/p\over p}\right) &\text{if $p\vert b_{jj}$,}\\
\left({db_{ii}b_{jj}/p\over p}\right)
&\text{otherwise}
\endcases
$$
for odd prime $p\vert d$. 

(6) If $\alpha_{ij}=1$, $1\le i<j\le k$, then $d\equiv 2\mod 3$ and  
$(-1)^{\eta_p}=\left({3d\over p}\right)$ for odd prime $p\vert d$.

(7) If $\alpha_{ij}=2$, $1\le i<j\le k$, then the invariant 
$\eta$ of $S$ is defined by $(-1)^{\eta_p}=\left({d\over p}\right)$ 
for odd prime $p\vert d$. 

(8) If $\alpha_{ij}=3$, $1\le i<j\le k$, then $\eta_3=0$ and 
$(-1)^{\eta_p}=\left({d/3\over p}\right)$ for odd prime $p\vert (d/3)$.
 
(9) If $\alpha_{ij}=4$, $1\le i<j\le k$, then  
$(-1)^{\eta_p}=\left({-d\over p}\right)$ for odd prime $p\vert d$.

(10) $h=hnr(S)=h(d,\eta)=0$ or $2$. }

Here we use that $\det(K_1)=\det(K)[K:K_1]^2$ for a sublattice 
$K_1\subset K$ of finite index, and 
$\det(K(\lambda))=\lambda^{\rk K}\det K$ where $K(\lambda)$ is a lattice 
$K$ with the form multiplied by $\lambda\in \bq$. Moreover, we use that 
$d$ is odd if the lattice $S$ is odd. 
(See similar considerations for the type (II0) in Sect. 5.1.)    

In Sect. 8.5.1 (the type (AI1) of the narrow part),  
$1\le \alpha_{12}\le \alpha_{23}\le 4$. We can use $i=1$, $j=2$ in  
(6)--(9) to find the invariant $\eta$ of $S$. Moreover, $\lambda\le 2$,  
if $\alpha_{12}\le 3$. In Appendix 2, we give the Program 22: funda11.main
which enumerates all possible matrices $\Gamma$ satisfying 
conditions of Sect. 8.5.1 and the additional conditions (1) --- (10), 
and their invariants $(d,\eta,h)$. There are 28 these triplets given 
in Table 7. 

In Sect. 8.5.2 (the type (AI0) of the narrow part),  
$\alpha_{12}=0$. Then $\lambda \le 2$, and we can use $i=1$, $j=2$ in 
(5) to find the invariant $\eta$ of $S$. In Appendix 2, we give the 
Program 23: funda10.main which enumerates all possible matrices 
$\Gamma$ satisfying Sect. 8.5.2 and (1) --- (10), and their invariants 
$(d,\eta,h)$. There are 86 these triplets given in Table 7. 

In Sect. 8.5.3 (the type (AII1) of the narrow part), 
$\alpha_{12}=0$. Then $\lambda \le 2$ and we can use $i=1$, $j=2$ in 
(5) to find the invariant $\eta$ of $S$. In Appendix 2 we give the 
Program 24: funda21.main which enumerates all possible matrices $\Gamma$ 
satisfying Sect. 8.5.3 and (1) --- (10), 
and their invariants $(d,\eta,h)$. There are 44 these triplets 
given in Table 7. 

In Sect. 8.5.4 (the type (AII0) of the narrow part), 
$\alpha_{12}=0$. Then $\lambda \le 2$ and we can use $i=1$, $j=2$ in 
(5) to find the invariant $\eta$ of $S$. In Appendix 2 we give the 
Program 25: funda20.main which enumerates all possible matrices $\Gamma$ 
satisfying Sect. 8.5.4 and (1) --- (10), and their invariants 
$(d,\eta,h)$. There are 164 these triplets given in Table 7.

 In Sect. 8.5.5 (the type (AIII) of the narrow part), 
$\alpha_{12}=0$. Then $\lambda \le 2$ and we can use $i=1$, $j=2$ in 
(5) to find the invariant $\eta$ of $S$. In Appendix 2 we give the 
Program 26: funda3.main which enumerates all possible matrices $\Gamma$ 
satisfying Sect. 8.5.5 and (1) --- (10), and their invariants 
$(d,\eta,h)$. There are 69 these triplets given in Table 7.

In Sect. 8.5.6 (the type (BI) of the narrow part), 
$1\le \alpha_{12}=\le 4$. We can use $i=1$, $j=2$ in 
(6)---(9) to find the invariant $\eta$ of $S$. Moreover, $\lambda\le 2$ 
if $\alpha_{12}\le 3$. In Appendix 2 we give the 
Program 27: fundb1.main which enumerates all possible matrices $\Gamma$ 
satisfying Sect. 8.5.6 and (1) --- (10), and their invariants 
$(d,\eta,h)$. There are 155 these triplets given in Table 7.

In Sect. 8.5.7 (the type (BII$_1$) of the narrow part), 
$\alpha_{12}=0$. Then $\lambda \le 2$ and we can use $i=1$, $j=2$ in 
(5) to find the invariant $\eta$ of $S$. In Appendix 2 we give the 
Program 28: fundb21.main which enumerates all possible matrices $\Gamma$ 
satisfying Sect. 8.5.7 and (1) --- (10), and their invariants 
$(d,\eta,h)$. There are 234 these triplets given in Table 7.

In Sect. 8.5.8 (the type (BII$_2$) of the narrow part), 
$\alpha_{12}=0$. Then $\lambda \le 2$ and we can use $i=1$, $j=2$ in 
(5) to find the invariant $\eta$ of $S$. In Appendix 2 we give the 
Program 29: fundb22.main which enumerates all possible matrices $\Gamma$ 
satisfying Sect. 8.5.8 and (1) --- (10), and their invariants 
$(d,\eta,h)$. There are 223 these triplets given in Table 7.

In Sect. 8.5.9 (the type (BIII) of the narrow part), 
$\alpha_{12}=0$. Then $\lambda \le 2$ and we can use $i=1$, $j=2$ in 
(5) to find the invariant $\eta$ of $S$. In Appendix 2 we give the 
Program 30: fundb3.main which enumerates all possible matrices $\Gamma$ 
satisfying Sect. 8.5.9 and (1) --- (10), and their invariants 
$(d,\eta,h)$. There are 216 these triplets given in Table 7.

In Table 7, for each triplet $(d,\eta,h)$ of Table 3 with $h=0$ and 
Table 6 we give the types (AI1 --- BIII) of narrow parts which give 
these triplets. All these triplets have $d\le 100000$ and are contained 
in Tables 3 and 6. 

This finishes the proof of Theorem 9.1.1. 
\enddemo

\centerline{\bf Table 7}
\centerline{Possible types of narrow parts from Theorem 8.3.4} 
\centerline{of main hyperbolic lattices of the rank $3$} 
\centerline{with square-free determinant $d$ and $h=0$ or $2$} 

\vskip20pt 
  
{\settabs 12 \columns
\+$n=1$ $d=1$ $\eta=0$ $h=0$ &&&&     
AI1 & AI0 &AII1& AII0 & AIII & BI & BII$_1$ & BII$_2$ & BIII\cr
\+$n=2$ $d=2$ $\eta=0$ $h=0$ &&&&    
AI1 & AI0 &   &&&        BI   &&      BII$_2$&\cr
\+$n=3$ $d=3$ $\eta=0$ $h=0$&&&&      
AI1& AI0& AII1& AII0& AIII& BI& BII$_1$& BII$_2$& BIII\cr
\+$n=4$ $d=3$ $\eta=1$ $h=0$&&&&      
AI1& AI0& AII1& AII0& AIII& BI& BII$_1$& BII$_2$& BIII\cr
\+$n=5$ $d=5$ $\eta=0$ $h=0$&&&&      
AI1& AI0& AII1& AII0& AIII& BI& BII$_1$& BII$_2$& BIII\cr
\+$n=6$ $d=5$ $\eta=1$ $h=0$&&&&      
AI1& AI0& AII1& AII0& AIII& BI& BII$_1$& BII$_2$& BIII\cr
\+$n=7$ $d=6$ $\eta=0$ $h=0$&&&&      
AI1& AI0& AII1&     &     & BI& BII$_1$& BII$_2$& BIII\cr
\+$n=9$ $d=7$ $\eta=1$ $h=0$&&&&      
AI1& AI0& AII1& AII0&     & BI& BII$_1$& BII$_2$& BIII\cr
\+$n=10$ $d=10$ $\eta=1$ $h=0$&&&&    
AI1& AI0& AII1& AII0& AIII& BI& BII$_1$& BII$_2$& BIII\cr
\+$n=11$ $d=11$ $\eta=0$ $h=0$&&&&    
AI1& AI0& AII1& AII0&     & BI& BII$_1$& BII$_2$& BIII\cr
\+$n=15$ $d=14$ $\eta=1$ $h=0$&&&&    
AI1& AI0&     &     &     & BI& BII$_1$& BII$_2$&\cr
\+$n=16$ $d=15$ $\eta=0$ $h=0$&&&&    
AI1& AI0& AII1& AII0& AIII& BI& BII$_1$& BII$_2$& BIII\cr
\+$n=17$ $d=15$ $\eta=1$ $h=0$&&&&      
 &    &     & AII0& AIII&   & BII$_1$& BII$_2$& BIII\cr
\+$n=18$ $d=15$ $\eta=2$ $h=0$&&&&  
     & AI0& AII1& AII0& AIII& BI& BII$_1$& BII$_2$& BIII\cr
\+$n=19$ $d=15$ $\eta=3$ $h=0$&&&&  
  AI1& AI0& AII1& AII0& AIII& BI& BII$_1$& BII$_2$& BIII\cr
\+$n=20$ $d=17$ $\eta=0$ $h=0$&&&&  
     & AI0&     & AII0& AIII& BI& BII$_1$& BII$_2$& BIII\cr
\+$n=22$ $d=21$ $\eta=0$ $h=0$&&&&  
  AI1& AI0& AII1& AII0& AIII& BI& BII$_1$& BII$_2$& BIII\cr
\+$n=24$ $d=21$ $\eta=2$ $h=0$&&&& 
   AI1& AI0& AII1& AII0& AIII& BI& BII$_1$& BII$_2$& BIII\cr
\+$n=25$ $d=21$ $\eta=3$ $h=0$&&&&       &    &     & 
AII0& AIII&   & BII$_1$& BII$_2$& BIII\cr
\+$n=28$ $d=26$ $\eta=1$ $h=0$&&&&    
AI1& AI0& AII1&     &     & BI& BII$_1$& BII$_2$& BIII\cr
\+$n=32$ $d=30$ $\eta=2$ $h=0$&&&&    
AI1& AI0&     & AII0& AIII& BI& BII$_1$& BII$_2$& BIII\cr
\+$n=34$ $d=33$ $\eta=1$ $h=0$&&&&    
   & AI0&     & AII0& AIII& BI& BII$_1$& BII$_2$& BIII\cr
\+$n=38$ $d=35$ $\eta=0$ $h=0$&&&&    
AI1& AI0&     & AII0&     & BI& BII$_1$& BII$_2$& BIII\cr
\+$n=40$ $d=35$ $\eta=3$ $h=0$&&&&    
AI1& AI0& AII1& AII0& AIII& BI& BII$_1$& BII$_2$& BIII\cr
\+$n=42$ $d=39$ $\eta=0$ $h=0$&&&&    
AI1& AI0&     & AII0& AIII& BI& BII$_1$& BII$_2$& BIII\cr
\+$n=46$ $d=42$ $\eta=0$ $h=0$&&&&    
   & AI0& AII1& AII0&     & BI& BII$_1$& BII$_2$& BIII\cr
\+$n=48$ $d=51$ $\eta=0$ $h=0$&&&&    
AI1& AI0& AII1& AII0& AIII& BI& BII$_1$& BII$_2$& BIII\cr
\+$n=56$ $d=65$ $\eta=0$ $h=0$&&&&    
   &    &     & AII0&     &   & BII$_1$& BII$_2$& BIII\cr
\+$n=61$ $d=69$ $\eta=3$ $h=0$&&&&    
   &    &     & AII0& AIII&   & BII$_1$& BII$_2$& BIII\cr
\+$n=64$ $d=74$ $\eta=1$ $h=0$&&&&    
   & AI0&     &     &     & BI& BII$_1$& BII$_2$&\cr
\+$n=65$ $d=77$ $\eta=1$ $h=0$&&&&   
    &    &     & AII0&     & BI& BII$_1$& BII$_2$& BIII\cr
\+$n=66$ $d=77$ $\eta=3$ $h=0$&&&&    
AI1& AI0&     & AII0& AIII& BI& BII$_1$& BII$_2$& BIII\cr
\+$n=67$ $d=78$ $\eta=2$ $h=0$&&&&    
AI1& AI0& AII1& AII0&     & BI& BII$_1$& BII$_2$& BIII\cr
\+$n=74$ $d=95$ $\eta=0$ $h=0$&&&&     
  & AI0&     &     &     & BI& BII$_1$& BII$_2$& BIII\cr
\+$n=78$ $d=105$ $\eta=1$ $h=0$&&&&    
  & AI0&     & AII0&     & BI& BII$_1$& BII$_2$& BIII\cr
\+$n=80$ $d=105$ $\eta=3$ $h=0$&&&&    
  &    &     & AII0& AIII&   & BII$_1$& BII$_2$& BIII\cr
\+$n=82$ $d=105$ $\eta=5$ $h=0$&&&&    
  &    &     & AII0& AIII&   & BII$_1$& BII$_2$& BIII\cr
\+$n=83$ $d=105$ $\eta=6$ $h=0$&&&&   
AI1& AI0& AII1& AII0& AIII& BI& BII$_1$& BII$_2$& BIII\cr
\+$n=84$ $d=105$ $\eta=7$ $h=0$&&&&   
   &    &     & AII0& AIII&   & BII$_1$& BII$_2$& BIII\cr
\+$n=87$ $d=114$ $\eta=2$ $h=0$&&&&   
   & AI0&     &     &     & BI& BII$_1$& BII$_2$&\cr
\+$n=99$ $d=165$ $\eta=0$ $h=0$&&&&   
AI1& AI0& AII1& AII0& AIII& BI& BII$_1$& BII$_2$& BIII\cr
\+$n=104$ $d=165$ $\eta=6$ $h=0$&&&&   
  &    &     & AII0& AIII&   & BII$_1$& BII$_2$& BIII\cr
\+$n=106$ $d=170$ $\eta=2$ $h=0$&&&&  
   & AI0&     &     &     & BI& BII$_1$& BII$_2$&\cr
\+$n=108$ $d=186$ $\eta=0$ $h=0$&&&&  
   & AI0& AII1&     &     & BI& BII$_1$& BII$_2$& BIII\cr
\+$n=110$ $d=195$ $\eta=3$ $h=0$&&&&   
  &    &     & AII0& AIII&   & BII$_1$& BII$_2$& BIII\cr
\+$n=112$ $d=195$ $\eta=6$ $h=0$&&&&  
   & AI0& AII1& AII0& AIII& BI& BII$_1$& BII$_2$& BIII\cr
\+$n=115$ $d=210$ $\eta=2$ $h=0$&&&&  
   &    &     & AII0& AIII&   & BII$_1$& BII$_2$& BIII\cr
\+$n=118$ $d=219$ $\eta=0$ $h=0$&&&&  
AI1& AI0&     & AII0&     & BI& BII$_1$& BII$_2$& BIII\cr
\+$n=120$ $d=231$ $\eta=2$ $h=0$&&&&  
   &    &     & AII0&     &   & BII$_1$& BII$_2$& BIII\cr
\+$n=121$ $d=231$ $\eta=7$ $h=0$&&&&   
  & AI0&     & AII0& AIII& BI& BII$_1$& BII$_2$& BIII\cr
\+$n=125$ $d=255$ $\eta=6$ $h=0$&&&&  
   &    & AII1& AII0& AIII& BI& BII$_1$& BII$_2$& BIII\cr
\+$n=133$ $d=285$ $\eta=5$ $h=0$&&&&   
  &    &     & AII0& AIII&   & BII$_1$& BII$_2$& BIII\cr
\+$n=139$ $d=330$ $\eta=3$ $h=0$&&&&  
   &    &     & AII0& AIII&   & BII$_1$& BII$_2$& BIII\cr
\+$n=141$ $d=345$ $\eta=6$ $h=0$&&&&  
   & AI0& AII1& AII0& AIII& BI& BII$_1$& BII$_2$& BIII\cr
\+$n=152$ $d=399$ $\eta=4$ $h=0$&&&&  
   & AI0& AII1& AII0& AIII& BI& BII$_1$& BII$_2$& BIII\cr
\+$n=155$ $d=435$ $\eta=0$ $h=0$&&&&   
  & AI0& AII1& AII0&     & BI& BII$_1$& BII$_2$& BIII\cr
\+$n=164$ $d=546$ $\eta=2$ $h=0$&&&&  
   &    &     &     &     & BI&        &        &\cr
\+$n=165$ $d=546$ $\eta=4$ $h=0$&&&&   
  &    &     & AII0&     &   & BII$_1$& BII$_2$& BIII\cr
\+$n=187$ $d=1155$ $\eta=2$ $h=0$&&&& 
   &    &     & AII0& AIII&   & BII$_1$& BII$_2$& BIII\cr
\+$n=190$ $d=1155$ $\eta=14$ $h=0$&&&&  
 & AI0& AII1& AII0& AIII& BI& BII$_1$& BII$_2$& BIII\cr
\+$n=202$ $d=2730$ $\eta=6$ $h=0$&&&&  
  &    &     &     &     &   & BII$_1$& BII$_2$& BIII\cr
\+&&&&&&&&&&&&\cr
\+$n=1$ $d=17$ $\eta=1$ $h=2$&&&&
     &     &     &     &     &BI &&&\cr                        
\+$n=2$ $d=19$ $\eta=1$ $h=2$&&&&  
     &  AI0&     &     &     & BI& BII$_1$&  BII$_2$&\cr   
\+$n=3$ $d=23$ $\eta=1$ $h=2$&&&&   
   &  AI0&     & AII0& AIII& BI& BII$_1$&  BII$_2$&  BIII\cr 
\+$n=4$ $d=31$ $\eta=0$ $h=2$&&&&   
   &  AI0&     &     &     & BI& BII$_1$&  BII$_2$&\cr  
\+$n=5$ $d=35$ $\eta=1$ $h=2$&&&&  
    &  AI0& AII1& AII0& AIII& BI& BII$_1$&  BII$_2$&  BIII\cr 
\+$n=6$ $d=37$ $\eta=0$ $h=2$&&&&  
    &  AI0&     &     &     & BI&&&\cr
\+$n=7$ $d=37$ $\eta=1$ $h=2$&&&&   
   &     &     & AII0&     &   & BII$_1$&  BII$_2$&  BIII\cr 
\+$n=8$ $d=39$ $\eta=1$ $h=2$&&&&    
  &  AI0&     & AII0& AIII& BI& BII$_1$&  BII$_2$&  BIII\cr 
\+$n=9$ $d=43$ $\eta=0$ $h=2$&&&&      
&     &     &     &     & BI& BII$_1$&  BII$_2$&  BIII\cr 
\+$n=10$ $d=46$ $\eta=1$ $h=2$&&&&     
&  AI0&     &     &     & BI& BII$_1$&  BII$_2$&\cr 
\+$n=11$ $d=47$ $\eta=0$ $h=2$&&&&     
&  AI0&     &     &     & BI&&&\cr
\+$n=12$ $d=51$ $\eta=2$ $h=2$&&&&     
&  AI0&     &     &     & BI& BII$_1$&&\cr 
\+$n=13$ $d=53$ $\eta=0$ $h=2$&&&&     
&  AI0&     &     &     & BI&&&\cr
\+$n=14$ $d=53$ $\eta=1$ $h=2$&&&&     
&  AI0&     & AII0& AIII& BI& BII$_1$&  BII$_2$&  BIII\cr 
\+$n=15$ $d=55$ $\eta=1$ $h=2$&&&&     
&     &     & AII0& AIII&   & BII$_1$&  BII$_2$&  BIII\cr 
\+$n=16$ $d=57$ $\eta=3$ $h=2$&&&&     
&     &     & AII0&     &   & BII$_1$&  BII$_2$&  BIII\cr 
\+$n=17$ $d=58$ $\eta=1$ $h=2$&&&&     
&  AI0&     &     &     & BI& BII$_1$&  BII$_2$&  BIII\cr 
\+$n=18$ $d=59$ $\eta=0$ $h=2$&&&&  
AI1&     & AII1&     &     & BI&        &         &  BIII\cr 
\+$n=19$ $d=62$ $\eta=1$ $h=2$&&&&     
&  AI0&     &     &     & BI&        &  BII$_2$&\cr 
\+$n=20$ $d=65$ $\eta=1$ $h=2$&&&&     
&     &     & AII0& AIII&   & BII$_1$&  BII$_2$&  BIII\cr 
\+$n=21$ $d=69$ $\eta=0$ $h=2$&&&&     
&  AI0&     & AII0&     & BI& BII$_1$&  BII$_2$&  BIII\cr 
\+$n=22$ $d=69$ $\eta=2$ $h=2$&&&&     
&  AI0&     & AII0&     & BI& BII$_1$&         &  BIII\cr 
\+$n=23$ $d=70$ $\eta=0$ $h=2$&&&&     
&  AI0&     &     &     & BI&&&\cr
\+$n=24$ $d=77$ $\eta=0$ $h=2$&&&&     
&     &     & AII0& AIII&   & BII$_1$&  BII$_2$&  BIII\cr 
\+$n=25$ $d=79$ $\eta=0$ $h=2$&&&&     
&  AI0&     &     &     & BI& BII$_1$&&\cr
\+$n=26$ $d=85$ $\eta=2$ $h=2$&&&&     
&     &     & AII0&     &   & BII$_1$&  BII$_2$&  BIII\cr 
\+$n=27$ $d=85$ $\eta=3$ $h=2$&&&&     
&     &     &     &     & BI&        &  BII$_2$&\cr 
\+$n=28$ $d=87$ $\eta=0$ $h=2$&&&&     
&     &     & AII0& AIII& BI& BII$_1$&  BII$_2$&  BIII\cr 
\+$n=29$ $d=91$ $\eta=0$ $h=2$&&&&     
&     &     & AII0& AIII&   & BII$_1$&  BII$_2$&  BIII\cr 
\+$n=30$ $d=91$ $\eta=1$ $h=2$&&&&     
&     &     & AII0&     &   & BII$_1$&  BII$_2$&  BIII\cr 
\+$n=31$ $d=93$ $\eta=3$ $h=2$&&&&     
&     &     & AII0& AIII&   & BII$_1$&  BII$_2$&  BIII\cr 
\+$n=32$ $d=106$ $\eta=1$ $h=2$&&&&    
&  AI0&     &     &     & BI& BII$_1$&  BII$_2$&\cr 
\+$n=33$ $d=107$ $\eta=0$ $h=2$&&&&    
&     & AII1&     &     & BI& BII$_1$&         &  BIII\cr 
\+$n=34$ $d=113$ $\eta=0$ $h=2$&&&&    
&     &     & AII0&     & BI& BII$_1$&  BII$_2$&  BIII\cr 
\+$n=35$ $d=114$ $\eta=1$ $h=2$&&&&    
&     &     &     &     &   & BII$_1$&  BII$_2$&  BIII\cr 
\+$n=36$ $d=115$ $\eta=0$ $h=2$&&&&    
&     &     & AII0&     &   & BII$_1$&  BII$_2$&  BIII\cr 
\+$n=37$ $d=115$ $\eta=3$ $h=2$&&&&    
&  AI0&     & AII0& AIII& BI& BII$_1$&  BII$_2$&  BIII\cr 
\+$n=38$ $d=119$ $\eta=0$ $h=2$&&&&    
&  AI0&     & AII0&     & BI& BII$_1$&  BII$_2$&  BIII\cr 
\+$n=39$ $d=122$ $\eta=1$ $h=2$&&&&    
&  AI0&     &     &     & BI& BII$_1$&  BII$_2$&\cr 
\+$n=40$ $d=123$ $\eta=0$ $h=2$&&&&    
&  AI0&     & AII0&     & BI& BII$_1$&  BII$_2$&  BIII\cr 
\+$n=41$ $d=123$ $\eta=1$ $h=2$&&&&    
&     &     & AII0&     &   & BII$_1$&  BII$_2$&  BIII\cr 
\+$n=42$ $d=123$ $\eta=3$ $h=2$&&&&    
&     & AII1& AII0&     & BI& BII$_1$&  BII$_2$&  BIII\cr 
\+$n=43$ $d=129$ $\eta=0$ $h=2$&&&&    
&     & AII1&     &     & BI& BII$_1$&  BII$_2$&  BIII\cr 
\+$n=44$ $d=129$ $\eta=2$ $h=2$&&&&    
&     & AII1& AII0&     & BI& BII$_1$&  BII$_2$&  BIII\cr 
\+$n=45$ $d=133$ $\eta=1$ $h=2$&&&&    
&     &     & AII0&     & BI&        &         &  BIII\cr 
\+$n=46$ $d=133$ $\eta=3$ $h=2$&&&&    
&     &     & AII0&     &   & BII$_1$&  BII$_2$&  BIII\cr 
\+$n=47$ $d=138$ $\eta=0$ $h=2$&&&&    
&     &     &     &     & BI& BII$_1$&  BII$_2$&  BIII\cr 
\+$n=48$ $d=141$ $\eta=0$ $h=2$&&&&    
&     &     &     &     & BI& BII$_1$&&\cr
\+$n=49$ $d=143$ $\eta=3$ $h=2$&&&& AI1
&  AI0& AII1& AII0&     & BI& BII$_1$&  BII$_2$&  BIII\cr 
\+$n=50$ $d=145$ $\eta=3$ $h=2$&&&&    
&     &     & AII0&     &   & BII$_1$&  BII$_2$&  BIII\cr 
\+$n=51$ $d=155$ $\eta=0$ $h=2$&&&&    
&     &     &     &     & BI& BII$_1$&  BII$_2$&\cr 
\+$n=52$ $d=159$ $\eta=0$ $h=2$&&&&    
&  AI0&     &     &     &   & BII$_1$&         &  BIII\cr 
\+$n=53$ $d=159$ $\eta=3$ $h=2$&&&&    
&  AI0&     & AII0&     & BI& BII$_1$&         &  BIII\cr 
\+$n=54$ $d=165$ $\eta=2$ $h=2$&&&&    
&     & AII1& AII0&     & BI& BII$_1$&  BII$_2$&  BIII\cr 
\+$n=55$ $d=174$ $\eta=2$ $h=2$&&&&    
&     &     & AII0&     & BI& BII$_1$&  BII$_2$&  BIII\cr 
\+$n=56$ $d=177$ $\eta=1$ $h=2$&&&&    
&     &     &     &     & BI& BII$_1$&  BII$_2$&  BIII\cr 
\+$n=57$ $d=182$ $\eta=0$ $h=2$&&&&    
&  AI0&     &     &     & BI&        &  BII$_2$&\cr 
\+$n=58$ $d=183$ $\eta=0$ $h=2$&&&&    
&  AI0& AII1& AII0&     & BI& BII$_1$&  BII$_2$&  BIII\cr 
\+$n=59$ $d=185$ $\eta=0$ $h=2$&&&&    
&     &     & AII0&     &   & BII$_1$&  BII$_2$&  BIII\cr 
\+$n=60$ $d=191$ $\eta=0$ $h=2$&&&&    
&     &     &     &     & BI&&&\cr
\+$n=61$ $d=195$ $\eta=0$ $h=2$&&&&    
&     &     & AII0&     &   & BII$_1$&  BII$_2$&  BIII\cr 
\+$n=62$ $d=195$ $\eta=1$ $h=2$&&&&    
&     &     & AII0& AIII&   & BII$_1$&  BII$_2$&  BIII\cr 
\+$n=63$ $d=195$ $\eta=7$ $h=2$&&&&    
&     &     & AII0&     &   & BII$_1$&  BII$_2$&  BIII\cr 
\+$n=64$ $d=205$ $\eta=0$ $h=2$&&&&    
&  AI0&     &     &     &   & BII$_1$&  BII$_2$&  BIII\cr 
\+$n=65$ $d=210$ $\eta=1$ $h=2$&&&&    
&     &     &     &     &   & BII$_1$&  BII$_2$&  BIII\cr 
\+$n=66$ $d=210$ $\eta=7$ $h=2$&&&&    
&  AI0&     & AII0& AIII& BI& BII$_1$&  BII$_2$&  BIII\cr 
\+$n=67$ $d=213$ $\eta=3$ $h=2$&&&&    
&     &     & AII0&     &   & BII$_1$&  BII$_2$&  BIII\cr 
\+$n=68$ $d=217$ $\eta=3$ $h=2$&&&&    
&     &     & AII0&     &   & BII$_1$&  BII$_2$&  BIII\cr 
\+$n=69$ $d=218$ $\eta=1$ $h=2$&&&&    
&  AI0&     &     &     & BI&        &  BII$_2$&\cr  
\+$n=70$ $d=221$ $\eta=0$ $h=2$&&&&    
&     &     &     &     &   &        &  BII$_2$&\cr 
\+$n=71$ $d=221$ $\eta=1$ $h=2$&&&&    
&     &     & AII0&     &   & BII$_1$&  BII$_2$&  BIII\cr 
\+$n=72$ $d=222$ $\eta=2$ $h=2$&&&&    
&  AI0&     &     &     & BI& BII$_1$&  BII$_2$&  BIII\cr 
\+$n=73$ $d=230$ $\eta=0$ $h=2$&&&&    
&  AI0&     &     &     & BI&&&\cr
\+$n=74$ $d=231$ $\eta=3$ $h=2$&&&&    
&     &     & AII0& AIII&   & BII$_1$&  BII$_2$&  BIII\cr 
\+$n=75$ $d=231$ $\eta=4$ $h=2$&&&&    
&  AI0&     & AII0&     & BI& BII$_1$&  BII$_2$&  BIII\cr 
\+$n=76$ $d=235$ $\eta=0$ $h=2$&&&&    
&     &     & AII0&     &   & BII$_1$&         &  BIII\cr 
\+$n=77$ $d=239$ $\eta=0$ $h=2$&&&&    
&     &     &     &     & BI&&&\cr
\+$n=78$ $d=255$ $\eta=2$ $h=2$&&&&    
&     &     & AII0& AIII&   & BII$_1$&  BII$_2$&  BIII\cr 
\+$n=79$ $d=255$ $\eta=3$ $h=2$&&&&    
&     &     & AII0& AIII&   & BII$_1$&  BII$_2$&  BIII\cr 
\+$n=80$ $d=258$ $\eta=2$ $h=2$&&&&    
&  AI0&     &     &     & BI&        &  BII$_2$&\cr 
\+$n=81$ $d=273$ $\eta=2$ $h=2$&&&&    
&     &     &     &     &   & BII$_1$&  BII$_2$&  BIII\cr 
\+$n=82$ $d=282$ $\eta=0$ $h=2$&&&&    
&  AI0&     &     &     & BI& BII$_1$&  BII$_2$&\cr 
\+$n=83$ $d=287$ $\eta=3$ $h=2$&&&&    
&     &     &     &     & BI&&&\cr
\+$n=84$ $d=290$ $\eta=0$ $h=2$&&&&    
&     &     &     &     & BI&&&\cr
\+$n=85$ $d=290$ $\eta=3$ $h=2$&&&&    
&     &     & AII0&     & BI& BII$_1$&  BII$_2$&  BIII\cr 
\+$n=86$ $d=301$ $\eta=1$ $h=2$&&&&    
&     &     &     &     &   &&&\cr  
\+$n=87$ $d=301$ $\eta=3$ $h=2$&&&&    
&     &     &     &     &   & BII$_1$&  BII$_2$&  BIII\cr 
\+$n=88$ $d=303$ $\eta=3$ $h=2$&&&&    
&     &     & AII0&     & BI& BII$_1$&  BII$_2$&  BIII\cr 
\+$n=89$ $d=309$ $\eta=0$ $h=2$&&&&    
&     &     &     &     & BI& BII$_1$&  BII$_2$&  BIII\cr 
\+$n=90$ $d=314$ $\eta=1$ $h=2$&&&&    
&     &     &     &     & BI&&&\cr
\+$n=91$ $d=319$ $\eta=3$ $h=2$&&&&    
&     &     &     &     & BI& BII$_1$&  BII$_2$&  BIII\cr 
\+$n=92$ $d=323$ $\eta=0$ $h=2$&&&&    
&  AI0&     & AII0& AIII& BI& BII$_1$&  BII$_2$&  BIII\cr 
\+$n=93$ $d=329$ $\eta=3$ $h=2$&&&&    
&     &     & AII0&     &   & BII$_1$&  BII$_2$&  BIII\cr 
\+$n=94$ $d=330$ $\eta=6$ $h=2$&&&&    
&  AI0&     & AII0&     & BI& BII$_1$&  BII$_2$&  BIII\cr 
\+$n=95$ $d=345$ $\eta=2$ $h=2$&&&&    
&     &     & AII0&     &   & BII$_1$&  BII$_2$&  BIII\cr 
\+$n=96$ $d=345$ $\eta=4$ $h=2$&&&&    
&     &     & AII0&     & BI& BII$_1$&  BII$_2$&  BIII\cr 
\+$n=97$ $d=345$ $\eta=5$ $h=2$&&&&    
&     &     & AII0&     &   & BII$_1$&  BII$_2$&  BIII\cr 
\+$n=98$ $d=357$ $\eta=6$ $h=2$&&&&    
&     &     & AII0&     &   & BII$_1$&         &  BIII\cr 
\+$n=99$ $d=366$ $\eta=2$ $h=2$&&&&    
&     &     &     &     & BI& BII$_1$&  BII$_2$&  BIII\cr 
\+$n=100$ $d=374$ $\eta=0$ $h=2$&&&&   
&     &     &     &     & BI&        &  BII$_2$&\cr 
\+$n=101$ $d=390$ $\eta=0$ $h=2$&&&&   
&  AI0&     &     &     & BI& BII$_1$&  BII$_2$&  BIII\cr 
\+$n=102$ $d=402$ $\eta=2$ $h=2$&&&&   
&     &     &     &     & BI&&&\cr
\+$n=103$ $d=406$ $\eta=3$ $h=2$&&&&   
&     &     & AII0&     & BI& BII$_1$&         &  BIII\cr 
\+$n=104$ $d=410$ $\eta=2$ $h=2$&&&&   
&     &     &     &     & BI&&&\cr
\+$n=105$ $d=426$ $\eta=0$ $h=2$&&&&   
&     &     &     &     & BI& BII$_1$&  BII$_2$&\cr 
\+$n=106$ $d=429$ $\eta=0$ $h=2$&&&&   
&     &     &     &     &   & BII$_1$&  BII$_2$&\cr 
\+$n=107$ $d=429$ $\eta=5$ $h=2$&&&&   
&     &     & AII0& AIII&   & BII$_1$&  BII$_2$&  BIII\cr 
\+$n=108$ $d=429$ $\eta=6$ $h=2$&&&&   
&     & AII1& AII0&     & BI& BII$_1$&  BII$_2$&\cr 
\+$n=109$ $d=431$ $\eta=0$ $h=2$&&&&   
&     &     &     &     & BI&&&\cr
\+$n=110$ $d=434$ $\eta=1$ $h=2$&&&&   
&     &     &     &     &   &        &  BII$_2$&\cr 
\+$n=111$ $d=438$ $\eta=0$ $h=2$&&&&   
&     &     &     &     & BI& BII$_1$&  BII$_2$&  BIII\cr 
\+$n=112$ $d=455$ $\eta=3$ $h=2$&&&&   
&     &     & AII0&     & BI& BII$_1$&         &  BIII\cr 
\+$n=113$ $d=455$ $\eta=6$ $h=2$&&&&   
&     &     &     &     &   &        &  BII$_2$&\cr 
\+$n=114$ $d=462$ $\eta=0$ $h=2$&&&&   
&     & AII1& AII0&     & BI& BII$_1$&  BII$_2$&  BIII\cr 
\+$n=115$ $d=462$ $\eta=6$ $h=2$&&&&   
&     &     & AII0&     &   & BII$_1$&  BII$_2$&  BIII\cr 
\+$n=116$ $d=465$ $\eta=1$ $h=2$&&&&   
&     &     & AII0&     &   & BII$_1$&  BII$_2$&  BIII\cr 
\+$n=117$ $d=465$ $\eta=6$ $h=2$&&&&   
&     &     &     &     &   & BII$_1$&  BII$_2$&  BIII\cr 
\+$n=118$ $d=469$ $\eta=1$ $h=2$&&&&   
&     &     &     &     &   & BII$_1$&         &  BIII\cr 
\+$n=119$ $d=469$ $\eta=3$ $h=2$&&&&   
&     &     & AII0& AIII&   & BII$_1$&         &  BIII\cr 
\+$n=120$ $d=470$ $\eta=0$ $h=2$&&&&   
&     &     &     &     & BI&&&\cr
\+$n=121$ $d=471$ $\eta=0$ $h=2$&&&&   
&     &     &     &     & BI&        &         &  BIII\cr 
\+$n=122$ $d=474$ $\eta=0$ $h=2$&&&&   
&     &     &     &     & BI&&&\cr
\+$n=123$ $d=483$ $\eta=2$ $h=2$&&&&   
&     &     & AII0& AIII&   & BII$_1$&  BII$_2$&  BIII\cr 
\+$n=124$ $d=515$ $\eta=0$ $h=2$&&&&   
&     & AII1&     &     & BI&        &         &  BIII\cr 
\+$n=125$ $d=518$ $\eta=3$ $h=2$&&&&   
&     &     &     &     & BI& BII$_1$&  BII$_2$&\cr 
\+$n=126$ $d=527$ $\eta=3$ $h=2$&&&&   
&     &     &     &     & BI&&&\cr
\+$n=127$ $d=554$ $\eta=1$ $h=2$&&&&   
&     &     &     &     & BI&&&\cr
\+$n=128$ $d=555$ $\eta=3$ $h=2$&&&&   
&     &     & AII0& AIII&   & BII$_1$&  BII$_2$&  BIII\cr 
\+$n=129$ $d=555$ $\eta=5$ $h=2$&&&&   
&     &     &     &     & BI& BII$_1$&  BII$_2$&  BIII\cr 
\+$n=130$ $d=555$ $\eta=6$ $h=2$&&&&   
&     &     & AII0&     & BI& BII$_1$&  BII$_2$&  BIII\cr 
\+$n=131$ $d=561$ $\eta=6$ $h=2$&&&&   
&  AI0& AII1& AII0&     & BI& BII$_1$&  BII$_2$&  BIII\cr 
\+$n=132$ $d=582$ $\eta=0$ $h=2$&&&&   
&     &     &     &     & BI& BII$_1$&  BII$_2$&  BIII\cr 
\+$n=133$ $d=595$ $\eta=3$ $h=2$&&&&   
&     &     & AII0&     &   & BII$_1$&  BII$_2$&  BIII\cr 
\+$n=134$ $d=609$ $\eta=3$ $h=2$&&&&   
&     &     & AII0& AIII&   & BII$_1$&  BII$_2$&  BIII\cr 
\+$n=135$ $d=609$ $\eta=6$ $h=2$&&&&   
&     &     & AII0&     &   & BII$_1$&  BII$_2$&  BIII\cr 
\+$n=136$ $d=618$ $\eta=0$ $h=2$&&&&   
&     &     &     &     & BI&&&\cr
\+$n=137$ $d=623$ $\eta=3$ $h=2$&&&&   
&     &     &     &     & BI&&&\cr
\+$n=138$ $d=627$ $\eta=1$ $h=2$&&&&   
&     &     &     &     &   & BII$_1$&  BII$_2$&  BIII\cr 
\+$n=139$ $d=627$ $\eta=3$ $h=2$&&&&   
&     &     &     &     &   & BII$_1$&  BII$_2$&  BIII\cr 
\+$n=140$ $d=663$ $\eta=3$ $h=2$&&&&   
&     &     & AII0&     &   & BII$_1$&  BII$_2$&  BIII\cr 
\+$n=141$ $d=665$ $\eta=0$ $h=2$&&&&   
&     &     &     &     &   & BII$_1$&  BII$_2$&  BIII\cr 
\+$n=142$ $d=690$ $\eta=1$ $h=2$&&&&   
&     &     &     &     &   & BII$_1$&  BII$_2$&  BIII\cr 
\+$n=143$ $d=690$ $\eta=2$ $h=2$&&&&   
&     &     &     &     & BI& BII$_1$&  BII$_2$&  BIII\cr 
\+$n=144$ $d=690$ $\eta=7$ $h=2$&&&&   
&     &     & AII0&     &   & BII$_1$&  BII$_2$&  BIII\cr 
\+$n=145$ $d=705$ $\eta=6$ $h=2$&&&&   
&     & AII1& AII0&     & BI& BII$_1$&  BII$_2$&  BIII\cr 
\+$n=146$ $d=715$ $\eta=0$ $h=2$&&&&   
&     &     & AII0&     &   & BII$_1$&  BII$_2$&  BIII\cr 
\+$n=147$ $d=715$ $\eta=6$ $h=2$&&&&   
&     &     & AII0&     &   & BII$_1$&  BII$_2$&  BIII\cr 
\+$n=148$ $d=741$ $\eta=0$ $h=2$&&&&   
&     &     & AII0&     &   & BII$_1$&  BII$_2$&  BIII\cr 
\+$n=149$ $d=741$ $\eta=6$ $h=2$&&&&   
&  AI0&     &     &     & BI& BII$_1$&  BII$_2$&\cr 
\+$n=150$ $d=770$ $\eta=2$ $h=2$&&&&   
&     &     &     &     &   & BII$_1$&  BII$_2$&\cr 
\+$n=151$ $d=770$ $\eta=7$ $h=2$&&&&   
&     &     & AII0&     &   & BII$_1$&  BII$_2$&  BIII\cr 
\+$n=152$ $d=791$ $\eta=3$ $h=2$&&&&   
&     &     &     &     & BI&&&\cr
\+$n=153$ $d=794$ $\eta=1$ $h=2$&&&&   
&&&&&&&&\cr
\+$n=154$ $d=798$ $\eta=0$ $h=2$&&&&   
&     &     & AII0&     &   & BII$_1$&  BII$_2$&  BIII\cr 
\+$n=155$ $d=798$ $\eta=6$ $h=2$&&&&   
&     &     &     &     &   & BII$_1$&  BII$_2$&\cr 
\+$n=156$ $d=806$ $\eta=0$ $h=2$&&&&   
&&&&&&&&\cr
\+$n=157$ $d=834$ $\eta=2$ $h=2$&&&&   
&     &     &     &     & BI&&&\cr
\+$n=158$ $d=858$ $\eta=0$ $h=2$&&&&   
&     &     &     &     & BI& BII$_1$&&\cr 
\+$n=159$ $d=890$ $\eta=2$ $h=2$&&&&   
&     &     &     &     & BI&&&\cr
\+$n=160$ $d=897$ $\eta=6$ $h=2$&&&&   
&     &     & AII0&     &   & BII$_1$&  BII$_2$&  BIII\cr 
\+$n=161$ $d=903$ $\eta=2$ $h=2$&&&&   
&  AI0& AII1& AII0&     & BI& BII$_1$&  BII$_2$&  BIII\cr 
\+$n=162$ $d=906$ $\eta=0$ $h=2$&&&&   
&     &     &     &     &   &        &  BII$_2$&\cr 
\+$n=163$ $d=915$ $\eta=5$ $h=2$&&&&   
&     &     & AII0&     &   & BII$_1$&  BII$_2$&  BIII\cr 
\+$n=164$ $d=935$ $\eta=0$ $h=2$&&&&   
&     &     &     &     &   & BII$_1$&  BII$_2$&  BIII\cr 
\+$n=165$ $d=959$ $\eta=3$ $h=2$&&&&   
&&&&&&&&\cr
\+$n=166$ $d=966$ $\eta=4$ $h=2$&&&&   
&     &     &     &     &   & BII$_1$&  BII$_2$&  BIII\cr 
\+$n=167$ $d=986$ $\eta=1$ $h=2$&&&&   
&&&&&&&&\cr
\+$n=168$ $d=987$ $\eta=7$ $h=2$&&&&   
&     &     & AII0&     &   & BII$_1$&         &  BIII\cr 
\+$n=169$ $d=1001$ $\eta=5$ $h=2$&&&&  
&     &     & AII0&     &   & BII$_1$&  BII$_2$&  BIII\cr 
\+$n=170$ $d=1010$ $\eta=0$ $h=2$&&&&  
&     &     &     &     & BI&&&\cr
\+$n=171$ $d=1015$ $\eta=6$ $h=2$&&&&  
&     &     &     &     &   & BII$_1$&  BII$_2$&  BIII\cr 
\+$n=172$ $d=1023$ $\eta=1$ $h=2$&&&&  
&  AI0&     & AII0&     & BI& BII$_1$&  BII$_2$&  BIII\cr 
\+$n=173$ $d=1066$ $\eta=2$ $h=2$&&&&  
&&&&&&&&\cr
\+$n=174$ $d=1085$ $\eta=6$ $h=2$&&&&  
&     &     & AII0& AIII&   & BII$_1$&  BII$_2$&  BIII\cr 
\+$n=175$ $d=1095$ $\eta=3$ $h=2$&&&&  
&     &     & AII0&     &   & BII$_1$&  BII$_2$&  BIII\cr 
\+$n=176$ $d=1110$ $\eta=0$ $h=2$&&&&  
&     &     &     &     & BI& BII$_1$&  BII$_2$&  BIII\cr 
\+$n=177$ $d=1122$ $\eta=1$ $h=2$&&&&  
&     &     &     &     &   & BII$_1$&  BII$_2$&  BIII\cr 
\+$n=178$ $d=1130$ $\eta=2$ $h=2$&&&&  
&&&&&&&&\cr
\+$n=179$ $d=1155$ $\eta=10$ $h=2$&&&& 
&     &     & AII0&     &   & BII$_1$&  BII$_2$&  BIII\cr 
\+$n=180$ $d=1173$ $\eta=6$ $h=2$&&&&  
&     &     &     &     & BI& BII$_1$&  BII$_2$&  BIII\cr 
\+$n=181$ $d=1194$ $\eta=0$ $h=2$&&&&  
&&&&&&&&\cr
\+$n=182$ $d=1218$ $\eta=2$ $h=2$&&&&  
&     &     &     &     &   & BII$_1$&  BII$_2$&\cr 
\+$n=183$ $d=1235$ $\eta=6$ $h=2$&&&&  
&     &     & AII0&     &   & BII$_1$&  BII$_2$&  BIII\cr 
\+$n=184$ $d=1245$ $\eta=6$ $h=2$&&&&  
&     &     & AII0& AIII& BI& BII$_1$&  BII$_2$&  BIII\cr 
\+$n=185$ $d=1254$ $\eta=2$ $h=2$&&&&  
&     &     &     &     &   & BII$_1$&         &  BIII\cr 
\+$n=186$ $d=1271$ $\eta=3$ $h=2$&&&&  
&&&&&&&&\cr
\+$n=187$ $d=1295$ $\eta=3$ $h=2$&&&&  
&     &     & AII0&     &   & BII$_1$&         &  BIII\cr 
\+$n=188$ $d=1302$ $\eta=2$ $h=2$&&&&  
&     &     &     &     & BI& BII$_1$&&\cr
\+$n=189$ $d=1302$ $\eta=4$ $h=2$&&&&  
&     &     &     &     &   & BII$_1$&&           BIII\cr 
\+$n=190$ $d=1338$ $\eta=0$ $h=2$&&&&  
&&&&&&&&\cr
\+$n=191$ $d=1365$ $\eta=9$ $h=2$&&&&  
&     &     & AII0& AIII&   & BII$_1$&  BII$_2$&  BIII\cr 
\+$n=192$ $d=1365$ $\eta=10$ $h=2$&&&& 
&     &     & AII0&     &   & BII$_1$&  BII$_2$&  BIII\cr 
\+$n=193$ $d=1365$ $\eta=12$ $h=2$&&&& 
&     &     & AII0&     &   & BII$_1$&  BII$_2$&  BIII\cr 
\+$n=194$ $d=1365$ $\eta=13$ $h=2$&&&& 
&     &     & AII0&     &   & BII$_1$&  BII$_2$&  BIII\cr 
\+$n=195$ $d=1370$ $\eta=2$ $h=2$&&&&  
&&&&&&&&\cr
\+$n=196$ $d=1446$ $\eta=0$ $h=2$&&&&  
&     &     &     &     & BI&&&\cr
\+$n=197$ $d=1463$ $\eta=1$ $h=2$&&&&  
&     &     &     &     &   &        &  BII$_2$&\cr 
\+$n=198$ $d=1479$ $\eta=6$ $h=2$&&&&  
&     &     &     &     & BI&        &         &  BIII\cr 
\+$n=199$ $d=1482$ $\eta=0$ $h=2$&&&&  
&&&&&&&&\cr
\+$n=200$ $d=1482$ $\eta=6$ $h=2$&&&&  
&     &     & AII0&     &   &        &  BII$_2$&  BIII\cr 
\+$n=201$ $d=1495$ $\eta=5$ $h=2$&&&&  
&     &     &     &     &   &        &  BII$_2$&  BIII\cr 
\+$n=202$ $d=1526$ $\eta=3$ $h=2$&&&&  
&&&&&&&&\cr
\+$n=203$ $d=1551$ $\eta=4$ $h=2$&&&&  
&     &     &     &     & BI& BII$_1$&  BII$_2$&  BIII\cr 
\+$n=204$ $d=1554$ $\eta=2$ $h=2$&&&&  
&&&&&&&&\cr
\+$n=205$ $d=1554$ $\eta=4$ $h=2$&&&&  
&  AI0&     & AII0&     &   & BII$_1$&  BII$_2$&  BIII\cr 
\+$n=206$ $d=1586$ $\eta=0$ $h=2$&&&&  
&&&&&&&&\cr
\+$n=207$ $d=1586$ $\eta=3$ $h=2$&&&&  
&     &     &     &     & BI&        &         &  BIII\cr 
\+$n=208$ $d=1605$ $\eta=6$ $h=2$&&&&  
&     &     &     & AIII& BI& BII$_1$&  BII$_2$&  BIII\cr 
\+$n=209$ $d=1626$ $\eta=0$ $h=2$&&&&  
&&&&&&&&\cr
\+$n=210$ $d=1635$ $\eta=5$ $h=2$&&&&  
&     &     &     &     &   & BII$_1$&  BII$_2$&  BIII\cr 
\+$n=211$ $d=1677$ $\eta=6$ $h=2$&&&&  
&     &     & AII0& AIII&   & BII$_1$&  BII$_2$&  BIII\cr 
\+$n=212$ $d=1751$ $\eta=3$ $h=2$&&&&  
&     &     &     &     & BI&&&\cr
\+$n=213$ $d=1771$ $\eta=7$ $h=2$&&&&  
&     &     & AII0&     &   & BII$_1$&&\cr 
\+$n=214$ $d=1785$ $\eta=10$ $h=2$&&&& 
&     &     & AII0&     &   & BII$_1$&  BII$_2$&  BIII\cr 
\+$n=215$ $d=1794$ $\eta=4$ $h=2$&&&&  
&&&&&&&&\cr
\+$n=216$ $d=1806$ $\eta=6$ $h=2$&&&&  
&     &     &     &     &   & BII$_1$&         &  BIII\cr 
\+$n=217$ $d=1898$ $\eta=2$ $h=2$&&&&  
&&&&&&&&\cr
\+$n=218$ $d=1986$ $\eta=2$ $h=2$&&&&  
&&&&&&&&\cr
\+$n=219$ $d=1995$ $\eta=7$ $h=2$&&&&  
&     &     & AII0& AIII&   & BII$_1$&  BII$_2$&  BIII\cr 
\+$n=220$ $d=1995$ $\eta=14$ $h=2$&&&& 
&     &     & AII0&     &   & BII$_1$&  BII$_2$&  BIII\cr 
\+$n=221$ $d=2015$ $\eta=6$ $h=2$&&&&  
&     &     &     &     &   & BII$_1$&  BII$_2$&  BIII\cr 
\+$n=222$ $d=2046$ $\eta=0$ $h=2$&&&&  
&     &     & AII0&     &   & BII$_1$&         &  BIII\cr 
\+$n=223$ $d=2090$ $\eta=0$ $h=2$&&&&  
&&&&&&&&\cr
\+$n=224$ $d=2145$ $\eta=14$ $h=2$&&&& 
&     &     & AII0&     &   & BII$_1$&  BII$_2$&  BIII\cr 
\+$n=225$ $d=2226$ $\eta=2$ $h=2$&&&&  
&&&&&&&&\cr
\+$n=226$ $d=2415$ $\eta=1$ $h=2$&&&&  
&     &     & AII0&     &   & BII$_1$&  BII$_2$&  BIII\cr 
\+$n=227$ $d=2415$ $\eta=4$ $h=2$&&&&  
&  AI0& AII1& AII0& AIII&   & BII$_1$&  BII$_2$&  BIII\cr 
\+$n=228$ $d=2415$ $\eta=8$ $h=2$&&&&  
&     &     & AII0&     &   &BII$_1$&   BII$_2$&  BIII\cr 
\+$n=229$ $d=2454$ $\eta=0$ $h=2$&&&&  
&&&&&&&&\cr
\+$n=230$ $d=2562$ $\eta=2$ $h=2$&&&&  
&&&&&&&&\cr
\+$n=231$ $d=2570$ $\eta=2$ $h=2$&&&&  
&&&&&&&&\cr
\+$n=232$ $d=2639$ $\eta=5$ $h=2$&&&&  
&      &    &     &     &   & BII$_1$&  BII$_2$&\cr
\+$n=233$ $d=2730$ $\eta=0$ $h=2$&&&&  
&      &    & AII0&     &   & BII$_1$&         &  BIII\cr  
\+$n=234$ $d=2730$ $\eta=10$ $h=2$&&&& 
&      &    & AII0&     & BI& BII$_1$&  BII$_2$&  BIII\cr 
\+$n=235$ $d=2805$ $\eta=3$ $h=2$&&&&  
&      &    & AII0& AIII&   & BII$_1$&  BII$_2$&  BIII\cr 
\+$n=236$ $d=2829$ $\eta=6$ $h=2$&&&&  
&      &    &     &     &   & BII$_1$&         &  BIII\cr 
\+$n=237$ $d=2886$ $\eta=0$ $h=2$&&&&  
&      &    &     &     & BI&&&\cr
\+$n=238$ $d=3003$ $\eta=11$ $h=2$&&&& 
&      &    & AII0&     &   & BII$_1$&  BII$_2$&  BIII\cr 
\+$n=239$ $d=3045$ $\eta=10$ $h=2$&&&& 
&      &    & AII0&     &   & BII$_1$&  BII$_2$&  BIII\cr 
\+$n=240$ $d=3066$ $\eta=6$ $h=2$&&&&  
&&&&&&&&\cr
\+$n=241$ $d=3135$ $\eta=2$ $h=2$&&&&  
&      &    & AII0&     &   & BII$_1$&  BII$_2$&  BIII\cr 
\+$n=242$ $d=3135$ $\eta=8$ $h=2$&&&&  
&      &    & AII0& AIII&   & BII$_1$&  BII$_2$&  BIII\cr 
\+$n=243$ $d=3315$ $\eta=0$ $h=2$&&&&  
&      &    & AII0&     &   & BII$_1$&  BII$_2$&  BIII\cr 
\+$n=244$ $d=3318$ $\eta=2$ $h=2$&&&&  
&      &    &     &     &   &        &         &  BIII\cr 
\+$n=245$ $d=3354$ $\eta=0$ $h=2$&&&&  
&&&&&&&&\cr
\+$n=246$ $d=3354$ $\eta=6$ $h=2$&&&&  
&&&&&&&&\cr
\+$n=247$ $d=3410$ $\eta=4$ $h=2$&&&&  
&&&&&&&&\cr
\+$n=248$ $d=4026$ $\eta=0$ $h=2$&&&&  
&&&&&&&&\cr
\+$n=249$ $d=4074$ $\eta=6$ $h=2$&&&&  
&&&&&&&&\cr
\+$n=250$ $d=4290$ $\eta=1$ $h=2$&&&&  
&      &    &     &     &   & BII$_1$&  BII$_2$&  BIII\cr 
\+$n=251$ $d=4326$ $\eta=2$ $h=2$&&&&  
&&&&&&&&\cr
\+$n=252$ $d=4902$ $\eta=4$ $h=2$&&&&  
&&&&&&&&\cr
\+$n=253$ $d=4991$ $\eta=7$ $h=2$&&&&  
&&&&&&&&\cr
\+$n=254$ $d=5226$ $\eta=0$ $h=2$&&&&  
&      &    &     &     & BI&&&\cr
\+$n=255$ $d=5334$ $\eta=2$ $h=2$&&&&  
&&&&&&&&\cr
\+$n=256$ $d=6006$ $\eta=2$ $h=2$&&&&  
&      &    &     &     &   & BII$_1$&  BII$_2$&\cr
\+$n=257$ $d=7590$ $\eta=8$ $h=2$&&&&  
&      &    &     &     &   & BII$_1$&         &  BIII\cr 
\+$n=258$ $d=10374$ $\eta=2$ $h=2$&&&& 
&&&&&&&&\cr
\+$n=259$ $d=29526$ $\eta=2$ $h=2$&&&& 
&&&&&&&&\cr}

\vskip30pt 

By Theorem 9.1.1, to prove Theorem 7.1, we now should check reflective 
type of all lattices of Table 3 with $h=0$ and Table 6. 
To find reflective type of these lattices, we should write down these 
lattices exactly. Suppose that the lattice $S$ has invariants 
$(d,\,\eta)$. Any lattice $S$ of Table 3   
has one of forms \thetag{9.1.1} or \thetag{9.1.3} below 
(see \thetag{5.1.18} --- \thetag{5.1.21}): 
$$
S=U\oplus \langle -d \rangle,
\tag{9.1.1}
$$
where 
$$
(-1)^{\eta_p}=\left({-d/p\over p}\right),\
\text{for any odd\ } p\vert d,
\tag{9.1.2}
$$
or 
$$
S=\langle n_1\rangle \oplus \langle -n_2 \rangle\oplus
\langle -n_3 \rangle (\epsilon_1/2,\epsilon_2/2,\epsilon_3/2)
\tag{9.1.3}
$$
where $\epsilon_1,\epsilon_2,\epsilon_3 \in \{0,\,1\}$,
we have $d=n_1n_2n_3$ if
$\epsilon_1=\epsilon_2=\epsilon_3=0$, and
$d=n_1n_2n_3/4$ otherwise. Moreover, $n_i\equiv 0\mod 2$ if
$\epsilon_i=1$;
 $(n_1\epsilon_1-n_2\epsilon_2-n_3\epsilon_3)/4\in \bz$;
$(n_1\epsilon_1-n_2\epsilon_2-n_3\epsilon_3)/4\in 2\bz$ if
$d\equiv 0\mod 2$; for any odd $p\vert d$ we have
$$
(-1)^{\eta_p}=
\cases
\left({n_1/p\over p}\right) &\text{if $p\vert n_1$,}\\
\left({-n_2/p\over p}\right) &\text{if $p\vert n_2$,}\\
\left({-n_3/p\over p}\right) &\text{if $p\vert n_3$.}
\endcases
\tag{9.1.4}
$$
Not all lattices of Table 6 have these forms. Some of them 
have the form  
$$
\langle 3d \rangle \oplus A_2(1/3,1/3,-1/3),\ \ \ 
A_2=\pmatrix-2&1\cr1&-2\endpmatrix, 
\tag{9.1.5}
$$
where $d\equiv 2\mod 3$. Then 
$$
(-1)^{\eta_p}=\left({3d\over p}\right)
\tag{9.1.6}
$$
for any odd prime $p\vert d$. 

We have 
\proclaim{Lemma 9.1.3} Any reflective main hyperbolic lattice $S$ of 
rank $3$ and with square-free determinant $d$ has one of forms 
\thetag{9.1.1}, \thetag{9.1.3} or \thetag{9.1.5}.
\endproclaim

\demo{Proof} Let us consider a vertex of a fundamental polygon $\M$ of 
$W(S)$. If the lattice $S$ is reflective, this vertex does exist. 
Let $r_1$ and $r_2$ are primitive roots of $S$ orthogonal to sides of 
$\M$ in this vertex. We have 
$$
0\le \alpha_{12}={4(r_1,\,r_2)^2\over r_1^2r_2^2}\le 4
$$
where $\alpha_{12}=4(\cos{\pi/k})^2$, $k=2,\,3,\,4,\,6$ or $+\infty$. 
We apply Lemma 9.1.2 (see also its proof). If $k=+\infty$, 
then $S$ has the form \thetag{9.1.1}. If $k=3$, then $S$ has the form 
\thetag{9.1.5}. If $k=2,\,4,\,6$, then $S$ has the form \thetag{9.1.3}.

This finishes the proof.    
\enddemo 

By Lemma 9.1.3, if a lattice $S$ of Table 6 does not have one of forms 
\thetag{9.1.1}, \thetag{9.1.3} or \thetag{9.1.5}, then $S$ is not 
reflective. This is the case for $n=204$, $206$, $215$, $245$, $247$ of 
Table 6.

All other lattices of Tables 3 and 6 have one of forms \thetag{9.1.1}, 
\thetag{9.1.3} or \thetag{9.1.5}. We apply to them Vinberg's algorithm 
(see Sect. 5). For lattices $S$ of forms \thetag{9.1.1} and 
\thetag{9.1.3}, we have described this algorithm in Sect. 5. 
For lattices $S$ of form \thetag{9.1.5}, we take $c=(1,0,0)$ as the 
center of Vinberg's algorithm, and roots $v_1=(0,1,0)$ and $v_2=(0,0,1)$ 
as roots of the height $0$. In Appendix 2, Program 31: refl0.14 we 
give Vinberg's algorithm for lattices of the form \thetag{9.1.5}. 

In Sect. 5 we gave several examples of calculations using Vinberg's 
algorithm for lattices of Table 3. Similar calculations we should do for all lattices of Table 6.  As the result of these calculations, 
we get that only 66 lattices  of Tables 3 and 6 give hyperbolically 
reflective lattices. 5 of them belong to Table 3 and 
have $h=0$, the rest 61 belong to Table 6 and have $h=2$.   
All these 66 lattices are given in Table 4.  

This finishes the proof of Basic Theorem 7.1.

We remark that for some lattices of Table 6 we can prove that they 
are not reflective by considering Table 7. Lattices of Table 7 which 
don't have any type AI1 --- BIII of narrow part are not reflective: 
their fundamental polygon $\M$ does not have a narrow part satisfying 
Theorem 8.3.4. This is true for 36 triplets: 
$$
\split
&(d,\eta,h)=\\
&(301,1,2),\ (794,1,2),\ (806,0,2),\, (959,3,2),\ (986,1,2),(1066,2,2),\\
&(1130,2,2),\ (1194,0,2),\ (1271,3,2),\ (1338,0,2),\ (1370,2,2),\ 
(1482,0,2),\\
&(1526,3,2),\ (1554,2,2),\ (1586,0,2),\ (1626,0,2),\ (1794,4,2),
\ (1898,2,2),\\ &(1986,2,2),\ (2090,0,2),\ (2226,2,2),\ (2454,0,2),
\ (2562,2,2),\ (2570,2,2),\\  &(3066,6,2),\ (3354,0,2),\ (3354,6,2),
\ (3410,4,2),\ (4026,0,2),\ (4074,6,2),\\ &(4326,2,2),\ (4902,4,2),
\ (4991,7,2),\ (5334,2,2),\ (10374,2,2),\ (29526,2,2).
\endsplit 
\tag{9.1.7} 
$$          
Thus, all triplets \thetag{9.1.7} give non-reflective lattices. 
\subhead
9.2. Proof of Theorem 7.2
\endsubhead 

We use notations and results of Sect. 2.2. Especially see Proposition 
2.2.6. Assume that 
$\widetilde{S}$ is hyperbolically reflective non-main hyperbolic lattice 
of determinant $2d$ where $S$ is main hyperbolic lattice of odd determinant 
$d$ where $d$ is square-free. Since $O(\widetilde{S})\subset O(S)$, 
reflective type of $\widetilde{S}$ is dominated by reflective type of $S$. 
In particular, $S$ is reflective. 

First, assume that $S$ is elliptically or 
parabolically reflective. Then $S$ belongs to Table 1 of 
Basic Theorem 2.3.2.1. Then we argue like in Sect. 5.2. The 
lattice $\widetilde{S}$ also has one of forms \thetag{9.1.1} or 
\thetag{9.1.3}, and 
we can calculate its reflective type using  Vinberg's algorithm. 
These  calculations are of the same difficulty as for main lattices. 
If by these calculations a lattice $\widetilde{S}$ is elliptically or 
parabolically reflective, it belongs to Table 2 of Theorem 2.3.3.1. 
If $\widetilde{S}$ is hyperbolically reflective, we give it in  
Table 5 of Theorem 7.2. There are 10 these cases corresponding to lattices 
$S$ with $h=0$ or $h=1$ which are not equivariantly equivalent to   
$\widetilde{S}$.   

Now suppose that the lattice $S$ is hyperbolically reflective. 
Thus, $S$ belongs to the Table 4 of Theorem 7.1. Let us show that 
$\widetilde{S}$ is also hyperbolically reflective if and only if 
$S$ and $\widetilde{S}$ are equivariantly equivalent, i. e. 
$$
\sum_{p\mid d}
{(1-p+4\eta_p+4\omega(d)_p)}\equiv \pm 1-1\equiv 0\ \text{or}\  6 \mod 8,
\tag{9.2.1}
$$
see Proposition 2.2.6. Suppose that \thetag{9.2.1} is valid. Then 
$O(\widetilde S)=O(S)$ (see Proposition 2.2.6), and the lattice 
$\widetilde{S}$ is hyperbolically reflective if $S$ does. All 
these cases are given in Table 5 as equivariantly equivalent cases.  

Suppose that 
\thetag{9.2.1} is not valid but $\widetilde{S}$ is hyperbolically reflective. 
Let $\widetilde{\M}$ be a fundamental polygon of $W(\widetilde{S})$. 
Since $W(\widetilde{S})\subset W(S)$, we can choose a fundamental polygon 
$\M$ for $W(S)$ such that $\M\subset \widetilde{\M}$. 
Suppose that $\M\not=\widetilde{\M}$ (equivalently 
$W(\widetilde{S})\not=W(S)$). 
Then there exists a reflection $s_r\in W(S)$ with respect to a 
side of $\M$ such that $s_r\notin W(\widetilde{S})$. 
Then $s_r(\M)$ is also contained in $\widetilde{\M}$. Let $l$ be 
the axis of $\M$. For all 66 cases of Theorem 7.1 the axis $l$ is 
contained inside 
$\M$. Then lines $l$ 
and $s_r(l)$ are contained in $\widetilde{\M}$ and are different. 
This is impossible for a polygon $\widetilde{\M}$ of 
restricted hyperbolic type. Only its axis has this property, 
and it is unique. Thus, we have proved that $W(\widetilde{S})=W(S)$. 
We remind that $S(2)\subset \widetilde{S}$ is an overlattice of $S(2)$ of 
index 2, and $S(2)$ is the maximal even sublattice of $\widetilde{S}$. 
If \thetag{9.2.1} is not valid, then the overlattice $\widetilde{S}$ is not 
unique, there are 3 of them, but all of them are conjugate by 
$O(S)=O(S(2))$. Since $W(\widetilde{S})=W(S)$, then all these three 
overlattices are invariant with respect to $W(S)$. Simple 
considerations over $\bz_2$ show that this is true for a 
reflection $s_r$ in a primitive root $r\in S$ if and only 
if $r^2\equiv 1\mod 2$. 
Thus, the equality $W(\widetilde{S})=W(S)$ implies that all roots $r\in S$ 
have odd square $r^2$. In Table 4 we are describing all primitive 
roots of $S$. There always exists a primitive root $r\in S$ with even $r^2$. 
We get a contradiction. 

This finishes the proof of Theorem 7.2.    

\subhead
9.3. Proof of Theorem 7.4 
\endsubhead 

The proof is the same as in Sect. 9.2 above where we proved 
that $W(\widetilde{S})=W(S)$ if $\widetilde{S}$ is hyperbolically 
reflective.

\head
10. Appendix 2: Programs for GP/PARI calculator
\endhead

\centerline{{\bf Program 13:} funda11.gen} 

\noindent
$\backslash$$\backslash$funda11.gen 
\newline
$\backslash$$\backslash$hyperbolic type, case AI1.
\newline
$\backslash$$\backslash$cosh(eta/2)=3/2
\newline
$\backslash$l;$\backslash$
\newline
epsilon=1;epsilon1=1;epsilon2=1;n=0;$\backslash$
\newline
rh1=[12.25000000,13.37793021,14.23012096,14.94097150;$\backslash$
\newline
 13.37793021,14.57106781,15.47225159,16.22381115;$\backslash$
\newline
 14.23012096,15.47225159,16.41025403,17.19241152;$\backslash$
\newline
 14.94097150,16.22381115,17.19241152,18.00000000];$\backslash$
\newline
for(alpha12=1,4,for(alpha23=alpha12,4,u=(rh1[alpha12,alpha23]+0.1)$
^\wedge$2;$\backslash$
\newline
for(alpha13=5,u,$\backslash$
\newline
if(alpha23==0$||$issquare(alpha12$\ast$alpha23$\ast$alpha13)!=1$||$$
\backslash$
\newline
--8+2$\ast$isqrt(alpha12$\ast$alpha23$\ast$alpha13)+2$\ast$alpha12+2$
\ast$alpha13+2$\ast$alpha23$<$=0,,$\backslash$
\newline
alpha=4$\ast$idmat(3);alpha[1,2]=alpha12;alpha[2,1]=alpha12;$\backslash$
\newline
alpha[2,3]=alpha23;alpha[3,2]=alpha23;alpha[1,3]=alpha13;$\backslash$
\newline
alpha[3,1]=alpha13;dalpha=$\backslash$
\newline
--8+2$\ast$isqrt(alpha12$\ast$alpha23$\ast$alpha13)+2$\ast$alpha12+2$
\ast$alpha13+2$\ast$alpha23;$\backslash$
\newline
a=--2$\ast$idmat(3);$\backslash$
\newline
fordiv(alpha[1,2],a12,fordiv(alpha[2,3],a23,fordiv(alpha[1,3],a13,
$\backslash$
\newline
a21=alpha[1,2]/a12;a32=alpha[2,3]/a23;a31=alpha[1,3]/a13;$\backslash$
\newline
if(a12$\ast$a23$\ast$a31!=a21$\ast$a13$\ast$a32,,$\backslash$
\newline
a[1,2]=a12;a[2,1]=a21;a[2,3]=a23;a[3,2]=a32;a[1,3]=a13;a[3,1]=a31;$
\backslash$
\newline
dd=idmat(3);dd[1,1]=a13$\ast$a32;dd[2,2]=a23$\ast$a31;dd[3,3]=a31$
\ast$a32;$\backslash$
\newline
b=a$\ast$dd;b=b/content(b);n=n+1;$\backslash$
\newline
db=smith(b);r=db[1];$\backslash$
\newline
if(r$>$epsilon,epsilon=r,);$\backslash$
\newline
fr=factor(r);tfr=matsize(fr)[1];$\backslash$
\newline
r1=1;for(j=1,tfr,r1=r1$\ast$fr[j,1]);$\backslash$
\newline
if(type(r1/2)==1,r1=r1/2,);$\backslash$
\newline
if(r1$>$epsilon1,epsilon1=r1,);$\backslash$
\newline
if(tfr$<$=0,,if(fr[tfr,1]$>$epsilon2,epsilon2=fr[tfr,1],));$\backslash$
\newline
))))))));pprint("nAI1=",n);pprint("aAI1=",epsilon);$\backslash$
\newline
pprint("aAI1\_1=",epsilon1);pprint("aAI1\_2=",epsilon2);

\vskip10pt

\centerline{{\bf Program 14:} funda10.gen}

\noindent
$\backslash$$\backslash$funda10.gen
\newline
$\backslash$$\backslash$hyperbolic type, case A10
\newline
$\backslash$$\backslash$cosh(eta/2)=3/2
\newline
$\backslash$l;$\backslash$
\newline
epsilon=1;epsilon1=1;epsilon2=1;n=0;$\backslash$
\newline
rh1=[9.412375826,10.37390342,11.10113930,11.70820393];$\backslash$
\newline
alpha12=0;for(alpha23=1,4,w=(rh1[alpha23]+0.1)$^\wedge$2;
for(alpha13=5,w,$\backslash$
\newline
if(alpha23==0$||$issquare(alpha12$\ast$alpha23$\ast$alpha13)!=1$
||$$\backslash$
\newline
--8+2$\ast$isqrt(alpha12$\ast$alpha23$\ast$alpha13)+2$
\ast$alpha12+2$\ast$alpha13+2$\ast$alpha23$<$=0,,$\backslash$
\newline
alpha=4$\ast$idmat(3);alpha[1,2]=alpha12;alpha[2,1]=alpha12;$\backslash$
\newline
alpha[2,3]=alpha23;alpha[3,2]=alpha23;alpha[1,3]=alpha13;$\backslash$
\newline
alpha[3,1]=alpha13;dalpha=$\backslash$
\newline
--8+2$\ast$isqrt(alpha12$\ast$alpha23$\ast$alpha13)+2$
\ast$alpha12+2$\ast$alpha13+2$\ast$alpha23;$\backslash$
\newline
a=--2$\ast$idmat(3);a12=0;a21=0;$\backslash$
\newline
fordiv(alpha[2,3],a23,fordiv(alpha[1,3],a13,$\backslash$
\newline
a32=alpha[2,3]/a23;a31=alpha[1,3]/a13;$\backslash$
\newline
if(a12$\ast$a23$\ast$a31!=a21$\ast$a13$\ast$a32,,$\backslash$
\newline
a[1,2]=a12;a[2,1]=a21;a[2,3]=a23;a[3,2]=a32;a[1,3]=a13;
a[3,1]=a31;$\backslash$
\newline
dd=idmat(3);dd[1,1]=a13$\ast$a32;dd[2,2]=a23$\ast$a31;
dd[3,3]=a31$\ast$a32;$\backslash$
\newline
b=a$\ast$dd;b=b/content(b);n=n+1;$\backslash$
\newline
db=smith(b);r=db[1];$\backslash$
\newline
if(r$>$epsilon,epsilon=r,);$\backslash$
\newline
fr=factor(r);tfr=matsize(fr)[1];$\backslash$
\newline
r1=1;for(j=1,tfr,r1=r1$\ast$fr[j,1]);$\backslash$
\newline
if(type(r1/2)==1,r1=r1/2,);$\backslash$
\newline
if(r1$>$epsilon1,epsilon1=r1,);$\backslash$
\newline
if(tfr$<$=0,,if(fr[tfr,1]$>$epsilon2,epsilon2=fr[tfr,1],));$\backslash$
\newline
))))));pprint("nAI0=",n);pprint("aAI0=",epsilon);$\backslash$
\newline
pprint("aAI0\_1=",epsilon1);pprint("aAI0\_2=",epsilon2);

\vskip10pt 

\centerline{{\bf Program 15:} funda21.gen}

\noindent
$\backslash$$\backslash$funda21.gen
\newline
$\backslash$$\backslash$hyperbolic type, case AII1
\newline
$\backslash$$\backslash$cosh(eta/2)=3/2
\newline
$\backslash$l;$\backslash$
\newline
epsilon=1;epsilon1=1;epsilon2=1;n=0;$\backslash$
\newline
rh1=[9.412375826,10.37390342,11.10113930,11.70820393];$\backslash$
\newline
rh2=[38.68043607,41.73090517,44.05297726,46];alpha12=alpha23=0;$
\backslash$
\newline
for(alpha34=1,4,w=(rh2[alpha34]+0.1)$^\wedge$2;$\backslash$
\newline
for(alpha14=5,w,for(alpha13=5,(7+0.1)$^\wedge$2,$\backslash$
\newline
if(issquare(u=alpha13$\ast$alpha34$\ast$alpha14)!=1,,$
\backslash$
\newline
if(type(alpha24=4+4$\ast$(alpha14+alpha34+isqrt(u))/(alpha13--4))!=1$
||$$\backslash$
\newline
alpha24$<$=(rh1[alpha34]--0.1)$^\wedge$2,,alpha=4$\ast$idmat(4);$
\backslash$
\newline
alpha[1,3]=alpha13;alpha[3,1]=alpha13;alpha[3,4]=alpha34;$
\backslash$
\newline
alpha[4,3]=alpha34;alpha[2,4]=alpha24;alpha[4,2]=alpha24;
$\backslash$
\newline
alpha[1,4]=alpha14;alpha[4,1]=alpha14;$\backslash$
\newline
fordiv(alpha[3,4],a34,fordiv(alpha[1,3],a13,fordiv(alpha[2,4],a24,$
\backslash$
\newline
if(type(a41=isqrt(u)/(a13$\ast$a34))!=1,,a=--2$\ast$idmat(4);$
\backslash$
\newline
a[3,4]=a34;a[4,3]=alpha[3,4]/a34;a[2,4]=a24;a[4,2]=alpha[2,4]/a24;$
\backslash$
\newline
a[1,3]=a13;a[3,1]=alpha[1,3]/a13;$\backslash$
\newline
a[4,1]=a41;if(a41==0,a[1,4]=0,a[1,4]=alpha[1,4]/a41);$\backslash$
\newline
if(type(a[1,4])!=1,,diag=idmat(4);diag[1,1]=a[1,3]$\ast$a[3,4]
$\ast$a[4,2];$\backslash$
\newline
diag[2,2]=a[3,1]$\ast$a[4,3]$\ast$a[2,4];diag[3,3]=a[3,1]$\ast$a[3,4]$
\ast$a[4,2];$\backslash$
\newline
diag[4,4]=a[3,1]$\ast$a[4,3]$\ast$a[4,2];b=a$\ast
$diag;b=b/content(b);n=n+1;$\backslash$
\newline
db=smith(b);r=db[2];if(r$>$epsilon,epsilon=r,);$\backslash$
\newline
fr=factor(r);tfr=matsize(fr)[1];r1=1;for(j=1,tfr,r1=r1
$\ast$fr[j,1]);$\backslash$
\newline
if(type(r1/2)==1,r1=r1/2,);if(r1$>$epsilon1,epsilon1=r1,);
$\backslash$
\newline
if(tfr$<$=0,,if(fr[tfr,1]$>$epsilon2,epsilon2=fr[tfr,1],));
$\backslash$
\newline
))))))))));pprint("nAII1=",n);pprint("aAII1=",epsilon);
$\backslash$
\newline
pprint("aAII1\_1=",epsilon1);pprint("aAII1\_2=",epsilon2);

\vskip10pt 

\centerline{{\bf Program 16:} funda20.gen}

\noindent
$\backslash$$\backslash$program funda20.gen
\newline
$\backslash$$\backslash$hyperbolic type, case AII0 
\newline
$\backslash$$\backslash$cosh(eta/2)=3/2
\newline
$\backslash$l;$\backslash$
\newline
epsilon=1;epsilon1=1;epsilon2=1;n=0;alpha12=0;alpha23=0;alpha34=0;$
\backslash$
\newline
for(alpha14=5,(31.15549442+0.1)$^\wedge$2,fordiv(4$\ast$alpha14,aa,$
\backslash$
\newline
if(aa$^\wedge$2$>$4$\ast$alpha14,,alpha13=4+aa;alpha24=4$\ast
$alpha14/aa+4;$\backslash$
\newline
if(alpha13$>$(7+0.1)$^\wedge$2,,alpha=4$\ast$idmat(4);$\backslash$
\newline
alpha[1,3]=alpha13;alpha[3,1]=alpha13;alpha[3,4]=alpha34;$\backslash$
\newline
alpha[4,3]=alpha34;alpha[2,4]=alpha24;alpha[4,2]=alpha24;$\backslash$
\newline
alpha[1,4]=alpha14;alpha[4,1]=alpha14;$\backslash$
\newline
fordiv(alpha[1,3],a13,fordiv(alpha[2,4],a24,fordiv(alpha[1,4],a14,$
\backslash$
\newline
a=-2$\ast$idmat(4);a[1,3]=a13;a[3,1]=alpha[1,3]/a13;a[1,4]=a14;$
\backslash$
\newline
a[4,1]=alpha[1,4]/a14;a[1,3]=a13;a[3,1]=alpha[1,3]/a13;$\backslash$
\newline
a[2,4]=a24;a[4,2]=alpha[2,4]/a24;$\backslash$
\newline
diag=idmat(4);diag[1,1]=a[1,3]$\ast$a[1,4]$\ast$a[4,2];$\backslash$
\newline
diag[2,2]=a[1,3]$\ast$a[4,1]$\ast$a[2,4];diag[3,3]=a[3,1]$\ast$a[1,4]
$\ast$a[4,2];$\backslash$
\newline
diag[4,4]=a[1,3]$\ast$a[4,2]$\ast$a[4,1];b=a$\ast
$diag;b=b/content(b);n=n+1;$\backslash$
\newline
db=smith(b);r=db[2];if(r$>$epsilon,epsilon=r,);$\backslash$
\newline
fr=factor(r);tfr=matsize(fr)[1];r1=1;for(j=1,tfr,r1=r1$\ast$fr[j,1]);$
\backslash$
\newline
if(type(r1/2)==1,r1=r1/2,);if(r1$>$epsilon1,epsilon1=r1,);$\backslash$
\newline
if(tfr$<$=0,,if(fr[tfr,1]$>$epsilon2,epsilon2=fr[tfr,1],));$\backslash$
\newline
)))))));pprint("nAII0=",n);pprint("aAII0=",epsilon);$\backslash$
\newline
pprint("aAII0$\_$1=",epsilon1);pprint("aAII0$\_$2=",epsilon2);

\vskip10pt 

\centerline{{\bf Program 17:} funda3.gen}

\noindent
$\backslash$$\backslash$funda3.gen
\newline
$\backslash$$\backslash$hyperbolic type, case AIII
\newline
$\backslash$$\backslash$cosh(eta/2)=3/2
\newline
$\backslash$l;$\backslash$
\newline
epsilon=1;epsilon1=1;epsilon2=1;n=0;al12=al23=al34=al45=0;$\backslash$
\newline
for(al15=5,(68.1815011826+0.1)$^\wedge$2,for(al13=5,(7+0.1)$
^\wedge$2,$\backslash$
\newline
for(al35=al13,(7+0.1)$^\wedge$2,$\backslash$
\newline
if(issquare(q=al13$\ast$al35$\ast
$al15)==0,,d=(al13+al35+al15--4+isqrt(q))$\ast$4;$\backslash$
\newline
if(type(al14=d/(al35--4))!=1$||$al14$<$=(31.15549442--0.1)$
^\wedge$2$\backslash$
\newline
$||$type(al25=d/(al13--4))!=1$||$al25$<$=(31.15549442--0.1)$
^\wedge$2$||$$\backslash$
\newline
type(al24=(al13$\ast$al35+4$\ast$al15+4$\ast$isqrt(q))$\ast
$4/((al35--4)$\ast$(al13--4)))!=1$||$$\backslash$
\newline
issquare(q1=al13$\ast$al35$\ast$al25$\ast$al24$\ast$al14)==0,,$
\backslash$
\newline
al=idmat(5)$\ast$4;al[1,5]=al15;al[5,1]=al15;al[1,3]=al13;al[3,1]=al13;$
\backslash$
\newline
al[1,4]=al14;al[4,1]=al14;al[2,4]=al24;al[4,2]=al24;$\backslash$
\newline
al[2,5]=al25;al[5,2]=al25;al[3,5]=al35;al[5,3]=al35;$\backslash$
\newline
fordiv(al13,a13,fordiv(al35,a35,a51=isqrt(q)/a13/a35;$\backslash$
\newline
if(a51==0,a15==0,if(type(a15=al15/a51)!=1,,fordiv(al14,a14,$
\backslash$
\newline
fordiv(al24,a24,a52=isqrt(q1)/a13/a35/a24$\ast$a14/al14;$\backslash$
\newline
if(type(a25=al25/a52)!=1,,a=idmat(5)$\ast$--2;$\backslash$
\newline
a[1,3]=a13;a[3,1]=al13/a13;a[1,4]=a14;a[4,1]=al14/a14;$\backslash$
\newline
a[1,5]=a15;a[5,1]=a51;a[2,4]=a24;a[4,2]=al24/a24;$\backslash$
\newline
a[2,5]=a25;a[5,2]=a52;a[3,5]=a35;a[5,3]=al35/a35;$\backslash$
\newline
diag=idmat(5);diag[1,1]=a[1,4]$\ast$a[1,3]$\ast$a[4,2]$\ast$a[2,5];$
\backslash$
\newline
diag[2,2]=a[4,1]$\ast$a[1,3]$\ast$a[2,4]$\ast$a[2,5];$\backslash$
\newline
diag[3,3]=a[1,4]$\ast$a[3,1]$\ast$a[4,2]$\ast$a[2,5];$\backslash$
\newline
diag[4,4]=a[4,1]$\ast$a[1,3]$\ast$a[4,2]$\ast$a[2,5];$\backslash$
\newline
diag[5,5]=a[4,1]$\ast$a[1,3]$\ast$a[2,4]$\ast$a[5,2];$\backslash$
\newline
b=a$\ast$diag;b=b/content(b);n=n+1;db=smith(b);r=db[3];$\backslash$
\newline
if(r$>$epsilon,epsilon=r,);fr=factor(r);tfr=matsize(fr)[1];$\backslash$
\newline
r1=1;for(j=1,tfr,r1=r1$\ast$fr[j,1]);if(type(r1/2)==1,r1=r1/2,);$\backslash$
\newline
if(r1$>$epsilon1,epsilon1=r1,);$\backslash$
\newline
if(fr[tfr,1]$>$epsilon2,epsilon2=fr[tfr,1],);$\backslash$
\newline
))))))))))));pprint("nAIII=",n);$\backslash$
\newline
pprint("aAIII=",epsilon);pprint("aAIII$\_$1=",epsilon1);$\backslash$
\newline
pprint("aAIII$\_$2=",epsilon2);

\vskip10pt 

\centerline{{\bf Program 18:} fundb1.gen}

\noindent
$\backslash$$\backslash$fundb1.gen
\newline
$\backslash$$\backslash$hyperbolic type, case BI
\newline
$\backslash$$\backslash$cosh(eta/2)=3/2
\newline
$\backslash$l;$\backslash$
\newline
epsilon=1;epsilon1=1;epsilon2=1;n=0;$\backslash$
\newline
rh4=[1.191398091,1.095286061,1.022736500,0.962423650];$\backslash$
\newline
for(alpha12=1,4,for(alpha23=5,22$^\wedge$2,for(alpha13=5,alpha23,$
\backslash$
\newline
if(sqrt(alpha13)/8--5/4$<$1+0.000001,ach1=0,$\backslash$
\newline
ach1=acosh(sqrt(alpha13)/8--5/4));$\backslash$
\newline
if(sqrt(alpha23)/8--5/4$<$1+0.000001,ach2=0,$\backslash$
\newline
ach2=acosh(sqrt(alpha23)/8--5/4));$\backslash$
\newline
if(ach1+ach2$>$rh4[alpha12]+0.1$||$alpha23==0$||$$\backslash$
\newline
issquare(alpha12$\ast$alpha23$\ast$alpha13)!=1$||$$\backslash$
\newline
--8+2$\ast$isqrt(alpha12$\ast$alpha23$\ast$alpha13)+2$\ast$alpha12+2$
\ast$alpha13+2$\ast$alpha23$<$=0,,$\backslash$
\newline
alpha=4$\ast$idmat(3);alpha[1,2]=alpha12;alpha[2,1]=alpha12;$\backslash$
\newline
alpha[2,3]=alpha23;alpha[3,2]=alpha23;alpha[1,3]=alpha13;$\backslash$
\newline
alpha[3,1]=alpha13;dalpha=$\backslash$
\newline
--8+2$\ast$isqrt(alpha12$\ast$alpha23$\ast$alpha13)+2$\ast$alpha12+2$
\ast$alpha13+2$\ast$alpha23;$\backslash$
\newline
a=--2$\ast$idmat(3);$\backslash$
\newline
fordiv(alpha[1,2],a12,fordiv(alpha[2,3],a23,fordiv(alpha[1,3],a13,$
\backslash$
\newline
a21=alpha[1,2]/a12;a32=alpha[2,3]/a23;a31=alpha[1,3]/a13;$\backslash$
\newline
if(a12$\ast$a23$\ast$a31!=a21$\ast$a13$\ast$a32,,$\backslash$
\newline
a[1,2]=a12;a[2,1]=a21;a[2,3]=a23;a[3,2]=a32;a[1,3]=a13;a[3,1]=a31;$
\backslash$
\newline
dd=idmat(3);dd[1,1]=a13$\ast$a32;dd[2,2]=a23$\ast$a31;dd[3,3]=a31$
\ast$a32;$\backslash$
\newline
b=a$\ast$dd;b=b/content(b);n=n+1;db=smith(b);r=db[1];$\backslash$
\newline
if(r$>$epsilon,epsilon=r,);fr=factor(r);tfr=matsize(fr)[1];$\backslash$
\newline
r1=1;for(j=1,tfr,r1=r1$\ast$fr[j,1]);if(type(r1/2)==1,r1=r1/2,);$
\backslash$
\newline
if(r1$>$epsilon1,epsilon1=r1,);$\backslash$
\newline
if(tfr$<$=0,,if(fr[tfr,1]$>$epsilon2,epsilon2=fr[tfr,1],));$\backslash$
\newline
))))))));pprint("nBI=",n);pprint("aBI1=",epsilon);$\backslash$
\newline
pprint("aBI\_1=",epsilon1);pprint("aB1\_2=",epsilon2); 

\vskip10pt 

\centerline{{\bf Program 19:} fundb21.gen}

\noindent
$\backslash$$\backslash$fundb21.gen
\newline
$\backslash$$\backslash$hyperbolic type, case BII$\_$1
\newline
$\backslash$$\backslash$cosh(eta/2)=3/2
\newline
$\backslash$l;$\backslash$
\newline
epsilon=1;epsilon1=1;epsilon2=1;n=0;alpha12=alpha23=0;$
\backslash$
\newline
for(alpha24=5,22$^\wedge$2,if(sqrt(alpha24)/8--5/4$<
$1+0.000001,ach2=0,$\backslash$
\newline
ach2=acosh(sqrt(alpha24)/8--5/4));$\backslash$
\newline
vv=min((8$\ast$(cosh(1.429914377--ach2)+5/4))$^\wedge$2+0.1,22$^\wedge
$2);$\backslash$
\newline
for(alpha14=5,vv,for(alpha34=alpha14,vv,$\backslash$
\newline
if(issquare(discrim=alpha14$\ast$alpha34$\ast$(alpha14$\ast
$alpha34+(alpha24--4)$\ast$$\backslash$
\newline
(alpha14+alpha24+alpha34--4)))!=1,,$\backslash$
\newline
if(type(u=(2$\ast$alpha14$\ast$alpha34+2$\ast
$isqrt(discrim))/(alpha24--4))!=1,,$\backslash$
\newline
if(type(alpha13=u$^\wedge$2/(alpha14$\ast$alpha34))!=1,,alpha=4$
\ast$idmat(4);$\backslash$
\newline
alpha[1,3]=alpha13;alpha[3,1]=alpha13;alpha[3,4]=alpha34;alpha[4,3]=alpha34;$
\backslash$
\newline
alpha[2,4]=alpha24;alpha[4,2]=alpha24;alpha[1,4]=alpha14;alpha[4,1]=alpha14;$
\backslash$
\newline
fordiv(alpha[3,4],a34,fordiv(alpha[1,3],a13,fordiv(alpha[2,4],a24,$
\backslash$
\newline
if(type(a41=u/(a13$\ast$a34))!=1,,a=--2$\ast$idmat(4);$
\backslash$
\newline
a[3,4]=a34;a[4,3]=alpha[3,4]/a34;a[2,4]=a24;a[4,2]=alpha[2,4]/a24;
$\backslash$
\newline
a[1,3]=a13;a[3,1]=alpha[1,3]/a13;$\backslash$
\newline
a[4,1]=a41;if(a41==0,a[1,4]=0,a[1,4]=alpha[1,4]/a41);if(type(a[1,4])!=1,,$
\backslash$
\newline
diag=idmat(4);diag[1,1]=a[1,3]$\ast$a[3,4]$\ast$a[4,2];$
\backslash$
\newline
diag[2,2]=a[3,1]$\ast$a[4,3]$\ast$a[2,4];diag[3,3]=a[3,1]$
\ast$a[3,4]$\ast$a[4,2];$\backslash$
\newline
diag[4,4]=a[3,1]$\ast$a[4,3]$\ast$a[4,2];b=a$\ast
$diag;b=b/content(b);n=n+1;$\backslash$
\newline
db=smith(b);r=db[2];if(r$>$epsilon,epsilon=r,);$\backslash$
\newline
fr=factor(r);tfr=matsize(fr)[1];r1=1;for(j=1,tfr,r1=r1$\ast
$fr[j,1]);$\backslash$
\newline
if(type(r1/2)==1,r1=r1/2,);if(r1$>$epsilon1,epsilon1=r1,);$
\backslash$
\newline
if(tfr$<$=0,,if(fr[tfr,1]$>$epsilon2,epsilon2=fr[tfr,1],));$
\backslash$
\newline
)))))))))));pprint("nBII1=",n);pprint("aBII1=",epsilon);$\backslash$
\newline
pprint("aBII1$\_$1=",epsilon1);pprint("aBII1$\_$2=",epsilon2);

\vskip10pt 

\centerline{{\bf Program 20:} fundb22.gen}

\noindent
$\backslash$$\backslash$fundb22.gen
\newline
$\backslash$$\backslash$hyperbolic type, case BII\_2
\newline
$\backslash$$\backslash$cos(eta/2)=3/2
\newline
$\backslash$l;$\backslash$
\newline
epsilon=1;epsilon1=1;epsilon2=1;n=0;alpha12=alpha34=0;$\backslash$
\newline
for(alpha23=5,22$^\wedge$2,if(sqrt(alpha23)/8--5/4$<$1+0.000001,ach2=0,$
\backslash$
\newline
ach2=acosh(sqrt(alpha23)/8--5/4));$\backslash$
\newline
vv=min((8$\ast$(cosh(1.429914377--ach2)+5/4))$^\wedge$2+0.1,22$^\wedge$2);$
\backslash$
\newline
for(alpha13=5,vv,for(alpha24=alpha13,vv,$\backslash$
\newline
if(issquare(discrim=alpha13$\ast$alpha23$\ast$alpha24$\ast$$\backslash$
\newline
(alpha13$\ast$alpha23$\ast$alpha24--(alpha23--4)$\ast$((alpha13--4)$\ast
$(alpha24--4)--$\backslash$
\newline
4$\ast$alpha23)))!=1,,forstep(tau=--1,1,2,$\backslash$
\newline
if(type(u=(alpha13$\ast$alpha23$\ast$alpha24+tau$\ast
$isqrt(discrim))/(alpha23--4))!=1,,$\backslash$
\newline
if(type(alpha14=u$^\wedge$2/(alpha13$\ast$alpha23$\ast$alpha24))!=1$
||$alpha14$<$=4,,$\backslash$
\newline
alpha=4$\ast$idmat(4);$\backslash$
\newline
alpha[1,3]=alpha13;alpha[3,1]=alpha13;alpha[1,4]=alpha14;alpha[4,1]=alpha14;$
\backslash$
\newline
alpha[2,3]=alpha23;alpha[3,2]=alpha23;alpha[2,4]=alpha24;alpha[4,2]=alpha24;$
\backslash$
\newline
fordiv(alpha[1,3],a13,fordiv(alpha[2,3],a32,fordiv(alpha[2,4],a24,$
\backslash$
\newline
if(type(a41=u/(a13$\ast$a32$\ast$a24))!=1,,a=--2$\ast$idmat(4);$
\backslash$
\newline
a[1,3]=a13;a[3,1]=alpha[1,3]/a13;a[3,2]=a32;a[2,3]=alpha[2,3]/a32;$
\backslash$
\newline
a[2,4]=a24;a[4,2]=alpha[2,4]/a24;a[4,1]=a41;a[1,4]=alpha[1,4]/a41;$
\backslash$
\newline
if(type(a[1,4])!=1,,diag=idmat(4);diag[1,1]=a[1,3]$\ast$a[1,4]$\ast
$a[3,2];$\backslash$
\newline
diag[2,2]=a[3,1]$\ast$a[1,4]$\ast$a[2,3];diag[3,3]=a[3,1]$\ast$a[1,4]$\ast
$a[3,2];$\backslash$
\newline
diag[4,4]=a[1,3]$\ast$a[4,1]$\ast$a[3,2];b=a$\ast
$diag;b=b/content(b);n=n+1;$\backslash$
\newline
db=smith(b);r=db[2];if(r$>$epsilon,epsilon=r,);$\backslash$
\newline
fr=factor(r);tfr=matsize(fr)[1];r1=1;for(j=1,tfr,r1=r1$\ast$fr[j,1]);$
\backslash$
\newline
if(type(r1/2)==1,r1=r1/2,);if(r1$>$epsilon1,epsilon1=r1,);$\backslash$
\newline
if(tfr$<$=0,,if(fr[tfr,1]$>$epsilon2,epsilon2=fr[tfr,1],));$\backslash$
\newline
))))))))))));pprint("nBII2=",n);pprint("aBII2=",epsilon);$\backslash$
\newline
pprint("aBII2\_1=",epsilon1);pprint("aBII2\_2=",epsilon2);

\vskip10pt 

\centerline{{\bf Program 21:} fundb3.gen}

\noindent
$\backslash$$\backslash$fundb3.gen
\newline
$\backslash$$\backslash$hyperbolic type, case BIII
\newline
$\backslash$$\backslash$cosh(eta/2)=3/2
\newline
$\backslash$l;$\backslash$
\newline
epsilon=1;epsilon1=1;epsilon2=1;n=0;alpha12=alpha23=0;$\backslash$
\newline
for(alpha14=5,(31.70820393)$^\wedge$2+0.1,$\backslash$
\newline
if(sqrt(alpha14)/8--5/4$<$1+0.000001,ach1=0,$\backslash$
\newline
ach1=acosh(sqrt(alpha14)/8--5/4));$\backslash$
\newline
vv=min((8$\ast$(cosh(1.9248473002--ach1)+5/4))$^\wedge$2+0.1,(31.70820393)$
^\wedge$2+0.1);$\backslash$
\newline
for(alpha34=alpha14,vv,for(alpha13=5,(7+0.1)$^\wedge$2,$\backslash$
\newline
if(issquare(u=alpha13$\ast$alpha34$\ast$alpha14)!=1,,$\backslash$
\newline
if(type(alpha24=4+4$\ast$(alpha14+alpha34+isqrt(u))/(alpha13--4))!=1,,$
\backslash$
\newline
alpha=4$\ast$idmat(4);alpha[1,3]=alpha13;alpha[3,1]=alpha13;$\backslash$
\newline
alpha[3,4]=alpha34;alpha[4,3]=alpha34;alpha[2,4]=alpha24;alpha[4,2]=alpha24;$
\backslash$
\newline
alpha[1,4]=alpha14;alpha[4,1]=alpha14;fordiv(alpha[3,4],a34,$\backslash$
\newline
fordiv(alpha[1,3],a13,fordiv(alpha[2,4],a24,$\backslash$
\newline
if(type(a41=isqrt(u)/(a13$\ast$a34))!=1,,a=--2$\ast$idmat(4);$\backslash$
\newline
a[3,4]=a34;a[4,3]=alpha[3,4]/a34;a[2,4]=a24;a[4,2]=alpha[2,4]/a24;$
\backslash$
\newline
a[1,3]=a13;a[3,1]=alpha[1,3]/a13;a[4,1]=a41;$\backslash$
\newline
if(a41==0,a[1,4]=0,a[1,4]=alpha[1,4]/a41);if(type(a[1,4])!=1,,$\backslash$
\newline
diag=idmat(4);diag[1,1]=a[1,3]$\ast$a[3,4]$\ast$a[4,2];$\backslash$
\newline
diag[2,2]=a[3,1]$\ast$a[4,3]$\ast$a[2,4];diag[3,3]=a[3,1]$\ast$a[3,4]$
\ast$a[4,2];$\backslash$
\newline
diag[4,4]=a[3,1]$\ast$a[4,3]$\ast$a[4,2];b=a$\ast
$diag;b=b/content(b);n=n+1;$\backslash$
\newline
db=smith(b);r=db[2];if(r$>$epsilon,epsilon=r,);$\backslash$
\newline
fr=factor(r);tfr=matsize(fr)[1];r1=1;for(j=1,tfr,r1=r1$\ast$fr[j,1]);$
\backslash$
\newline
if(type(r1/2)==1,r1=r1/2,);if(r1$>$epsilon1,epsilon1=r1,);$\backslash$
\newline
if(tfr$<$=0,,if(fr[tfr,1]$>$epsilon2,epsilon2=fr[tfr,1],));$\backslash$
\newline
))))))))));pprint("nBIII=",n);pprint("aBIII=",epsilon);$\backslash$
\newline
pprint("aBIII\_1=",epsilon1);pprint("aBIII\_2=",epsilon2);

\vskip10pt  

\centerline{{\bf Program 22:} funda11.main}

\noindent
$\backslash$$\backslash$funda11.mai
\newline
$\backslash$$\backslash$hyperbolic type, case AI1
\newline
$\backslash$$\backslash$cosh(eta/2)=3/2
\newline
$\backslash$r h2
\newline
$\backslash$r h3
\newline
$\backslash$l;$\backslash$
\newline
n=0;$\backslash$
\newline
rh1=$\backslash$
\newline
[12.25000000,13.37793021,14.23012096,14.94097150;$\backslash$
\newline
 13.37793021,14.57106781,15.47225159,16.22381115;$\backslash$
\newline
 14.23012096,15.47225159,16.41025403,17.19241152;$\backslash$
\newline
 14.94097150,16.22381115,17.19241152,18.00000000];$\backslash$
\newline
for(alpha12=1,4,for(alpha23=alpha12,4,u1u=(rh1[alpha12,alpha23]+0.1)$
^\wedge$2;$\backslash$
\newline
for(alpha13=5,u1u,$\backslash$
\newline
if(alpha23==0$||$issquare(alpha12$\ast$alpha23$\ast$alpha13)!=1$||
$$\backslash$
\newline
--8+2$\ast$isqrt(alpha12$\ast$alpha23$\ast$alpha13)+2$\ast$alpha12+2$\ast
$alpha13+2$\ast$alpha23$<$=0,,$\backslash$
\newline
alpha=4$\ast$idmat(3);alpha[1,2]=alpha12;alpha[2,1]=alpha12;$\backslash$
\newline
alpha[2,3]=alpha23;alpha[3,2]=alpha23;alpha[1,3]=alpha13;$\backslash$
\newline
alpha[3,1]=alpha13;dalpha=$\backslash$
\newline
--8+2$\ast$isqrt(alpha12$\ast$alpha23$\ast$alpha13)+2$\ast$alpha12+2$
\ast$alpha13+2$\ast$alpha23;$\backslash$
\newline
a=--2$\ast$idmat(3);$\backslash$
\newline
fordiv(alpha[1,2],a12,fordiv(alpha[2,3],a23,fordiv(alpha[1,3],a13,$
\backslash$
\newline
a21=alpha[1,2]/a12;a32=alpha[2,3]/a23;a31=alpha[1,3]/a13;$\backslash$
\newline
if(a12$\ast$a23$\ast$a31!=a21$\ast$a13$\ast$a32,,$\backslash$
\newline
a[1,2]=a12;a[2,1]=a21;a[2,3]=a23;a[3,2]=a32;a[1,3]=a13;a[3,1]=a31;$
\backslash$
\newline
dd=idmat(3);dd[1,1]=a13$\ast$a32;dd[2,2]=a23$\ast$a31;dd[3,3]=a31$\ast$a32;
$\backslash$
\newline
b=a$\ast$dd;b=b/content(b);db=smith(b);fdb=factor(db[1]);$\backslash$
\newline
detb=db[1]$\ast$db[2]$\ast$db[3];fdetb=factor(detb);$\backslash$
\newline
if(issqfree(b[1,1])==0$||$issqfree(b[2,2])==0$||$issqfree(b[3,3])==0$||
$$\backslash$
\newline
(content([--b[1,1]$\ast$--b[2,2]$\ast$--b[3,3],2])$>$1\&\&$\backslash$
\newline
content([--b[1,1]$\ast$--b[2,2]$\ast$--b[3,3],8])$<$8\&\&$\backslash$
\newline
fdetb[1,1]==2\&\&mod(fdetb[1,2],2)==mod(1,2))$||$$\backslash$
\newline
(alpha[1,2]$<$4\&\&content([--b[1,1],--b[2,2]])$>$2)$||$$\backslash$
\newline
(alpha[2,3]$<$4\&\&content([--b[2,2],--b[3,3]])$>$2),,$\backslash$
\newline
gam=0;for(j=1,matsize(fdetb)[1],if(fdetb[j,1]!=2\&\&$\backslash$
\newline
type(--b[1,1]$\ast$--b[2,2]$\ast$--b[3,3]/fdetb[j,1])==1\&\&$\backslash$
\newline
mod(fdetb[j,2],2)==mod(0,2),gam=1,));$\backslash$
\newline
if(gam==1,,if(alpha[1,2]$<$4,b=b;$\backslash$
\newline
d=1;for(k=1,matsize(fdetb)[1],$\backslash$
\newline
if(mod(fdetb[k,2],2)==mod(1,2),d=d$\ast$fdetb[k,1],));$\backslash$
\newline
if(alpha[1,2]==1\&\&mod(d,3)!=mod(2,3),,if(d$<$=2,et=0,$\backslash$
\newline
fd=factor(d);if(fd[1,1]==2,d1=d/2,d1=d);fd1=factor(d1);et=0;$\backslash$
\newline
for(k=1,matsize(fd1)[1],if(alpha[1,2]==1,$\backslash$
\newline
if(kro(3$\ast$d/fd1[k,1],fd1[k,1])==1,,et=et+2$^\wedge
$(k--1)),if(alpha[1,2]==2,$\backslash$
\newline
if(kro(d/fd1[k,1],fd1[k,1])==1,,et=et+2$^\wedge$(k--1)),if(alpha[1,2]==3,$
\backslash$
\newline
if(fd1[k,1]==3,,if(kro(3$\ast$d/fd1[k,1],fd1[k,1])==1,,et=et+2$
^\wedge$(k--1))),))));$\backslash$
\newline
hhh=hnr(d,et);if(hhh$>$2$||$hhh==1,,$\backslash$
\newline
n=n+1;pprint("n=",n);pprint(b);pprint(db);$\backslash$
\newline
pprint(fdb);pprint("d=",d," et=",et," h=",hhh);$\backslash$
\newline
if(n==1,u=matrix(1,3,j,k,0);u[1,]=[d,et,hhh],aaa=0;nnn=matsize(u)[1];$
\backslash$
\newline
for(j=1,nnn,if(u[j,]==[d,et,hhh],aaa=1,));$\backslash$
\newline
if(aaa==1,,u1=matrix(nnn+1,3,j,k,0);for(t=1,nnn,$\backslash$
\newline
if(u[t,1]$<$d$||$(u[t,1]==d\&\&u[t,2]$<$et),u1[t,]=u[t,],u1[t+1,]=u[t,]));$
\backslash$
\newline
for(t=1,nnn+1,if(u1[t,]==[0,0,0],u1[t,]=[d,et,hhh],));u=u1))))),$\backslash$
\newline
if(alpha[1,2]==4,d2=1;for(k=1,matsize(fdetb)[1],$\backslash$
\newline
if(mod(fdetb[k,2],2)==mod(0,2)\&\&fdetb[k,1]!=2,d2=d2$\ast$fdetb[k,1],));
$\backslash$
\newline
d2=2$\ast$d2/content([2$\ast$d2,--b[1,1]$\ast$--b[2,2]$\ast
$--b[3,3]]);fordiv(d2,t,$\backslash$
\newline
if(mod(d2,2)==mod(0,2)\&\&fdetb[1,1]==2\&\&mod(fdetb[1,2],2)==mod(1,2)\&\&$
\backslash$
\newline
mod(t,2)==mod(1,2),,bb=b$\ast$t;dbb=smith(bb);fdbb=factor(dbb[1]);$
\backslash$
\newline
detbb=dbb[1]$\ast$dbb[2]$\ast$dbb[3];fdetbb=factor(detbb);$\backslash$
\newline
d=1;for(k=1,matsize(fdetbb)[1],$\backslash$
\newline
if(mod(fdetbb[k,2],2)==mod(1,2),d=d$\ast$fdetbb[k,1],));if(d$<$=2,et=0,$
\backslash$
\newline
fd=factor(d);if(fd[1,1]==2,d1=d/2,d1=d);fd1=factor(d1);et=0;$\backslash$
\newline
for(k=1,matsize(fd1)[1],$\backslash$
\newline
if(kro(--d/fd1[k,1],fd1[k,1])==1,,et=et+2$^\wedge$(k--1))));$\backslash$
\newline
hhh=hnr(d,et);if(hhh$>$2$||$hhh==1,,$\backslash$
\newline
n=n+1;pprint("n=",n);pprint(bb);pprint(dbb);$\backslash$
\newline
pprint(fdbb);pprint("d=",d," et=",et," h=",hhh);$\backslash$
\newline
if(n==1,u=matrix(1,3,j,k,0);u[1,]=[d,et,hhh],aaa=0;nnn=matsize(u)[1];$
\backslash$
\newline
for(j=1,nnn,if(u[j,]==[d,et,hhh],aaa=1,));$\backslash$
\newline
if(aaa==1,,u1=matrix(nnn+1,3,j,k,0);for(t=1,nnn,$\backslash$
\newline
if(u[t,1]$<$d$||$(u[t,1]==d\&\&u[t,2]$<$et),u1[t,]=u[t,],$\backslash$
\newline
u1[t+1,]=u[t,]));$\backslash$
\newline
for(t=1,nnn+1,if(u1[t,]==[0,0,0],u1[t,]=[d,et,hhh],));u=u1))))),$
\backslash$
\newline
))))))))))));pprint("u=",u);pprint("matsize u=",matsize(u));

\vskip10pt 

\centerline{{\bf Program 23:} funda10.main}

\noindent
$\backslash$$\backslash$funda10.mai
\newline
$\backslash$$\backslash$hyperbolic type, case A10
\newline
$\backslash$$\backslash$cosh(eta/2)=3/2
\newline
$\backslash$r h2
\newline
$\backslash$r h3
\newline
$\backslash$l;$\backslash$
\newline
n=0;$\backslash$
\newline
rh1=[9.412375826,10.37390342,11.10113930,11.70820393];$\backslash$
\newline
alpha12=0;for(alpha23=1,4,w=(rh1[alpha23]+0.1)$^\wedge$2;
for(alpha13=5,w,$
\backslash$
\newline
if(alpha23==0$||$issquare(alpha12$\ast$alpha23$\ast$alpha13)!=1$||
$$\backslash$
\newline
--8+2$\ast$isqrt(alpha12$\ast$alpha23$\ast$alpha13)+2$\ast$alpha12+2$
\ast$alpha13+2$\ast$alpha23$<$=0,,$\backslash$
\newline
alpha=4$\ast$idmat(3);alpha[1,2]=alpha12;alpha[2,1]=alpha12;$\backslash$
\newline
alpha[2,3]=alpha23;alpha[3,2]=alpha23;alpha[1,3]=alpha13;$\backslash$
\newline
alpha[3,1]=alpha13;dalpha=$\backslash$
\newline
--8+2$\ast$isqrt(alpha12$\ast$alpha23$\ast$alpha13)+2$\ast$alpha12+2$
\ast$alpha13+2$\ast$alpha23;$\backslash$
\newline
a=--2$\ast$idmat(3);a12=0;a21=0;$\backslash$
\newline
fordiv(alpha[2,3],a23,fordiv(alpha[1,3],a13,$\backslash$
\newline
a32=alpha[2,3]/a23;a31=alpha[1,3]/a13;if(a12$\ast$a23$\ast$a31!=a21$
\ast$a13$\ast$a32,,$\backslash$
\newline
a[1,2]=a12;a[2,1]=a21;a[2,3]=a23;a[3,2]=a32;a[1,3]=a13;a[3,1]=a31;$
\backslash$
\newline
dd=idmat(3);dd[1,1]=a13$\ast$a32;dd[2,2]=a23$\ast$a31;dd[3,3]=a31$
\ast$a32;$\backslash$
\newline
b=a$\ast$dd;b=b/content(b);db=smith(b);fdb=factor(db[1]);$\backslash$
\newline
if(issqfree(b[1,1])==0$||$issqfree(b[2,2])==0$||$issqfree(b[3,3])==0$
||$$\backslash$
\newline
content([--b[1,1],--b[2,2]])$>$2$||$$\backslash$
\newline
(alpha[2,3]$<$4\&\&content([--b[2,2],--b[3,3]])$>
$2),,detb=det(b);(if(detb==1,,$\backslash$
\newline
fdetb=factor(detb);$\backslash$
\newline
if(content([--b[1,1]$\ast$--b[2,2]$\ast$--b[3,3],2])$>$1\&\&$\backslash$
\newline
content([--b[1,1]$\ast$--b[2,2]$\ast$--b[3,3],8])$<$8\&\&$\backslash$
\newline
fdetb[1,1]==2\&\&mod(fdetb[1,2],2)==mod(1,2),,$\backslash$
\newline
gam=0;for(j=1,matsize(fdetb)[1],if(fdetb[j,1]!=2\&\&$\backslash$
\newline
type(--b[1,1]$\ast$--b[2,2]$\ast$--b[3,3]/fdetb[j,1])==1\&\&$\backslash$
\newline
mod(fdetb[j,2],2)==mod(0,2),gam=1,));if(gam==1,,for(k=0,1,$\backslash$
\newline
if(k==0\&\&content(--[b[1,1]$\ast$--b[2,2]$\ast$--b[3,3],2])==1\&\&$
\backslash$
\newline
fdetb[1,1]==2\&\&mod(fdetb[1,2],2)==mod(1,2),,$\backslash$
\newline
if(k==1\&\&content(--[b[1,1]$\ast$--b[2,2]$\ast$--b[3,3],2])==1,$
\backslash$
\newline
b=2$\ast$b;db=smith(b);fdb=factor(db[1]);gam1=0,gam1=1);$\backslash$
\newline
if(k==1\&\&gam1==1,,$\backslash$
\newline
detb=db[1]$\ast$db[2]$\ast
$db[3];fdetb=factor(detb);d=1;for(k=1,matsize(fdetb)[1],$\backslash$
\newline
if(mod(fdetb[k,2],2)==mod(1,2),d=d$\ast$fdetb[k,1],));if(d$<
$=2,et=0,$\backslash$
\newline
fd=factor(d);if(fd[1,1]==2,d1=d/2,d1=d);fd1=factor(d1);$\backslash$
\newline
et=0;for(k=1,matsize(fd1)[1],if(type(b[1,1]/fd1[k,1])==1,$\backslash$
\newline
if(kro(b[1,1]/fd1[k,1],fd1[k,1])==1,,et=et+2$^\wedge$(k--1)),$
\backslash$
\newline
if(type(b[2,2]/fd1[k,1])==1,$\backslash$
\newline
if(kro(b[2,2]/fd1[k,1],fd1[k,1])==1,,et=et+2$^\wedge$(k--1)),$
\backslash$
\newline
if(kro(d$\ast$b[1,1]$\ast$b[2,2]/fd1[k,1],fd1[k,1])==1,,et=et+2$
^\wedge$(k--1))))));$\backslash$
\newline
hhh=hnr(d,et);if(hhh$>$2$||$hhh==1,,$\backslash$
\newline
n=n+1;pprint("n=",n);pprint(b);pprint(db);$\backslash$
\newline
pprint(fdb);pprint("d=",d," et=",et," h=",hhh);$\backslash$
\newline
if(n==1,u=matrix(1,3,j,k,0);u[1,]=[d,et,hhh],aaa=0;nnn=matsize(u)[1];$
\backslash$
\newline
for(j=1,nnn,if(u[j,]==[d,et,hhh],aaa=1,));$\backslash$
\newline
if(aaa==1,,u1=matrix(nnn+1,3,j,k,0);for(t=1,nnn,$\backslash$
\newline
if(u[t,1]$<$d$||$(u[t,1]==d\&\&u[t,2]$<
$et),u1[t,]=u[t,],u1[t+1,]=u[t,]));$\backslash$
\newline
for(t=1,nnn+1,if(u1[t,]==[0,0,0],u1[t,]=[d,et,hhh],));u=u1));$\backslash$
\newline
)))))))))))))));pprint("u=",u);pprint("matsize u=",matsize(u));

\vskip10pt 

\centerline{{\bf Program 24:} funda21.main}

\noindent
$\backslash$$\backslash$funda21.main
\newline
$\backslash$$\backslash$hyperbolic type, case A21
\newline
$\backslash$$\backslash$cosh(eta/2)=3/2
\newline
$\backslash$r h2
\newline
$\backslash$r h3
\newline
$\backslash$l;$\backslash$
\newline
n=0;$\backslash$
\newline
rh1=[9.412375826,10.37390342,11.10113930,11.70820393];$\backslash$
\newline
rh2=[38.68043607,41.73090517,44.05297726,46];alpha12=alpha23=0;$
\backslash$
\newline
for(alpha34=1,4,w=(rh2[alpha34]+0.1)$^\wedge$2;for(alpha14=5,w,$
\backslash$
\newline
for(alpha13=5,(7+0.1)$^\wedge$2,if(issquare(u=alpha13$\ast$alpha34$
\ast$alpha14)!=1,,$\backslash$
\newline
if(type(alpha24=4+4$\ast$(alpha14+alpha34+isqrt(u))/(alpha13--4))!=1$
||$$\backslash$
\newline
alpha24$<$=(rh1[alpha34]--0.1)$^\wedge$2,,alpha=4$\ast$idmat(4);$\backslash$
\newline
alpha[1,3]=alpha13;alpha[3,1]=alpha13;alpha[3,4]=alpha34;$\backslash$
\newline
alpha[4,3]=alpha34;alpha[2,4]=alpha24;alpha[4,2]=alpha24;$\backslash$
\newline
alpha[1,4]=alpha14;alpha[4,1]=alpha14;$\backslash$
\newline
fordiv(alpha[3,4],a34,fordiv(alpha[1,3],a13,$\backslash$
\newline
fordiv(alpha[2,4],a24,if(type(a41=isqrt(u)/(a13$\ast$a34))!=1,,$\backslash$
\newline
a=--2$\ast$idmat(4);a[3,4]=a34;a[4,3]=alpha[3,4]/a34;a[2,4]=a24;$\backslash$
\newline
a[4,2]=alpha[2,4]/a24;a[1,3]=a13;a[3,1]=alpha[1,3]/a13;$\backslash$
\newline
a[4,1]=a41;if(a41==0,a[1,4]=0,a[1,4]=alpha[1,4]/a41);$\backslash$
\newline
if(type(a[1,4])!=1,,diag=idmat(4);diag[1,1]=a[1,3]$\ast$a[3,4]$
\ast$a[4,2];$\backslash$
\newline
diag[2,2]=a[3,1]$\ast$a[4,3]$\ast$a[2,4];diag[3,3]=a[3,1]$\ast$a[3,4]$
\ast$a[4,2];$\backslash$
\newline
diag[4,4]=a[3,1]$\ast$a[4,3]$\ast$a[4,2];b=a$\ast
$diag;b=b/content(b);$\backslash$
\newline
db=smith(b);dbb=db[2];fdbb=factor(dbb);$\backslash$
\newline
if(issqfree(b[1,1])==0$||$issqfree(b[2,2])==0$||$issqfree(b[3,3])==0$
||$$\backslash$
\newline
issqfree(b[4,4])==0$||$content([--b[1,1],--b[2,2]])$>$2$||$$\backslash$
\newline
content([--b[2,2],--b[3,3]])$>$2$||$$\backslash$
\newline
(alpha[3,4]$<$4\&\&content([--b[3,3],--b[4,4]])$>$2),,$\backslash$
\newline
detb=db[2]$\ast$db[3]$\ast$db[4];if(detb==1,,fdetb=factor(detb);
$\backslash$
\newline
if(content([--b[1,1]$\ast$--b[2,2]$\ast$--b[3,3]$\ast$--b[4,4],16])$
>$1\&\&$\backslash$
\newline
content([--b[1,1]$\ast$--b[2,2]$\ast$--b[3,3]$\ast$--b[4,4],16])$
<$16\&\&$\backslash$
\newline
fdetb[1,1]==2\&\&mod(fdetb[1,2],2)==mod(1,2),,$\backslash$
\newline
gam=0;for(j=1,matsize(fdetb)[1],if(fdetb[j,1]!=2\&\&$\backslash$
\newline
type(--b[1,1]$\ast$--b[2,2]$\ast$--b[3,3]$\ast
$--b[4,4]/fdetb[j,1])==1\&\&$\backslash$
\newline
mod(fdetb[j,2],2)==mod(0,2),gam=1,));if(gam==1,,for(k=0,1,$\backslash$
\newline
if(k==0\&\&content(--[b[1,1]$\ast$--b[2,2]$\ast$--b[3,3]$\ast
$--b[4,4],2])==1\&\&$\backslash$
\newline
fdetb[1,1]==2\&\&mod(fdetb[1,2],2)==mod(1,2),,$\backslash$
\newline
if(k==1\&\&content(--[b[1,1]$\ast$--b[2,2]$\ast$--b[3,3]$\ast
$--b[4,4],2])==1,$\backslash$
\newline
b=2$\ast$b;db=smith(b);dbb=db[2];fdbb=factor(dbb);gam1=0,gam1=1);$\backslash$
\newline
if(k==1\&\&gam1==1,,detb=db[2]$\ast$db[3]$\ast$db[4];$\backslash$
\newline
fdetb=factor(detb);d=1;for(k=1,matsize(fdetb)[1],$\backslash$
\newline
if(mod(fdetb[k,2],2)==mod(1,2),d=d$\ast$fdetb[k,1],));if(d$<
$=2,et=0,$\backslash$
\newline
fd=factor(d);if(fd[1,1]==2,d1=d/2,d1=d);fd1=factor(d1);$\backslash$
\newline
et=0;for(k=1,matsize(fd1)[1],if(type(b[1,1]/fd1[k,1])==1,$\backslash$
\newline
if(kro(b[1,1]/fd1[k,1],fd1[k,1])==1,,et=et+2$^\wedge$(k--1)),$\backslash$
\newline
if(type(b[2,2]/fd1[k,1])==1,$\backslash$
\newline
if(kro(b[2,2]/fd1[k,1],fd1[k,1])==1,,et=et+2$^\wedge$(k--1)),$\backslash$
\newline
if(kro(d$\ast$b[1,1]$\ast$b[2,2]/fd1[k,1],fd1[k,1])==1,,et=et+2$
^\wedge$(k--1))))));$\backslash$
\newline
hhh=hnr(d,et);if(hhh$>$2$||$hhh==1,,$\backslash$
\newline
n=n+1;pprint("n=",n);pprint(a);pprint(b);pprint(db);$\backslash$
\newline
pprint(fdbb);pprint("d=",d," et=",et," h=",hhh);$\backslash$
\newline
if(n==1,uuu=matrix(1,3,j,k,0);uuu[1,]=[d,et,hhh],$\backslash$
\newline
aaa=0;nnn=matsize(uuu)[1];for(j=1,nnn,if(uuu[j,]==[d,et,hhh],aaa=1,));$
\backslash$
\newline
if(aaa==1,,uuu1=matrix(nnn+1,3,j,k,0);for(t=1,nnn,$\backslash$
\newline
if(uuu[t,1]$<$d$||$(uuu[t,1]==d\&\&uuu[t,2]$<
$et),uuu1[t,]=uuu[t,],$\backslash$
\newline
uuu1[t+1,]=uuu[t,]));for(t=1,nnn+1,if(uuu1[t,]==[0,0,0],$\backslash$
\newline
uuu1[t,]=[d,et,hhh],));uuu=uuu1));))))))))))))))))));$\backslash$
\newline
pprint("uuu=",uuu);pprint("matsize uuu=",matsize(uuu));

\vskip10pt 

\centerline{{\bf Program 25:} funda20.main}

\noindent
$\backslash$$\backslash$program funda20.main
\newline
$\backslash$$\backslash$hyperbolic type, case AII0
\newline
$\backslash$$\backslash$cosh(eta/2)=3/2
\newline
$\backslash$r h2
\newline
$\backslash$r h3
\newline
$\backslash$l;$\backslash$
\newline
n=0;alpha12=0;alpha23=0;alpha34=0;for(alpha14=5,(31.15549442+0.1)$^\wedge
$2,$\backslash$
\newline
fordiv(4$\ast$alpha14,aa,$\backslash$
\newline
if(aa$^\wedge$2$>$4$\ast$alpha14,,alpha13=4+aa;alpha24=4$\ast$alpha14/aa+4;$
\backslash$
\newline
if(alpha13$>$(7+0.1)$^\wedge$2,,alpha=4$\ast$idmat(4);$\backslash$
\newline
alpha[1,3]=alpha13;alpha[3,1]=alpha13;alpha[3,4]=alpha34;$\backslash$
\newline
alpha[4,3]=alpha34;alpha[2,4]=alpha24;alpha[4,2]=alpha24;$\backslash$
\newline
alpha[1,4]=alpha14;alpha[4,1]=alpha14;$\backslash$
\newline
fordiv(alpha[1,3],a13,fordiv(alpha[2,4],a24,fordiv(alpha[1,4],a14,$
\backslash$
\newline
a=--2$\ast$idmat(4);a[1,3]=a13;a[3,1]=alpha[1,3]/a13;a[1,4]=a14;$
\backslash$
\newline
a[4,1]=alpha[1,4]/a14;a[1,3]=a13;a[3,1]=alpha[1,3]/a13;$\backslash$
\newline
a[2,4]=a24;a[4,2]=alpha[2,4]/a24;$\backslash$
\newline
diag=idmat(4);diag[1,1]=a[1,3]$\ast$a[1,4]$\ast$a[4,2];$\backslash$
\newline
diag[2,2]=a[1,3]$\ast$a[4,1]$\ast$a[2,4];diag[3,3]=a[3,1]$\ast$a[1,4]$
\ast$a[4,2];$\backslash$
\newline
diag[4,4]=a[1,3]$\ast$a[4,2]$\ast$a[4,1];b=a$\ast$diag;b=b/content(b);$
\backslash$
\newline
m=m+1;db=smith(b);dbb=db[2];fdbb=factor(dbb);$\backslash$
\newline
if(issqfree(b[1,1])==0$||$issqfree(b[2,2])==0$||$issqfree(b[3,3])==0$||$$
\backslash$
\newline
issqfree(b[4,4])==0$||$content([--b[1,1],--b[2,2]])$>$2$||$$\backslash$
\newline
content([--b[2,2],--b[3,3]])$>$2$||$content([--b[3,3],--b[4,4]])$>$2,,$
\backslash$
\newline
detb=db[2]$\ast$db[3]$\ast$db[4];if(detb==1,,fdetb=factor(detb);$\backslash$
\newline
if(content([--b[1,1]$\ast$--b[2,2]$\ast$--b[3,3]$\ast$--b[4,4],16])$>$1\&\&$
\backslash$
\newline
content([--b[1,1]$\ast$--b[2,2]$\ast$--b[3,3]$\ast$--b[4,4],16])$<$16\&\&$
\backslash$
\newline
fdetb[1,1]==2\&\&mod(fdetb[1,2],2)==mod(1,2),,$\backslash$
\newline
gam=0;for(j=1,matsize(fdetb)[1],if(fdetb[j,1]!=2\&\&$\backslash$
\newline
type(--b[1,1]$\ast$--b[2,2]$\ast$--b[3,3]$\ast$--b[4,4]/fdetb[j,1])==1\&\&$
\backslash$
\newline
mod(fdetb[j,2],2)==mod(0,2),gam=1,));if(gam==1,,for(k=0,1,$\backslash$
\newline
if(k==0\&\&content(--[b[1,1]$\ast$--b[2,2]$\ast
$--b[3,3]$\ast$--b[4,4],2])==1\&\&$
\backslash$
\newline
fdetb[1,1]==2\&\&mod(fdetb[1,2],2)==mod(1,2),,$\backslash$
\newline
if(k==1\&\&content(--[b[1,1]$\ast$--b[2,2]$\ast$--b[3,3]$
\ast$--b[4,4],2])==1,$
\backslash$
\newline
b=2$\ast$b;db=smith(b);dbb=db[2];fdbb=factor(dbb);gam1=0,gam1=1);$\backslash$
\newline
if(k==1\&\&gam1=1,,detb=db[2]$\ast$db[3]$\ast$db[4];fdetb=factor(detb);$
\backslash$
\newline
d=1;for(k=1,matsize(fdetb)[1],$\backslash$
\newline
if(mod(fdetb[k,2],2)==mod(1,2),d=d$\ast$fdetb[k,1],));if(d$<$=2,et=0,$
\backslash$
\newline
fd=factor(d);if(fd[1,1]==2,d1=d/2,d1=d);fd1=factor(d1);$\backslash$
\newline
et=0;for(k=1,matsize(fd1)[1],if(type(b[1,1]/fd1[k,1])==1,$\backslash$
\newline
if(kro(b[1,1]/fd1[k,1],fd1[k,1])==1,,et=et+2$^\wedge$(k--1)),$\backslash$
\newline
if(type(b[2,2]/fd1[k,1])==1,$\backslash$
\newline
if(kro(b[2,2]/fd1[k,1],fd1[k,1])==1,,et=et+2$^\wedge$(k--1)),$\backslash$
\newline
if(kro(d$\ast$b[1,1]$\ast$b[2,2]/fd1[k,1],fd1[k,1])==1,,et=et+2$^\wedge
$(k--1))))));$\backslash$
\newline
hhh=hnr(d,et);if(hhh$>$2$||$hhh==1,,$\backslash$
\newline
n=n+1;pprint("n=",n);pprint(a);pprint(b);pprint(db);$\backslash$
\newline
pprint(fdbb);pprint("d=",d," et=",et," h=",hhh);$\backslash$
\newline
if(n==1,u=matrix(1,3,j,k,0);u[1,]=[d,et,hhh],aaa=0;$\backslash$
\newline
nnn=matsize(u)[1];for(j=1,nnn,if(u[j,]==[d,et,hhh],aaa=1,));$\backslash$
\newline
if(aaa==1,,u1=matrix(nnn+1,3,j,k,0);for(t=1,nnn,$\backslash$
\newline
if(u[t,1]$<$d$||$(u[t,1]==d\&\&u[t,2]$<$et),u1[t,]=u[t,],u1[t+1,]=u[t,]));$
\backslash$
\newline
for(t=1,nnn+1,if(u1[t,]==[0,0,0],u1[t,]=[d,et,hhh],));u=u1));$\backslash$
\newline
)))))))))))))));pprint("u=",u);pprint("matsize u=",matsize(u));

\vskip10pt 

\centerline{{\bf Program 26:} funda3.main}

\noindent
$\backslash$$\backslash$funda3.mai
\newline
$\backslash$$\backslash$hyperbolic type, case AIII
\newline
$\backslash$$\backslash$cosh(eta/2)=3/2
\newline
$\backslash$r h2
\newline
$\backslash$r h3
\newline
$\backslash$l;$\backslash$
\newline
n=0;al12=al23=al34=al45=0;for(al15=5,(68.1815011826+0.1)$^\wedge$2,$
\backslash$
\newline
for(al13=5,(7+0.1)$^\wedge$2,for(al35=al13,(7+0.1)$^\wedge$2,$
\backslash$
\newline
if(issquare(q=al13$\ast$al35$\ast$al15)==0,,d=(al13+al35+al15--4+isqrt(q))$
\ast$4;$\backslash$
\newline
if(type(al14=d/(al35--4))!=1$||$al14$<$=(31.15549442--0.1)$
^\wedge$2$\backslash$
\newline
$||$type(al25=d/(al13--4))!=1$||$al25$<$=(31.15549442--0.1)$
^\wedge$2$||$$\backslash$
\newline
type(al24=(al13$\ast$al35+4$\ast$al15+4$\ast$isqrt(q))$\ast$4/((al35--4)$
\ast$(al13--4)))!=1$||$$\backslash$
\newline
issquare(q1=al13$\ast$al35$\ast$al25$\ast$al24$\ast$al14)==0,,$
\backslash$
\newline
al=idmat(5)$\ast$4;al[1,5]=al15;al[5,1]=al15;al[1,3]=al13;al[3,1]=al13;$
\backslash$
\newline
al[1,4]=al14;al[4,1]=al14;al[2,4]=al24;al[4,2]=al24;$\backslash$
\newline
al[2,5]=al25;al[5,2]=al25;al[3,5]=al35;al[5,3]=al35;$\backslash$
\newline
fordiv(al13,a13,fordiv(al35,a35,a51=isqrt(q)/a13/a35;$\backslash$
\newline
if(a51==0,a15==0,if(type(a15=al15/a51)!=1,,fordiv(al14,a14,$\backslash$
\newline
fordiv(al24,a24,a52=isqrt(q1)/a13/a35/a24$\ast$a14/al14;$\backslash$
\newline
if(type(a25=al25/a52)!=1,,a=idmat(5)$\ast$--2;$\backslash$
\newline
a[1,3]=a13;a[3,1]=al13/a13;a[1,4]=a14;a[4,1]=al14/a14;$\backslash$
\newline
a[1,5]=a15;a[5,1]=a51;a[2,4]=a24;a[4,2]=al24/a24;$\backslash$
\newline
a[2,5]=a25;a[5,2]=a52;a[3,5]=a35;a[5,3]=al35/a35;$\backslash$
\newline
diag=idmat(5);diag[1,1]=a[1,4]$\ast$a[1,3]$\ast$a[4,2]$\ast$a[2,5];$
\backslash$
\newline
diag[2,2]=a[4,1]$\ast$a[1,3]$\ast$a[2,4]$\ast$a[2,5];$\backslash$
\newline
diag[3,3]=a[1,4]$\ast$a[3,1]$\ast$a[4,2]$\ast$a[2,5];$\backslash$
\newline
diag[4,4]=a[4,1]$\ast$a[1,3]$\ast$a[4,2]$\ast$a[2,5];$\backslash$
\newline
diag[5,5]=a[4,1]$\ast$a[1,3]$\ast$a[2,4]$\ast$a[5,2];$\backslash$
\newline
b=a$\ast$diag;b=b/content(b);db=smith(b);dbb=db[3];fdbb=factor(dbb);$
\backslash$
\newline
if(issqfree(b[1,1])==0$||$issqfree(b[2,2])==0$||$issqfree(b[3,3])==0$
||$$\backslash$
\newline
issqfree(b[4,4])==0$||$issqfree(b[5,5])==0$||$$\backslash$
\newline
content([--b[1,1],--b[2,2]])$>$2$||$content([--b[2,2],--b[3,3]])$
>$2$||$$\backslash$
\newline
content([--b[3,3],--b[4,4]])$>$2$||$content([--b[4,4],--b[5,5]])$
>$2,,$\backslash$
\newline
detb=db[3]$\ast$db[4]$\ast$db[5];if(detb==1,,fdetb=factor(detb);$
\backslash$
\newline
if(content([--b[1,1]$\ast$--b[2,2]$\ast$--b[3,3]$\ast$--b[4,4]$
\ast$--b[5,5],2])$>$1\&\&$\backslash$
\newline
content([--b[1,1]$\ast$--b[2,2]$\ast$--b[3,3]$\ast$--b[4,4]$
\ast$--b[5,5],32])$<$32\&\&$\backslash$
\newline
fdetb[1,1]==2\&\&mod(fdetb[1,2],2)==mod(1,2),,$\backslash$
\newline
gam=0;for(j=1,matsize(fdetb)[1],if(fdetb[j,1]!=2\&\&$\backslash$
\newline
type(--b[1,1]$\ast$--b[2,2]$\ast$--b[3,3]$\ast$--b[4,4]$
\ast$--b[5,5]/fdetb[j,1])==1\&\&$\backslash$
\newline
mod(fdetb[j,2],2)==mod(0,2),gam=1,));if(gam==1,,for(k=0,1,$\backslash$
\newline
if(k==0\&\&content(--[b[1,1]$\ast$--b[2,2]$\ast$--b[3,3]$\ast$--b[4,4]$
\ast$--b[5,5],2])==1\&\&$\backslash$
\newline
fdetb[1,1]==2\&\&mod(fdetb[1,2],2)==mod(1,2),,$\backslash$
\newline
if(k==1\&\&content(--[b[1,1]$\ast$--b[2,2]$\ast$--b[3,3]$\ast$--b[4,4]$
\ast$--b[5,5],2])==1,$\backslash$
\newline
b=2$\ast$b;db=smith(b);dbb=db[3];fdbb=factor(dbb);gam1=0,gam1=1);$
\backslash$
\newline
if(k==1\&\&gam1==1,,detb=db[3]$\ast$db[4]$\ast$db[5];fdetb=factor(detb);$
\backslash$
\newline
d=1;for(k=1,matsize(fdetb)[1],$\backslash$
\newline
if(mod(fdetb[k,2],2)==mod(1,2),d=d$\ast$fdetb[k,1],));if(d$
<$=2,et=0,$\backslash$
\newline
fd=factor(d);if(fd[1,1]==2,d1=d/2,d1=d);fd1=factor(d1);$\backslash$
\newline
et=0;for(k=1,matsize(fd1)[1],if(type(b[1,1]/fd1[k,1])==1,$
\backslash$
\newline
if(kro(b[1,1]/fd1[k,1],fd1[k,1])==1,,et=et+2$^\wedge$(k--1)),$
\backslash$
\newline
if(type(b[2,2]/fd1[k,1])==1,$\backslash$
\newline
if(kro(b[2,2]/fd1[k,1],fd1[k,1])==1,,et=et+2$^\wedge$(k--1)),$\backslash$
\newline
if(kro(d$\ast$b[1,1]$\ast$b[2,2]/fd1[k,1],fd1[k,1])==1,,et=et+2$
^\wedge$(k--1))))));$\backslash$
\newline
hhh=hnr(d,et);if(hhh$>$2$||$hhh==1,,$\backslash$
\newline
n=n+1;pprint("n=",n);pprint(a);pprint(b);pprint(db);$\backslash$
\newline
pprint(fdbb);pprint("d=",d," et=",et," h=",hhh);$\backslash$
\newline
if(n==1,u=matrix(1,3,j,k,0);u[1,]=[d,et,hhh],aaa=0;nnn=matsize(u)[1];$
\backslash$
\newline
for(j=1,nnn,if(u[j,]==[d,et,hhh],aaa=1,));$\backslash$
\newline
if(aaa==1,,u1=matrix(nnn+1,3,j,k,0);for(t=1,nnn,$\backslash$
\newline
if(u[t,1]$<$d$||$(u[t,1]==d\&\&u[t,2]$<$et),u1[t,]=u[t,],u1[t+1,]=u[t,]));$
\backslash$
\newline
for(t=1,nnn+1,if(u1[t,]==[0,0,0],u1[t,]=[d,et,hhh],));u=u1));$\backslash$
\newline
))))))))))))))))))));pprint("u=",u);pprint("matsize u=",matsize(u));

\vskip10pt 

\centerline{{\bf Program 27:} fundb1.main}

\noindent
$\backslash$$\backslash$fundb1.mai
\newline
$\backslash$$\backslash$hyperbolic type, case BI
\newline
$\backslash$$\backslash$cosh(eta/2)=3/2
\newline
$\backslash$r h2
\newline
$\backslash$r h3
\newline
$\backslash$l;$\backslash$
\newline
n=0;$\backslash$
\newline
rh4=[1.191398091,1.095286061,1.022736500,0.962423650];$\backslash$
\newline
for(alpha12=1,4,for(alpha23=5,22$^\wedge$2,for(alpha13=5,alpha23,$
\backslash$
\newline
if(sqrt(alpha13)/8--5/4$<$1+0.000001,ach1=0,ach1=acosh(sqrt(alpha13)/8--5/4));$
\backslash$
\newline
if(sqrt(alpha23)/8--5/4$<
$1+0.000001,ach2=0,ach2=acosh(sqrt(alpha23)/8--5/4));$\backslash$
\newline
if(ach1+ach2$>$rh4[alpha12]+0.1$||$alpha23==0$||$$\backslash$
\newline
issquare(alpha12$\ast$alpha23$\ast$alpha13)!=1$||$$\backslash$
\newline
--8+2$\ast$isqrt(alpha12$\ast$alpha23$\ast$alpha13)+2$\ast$alpha12+2$
\ast$alpha13+2$\ast$alpha23$<$=0,,$\backslash$
\newline
alpha=4$\ast$idmat(3);alpha[1,2]=alpha12;alpha[2,1]=alpha12;$\backslash$
\newline
alpha[2,3]=alpha23;alpha[3,2]=alpha23;alpha[1,3]=alpha13;$\backslash$
\newline
alpha[3,1]=alpha13;dalpha=$\backslash$
\newline
--8+2$\ast$isqrt(alpha12$\ast$alpha23$\ast$alpha13)+2$\ast$alpha12+2$
\ast$alpha13+2$\ast$alpha23;$\backslash$
\newline
a=--2$\ast$idmat(3);$\backslash$
\newline
fordiv(alpha[1,2],a12,fordiv(alpha[2,3],a23,fordiv(alpha[1,3],a13,$
\backslash$
\newline
a21=alpha[1,2]/a12;a32=alpha[2,3]/a23;a31=alpha[1,3]/a13;$\backslash$
\newline
if(a12$\ast$a23$\ast$a31!=a21$\ast$a13$\ast$a32,,$\backslash$
\newline
a[1,2]=a12;a[2,1]=a21;a[2,3]=a23;a[3,2]=a32;a[1,3]=a13;a[3,1]=a31;$
\backslash$
\newline
dd=idmat(3);dd[1,1]=a13$\ast$a32;dd[2,2]=a23$\ast$a31;dd[3,3]=a31$\ast$a32;$
\backslash$
\newline
b=a$\ast$dd;b=b/content(b);db=smith(b);fdb=factor(db[1]);$\backslash$
\newline
detb=db[1]$\ast$db[2]$\ast$db[3];fdetb=factor(detb);$\backslash$
\newline
if(issqfree(b[1,1])==0$||$issqfree(b[2,2])==0$||$issqfree(b[3,3])==0$
||$$\backslash$
\newline
(content([--b[1,1]$\ast$--b[2,2]$\ast$--b[3,3],2])$>$1\&\&$\backslash$
\newline
content([--b[1,1]$\ast$--b[2,2]$\ast$--b[3,3],8])$<$8\&\&$\backslash$
\newline
fdetb[1,1]==2\&\&mod(fdetb[1,2],2)==mod(1,2))$||$$\backslash$
\newline
(alpha[1,2]$<$4\&\&content([--b[1,1],--b[2,2]])$>$2),,$\backslash$
\newline
gam=0;for(j=1,matsize(fdetb)[1],if(fdetb[j,1]!=2\&\&$\backslash$
\newline
type(--b[1,1]$\ast$--b[2,2]$\ast$--b[3,3]/fdetb[j,1])==1\&\&$\backslash$
\newline
mod(fdetb[j,2],2)==mod(0,2),gam=1,));if(gam==1,,$\backslash$
\newline
if(alpha[1,2]$<$4,b=b;d=1;for(k=1,matsize(fdetb)[1],$\backslash$
\newline
if(mod(fdetb[k,2],2)==mod(1,2),d=d$\ast$fdetb[k,1],));$\backslash$
\newline
if(alpha[1,2]==1\&\&mod(d,3)!=mod(2,3),,if(d$<$=2,et=0,$\backslash$
\newline
fd=factor(d);if(fd[1,1]==2,d1=d/2,d1=d);fd1=factor(d1);$\backslash$
\newline
et=0;for(k=1,matsize(fd1)[1],if(alpha[1,2]==1,$\backslash$
\newline
if(kro(3$\ast$d/fd1[k,1],fd1[k,1])==1,,et=et+2$^\wedge
$(k--1)),if(alpha[1,2]==2,$\backslash$
\newline
if(kro(d/fd1[k,1],fd1[k,1])==1,,et=et+2$^\wedge$(k--1)),$
\backslash$
\newline
if(alpha[1,2]==3,if(fd1[k,1]==3,,$\backslash$
\newline
if(kro(3$\ast$d/fd1[k,1],fd1[k,1])==1,,et=et+2$^\wedge$(k--1))),))));$
\backslash$
\newline
hhh=hnr(d,et);if(hhh$>$2$||$hhh==1,,n=n+1;pprint("n=",n);pprint(b);$
\backslash$
\newline
pprint(db);pprint(fdb);pprint("d=",d," et=",et," h=",hhh);$\backslash$
\newline
if(n==1,u=matrix(1,3,j,k,0);u[1,]=[d,et,hhh],aaa=0;nnn=matsize(u)[1];$
\backslash$
\newline
for(j=1,nnn,if(u[j,]==[d,et,hhh],aaa=1,));$\backslash$
\newline
if(aaa==1,,u1=matrix(nnn+1,3,j,k,0);for(t=1,nnn,$\backslash$
\newline
if(u[t,1]$<$d$||$(u[t,1]==d\&\&u[t,2]$<$et),u1[t,]=u[t,],u1[t+1,]=u[t,]));$
\backslash$
\newline
for(t=1,nnn+1,if(u1[t,]==[0,0,0],u1[t,]=[d,et,hhh],));u=u1))))),$
\backslash$
\newline
if(alpha[1,2]==4,d2=1;for(k=1,matsize(fdetb)[1],$\backslash$
\newline
if(mod(fdetb[k,2],2)==mod(0,2)\&\&fdetb[k,1]!=2,d2=d2$\ast$fdetb[k,1],));$
\backslash$
\newline
d2=2$\ast$d2/content([2$\ast$d2,--b[1,1]$\ast$--b[2,2]$\ast
$--b[3,3]]);fordiv(d2,t,$\backslash$
\newline
if(mod(d2,2)==mod(0,2)\&\&fdetb[1,1]==2\&\&mod(fdetb[1,2],2)==mod(1,2)\&\&$
\backslash$
\newline
mod(t,2)==mod(1,2),,bb=b$\ast$t;dbb=smith(bb);fdbb=factor(dbb[1]);$\backslash$
\newline
detbb=dbb[1]$\ast$dbb[2]$\ast$dbb[3];fdetbb=factor(detbb);d=1;$\backslash$
\newline
for(k=1,matsize(fdetbb)[1],$\backslash$
\newline
if(mod(fdetbb[k,2],2)==mod(1,2),d=d$\ast$fdetbb[k,1],));if(d$<$=2,et=0,$
\backslash$
\newline
fd=factor(d);if(fd[1,1]==2,d1=d/2,d1=d);fd1=factor(d1);et=0;$\backslash$
\newline
for(k=1,matsize(fd1)[1],if(kro(--d/fd1[k,1],fd1[k,1])==1,,et=et+2$
^\wedge$(k--1))));$\backslash$
\newline
hhh=hnr(d,et);if(hhh$>$2$||$hhh==1,,$\backslash$
\newline
n=n+1;pprint("n=",n);pprint(bb);pprint(dbb);$\backslash$
\newline
pprint(fdbb);pprint("d=",d," et=",et," h=",hhh);$\backslash$
\newline
if(n==1,u=matrix(1,3,j,k,0);u[1,]=[d,et,hhh],aaa=0;nnn=matsize(u)[1];$
\backslash$
\newline
for(j=1,nnn,if(u[j,]==[d,et,hhh],aaa=1,));$\backslash$
\newline
if(aaa==1,,u1=matrix(nnn+1,3,j,k,0);for(t=1,nnn,$\backslash$
\newline
if(u[t,1]$<$d$||$(u[t,1]==d\&\&u[t,2]$<$et),u1[t,]=u[t,],u1[t+1,]=u[t,]));$
\backslash$
\newline
for(t=1,nnn+1,if(u1[t,]==[0,0,0],u1[t,]=[d,et,hhh],));u=u1))))),$\backslash$
\newline
))))))))))));$\backslash$
\newline
pprint("u=",u);pprint("matsize u=",matsize(u));

\vskip10pt 

\centerline{{\bf Program 28:} fundb21.main}

\noindent
$\backslash$$\backslash$fundb21.mai
\newline
$\backslash$$\backslash$hyperbolic type, case BII\_1
\newline
$\backslash$$\backslash$cosh(eta/2)=3/2
\newline
$\backslash$r h2;$\backslash$
\newline
$\backslash$r h3;$\backslash$
\newline
$\backslash$l;$\backslash$
\newline
n=0;alpha12=alpha23=0;for(alpha24=5,22$^\wedge$2,$\backslash$
\newline
if(sqrt(alpha24)/8--5/4$<
$1+0.000001,ach2=0,ach2=acosh(sqrt(alpha24)/8--5/4));$\backslash$
\newline
vv=min((8$\ast$(cosh(1.429914377--ach2)+5/4))$^\wedge$2+0.1,22$
^\wedge$2);$\backslash$
\newline
for(alpha14=5,vv,for(alpha34=alpha14,vv,$\backslash$
\newline
if(issquare(discrim=alpha14$\ast$alpha34$\ast$(alpha14$
\ast$alpha34+(alpha24--4)$\ast$$\backslash$
\newline
(alpha14+alpha24+alpha34--4)))!=1,,$\backslash$
\newline
if(type(u=(2$\ast$alpha14$\ast$alpha34+2$\ast
$isqrt(discrim))/(alpha24--4))!=1,,$\backslash$
\newline
if(type(alpha13=u$^\wedge$2/(alpha14$\ast$alpha34))!=1,,alpha=4$
\ast$idmat(4);$\backslash$
\newline
alpha[1,3]=alpha13;alpha[3,1]=alpha13;alpha[3,4]=alpha34;alpha[4,3]=alpha34;$
\backslash$
\newline
alpha[2,4]=alpha24;alpha[4,2]=alpha24;alpha[1,4]=alpha14;alpha[4,1]=alpha14;$
\backslash$
\newline
fordiv(alpha[3,4],a34,fordiv(alpha[1,3],a13,fordiv(alpha[2,4],a24,$
\backslash$
\newline
if(type(a41=u/(a13$\ast$a34))!=1,,a=--2$\ast$idmat(4);$\backslash$
\newline
a[3,4]=a34;a[4,3]=alpha[3,4]/a34;a[2,4]=a24;a[4,2]=alpha[2,4]/a24;$
\backslash$
\newline
a[1,3]=a13;a[3,1]=alpha[1,3]/a13;a[4,1]=a41;$\backslash$
\newline
if(a41==0,a[1,4]=0,a[1,4]=alpha[1,4]/a41);if(type(a[1,4])!=1,,$
\backslash$
\newline
diag=idmat(4);diag[1,1]=a[1,3]$\ast$a[3,4]$\ast$a[4,2];$
\backslash$
\newline
diag[2,2]=a[3,1]$\ast$a[4,3]$\ast$a[2,4];diag[3,3]=a[3,1]$
\ast$a[3,4]$\ast$a[4,2];$\backslash$
\newline
diag[4,4]=a[3,1]$\ast$a[4,3]$\ast$a[4,2];b=a$\ast
$diag;b=b/content(b);$\backslash$
\newline
db=smith(b);dbb=db[2];fdbb=factor(dbb);$\backslash$
\newline
if(issqfree(b[1,1])==0$||$issqfree(b[2,2])==0$||$issqfree(b[3,3])==0$
||$$\backslash$
\newline
issqfree(b[4,4])==0$||$content([--b[1,1],--b[2,2]])$>$2$||$$\backslash$
\newline
content([--b[2,2],--b[3,3]])$>$2$||$$\backslash$
\newline
(alpha[3,4]==2\&\&content([--b[3,3],--b[4,4]])$>$1)$||$$\backslash$
\newline
(alpha[3,4]==3\&\&content([--b[3,3],--b[4,4]])$>$2),,$\backslash$
\newline
detb=db[2]$\ast$db[3]$\ast$db[4];if(detb==1,,fdetb=factor(detb);$
\backslash$
\newline
if(content([--b[1,1]$\ast$--b[2,2]$\ast$--b[3,3]$\ast$--b[4,4],16])$>
$1\&\&$\backslash$
\newline
content([--b[1,1]$\ast$--b[2,2]$\ast$--b[3,3]$\ast$--b[4,4],16])$<
$16\&\&$\backslash$
\newline
fdetb[1,1]==2\&\&mod(fdetb[1,2],2)==mod(1,2),,$\backslash$
\newline
gam=0;for(j=1,matsize(fdetb)[1],if(fdetb[j,1]!=2\&\&$\backslash$
\newline
type(--b[1,1]$\ast$--b[2,2]$\ast$--b[3,3]$\ast
$--b[4,4]/fdetb[j,1])==1\&\&$
\backslash$
\newline
mod(fdetb[j,2],2)==mod(0,2),gam=1,));if(gam==1,,for(k=0,1,$\backslash$
\newline
if(k==0\&\&content(--[b[1,1]$\ast$--b[2,2]$\ast$--b[3,3]$\ast
$--b[4,4],2])==1\&\&$\backslash$
\newline
fdetb[1,1]==2\&\&mod(fdetb[1,2],2)==mod(1,2),,$\backslash$
\newline
if(k==1\&\&content(--[b[1,1]$\ast$--b[2,2]$\ast$--b[3,3]$\ast
$--b[4,4],2])==1,$\backslash$
\newline
b=2$\ast$b;db=smith(b);dbb=db[2];fdbb=factor(dbb);gam1=0,gam1=1);$
\backslash$
\newline
if(k==1\&\&gam1==1,,detb=db[2]$\ast$db[3]$\ast$db[4];fdetb=factor(detb);$
\backslash$
\newline
d=1;for(k=1,matsize(fdetb)[1],if(mod(fdetb[k,2],2)==mod(1,2),$\backslash$
\newline
d=d$\ast$fdetb[k,1],));if(d$<$=2,et=0,$\backslash$
\newline
fd=factor(d);if(fd[1,1]==2,d1=d/2,d1=d);fd1=factor(d1);$\backslash$
\newline
et=0;for(k=1,matsize(fd1)[1],if(type(b[1,1]/fd1[k,1])==1,$\backslash$
\newline
if(kro(b[1,1]/fd1[k,1],fd1[k,1])==1,,et=et+2$^\wedge$(k--1)),$\backslash$
\newline
if(type(b[2,2]/fd1[k,1])==1,$\backslash$
\newline
if(kro(b[2,2]/fd1[k,1],fd1[k,1])==1,,et=et+2$^\wedge$(k--1)),$\backslash$
\newline
if(kro(d$\ast$b[1,1]$\ast$b[2,2]/fd1[k,1],fd1[k,1])==1,,et=et+2$
^\wedge$(k--1))))));$\backslash$
\newline
hhh=hnr(d,et);if(hhh$>$2$||$hhh==1,,$\backslash$
\newline
n=n+1;pprint("n=",n);pprint(a);pprint(b);pprint(db);$\backslash$
\newline
pprint(fdbb);pprint("d=",d," et=",et," h=",hhh);$\backslash$
\newline
if(n==1,uuu=matrix(1,3,j,k,0);uuu[1,]=[d,et,hhh],$\backslash$
\newline
aaa=0;nnn=matsize(uuu)[1];for(j=1,nnn,if(uuu[j,]==[d,et,hhh],aaa=1,));$
\backslash$
\newline
if(aaa==1,,uuu1=matrix(nnn+1,3,j,k,0);for(t=1,nnn,$\backslash$
\newline
if(uuu[t,1]$<$d$||$(uuu[t,1]==d\&\&uuu[t,2]$<$et),uuu1[t,]=uuu[t,],$
\backslash$
\newline
uuu1[t+1,]=uuu[t,]));for(t=1,nnn+1,if(uuu1[t,]==[0,0,0],$\backslash$
\newline
uuu1[t,]=[d,et,hhh],));uuu=uuu1));)))))))))))))))))));$\backslash$
\newline
pprint("uuu=",uuu);pprint("matsize uuu=",matsize(uuu));

\vskip10pt 

\centerline{{\bf Program 29:} fundb22.main}

\noindent
$\backslash$$\backslash$fundb22.mai
\newline
$\backslash$$\backslash$hyperbolic type, case BII\_2
\newline
$\backslash$$\backslash$cos(eta/2)=3/2
\newline
$\backslash$r h2;$\backslash$
\newline
$\backslash$r h3;$\backslash$
\newline
$\backslash$l;$\backslash$
\newline
n=0;alpha12=alpha34=0;for(alpha23=5,22$^\wedge$2,$\backslash$
\newline
if(sqrt(alpha23)/8--5/4$<$1+0.000001,ach2=0,$\backslash$
\newline
ach2=acosh(sqrt(alpha23)/8--5/4));$\backslash$
\newline
vv=min((8$\ast$(cosh(1.429914377--ach2)+5/4))$^\wedge$2+0.1,22$
^\wedge$2);$\backslash$
\newline
for(alpha13=5,vv,for(alpha24=alpha13,vv,$\backslash$
\newline
if(issquare(discrim=alpha13$\ast$alpha23$\ast$alpha24$\ast$$
\backslash$
\newline
(alpha13$\ast$alpha23$\ast$alpha24--(alpha23--4)$\ast$((alpha13--4)$
\ast$(alpha24--4)--$\backslash$
\newline
4$\ast$alpha23)))!=1,,forstep(tau=--1,1,2,$\backslash$
\newline
if(type(u=(alpha13$\ast$alpha23$\ast$alpha24+tau$\ast
$isqrt(discrim))/(alpha23--4))!=1,,$\backslash$
\newline
if(type(alpha14=u$^\wedge$2/(alpha13$\ast$alpha23$\ast$alpha24))!=1$
||$alpha14$<$=4,,$\backslash$
\newline
alpha=4$\ast$idmat(4);$\backslash$
\newline
alpha[1,3]=alpha13;alpha[3,1]=alpha13;alpha[1,4]=alpha14;alpha[4,1]=alpha14;$
\backslash$
\newline
alpha[2,3]=alpha23;alpha[3,2]=alpha23;alpha[2,4]=alpha24;alpha[4,2]=alpha24;$
\backslash$
\newline
fordiv(alpha[1,3],a13,fordiv(alpha[2,3],a32,fordiv(alpha[2,4],a24,$
\backslash$
\newline
if(type(a41=u/(a13$\ast$a32$\ast$a24))!=1,,a=--2$\ast$idmat(4);$
\backslash$
\newline
a[1,3]=a13;a[3,1]=alpha[1,3]/a13;a[3,2]=a32;a[2,3]=alpha[2,3]/a32;$
\backslash$
\newline
a[2,4]=a24;a[4,2]=alpha[2,4]/a24;a[4,1]=a41;a[1,4]=alpha[1,4]/a41;$
\backslash$
\newline
if(type(a[1,4])!=1,,diag=idmat(4);diag[1,1]=a[1,3]$\ast$a[1,4]$\ast
$a[3,2];$\backslash$
\newline
diag[2,2]=a[3,1]$\ast$a[1,4]$\ast$a[2,3];diag[3,3]=a[3,1]$\ast$a[1,4]$\ast
$a[3,2];$\backslash$
\newline
diag[4,4]=a[1,3]$\ast$a[4,1]$\ast$a[3,2];b=a$\ast$diag;b=b/content(b);$
\backslash$
\newline
db=smith(b);dbb=db[2];fdbb=factor(dbb);$\backslash$
\newline
if(issqfree(b[1,1])==0$||$issqfree(b[2,2])==0$||$issqfree(b[3,3])==0$
||$$\backslash$
\newline
issqfree(b[4,4])==0$||$content([--b[1,1],--b[2,2]])$>$2$||$$\backslash$
\newline
content([--b[3,3],--b[4,4]])$>$2,,detb=db[2]$\ast$db[3]$\ast$db[4];$
\backslash$
\newline
if(detb==1,,fdetb=factor(detb);$\backslash$
\newline
if(content([--b[1,1]$\ast$--b[2,2]$\ast$--b[3,3]$\ast$--b[4,4],16])$>
$1\&\&$\backslash$
\newline
content([--b[1,1]$\ast$--b[2,2]$\ast$--b[3,3]$\ast$--b[4,4],16])$<
$16\&\&$\backslash$
\newline
fdetb[1,1]==2\&\&mod(fdetb[1,2],2)==mod(1,2),,$\backslash$
\newline
gam=0;for(j=1,matsize(fdetb)[1],if(fdetb[j,1]!=2\&\&$\backslash$
\newline
type(--b[1,1]$\ast$--b[2,2]$\ast$--b[3,3]$\ast$--b[4,4]/fdetb[j,1])==1\&\&$
\backslash$
\newline
mod(fdetb[j,2],2)==mod(0,2),gam=1,));if(gam==1,,for(k=0,1,$\backslash$
\newline
if(k==0\&\&content(--[b[1,1]$\ast$--b[2,2]$\ast$--b[3,3]$\ast
$--b[4,4],2])==1\&\&$\backslash$
\newline
fdetb[1,1]==2\&\&mod(fdetb[1,2],2)==mod(1,2),,$\backslash$
\newline
if(k==1\&\&content(--[b[1,1]$\ast$--b[2,2]$\ast$--b[3,3]$\ast
$--b[4,4],2])==1,$\backslash$
\newline
b=2$\ast$b;db=smith(b);dbb=db[2];fdbb=factor(dbb);gam1=0,gam1=1);$
\backslash$
\newline
if(k==1\&\&gam1==1,,detb=db[2]$\ast$db[3]$\ast$db[4];fdetb=factor(detb);$
\backslash$
\newline
d=1;for(k=1,matsize(fdetb)[1],if(mod(fdetb[k,2],2)==mod(1,2),$
\backslash$
\newline
d=d$\ast$fdetb[k,1],));if(d$<$=2,et=0,fd=factor(d);$\backslash$
\newline
if(fd[1,1]==2,d1=d/2,d1=d);fd1=factor(d1);et=0;$\backslash$
\newline
for(k=1,matsize(fd1)[1],if(type(b[1,1]/fd1[k,1])==1,$\backslash$
\newline
if(kro(b[1,1]/fd1[k,1],fd1[k,1])==1,,et=et+2$^\wedge$(k--1)),$
\backslash$
\newline
if(type(b[2,2]/fd1[k,1])==1,$\backslash$
\newline
if(kro(b[2,2]/fd1[k,1],fd1[k,1])==1,,et=et+2$^\wedge$(k--1)),$
\backslash$
\newline
if(kro(d$\ast$b[1,1]$\ast$b[2,2]/fd1[k,1],fd1[k,1])==1,,et=et+2$
^\wedge$(k--1))))));$\backslash$
\newline
hhh=hnr(d,et);if(hhh$>$2$||$hhh==1,,$\backslash$
\newline
n=n+1;pprint("n=",n);pprint(a);pprint(b);pprint(db);$\backslash$
\newline
pprint(fdbb);pprint("d=",d," et=",et," h=",hhh);$\backslash$
\newline
if(n==1,uuu=matrix(1,3,j,k,0);uuu[1,]=[d,et,hhh],$\backslash$
\newline
aaa=0;nnn=matsize(uuu)[1];$\backslash$
\newline
for(j=1,nnn,if(uuu[j,]==[d,et,hhh],aaa=1,));$\backslash$
\newline
if(aaa==1,,uuu1=matrix(nnn+1,3,j,k,0);for(t=1,nnn,$\backslash$
\newline
if(uuu[t,1]$<$d$||$(uuu[t,1]==d\&\&uuu[t,2]$<$et),uuu1[t,]=uuu[t,],$
\backslash$
\newline
uuu1[t+1,]=uuu[t,]));for(t=1,nnn+1,if(uuu1[t,]==[0,0,0],$\backslash$
\newline
uuu1[t,]=[d,et,hhh],));uuu=uuu1));))))))))))))))))))));$\backslash$
\newline
pprint("uuu=",uuu);pprint("matsize uuu=",matsize(uuu));

\vskip10pt 

\centerline{{\bf Program 30:} fundb3.main} 

\noindent
$\backslash$$\backslash$fundb3.mai
\newline
$\backslash$$\backslash$hyperbolic type, case BIII
\newline
$\backslash$$\backslash$cosh(eta/2)=3/2
\newline
$\backslash$r h2
\newline
$\backslash$r h3
\newline
$\backslash$l;$\backslash$
\newline
n=0;alpha12=alpha23=0;for(alpha14=5,(31.70820393)$^\wedge$2+0.1,$
\backslash$
\newline
if(sqrt(alpha14)/8--5/4$<$1+0.000001,ach1=0,$\backslash$
\newline
ach1=acosh(sqrt(alpha14)/8--5/4));$\backslash$
\newline
vv=min((8$\ast$(cosh(1.9248473002--ach1)+5/4))$^\wedge$2+0.1,(31.70820393)$
^\wedge$2+0.1);$\backslash$
\newline
for(alpha34=alpha14,vv,for(alpha13=5,(7+0.1)$^\wedge$2,$\backslash$
\newline
if(issquare(u=alpha13$\ast$alpha34$\ast$alpha14)!=1,,$\backslash$
\newline
if(type(alpha24=4+4$\ast$(alpha14+alpha34+isqrt(u))/(alpha13--4))!=1,,$
\backslash$
\newline
alpha=4$\ast$idmat(4);alpha[1,3]=alpha13;alpha[3,1]=alpha13;$\backslash$
\newline
alpha[3,4]=alpha34;alpha[4,3]=alpha34;$\backslash$
\newline
alpha[2,4]=alpha24;alpha[4,2]=alpha24;$\backslash$
\newline
alpha[1,4]=alpha14;alpha[4,1]=alpha14;$\backslash$
\newline
fordiv(alpha[3,4],a34,fordiv(alpha[1,3],a13,fordiv(alpha[2,4],a24,$
\backslash$
\newline
if(type(a41=isqrt(u)/(a13$\ast$a34))!=1,,a=--2$\ast$idmat(4);$
\backslash$
\newline
a[3,4]=a34;a[4,3]=alpha[3,4]/a34;a[2,4]=a24;a[4,2]=alpha[2,4]/a24;$
\backslash$
\newline
a[1,3]=a13;a[3,1]=alpha[1,3]/a13;a[4,1]=a41;$\backslash$
\newline
if(a41==0,a[1,4]=0,a[1,4]=alpha[1,4]/a41);if(type(a[1,4])!=1,,$
\backslash$
\newline
diag=idmat(4);diag[1,1]=a[1,3]$\ast$a[3,4]$\ast$a[4,2];$\backslash$
\newline
diag[2,2]=a[3,1]$\ast$a[4,3]$\ast$a[2,4];diag[3,3]=a[3,1]$\ast$a[3,4]$
\ast$a[4,2];$\backslash$
\newline
diag[4,4]=a[3,1]$\ast$a[4,3]$\ast$a[4,2];b=a$\ast$diag;b=b/content(b);$
\backslash$
\newline
db=smith(b);dbb=db[2];fdbb=factor(dbb);$\backslash$
\newline
if(issqfree(b[1,1])==0$||$issqfree(b[2,2])==0$||$issqfree(b[3,3])==0$||$$
\backslash$
\newline
issqfree(b[4,4])==0$||$content([--b[1,1],--b[2,2]])$>$2$||$$\backslash$
\newline
content([--b[2,2],--b[3,3]])$>$2$||$$\backslash$
\newline
(alpha[3,4]==2\&\&content([--b[3,3],--b[4,4]])$>$1)$||$$\backslash$
\newline
(alpha[3,4]==3\&\&content([--b[3,3],--b[4,4]])$>$2),,$\backslash$
\newline
detb=db[2]$\ast$db[3]$\ast$db[4];if(detb==1,,fdetb=factor(detb);$\backslash$
\newline
if(content([--b[1,1]$\ast$--b[2,2]$\ast$--b[3,3]$\ast$--b[4,4],16])$>$1\&\&$
\backslash$
\newline
content([--b[1,1]$\ast$--b[2,2]$\ast$--b[3,3]$\ast$--b[4,4],16])$<$16\&\&$
\backslash$
\newline
fdetb[1,1]==2\&\&mod(fdetb[1,2],2)==mod(1,2),,$\backslash$
\newline
gam=0;for(j=1,matsize(fdetb)[1],if(fdetb[j,1]!=2\&\&$\backslash$
\newline
type(--b[1,1]$\ast$--b[2,2]$\ast$--b[3,3]$\ast$--b[4,4]/fdetb[j,1])==1\&\&$
\backslash$
\newline
mod(fdetb[j,2],2)==mod(0,2),gam=1,));if(gam==1,,for(k=0,1,$\backslash$
\newline
if(k==0\&\&content(--[b[1,1]$\ast$--b[2,2]$\ast$--b[3,3]$\ast
$--b[4,4],2])==1\&\&$\backslash$
\newline
fdetb[1,1]==2\&\&mod(fdetb[1,2],2)==mod(1,2),,$\backslash$
\newline
if(k==1\&\&content(--[b[1,1]$\ast$--b[2,2]$\ast$--b[3,3]$\ast
$--b[4,4],2])==1,$\backslash$
\newline
b=2$\ast$b;db=smith(b);dbb=db[2];fdbb=factor(dbb);gam1=0,gam1=1);$
\backslash$
\newline
if(k==1\&\&gam1==1,,detb=db[2]$\ast$db[3]$\ast$db[4];$\backslash$
\newline
fdetb=factor(detb);d=1;for(k=1,matsize(fdetb)[1],$\backslash$
\newline
if(mod(fdetb[k,2],2)==mod(1,2),d=d$\ast$fdetb[k,1],));if(d$<
$=2,et=0,$\backslash$
\newline
fd=factor(d);if(fd[1,1]==2,d1=d/2,d1=d);fd1=factor(d1);$\backslash$
\newline
et=0;for(k=1,matsize(fd1)[1],if(type(b[1,1]/fd1[k,1])==1,$
\backslash$
\newline
if(kro(b[1,1]/fd1[k,1],fd1[k,1])==1,,et=et+2$^\wedge$(k--1)),$
\backslash$
\newline
if(type(b[2,2]/fd1[k,1])==1,$\backslash$
\newline
if(kro(b[2,2]/fd1[k,1],fd1[k,1])==1,,et=et+2$^\wedge$(k--1)),$\backslash$
\newline
if(kro(d$\ast$b[1,1]$\ast$b[2,2]/fd1[k,1],fd1[k,1])==1,,et=et+2$
^\wedge$(k--1))))));$\backslash$
\newline
hhh=hnr(d,et);if(hhh$>$2$||$hhh==1,,$\backslash$
\newline
n=n+1;pprint("n=",n);pprint(a);pprint(b);pprint(db);$\backslash$
\newline
pprint(fdbb);pprint("d=",d," et=",et," h=",hhh);$\backslash$
\newline
if(n==1,uuu=matrix(1,3,j,k,0);uuu[1,]=[d,et,hhh],$\backslash$
\newline
aaa=0;nnn=matsize(uuu)[1];$\backslash$
\newline
for(j=1,nnn,if(uuu[j,]==[d,et,hhh],aaa=1,));$\backslash$
\newline
if(aaa==1,,uuu1=matrix(nnn+1,3,j,k,0);for(t=1,nnn,$\backslash$
\newline
if(uuu[t,1]$<$d$||$(uuu[t,1]==d\&\&uuu[t,2]$<$et),uuu1[t,]=uuu[t,],$
\backslash$
\newline
uuu1[t+1,]=uuu[t,]));$\backslash$
\newline
for(t=1,nnn+1,if(uuu1[t,]==[0,0,0],uuu1[t,]=[d,et,hhh],));uuu=uuu1));$
\backslash$
\newline
))))))))))))))))));pprint("uuu=",uuu);$\backslash$
\newline
pprint("matsize uuu=",matsize(uuu));

\vskip10pt 

\centerline{{\bf Program 31:} refl0.14} 

\noindent 
$\backslash$$\backslash$main(n,h) calculates for
\newline
$\backslash$$\backslash$$<$3n$>$$\backslash$oplus A$\_$2(1/3,1/3,--1/3) where 
\newline
$\backslash$$\backslash$where n$\backslash$equiv 2 mod 3
\newline
$\backslash$$\backslash$it is given by the matrix g
\newline
$\backslash$$\backslash$the vectors v of the height $
\backslash$le h from P($\backslash$M)$\_${$\backslash$pr}
\newline
$\backslash$$\backslash$and calculates chains e and f 
and their Gram matrices ge, gf
\newline
$\backslash$$\backslash$comment:
\newline
$\backslash$$\backslash$possible squares d of roots divide dd
\newline
$\backslash$$\backslash$y1t=product of a root with f=(1,0,0)
\newline
$\backslash$$\backslash$(y1,y2,y3) are coordinates of a root
\newline
$\backslash$$\backslash$we have y1$\backslash$ge 0 and y2,y3$
\backslash$le 0;z1=3y1;z2=--3y2;z3=--3y3;
\newline
refl(n,h,m,k,v1,y1,y2,y3,y1t,w1,w,z1,z2,z3,dd,h1,u,$\backslash$
\newline
u1,w,w1,alpha,m1,divdd,mdivdd,h11)=$\backslash$
\newline
g=[3$\ast$n,0,0;0,--2,1;0,1,--2];eps=[1/3,1/3,--1/3];$\backslash$
\newline
m=2;v=matrix(m,3,j,k,0);$\backslash$
\newline
v[1,]=[0,1,0];v[2,]=[0,0,1];$\backslash$
\newline
dd=2$\ast$n;divdd=divisors(dd);mdivdd=matsize(divdd)[2];$
\backslash$
\newline
forstep(h1=n,h,n,$\backslash$
\newline
for(j1=1,mdivdd,d=divdd[j1];h11=d$\ast$h1/2;$\backslash$
\newline
if(type(h11)!=1,,$\backslash$
\newline
if(issquare(h11)!=1,,y1t=isqrt(h11);$\backslash$
\newline
if(type(2$\ast$y1t/d)!=1,,y1=y1t/(3$\ast$n);z1=y1$\ast$3;$
\backslash$
\newline
if(type(z1)!=1,,$\backslash$
\newline
for(z3=0,floor(sqrt(w1=2$\ast$n$\ast$z1$^\wedge$2+6$\ast$d)+0.000001),$
\backslash$
\newline
if(issquare(w=3$\ast$w1--3$\ast$z3$^\wedge$2)!=1,,$\backslash$
\newline
z2=(z3+isqrt(w))/2;$\backslash$
\newline
if(type(z2)!=1,,$\backslash$
\newline
if(mod(z1,3)!=mod(--z2,3)$||$mod(z1,3)!=mod(z3,3),,$\backslash$
\newline
y2=--z2/3;y3=--z3/3;u=[y1,y2,y3];$\backslash$
\newline
if(type(2$\ast$[1,0,0]$\ast$g$\ast$u~/d)!=1$||$type(2$\ast$[0,1,0]$\ast
$g$\ast$u~/d)!=1$||$$\backslash$
\newline
type(2$\ast$[0,0,1]$\ast$g$\ast$u~/d)!=1$||$type(2$\ast$eps$\ast
$g$\ast$u~/d)!=1,,$\backslash$
\newline
alpha=1;m1=1;$\backslash$
\newline
while(alpha,if(m1$>$m,alpha=0,if(v[m1,]$\ast$g$\ast$u~$>
$=0,m1=m1+1,alpha=0));$\backslash$
\newline
if(m1$<$=m,,m=m1;v1=matrix(m,3,j,k,0);$\backslash$
\newline
for(j=1,m--1,v1[j,]=v[j,]);v1[m,]=[y1,y2,y3];v=v1;kill(v1)$
\backslash$
\newline
)))))))))))));v;

\vskip10pt

\noindent
$\backslash$$\backslash$$\ast$$\ast$
\newline
This part is the same as in Program 10: refl0.1 between $
\backslash
$$\backslash$$\ast$$\ast$ and $\backslash$$\backslash$$\ast
$$\ast$$\ast$

\noindent
$\backslash$$\backslash$$\ast$$\ast$$\ast$

\vskip10pt  

\noindent
$\backslash$l;$\backslash$
\newline
main(n,h)=$\backslash$
\newline
refl(n,h);$\backslash$
\newline
e1fromv(v);e2fromv(v,e1);efromv(v,e1,e2);$\backslash$
\newline
f1fromv(v,e);f2fromv(f1);ffromv(f1,f2);

\enddocument

\end